# Differential Equation over Banach Algebra

Aleks Kleyn

Aleks_Kleyn@MailAPS.org
http://AleksKleyn.dyndns-home.com:4080/
http://sites.google.com/site/AleksKleyn/
http://arxiv.org/a/kleyn_a_1
http://AleksKleyn.blogspot.com/

Abstract. In the book, I considered differential equations of order 1 over Banach $D$-algebra: differential equation solved with respect to the derivative; exact differential equation; linear homogeneous equation.

I considered examples of differential equations in quaternion algebra.

In order to study homogeneous system of linear differential equations, I considered vector space over division $D$-algebra, solving of linear equations over division $D$-algebra and the theory of eigenvalues in non commutative division $D$-algebra.

I considered example of homogeneous system of linear differential equations in quaternion algebra, for which initial value problem has infinitely many solutions.

# Contents

CHAPTER 1

# Preface

## 1.1. Intention of the Book

When I began to study noncommutative algebra, I could write a linear map only together with its argument. Tensor representation of linear map allowed me to simplify notation and expanded the field of my research. Now I am ready to study differential equations over Banach algebra.

Since the product is noncommutative, the set of equations which I can consider is more limited than in the commutative case. However the set of differential equations considered in this book is source of extremely interesting theory.

## 1.2. Preface to Version 1

Derivative of map in Banach spaces has different representations: if we consider a map $f$ of $D$-algebra $A$ represented using coordinates with respect to basis of $D$-module $A$, then derivative is represented by Jacoby matrix of the map $f$; if we consider coordinateless representation of the map $f$, then derivative is represented by $A\otimes A$- number. So I decided to write this paper such way that somebody can use it regardless of the presentation, and, in examples, I want to consider solving of differential equation in different representations.

I dedicated the chapter 6 to differential equation solved with respect to the derivative. To make the picture complete, I consider different forms of representation of differential equation and compare the results. Examples in the text give ability to see how calculations work in noncommutative algebra. At the same time, these examples are source of new ideas and they are an integral part of the book.

The statement that problem with initial condition for linear homogeneous equation has infinitely many solutions was the most interesting and unexpected statement.

January, 2018

## 1.3. Preface to Version 2

I was extremely interested in papers [16], [18]. The authors of these papers, as well as I, consider a system of differential equations for maps into quaternion algebra. However, maps in these papers depend on the real variable. This makes it easier to solve the problem and, at the same time, allows us to see new details.

So, without delay, I decided to consider a system of differential equations

$$\frac{dx^1}{dt} = a^1_1 x^1 + ... + a^1_n x^n$$
$$......$$
$$\frac{dx^n}{dt} = a^n_1 x^1 + ... + a^n_n x^n$$

in general, as did the authors of papers [16], [18], and then try to solve particular system of differential equations.

I did not see any problem with system of differential equations

$$\frac{dx^1}{dt} = x^2$$
$$\frac{dx^2}{dt} = x^1$$

This system of differential equations has unique solution which satisfies to initial condition

$$t = 0 \quad x^1 = 0 \quad x^2 = 1$$

and I immediately get the Euler's formula. However the system of differential equations

$$\frac{dx^1}{dt} = x^2$$
$$\frac{dx^2}{dt} = -x^1$$

has infinitely many solutions which satisfies to initial condition

$$t = 0 \quad x^1 = 0 \quad x^2 = 1$$

This statement alerted me; and I decided to return to differential equation

$$\frac{dy}{dx} = \frac{1}{2}(y \otimes 1 + 1 \otimes y) \tag{1.3.1}$$

with initial condition

$$y(0) = 1$$

Although, last year, I came to a conclusion that this problem also has infinitely many solutions, it was important for me to understand why Taylor series decomposition is unique, but I see infinitely many solutions.

When I started to study the differential equation (1.3.1), I was strongly influenced by the classical tradition: exponent is unique mapping which coincides with its derivative. However, this is not true in non commutative Banach algebra. When I realized this, I found a large set of maps, which are different from exponent and preserve their Taylor series decomposition after differentiation. Since behavior of these maps is similar to behavior of exponent, I called these maps quasiexponent.

As always, I have a lot of unanswered questions.

One of these questions relates to the statement that differential equation for exponent is not completely integrable in non commutative algebra.

I am interested in how it is possible that the system of differential equations

$$\frac{dx^1}{dt} = x^2$$
$$\frac{dx^2}{dt} = -x^1$$

has non unique solution, while Taylor series expansion is unique. This example also does not agree with the theorem [17]-3.1 on the uniqueness of the solution of the initial value problem.

The theory of eigenvalues is closely related to the theory of polynomials over non commutative algebra. This is a new and rapidly developing theory. According to the theorem [14]-7.2, in quaternion algebra, there exists quadratic equation which has 1 root. According to the theorem [14]-7.3, in quaternion algebra, there exists quadratic equation which has no root. Therefore, $2 \times 2$ matrix can have a unique eigenvalue or not have eigenvalues at all.

There is one more question: the question about covariance. We solve the system of differential equations with respect to selected basis in module over $D$-algebra. How does the system of differential equations change and how does its solution change, if we select another basis?

May, 2019

## 1.4. Preface to Version 3

Mathematics is deductive science. It means that the statement $B$ should follow from the statement $A$. It seems everything is simple. However, deduction is not enough when you try to solve the problem. Induction is a step into the unknown, step to where you have never been before. Or maybe it is terra incognita, where somebody has never walked.

In 2008, I started to study calculus over non commutative algebra. However it took me one more year to understand the structure of linear map. Almost immediately I realized that derivative of higher order is a polylinear map which is symmetric with respect to the order of arguments. This allowed me to solve the most simple differential equations. Thereafter I tried to understand what is the form of the differential equation for exponent.

It was clear that classical equation $y' = y$ is not good because on the left and on the right sides I see different types of expression. After a brief search, I discovered that the equation

$$\frac{dy}{dx} \circ h = \frac{1}{2}(yh + hy) \tag{1.4.1}$$

has the symmetry that I am looking for. Even more, Taylor series expansion of this solution was the same as for exponent. I was sure that I found proper solution.

The last year I started to write the book dedicated to differetial equations over non commutative algebra. In particular, I wrote down the condition when a differential equation has a solution. It turned out that the differential equation which I have written for exponent is wrong. At this time, the search of a differential equation took more time. I discovered the differential equation like

$$\frac{dy}{dx} \circ 1 = y$$

However it turned out that this equation has a lot of other solutions besides the exponent. At the beginning, these maps seemed to me monsters; then I put up with their existence. Suddenly I realized that these maps are very important. So I called new map quasiexponent.

Quasiexponent occurs in the process of differentiation of exponent. In commutative algebra, the differential equation $y' = y$ says that the rate of change of the map $y = e^x$ is equal $y$. Quasiexponent plays a similar role.

The research is not finished. I see a lot of questions that I need to answer.

June, 2019

## 1.5. Conventions

Convention 1.5.1. *We will use Einstein summation convention in which repeated index (one above and one below) implies summation with respect to repeated index. In this case we assume that we know the set of summation index and do not use summation symbol*

$$c^i v_i = \sum_{i \in I} c^i v_i$$

*If needed to clearly show a set of indices, I will do it.* □

Convention 1.5.2. *Let $A$ be free algebra with finite or countable basis. Considering expansion of element of algebra $A$ relative basis $\overline{\overline{e}}$ we use the same root letter to denote this element and its coordinates. In expression $a^2$, it is not clear whether this is component of expansion of element $a$ relative basis, or this is operation $a^2 = aa$. To make text clearer we use separate color for index of element of algebra. For instance,*

$$a = a^i e_i$$

□

Convention 1.5.3. *If free finite dimensional algebra has unit, then we identify the vector of basis $e_0$ with unit of algebra.* □

CHAPTER 2

# Preliminary Statements

## 2.1. Module of Linear Maps into $D$-Algebra

Let $A$ be $D$-module. Let $B$ be free finite dimensional associative $D$-algebra. According to the theorem [13]-6.4.8, left $B\otimes B$-module $\mathcal{L}(D;A\to B)$ has finite basis $\overline{\overline{F}}_{AB}$. However, it is obvious that in the general case the choice of this basis is arbitrary. At the same time, the choice of basis $\overline{\overline{F}}$ of left $B\otimes B$-module $\mathcal{L}(D;B\to B)$ is usually associated with a structure of $D$-algebra $B$. So the question arises whether there exists relation between bases $\overline{\overline{F}}_{AB}$ and $\overline{\overline{F}}$.

THEOREM 2.1.1. *Let $A$ be $D$-module, $n=\dim A$. Let $B$ be associative $D$-algebra, $m=\dim B$. Let $\overline{\overline{F}}$ be basis of left $B\otimes B$-module $\mathcal{L}(D;B\to B)$. Let $n\le m$. Let*

$$G:A\to B$$

*be linear map of maximal rank. The set*

$$\overline{\overline{F}}\circ G=\{F_k\circ G:F_k\in\overline{\overline{F}}\}\tag{2.1.1}$$

*generates left $B\otimes B$-module*[2.1] *$\mathcal{L}(D;A\to B)$.*

PROOF. Let

$$g:A\to B$$

be a linear map. Let $\overline{\overline{e}}_A$ be the basis of $D$-module $A$. Let $\overline{\overline{e}}_B$ be the basis of $D$-module $B$. According to the theorem [13]-4.2.2, the linear map $G$ has coordinates

$$G=\begin{pmatrix}G^1_1&...&G^1_n\\...&...&...\\G^m_1&...&G^m_n\end{pmatrix}$$

with respect to bases $\overline{\overline{e}}_A$, $\overline{\overline{e}}_B$ and the linear map $g$ has coordinates

$$g=\begin{pmatrix}g^1_1&...&g^1_n\\...&...&...\\g^m_1&...&g^m_n\end{pmatrix}$$

with respect to bases $\overline{\overline{e}}_A$, $\overline{\overline{e}}_B$. A row $G_i$ of the matrix $G$, as well a row $g_i$ of the matrix $g$, is coordinates of linear form $A\to D$. Since the matrix $G$ has maximal

[2.1] I do not claim that this set is a basis, because maps $F_i\circ G$, $i\in I$, can be linearly dependent.

rank, then rows of the matrix $G$ generate $D$-module $\mathcal{L}(D;A\to D)$ and rows of the matrix $g$ are linear combination of rows of the matrix $G$

$$g^i_k = C^i_j G^j_k \tag{2.1.2}$$

Since we can consider the matrix $C$ as coordinates of linear map

$$C : B \to B$$

then the equality

$$g^i_k = (c_{1.k.l} F_{k\cdot}{}^i_j c_{2.k.l}) G^j_k \tag{2.1.3}$$

follows from the equality (2.1.2) and from the equality

$$C = c_{1.k.l} F_k c_{2.k.l}$$

Since $G^j_k \in D$, then the equality

$$g^i_k = c_{1.k.l}(F_{k\cdot}{}^i_j G^j_k) c_{2.k.l} \tag{2.1.4}$$

follows from the equality (2.1.3). Therefore, the map $g$ belongs to linear span of the set of maps (2.1.1). □

Theorem 2.1.2. *Let* $n = \dim A$, $m = \dim B$. *Let* $\overline{\overline{F}}$ *be basis of left* $B\otimes B$-*module* $\mathcal{L}(D;B\to B)$. *Let* $n > m$. *Let*

$$G : A \to B$$

*be linear map of maximal rank. The set*

$$\overline{\overline{F}} \circ G = \{F_k \circ G : F_k \in \overline{\overline{F}}\}$$

*generates the set of maps*

$$\{g \in \mathcal{L}(D;A\to B) : \ker G \subseteq \ker g\} \tag{2.1.5}$$

Proof. The proof of the theorem is similar to the proof of the theorem 2.1.1. However, since number of rows of the matrix $G$ less then dimention of $D$-module $A$, then rows of the matrix $G$ do not generate $D$-module $\mathcal{L}(D;A\to D)$ and the map $G$ has non trivial kernel. In particular, rows of the matrix $g$ linearly depend on rows of the matrix $G$ iff $\ker G \subseteq \ker g$. □

From the theorem 2.1.2, it follows that choice of the map $G$ depends on the map $g$. It is easy to see that theorem 2.1.1 is particular case of the theorem 2.1.2, because, in the theorem 2.1.2, $\ker G = \emptyset$.

Theorem 2.1.3. *Let* $\overline{\overline{e}}$ *be the basis of* $D$-*module* $A$. *Then the set of maps*

$$P_i : d \in D \to de_i \in A$$

*is the basis of* $D$-*module* $\mathcal{L}(D;D\to A)$.

Proof. The theorem follows from the equality

$$f \circ t = (f^i t) e_i = P_i \circ (f^i t)$$

□

## 2.2. Direct Sum of $D$-modules

DEFINITION 2.2.1. *Let $\mathcal{A}$ be a category. Let $\{B_i, i \in I\}$ be the set of objects of $\mathcal{A}$. Object*

$$P = \coprod_{i \in I} B_i$$

*and set of morphisms*

$$\{f_i : B_i \to P, i \in I\}$$

*is called a* **coproduct of set of objects** $\{B_i, i \in I\}$ **in category** $\mathcal{A}$[2.2] *if for any object $R$ and set of morphisms*

$$\{g_i : B_i \to R, i \in I\}$$

*there exists a unique morphism*

$$h : P \to R$$

*such that diagram*

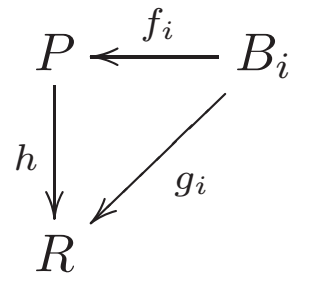

$$h \circ f_i = g_i$$

*is commutative for all $i \in I$.*

*If $|I| = n$, then we also will use notation*

$$P = \coprod_{i=1}^{n} B_i = B_1 \coprod ... \coprod B_n$$

*for coproduct of set of objects $\{B_i, i \in I\}$ in $\mathcal{A}$.* □

DEFINITION 2.2.2. *Coproduct in category of Abelian groups Ab is called* **direct sum.**[2.3] *We will use notation $A \oplus B$ for direct sum of Abelian groups $A$ and $B$.* □

THEOREM 2.2.3. *Let $\{A_i, i \in I\}$ be set of Abelian groups. Let*

$$A \subseteq \prod_{i \in I} A_i$$

*be such set that $(x_i, i \in I) \in A$, if $x_i \neq 0$ for finite number of indices $i$. Then*[2.4]

$$A = \bigoplus_{i \in I} A_i \tag{2.2.1}$$

PROOF. According to construction, $A$ is subgroup of Abelian group $\prod A_i$. The map

$$\lambda_j : A_j \to A$$

defined by the equality

$$\lambda_j(x) = (\delta^i_j x, i \in I) \tag{2.2.2}$$

is an injective homomorphism.

[2.2] I made definition according to [1], page 59.
[2.3] See also definition in [1], pages 36, 37.
[2.4] See also proposition [1]-7.1, page 37.

Let

$$\{f_i : A_i \to B, i \in I\}$$

be set of homomorphisms into Abelian group $B$. We define the map

$$f : A \to B$$

by the equality

$$f\left(\bigoplus_{i\in I} x_i\right) = \sum_{i\in I} f_i(x_i) \tag{2.2.3}$$

The sum in the right side of the equality (2.2.3) is finite, since all summands, except for a finite number, equal 0. From the equality

$$f_i(x_i + y_i) = f_i(x_i) + f_i(y_i)$$

and the equality (2.2.3), it follows that

$$\begin{aligned} f(\bigoplus_{i\in I}(x_i + y_i)) &= \sum_{i\in I} f_i(x_i + y_i) = \sum_{i\in I}(f_i(x_i) + f_i(y_i)) \\ &= \sum_{i\in I} f_i(x_i) + \sum_{i\in I} f_i(y_i) \\ &= f(\bigoplus_{i\in I} x_i) + f(\bigoplus_{i\in I} y_i) \end{aligned}$$

Therefore, the map $f$ is homomorphism of Abelian group. The equality

$$f \circ \lambda_j(x) = \sum_{i\in I} f_i(\delta^i_j x) = f_j(x)$$

follows from equalities (2.2.2), (2.2.3). Since the map $\lambda_i$ is injective, then the map $f$ is unique. Therefore, the theorem folows from definitions 2.2.1, 2.2.2. □

THEOREM 2.2.4. *Direct sum of Abelian groups* $A_1$, ..., $A_n$ *coincides with their Cartesian product*

$$A_1 \oplus ... \oplus A_n = A_1 \times ... \times A_n$$

PROOF. The theorem follows from the theorem 2.2.3. □

Let

$$A = A_1 \oplus ... \oplus A_n$$

be direct sum of Abelian groups $A_1$, ..., $A_n$. According to the proof of the theorem 2.2.3, any $A$-number $a$ has form $(a_1, ..., a_n)$ where $a_i \in A_i$. We also will use notation

$$a = a_1 \oplus ... \oplus a_n$$

DEFINITION 2.2.5. *Coproduct in category of* $D$*-modules is called* **direct sum**.[2.5] *We will use notation* $A \oplus B$ *for direct sum of* $D$*-modules* $A$ *and* $B$. □

THEOREM 2.2.6. *Let* $\{A_i, i \in I\}$ *be set of* $D$*-modules. Then the representation*

$$D \xrightarrow{*} \bigoplus_{i\in I} A_i \qquad d\left(\bigoplus_{i\in I} a_i\right) = \bigoplus_{i\in I} da_i$$

[2.5] See also the definition of direct sum of modules in [1], page 128. On the same page, Lang proves the existence of direct sum of modules.

*of the ring $D$ in direct sum of Abelian groups*

$$A=\bigoplus_{i\in I}A_i$$

*is direct sum of $D$-modules*

$$A=\bigoplus_{i\in I}A_i$$

Proof.
Let

$$\{f_i:A_i\to B,i\in I\}$$

be set of linear maps into $D$-module $B$. We define the map

$$f:A\to B$$

by the equality

$$f\left(\bigoplus_{i\in I}x_i\right)=\sum_{i\in I}f_i(x_i) \tag{2.2.4}$$

The sum in the right side of the equality (2.2.4) is finite, since all summands, except for a finite number, equal 0. From the equality

$$f_i(x_i+y_i)=f_i(x_i)+f_i(y_i)$$

and the equality (2.2.4), it follows that

$$\begin{aligned}f(\bigoplus_{i\in I}(x_i+y_i))&=\sum_{i\in I}f_i(x_i+y_i)=\sum_{i\in I}(f_i(x_i)+f_i(y_i))\\&=\sum_{i\in I}f_i(x_i)+\sum_{i\in I}f_i(y_i)\\&=f(\bigoplus_{i\in I}x_i)+f(\bigoplus_{i\in I}y_i)\end{aligned}$$

From the equality

$$f_i(dx_i)=df_i(x_i)$$

and the equality (2.2.4), it follows that

$$\begin{aligned}f((dx_i,i\in I))&=\sum_{i\in I}f_i(dx_i)=\sum_{i\in I}df_i(x_i)=d\sum_{i\in I}f_i(x_i)\\&=df((x_i,i\in I))\end{aligned}$$

Therefore, the map $f$ is linear map. The equality

$$f\circ\lambda_j(x)=\sum_{i\in I}f_i(\delta^i_jx)=f_j(x)$$

follows from equalities (2.2.2), (2.2.4). Since the map $\lambda_i$ is injective, then the map $f$ is unique. Therefore, the theorem folows from definitions 2.2.1, 2.2.2. □

Theorem 2.2.7. *Direct sum of $D$-modules $A_1$, ..., $A_n$ coincides with their Cartesian product*

$$A_1\oplus...\oplus A_n=A_1\times...\times A_n$$

Proof. The theorem follows from the theorem 2.2.6. □

Theorem 2.2.8. *Let* $A^1$, ..., $A^n$ *be* $D$*-modules and*

$$A = A^1 \oplus ... \oplus A^n$$

*Let us represent* $A$*-number*

$$a = a^1 \oplus ... \oplus a^n$$

*as column vector*

$$a = \begin{pmatrix} a^1 \\ ... \\ a^n \end{pmatrix}$$

*Let us represent a linear map*

$$f : A \to B$$

*as row vector*

$$f = \begin{pmatrix} f_1 & ... & f_n \end{pmatrix}$$

$$f_i : A^i \to B$$

*Then we can represent value of the map* $f$ *in* $A$*-number* $a$ *as product of matrices*

$$f \circ a = \begin{pmatrix} f_1 & ... & f_n \end{pmatrix} {}_\circ{}^\circ \begin{pmatrix} a^1 \\ ... \\ a^n \end{pmatrix} = f_i \circ a^i \tag{2.2.5}$$

Proof. The theorem follows from the definition (2.2.4). □

Theorem 2.2.9. *Let* $B^1$, ..., $B^m$ *be* $D$*-modules and*

$$B = B^1 \oplus ... \oplus B^m$$

*Let us represent* $B$*-number*

$$b = b^1 \oplus ... \oplus b^m$$

*as column vector*

$$b = \begin{pmatrix} b^1 \\ ... \\ b^m \end{pmatrix}$$

*Then the linear map*

$$f : A \to B$$

*has representation as column vector of maps*

$$f = \begin{pmatrix} f^1 \\ ... \\ f^m \end{pmatrix}$$

*such way that, if* $b = f \circ a$, *then*

$$\begin{pmatrix} b^1 \\ ... \\ b^m \end{pmatrix} = \begin{pmatrix} f^1 \\ ... \\ f^m \end{pmatrix} \circ a = \begin{pmatrix} f^1 \circ a \\ ... \\ f^m \circ a \end{pmatrix}$$

PROOF. The theorem follows from the theorem [13]-2.1.5. □

THEOREM 2.2.10. *Let* $A^1$, ..., $A^n$, $B^1$, ..., $B^m$ *be* $D$*-modules and*

$$A = A^1 \oplus ... \oplus A^n$$
$$B = B^1 \oplus ... \oplus B^m$$

*Let us represent* $A$*-number*

$$a = a^1 \oplus ... \oplus a^n$$

*as column vector*

$$a = \begin{pmatrix} a^1 \\ ... \\ a^n \end{pmatrix}$$

*Let us represent* $B$*-number*

$$b = b^1 \oplus ... \oplus b^m$$

*as column vector*

$$b = \begin{pmatrix} b^1 \\ ... \\ b^m \end{pmatrix}$$

*Then the linear map* $f$ *has representation as a matrix of maps*

$$f = \begin{pmatrix} f^1_1 & ... & f^1_n \\ ... & ... & ... \\ f^m_1 & ... & f^m_n \end{pmatrix}$$

*such way that, if* $b = f \circ a$*, then*

$$\begin{pmatrix} b^1 \\ ... \\ b^m \end{pmatrix} = \begin{pmatrix} f^1_1 & ... & f^1_n \\ ... & ... & ... \\ f^m_1 & ... & f^m_n \end{pmatrix} {}_\circ{}^\circ \begin{pmatrix} a^1 \\ ... \\ a^n \end{pmatrix} = \begin{pmatrix} f^1_i \circ a^i \\ ... \\ f^m_i \circ a^i \end{pmatrix} \tag{2.2.6}$$

*The map*

$$f^i_j : A^j \to B^i$$

*is a linear map and is called* **partial linear map**.

PROOF. According to the theorem 2.2.9, there exists the set of linear maps

$$f^i : A \to B^i$$

such that

$$\begin{pmatrix} b^1 \\ ... \\ b^m \end{pmatrix} = \begin{pmatrix} f^1 \\ ... \\ f^m \end{pmatrix} \circ a = \begin{pmatrix} f^1 \circ a \\ ... \\ f^m \circ a \end{pmatrix} \tag{2.2.7}$$

According to the theorem 2.2.8, for every $i$, there exists the set of linear maps

$$f^i_j : A^j \to B^i$$

such that

$$f^i \circ a = \begin{pmatrix} f^i_1 & ... & f^i_n \end{pmatrix} {}_\circ{}^\circ \begin{pmatrix} a^1 \\ ... \\ a^n \end{pmatrix} = f^i_j \circ a^j \tag{2.2.8}$$

If we identify matrices

$$\begin{pmatrix} \begin{pmatrix} f^1_1 & ... & f^1_n \end{pmatrix} \\ ... \\ \begin{pmatrix} f^m_1 & ... & f^m_n \end{pmatrix} \end{pmatrix} = \begin{pmatrix} f^1_1 & ... & f^1_n \\ ... & ... & ... \\ f^m_1 & ... & f^m_n \end{pmatrix}$$

then the equality (2.2.6) follows from equalities (2.2.7), (2.2.8). □

## 2.3. Direct Sum of Banach $D$-Modules

THEOREM 2.3.1. *Let $A^1$, ..., $A^n$ be Banach $D$-modules and*

$$A = A^1 \oplus ... \oplus A^n$$

*Then, in $D$-module $A$, we can introduce norm such that $D$-module $A$ is Banach $D$-module.*

PROOF. Let $\|a^i\|_i$ be norm in $D$-module $A^i$.

2.3.1.1: We introduce norm in $D$-module $A$ by the equality

$$\|b\| = \max(\|b^i\|_i, i = 1, ..., n)$$

where

$$b = b^1 \oplus ... \oplus b^n$$

2.3.1.2: Let $\{a_p\}$, $p = 1, ...,$ be fundamental sequence where

$$a_p = a^1_p \oplus ... \oplus a^n_p$$

2.3.1.3: Therefore, for any $\epsilon \in R$, $\epsilon > 0$, there exists $N$ such that for any $p$, $q > N$

$$\|a_p - a_q\| < \epsilon$$

2.3.1.4: According to statements 2.3.1.1, 2.3.1.2, 2.3.1.3,

$$\|a^i_p - a^i_q\|_i < \epsilon$$

for any $p$, $q > N$ and $i = 1, ..., n$.

2.3.1.5: Therefore, the sequence $\{a^i_p\}$, $i = 1, ..., n$, $p = 1, ...,$ is fundamental sequence in $D$-module $A^i$ and there exists limit

$$a^i = \lim_{p\to\infty} a^i_p$$

2.3.1.6: Let

$$a = a^1 \oplus ... \oplus a^n$$

2.3.1.7: According to the statement 2.3.1.5, for any $\epsilon \in R$, $\epsilon > 0$, there exists $N_i$ such that for any $p > N_i$

$$\|a^i - a^i_p\|_i < \epsilon$$

2.3.1.8: Let

$$N = \max(N_1, ..., N_n)$$

2.3.1.9: According to statements 2.3.1.6, 2.3.1.7, 2.3.1.8, for any $\epsilon \in R$, $\epsilon > 0$, there exists $N$ such that for any $p > N$

$$\|a - a_p\|_i < \epsilon$$

2.3.1.10: Therefore,

$$a = \lim_{p\to\infty} a_p$$

The theorem follows from statements 2.3.1.1, 2.3.1.2, 2.3.1.10. □

Using the theorem 2.3.1, we can consider the derivative of a map

$$f : A^1 \oplus ... \oplus A^n \to B^1 \oplus ... \oplus B^m$$

Theorem 2.3.2. *Let* $A^1, ..., A^n$, $B^1, ..., B^m$ *be Banach* $D$*-modules and*

$$A = A^1 \oplus ... \oplus A^n$$
$$B = B^1 \oplus ... \oplus B^m$$

*Let us represent differential*

$$dx = dx^1 \oplus ... \oplus dx^n$$

*as column vector*

$$dx = \begin{pmatrix} dx^1 \\ ... \\ dx^n \end{pmatrix}$$

*Let us represent differential*

$$dy = dy^1 \oplus ... \oplus dy^m$$

*as column vector*

$$dy = \begin{pmatrix} dy^1 \\ ... \\ dy^m \end{pmatrix}$$

*Then the derivative of the map*

$$f : A \to B$$
$$f = f^1 \oplus ... \oplus f^m$$

*has representation*

$$\frac{df}{dx} = \begin{pmatrix} \frac{\partial f^1}{\partial x^1} & ... & \frac{\partial f^1}{\partial x^n} \\ ... & ... & ... \\ \frac{\partial f^m}{\partial x^1} & ... & \frac{\partial f^m}{\partial x^n} \end{pmatrix}$$

*such way that*

$$\begin{pmatrix} dy^1 \\ ... \\ dy^m \end{pmatrix} = \begin{pmatrix} \frac{\partial f^1}{\partial x^1} & ... & \frac{\partial f^1}{\partial x^n} \\ ... & ... & ... \\ \frac{\partial f^m}{\partial x^1} & ... & \frac{\partial f^m}{\partial x^n} \end{pmatrix} {}_\circ{}^\circ \begin{pmatrix} dx^1 \\ ... \\ dx^n \end{pmatrix} = \begin{pmatrix} \frac{\partial f^1}{\partial x^i} \circ dx^i \\ ... \\ \frac{\partial f^m}{\partial x^i} \circ dx^i \end{pmatrix} \tag{2.3.1}$$

Statement 2.3.3. *The linear map* $\frac{\partial f^i}{\partial x^j}$ *is called* **partial derivative** *and this map is the derivative of map* $f^i$ *with respect to variable* $x^j$ *assuming that other coordinates of* $A$*-number* $x$ *are fixed.* ⊙

Proof. The equality (2.3.1) follows from the equality (2.2.6).

We can represent the map

$$f^i : A \to B^i$$

as

$$f^i(x) = f^i(x^1, ..., x^n)$$

The equality

$$\frac{df^i(x)}{dx} \circ dx = \frac{\partial f^i(x^1, ..., x^n)}{\partial x^j} \circ dx^j \tag{2.3.2}$$

follows from the equality (2.3.1). According to the theorem [15]-3.3.4,

$$\begin{aligned}
\frac{df^i(x)}{dx} \circ dx =& \lim_{t\to 0,\ t\in R}(t^{-1}(f^(x + tdx) - f(x)))\\
=& \lim_{t\to 0,\ t\in R}(t^{-1}(f^i(x^1 + tdx^1, x^2 + tdx^2, ..., x^n + tdx^n)\\
&- f^i(x^1, x^2 + tdx^2, ..., x^n + tdx^n)\\
&+ f^i(x^1, x^2 + tdx^2, ..., x^n + tdx^n) - ...\\
&- f^i(x^1, x^2, ..., x^n)))\\
=& \lim_{t\to 0,\ t\in R}(t^{-1}(f^i(x^1 + tdx^1, x^2 + tdx^2, ..., x^n + tdx^n)\\
&- f^i(x^1, x^2 + tdx^2, ..., x^n + tdx^n))) + ...\\
&+ \lim_{t\to 0,\ t\in R}(t^{-1}(f^i(x^1, ..., x^n + tdx^n) - f^i(x^1, ..., x^n)))\\
=& f^i_1 \circ dx^1 + ... + f^i_n \circ dx^n
\end{aligned} \tag{2.3.3}$$

where $f^i_j$ is the derivative of map $f^i$ with respect to variable $x^j$ assuming that other coordinates of $A$-number $x$ are fixed. The equality

$$f^i_j = \frac{\partial f^i(x^1, ..., x^n)}{\partial x^j} \tag{2.3.4}$$

follows from equalities (2.3.2), (2.3.3). The statement 2.3.3 follows from the equality (2.3.4). □

Example 2.3.4. *Consider map*

$$\begin{aligned}
y^1 &= f^1(x^1, x^2, x^3) = (x^1)^2 + x^2x^3\\
y^2 &= f^2(x^1, x^2, x^3) = x^1x^2 + (x^3)^2
\end{aligned} \tag{2.3.5}$$

*Therefore*

$$\frac{\partial y^1}{\partial x^1} = x^1 \otimes 1 + 1 \otimes x^1 \qquad \frac{\partial y^1}{\partial x^2} = 1 \otimes x^3 \qquad \frac{\partial y^1}{\partial x^3} = x^2 \otimes 1$$

$$\frac{\partial y^2}{\partial x^1} = 1 \otimes x^2 \qquad \frac{\partial y^2}{\partial x^2} = x^1 \otimes 1 \qquad \frac{\partial y^2}{\partial x^3} = x^3 \otimes 1 + 1 \otimes x^3$$

*and the derivative of the map* (2.3.5) *is*

$$\frac{df}{dx}=\begin{pmatrix}x^1\otimes 1+1\otimes x^1 & 1\otimes x^3 & x^2\otimes 1\\ 1\otimes x^2 & x^1\otimes 1 & x^3\otimes 1+1\otimes x^3\end{pmatrix} \tag{2.3.6}$$

*The equality*

$$\begin{aligned}dy^1&=(x^1\otimes 1+1\otimes x^1)\circ dx^1+(1\otimes x^3)\circ dx^2+(x^2\otimes 1)\circ dx^3\\ &=x^1dx^1+dx^1x^1+dx^2x^3+x^2dx^3\\ dy^2&=(1\otimes x^2)\circ dx^1+(x^1\otimes 1)\circ dx^2+(x^3\otimes 1+1\otimes x^3)\circ dx^3\\ &=dx^1x^2+x^1dx^2+x^3dx^3+dx^3x^3\end{aligned} \tag{2.3.7}$$

*follows from the equality* (2.3.6). *We also can get the expression* (2.3.7) *by direct calculation*

$$\begin{aligned}dy^1&=f^1(x+dx)-f^1(x)\\ &=(x^1+dx^1)^2+(x^2+dx^2)(x^3+dx^3)-(x^1)^2-x^2x^3\\ &=(x^1)^2+x^1dx^1+dx^1x^1+x^2x^3+dx^2x^3+x^2dx^3-(x^1)^2-x^2x^3\\ &=x^1dx^1+dx^1x^1+dx^2x^3+x^2dx^3\\ dy^2&=f^2(x+dx)-f^2(x)\\ &=(x^1+dx^1)(x^2+dx^2)+(x^3+dx^3)^2-x^1x^2-(x^3)^2\\ &=x^1x^2+dx^1x^2+x^1dx^2+(x^3)^2+x^3dx^3+dx^3x^3-x^1x^2-(x^3)^2\\ &=dx^1x^2+x^1dx^2+x^3dx^3+dx^3x^3\end{aligned} \tag{2.3.8}$$

□

CHAPTER 3

# Biring of Matrices

## 3.1. Concept of Generalized Index

Studying tensor calculus we start from studying univalent covariant and contravariant tensors. In spite on difference of properties both these objects are elements of respective vector spaces. Suppose we introduce a generalized index according to the rule $a^{i} = a^{i}$, $b^{i} = b^{\cdot}{}_{i}^{-}$. Then we see that these tensors have similar behavior. For instance, the transformation of a covariant tensor takes on the form

$$b'^{i} = b'^{\cdot}{}_{i}^{-} = f^{\cdot}{}_{i}^{-}{}_{\cdot}^{j}{}_{-} b^{\cdot}{}_{j}^{-} = f^{i}_{j} b^{j}$$

This similarity goes as far as we need because tensors also form a vector space.

These observations of the similarity between properties of covariant and contravariant tensors lead us to the concept of generalized index. We will use the symbol $\cdot$ in front of a generalized index when we need to describe its structure. I put the sign $'-'$ in place of the index whose position was changed. For instance, if an original term was $a_{ij}$ I will use notation $a_{i\cdot}{}^{j}_{-}$ instead of notation $a^{j}_{i}$.

Even though the structure of a generalized index is arbitrary we assume that there exists a one-to-one map of the interval of positive integers $\mathbf{1}$, ..., $\mathbf{n}$ to the range of index. Let $\mathbf{I}$ be the range of the index $i$. We denote the power of this set by symbol $|\mathbf{I}|$ and assume that $|\mathbf{I}| = \mathbf{n}$. If we want to enumerate elements $a_{i}$ we use notation $a_{\mathbf{1}}$, ..., $a_{\mathbf{n}}$.

Representation of coordinates of a vector as a matrix allows making a notation more compact. The question of the presentation of vector as a row or a column of the matrix is just a question of convention. We extend the concept of generalized index to entries of the matrix. A matrix is a two dimensional table, the rows and columns of which are enumerated by generalized indeices. Since we use generalized index, we cannot tell whether index $a$ of matrix enumerates rows or columns until we know the structure of index. However as we can see below the form of presentation of matrix is not important for us. To make sure that notation offered below is consistent with the traditional we will assume that the matrix is presented in the form

$$a = \begin{pmatrix} a^{\mathbf{1}}_{\mathbf{1}} & ... & a^{\mathbf{1}}_{\mathbf{n}} \\ ... & ... & ... \\ a^{\mathbf{m}}_{\mathbf{1}} & ... & a^{\mathbf{m}}_{\mathbf{n}} \end{pmatrix}$$

The upper index enumerates rows and the lower index enumerates columns.

Definition 3.1.1. *I use the following names and notation for different* **minor matrices** *of the matrix $a$*

- **column of matrix** *with the index* $i$

$$a_i = \begin{pmatrix} a_i^1 \\ ... \\ a_i^m \end{pmatrix}$$

*The upper index enumerates entries of the column and the lower index enumerates columns.*

$a_T$ : *the minor matrix obtained from the matrix* $a$ *by selecting columns with an index from the set* $T$

$a_{[i]}$ : *the minor matrix obtained from the matrix* $a$ *by deleting column* $a_i$

$a_{[T]}$ : *the minor matrix obtained from the matrix* $a$ *by deleting columns with an index from the set* $T$

- **row of matrix** *with the index* $j$

$$a^j = \begin{pmatrix} a_1^j & ... & a_n^j \end{pmatrix}$$

*The lower index enumerates entries of the row and the upper index enumerates rows.*

$a^S$ : *the minor matrix obtained from the matrix* $a$ *by selecting rows with an index from the set* $S$

$a^{[j]}$ : *the minor matrix obtained from the matrix* $a$ *by deleting row* $a^j$

$a^{[S]}$ : *the minor matrix obtained from the matrix* $a$ *by deleting rows with an index from the set* $S$

□

Remark 3.1.2. *We will combine the notation of indices. Thus* $a_i^j$ *is* $1 \times 1$ *minor matrix. The same time this is the notation for a matrix entry. This allows an identifying of* $1 \times 1$ *matrix and its entry. The index* $a$ *is number of column of matrix and the index* $b$ *is number of row of matrix.* □

Let $I$, $|I| = n$, be a set of indices. We introduce the **Kronecker symbol**

$$\delta_j^i = \begin{cases} 1 & i = j \\ 0 & i \neq j \end{cases} \qquad i, j \in I$$

## 3.2. Biring

We consider matrices whose entries belong to associative division $D$-algebrs $A$.

The product of matrices is associated with the product of homomorphisms of vector spaces over field. According to the custom the product of matrices $a$ and $b$ is defined as product of rows of the matrix $a$ and columns of the matrix $b$.

Example 3.2.1. *Let* $\overline{\overline{e}}$ *be basis of right vector space* $V$ *over* $D$*-algebra* $A$ *(see the definition 5.1.4 and the theorem 5.2.7). We represent the basis* $\overline{\overline{e}}$ *as row of matrix*

$$e = \begin{pmatrix} e_1 & ... & e_n \end{pmatrix}$$

*We represent coordinates of vector $v$ as vector column*

$$v = \begin{pmatrix} v^1 \\ ... \\ v^n \end{pmatrix}$$

*Therefore, we can represent the vector $v$ as product of matrices*

$$v = \begin{pmatrix} e_1 & ... & e_n \end{pmatrix} \begin{pmatrix} v^1 \\ ... \\ v^n \end{pmatrix} = e_i v^i$$

*We represent homomorphism of right vector space $V$ using matrix*

$$v'^i = f^i_j v^j \tag{3.2.1}$$

*The equality (3.2.1) expresses a traditional product of matrices $f$ and $v$.* □

Example 3.2.2. *Let $\overline{\overline{e}}$ be basis of left vector space $V$ over $D$-algebra $A$ (see the definition 5.1.3 and the theorem 5.2.7). We represent the basis $\overline{\overline{e}}$ as row of matrix*

$$e = \begin{pmatrix} e_1 & ... & e_n \end{pmatrix}$$

*We represent coordinates of vector $v$ as vector column*

$$v = \begin{pmatrix} v^1 \\ ... \\ v^n \end{pmatrix}$$

*However, we cannot represent the vector*

$$v = v^i e_i$$

*as product of matrices*

$$v = \begin{pmatrix} v^1 \\ ... \\ v^n \end{pmatrix} \qquad e = \begin{pmatrix} e_1 & ... & e_n \end{pmatrix}$$

*because this product is not defined. We represent homomorphism of left vector space $V$ using matrix*

$$v'^i = v^j f^i_j \tag{3.2.2}$$

*We cannot express the equality (3.2.2) as traditional product of matrices $v$ and $f$.* □

From examples 3.2.1, 3.2.2, it follows that we cannot confine ourselves to traditional product of matrices and we need to define two products of matrices. To distinguish between these products we introduced a new notation. In order to keep this notation consistent with the existing one we assume that we have in mind ${}_*{}^*$-product when no clear notation is present.

Definition 3.2.3. *Let the nubmer of columns of the matrix $a$ equal the number of rows of the matrix $b$.* ${}_*{}^*$**-product** *of matrices $a$ and $b$ has form*

$$\begin{cases} a{}_*{}^*b = \left(a_k^i b_j^k\right) \\ (a{}_*{}^*b)_j^i = a_k^i b_j^k \end{cases} \tag{3.2.3}$$

*and can be expressed as product of a row of matrix $a$ over a column of matrix $b$.*[3.1]

□

Definition 3.2.4. *Let the nubmer of rows of the matrix $a$ equal the number of columns of the matrix $b$.* ${}^*{}_*$**-product** *of matrices $a$ and $b$ has form*

$$\begin{cases} a{}^*{}_*b = \left(a_i^k b_k^j\right) \\ (a{}^*{}_*b)_j^i = a_i^k b_k^j \end{cases} \tag{3.2.4}$$

*and can be expressed as product of a column of matrix $a$ over a row of matrix $b$.*[3.2]

□

We also consider following operations on the set of matrices.

Definition 3.2.5. *The transpose $a^T$ of the matrix $a$ exchanges rows and columns*

$$(a^T)_j^i = a_i^j \tag{3.2.5}$$

□

Definition 3.2.6. *The sum of matrices $a$ and $b$ is defined by the equality*

$$(a+b)_j^i = a_j^i + b_j^i$$

□

Let

$$a = \begin{pmatrix} a_1^1 & ... & a_n^1 \\ ... & ... & ... \\ a_1^m & ... & a_n^m \end{pmatrix}$$

be a matrix of $A$-numbers. We call matrix[3.3]

$$\mathcal{H}a = \begin{pmatrix} (a_1^1)^{-1} & ... & (a_1^m)^{-1} \\ ... & ... & ... \\ (a_n^1)^{-1} & ... & (a_n^m)^{-1} \end{pmatrix} \tag{3.2.6}$$

[3.1] We will use symbol ${}_*{}^*$- in following terminology and notation. We will read symbol ${}_*{}^*$ as $rc$ -product or product of row over column. To draw symbol of product of row over column, we put two symbols of product in the place of index which participate in sum. For instance, if product of $A$-numbers has form $a\circ b$, then ${}_*{}^*$-product of matrices $a$ and $b$ has form $a{}_\circ{}^\circ b$.

[3.2] We will use symbol ${}^*{}_*$- in following terminology and notation. We will read symbol ${}^*{}_*$ as $cr$ -product or product of column over row. To draw symbol of product of column over row, we put two symbols of product in the place of index which participate in sum. For instance, if product of $A$-numbers has form $a\circ b$, then ${}^*{}_*$-product of matrices $a$ and $b$ has form $a{}^\circ{}_\circ b$.

[3.3] The notation in the equality (3.2.7) means that we exchange rows and columns in Hadamard inverse. We can formally write right site of the equality (3.2.7) in the following form

$$(a_i^j)^{-1} = \frac{1}{a_j^i}$$

$$(\mathcal{H}a)^i_j = ((a^T)^i_j)^{-1} = (a^j_i)^{-1} \tag{3.2.7}$$

**Hadamard inverse of matrix** $a$ ([5]-page 4).

Set of $n \times n$ matrices is closed relative to ${}_*{}^*$-product and ${}^*{}_*$-product as well relative to sum.

THEOREM 3.2.7.

$$(a_*{}^*b)^T = a^T{}^*{}_*b^T \tag{3.2.8}$$

PROOF. The chain of equalities

$$((a_*{}^*b)^T)^j_i = (a_*{}^*b)^i_j = a^i_k b^k_j = (a^T)^k_i (b^T)^j_k = ((a^T)^*{}_*(b^T))^j_i \tag{3.2.9}$$

follows from (3.2.5), (3.2.3) and (3.2.4). The equality (3.2.8) follows from (3.2.9). □

Matrix $\delta = (\delta^i_j)$ is identity for both products.

DEFINITION 3.2.8. *The set* $\mathcal{A}$ *is a* **biring** *if we defined on* $\mathcal{A}$ *an unary operation, say transpose, and three binary operations, say* ${}_*{}^*$*-product,* ${}^*{}_*$*-product and sum, such that*

- ${}_*{}^*$*-product and sum define structure of ring on* $\mathcal{A}$
- ${}^*{}_*$*-product and sum define structure of ring on* $\mathcal{A}$
- *both products have common identity* $\delta$
- *products satisfy equation*

$$(a_*{}^*b)^T = a^T{}^*{}_*b^T$$

- *transpose of identity is identity*

$$\delta^T = \delta \tag{3.2.10}$$

- *double transpose is original element*

$$(a^T)^T = a \tag{3.2.11}$$

□

THEOREM 3.2.9.

$$(a^*{}_*b)^T = (a^T)_*{}^*(b^T) \tag{3.2.12}$$

PROOF. We can prove (3.2.12) in case of matrices the same way as we proved (3.2.8). However it is more important for us to show that (3.2.12) follows directly from (3.2.8).

Applying (3.2.11) to each term in left side of (3.2.12) we get

$$(a^*{}_*b)^T = ((a^T)^T{}^*{}_*(b^T)^T)^T \tag{3.2.13}$$

From (3.2.13) and (3.2.8) it follows that

$$(a^*{}_*b)^T = ((a^T{}_*{}^*b^T)^T)^T \tag{3.2.14}$$

(3.2.12) follows from (3.2.14) and (3.2.11). □

DEFINITION 3.2.10. *We introduce* ${}_*{}^*$**-power** *of* $\mathcal{A}$*-number* $a$ *using recursive definition*

$$a^{0_*{}^*} = \delta \tag{3.2.15}$$

$$a^{n_*{}^*} = a^{n-1_*{}^*}{}_*{}^*a \tag{3.2.16}$$

□

Theorem 3.2.11.

$$a^{1_*{}^*} = a \tag{3.2.17}$$

Proof. The equality

$$a^{1_*{}^*} = a^{0_*{}^*}{}_*{}^*a = \delta_*{}^*a = a \tag{3.2.18}$$

follows from equalities (3.2.15), (3.2.16). The equality (3.2.17) follows from the equality (3.2.18). □

Definition 3.2.12. *We introduce* ${}^*{}_*$**-power** *of* $\mathcal{A}$*-number* $a$ *using recursive definition*

$$a^{0^*{}_*} = \delta \tag{3.2.19}$$

$$a^{n^*{}_*} = a^{n-1^*{}_*}{}^*{}_*a \tag{3.2.20}$$

□

Theorem 3.2.13.

$$a^{1^*{}_*} = a \tag{3.2.21}$$

Proof. The equality

$$a^{1^*{}_*} = a^{0^*{}_*}{}^*{}_*a = \delta^*{}_*a = a \tag{3.2.22}$$

follows from equalities (3.2.19), (3.2.20). The equality (3.2.21) follows from the equality (3.2.22). □

Theorem 3.2.14.

$$(a^T)^{n_*{}^*} = (a^{n^*{}_*})^T \tag{3.2.23}$$

$$(a^T)^{n^*{}_*} = (a^{n_*{}^*})^T \tag{3.2.24}$$

Proof. We proceed by induction on $n$.

For $n = 0$ the statement immediately follows from equalities (3.2.15), (3.2.19), and (3.2.10).

Suppose the statement of theorem holds when $n = k - 1$

$$(a^T)^{n-1_*{}^*} = (a^{n-1^*{}_*})^T \tag{3.2.25}$$

It follows from (3.2.16) that

$$(a^T)^{k_*{}^*} = (a^T)^{k-1_*{}^*}{}_*{}^*a^T \tag{3.2.26}$$

It follows from (3.2.26) and (3.2.25) that

$$(a^T)^{k_*{}^*} = (a^{k-1^*{}_*})^T{}_*{}^*a^T \tag{3.2.27}$$

It follows from (3.2.27) and (3.2.12) that

$$(a^T)^{k_*{}^*} = (a^{k-1^*{}_*}{}^*{}_*a)^T \tag{3.2.28}$$

(3.2.23) follows from (3.2.26) and (3.2.20).

We can prove (3.2.24) by similar way. □

Definition 3.2.15. *$\mathcal{A}$-number* $a^{-1_*{}^*}$ *is* ${}_*{}^*$**-inverse element** *of $\mathcal{A}$-number $a$ if*

$$a_*{}^*a^{-1_*{}^*} = \delta \tag{3.2.29}$$

*$\mathcal{A}$-number* $a^{-1^*{}_*}$ *is* ${}^*{}_*$**-inverse element** *of $\mathcal{A}$-number $a$ if*

$$a^*{}_*a^{-1^*{}_*} = \delta \tag{3.2.30}$$

□

Theorem 3.2.16. *Let $\mathcal{A}$-number $a$ have ${}_*{}^*$-inverse $\mathcal{A}$-number. Then transpose $\mathcal{A}$-number $a^T$ has ${}^*{}_*$-inverse $\mathcal{A}$-number and these $\mathcal{A}$-numbers satisfy the equality*

$$(a^T)^{-1^*{}_*} = (a^{-1_*{}^*})^T \tag{3.2.31}$$

*Let $\mathcal{A}$-number $a$ have ${}^*{}_*$-inverse $\mathcal{A}$-number. Then transpose $\mathcal{A}$-number $a^T$ has ${}_*{}^*$-inverse $\mathcal{A}$-number and these $\mathcal{A}$-numbers satisfy the equality*

$$(a^T)^{-1_*{}^*} = (a^{-1^*{}_*})^T \tag{3.2.32}$$

Proof. If we get transpose of both side (3.2.29) and apply (3.2.10) we get

$$(a_*{}^*a^{-1_*{}^*})^T = \delta^T = \delta$$

Applying (3.2.8) we get

$$\delta = a^T{}^*{}_*(a^{-1_*{}^*})^T \tag{3.2.33}$$

(3.2.31) follows from comparison (3.2.30) and (3.2.33).

We can prove (3.2.32) similar way. □

Theorems 3.2.7, 3.2.9, 3.2.14, and 3.2.16 show that some kind of duality exists between ${}_*{}^*$-product and ${}^*{}_*$-product. We can combine these statements.

Theorem 3.2.17 (**duality principle for biring**). *Let $\mathfrak{A}$ be true statement about biring $\mathcal{A}$. If we exchange the same time*

- *$a \in \mathcal{A}$ and $a^T$*
- *${}_*{}^*$-product and ${}^*{}_*$-product*

*then we soon get true statement.*

Theorem 3.2.18 (**duality principle for biring of matrices**). *Let $\mathcal{A}$ be biring of matrices. Let $\mathfrak{A}$ be true statement about matrices. If we exchange the same time*

- *columns and rows of all matrices*
- *${}_*{}^*$-product and ${}^*{}_*$-product*

*then we soon get true statement.*

Proof. This is the immediate consequence of the theorem 3.2.17. □

Remark 3.2.19. *We execute operations in expression*

$$a_*{}^*b_*{}^*c \tag{3.2.34}$$

*from left to right. However we can execute product from right to left. In such case the expression (3.2.34) gets the form*

$$c^*{}_*b^*{}_*a$$

*We follow the rule that to write power and indices from right of expression. For instance, let original expression be like*

$$a^{-1_*{}^*}{}_*{}^*b^i$$

*Then expression which we read from right to left is like*

$$b^{i*}{}_{*}a^{-1^*}{}_{*}$$

*Suppose we established the order in which we write indices. Then we state that we read an expression from top to bottom reading first upper indices, then lower ones. We can read this expression from bottom to top. We extend this rule stating that we read symbols of operation in the same order as indices. For instance, if we read expression*

$$a^{i}{}_{*}{}^{*}b^{-1_*{}^*}=c^{i}$$

*from bottom to top, then we can write this expression as*

$$a_{i}{}^{*}{}_{*}b^{-1^*{}_*}=c_{i}$$

*According to the duality principle if we can prove one statement then we can prove other as well.* □

Theorem 3.2.20. *Let matrix $a$ have ${}_*{}^*$-inverse matrix. Then for any matrices $b$ and $c$ the equality*

$$b=c \tag{3.2.35}$$

*follows from the equality*

$$b_*{}^*a=c_*{}^*a \tag{3.2.36}$$

Proof. The equality (3.2.35) follows from the equality (3.2.36) if we multiply both parts of the equality (3.2.36) over $a^{-1_*{}^*}$. □

## 3.3. Quasideterminant

Theorem 3.3.1. *Let $n\times n$ matrix $a$ have ${}_*{}^*$-inverse matrix.*[3.4] *Then $k\times k$ minor matrix of ${}_*{}^*$-inverse matrix satisfies to the equality*

$$\left((a^{-1_*{}^*})^{I}_{J}\right)^{-1_*{}^*}=a^{J}_{I}-a^{J}_{[I]}{}_*{}^*\left(a^{[J]}_{[I]}\right)^{-1_*{}^*}{}_*{}^*a^{[J]}_{I} \tag{3.3.1}$$

Proof. Consider minor matrix $(a^{-1_*{}^*})_{J}$ obtained from the matrix $a^{-1_*{}^*}$ by selecting columns with an index from the set $J$. Definition (3.2.29) of ${}_*{}^*$-inverse matrix leads to the system of linear equations

$$a^{[J]}{}_*{}^*(a^{-1_*{}^*})_{J}=0 \tag{3.3.2}$$

$$a^{J}{}_*{}^*(a^{-1_*{}^*})_{J}=\delta \tag{3.3.3}$$

where $a^{J}$ is minor matrix obtained from the matrix $a$ by selecting rows with an index from the set $J$ and $a^{[J]}$ is minor matrix obtained from the matrix $a$ by deleting rows with an index from the set $J$. The equality

$$a^{[J]}_{[I]}{}_*{}^*(a^{-1_*{}^*})^{[I]}_{J}+a^{[J]}_{I}{}_*{}^*(a^{-1_*{}^*})^{I}_{J}=0 \tag{3.3.4}$$

follows from the equality (3.3.2) if we consider the minor matrix $a^{[J]}$ as a union of minor matrices $a^{[J]}_{I}$ and $a^{[J]}_{[I]}$. The equality

$$a^{J}_{[I]}{}_*{}^*(a^{-1_*{}^*})^{[I]}_{J}+a^{J}_{I}{}_*{}^*(a^{-1_*{}^*})^{I}_{J}=\delta \tag{3.3.5}$$

[3.4]This statement and its proof is based on statement 1.2.1 from [4] (page 8) for matrix over free division ring.

follows from the equality (3.3.3) if we consider the minor matrix $a^J$ as a union of minor matrices $a_I^J$ and $a_{[I]}^J$. We multiply the equality (3.3.4) by $\left(a_{[I]}^{[J]}\right)^{-1_*{}^*}$

$$(a^{-1_*{}^*})_J^{[I]} + \left(a_{[I]}^{[J]}\right)^{-1_*{}^*} {}_*{}^* a_I^{[J]} {}_*{}^* (a^{-1_*{}^*})_J^I = 0 \tag{3.3.6}$$

The equality

$$-a_{[I]}^J {}_*{}^* \left(a_{[I]}^{[J]}\right)^{-1_*{}^*} {}_*{}^* a_I^{[J]} {}_*{}^* (a^{-1_*{}^*})_J^I + a_I^J {}_*{}^* (a^{-1_*{}^*})_J^I = \delta \tag{3.3.7}$$

follows from equalities (3.3.6), (3.3.5). The equality (3.3.1) follows from (3.3.7) when we both sides of the equality (3.3.7) multiply by $\left((a^{-1_*{}^*})_J^I\right)^{-1_*{}^*}$. □

THEOREM 3.3.2. *Let $n \times n$ matrix $a$ have a ${}_*{}^*$-inverse matrix. Then entries of the ${}_*{}^*$-inverse matrix satisfy to the equality*[3.5]

$$(a^{-1_*{}^*})_j^i = \left(a_i^j - a_{[i]}^j {}_*{}^* \left(a_{[i]}^{[j]}\right)^{-1_*{}^*} {}_*{}^* a_i^{[j]}\right)^{-1} \tag{3.3.8}$$

$$(\mathcal{H}a^{-1_*{}^*})_i^j = a_i^j - a_{[i]}^j {}_*{}^* \left(a_{[i]}^{[j]}\right)^{-1_*{}^*} {}_*{}^* a_i^{[j]} \tag{3.3.9}$$

PROOF. The equality (3.3.8) follows from the equality (3.3.1) when we assume $I = i$, $J = j$. The equality (3.3.9) follows from equalities (3.2.7), (3.3.8). □

According to [4], page 3 we do not have an appropriate definition of a determinant for a division algebra.[3.6] However, we can define a quasideterminant which finally gives a similar picture. In definition below we follow definition [4]-1.2.2.

DEFINITION 3.3.3. $\binom{j}{i}$**-${}_*{}^*$-quasideterminant** *of $n \times n$ matrix $a$ is formal expression*

$$\det({}_*{}^*)_i^j\, a = (\mathcal{H}a^{-1_*{}^*})_i^j \tag{3.3.10}$$

*According to the remark 3.1.2 we consider $\binom{j}{i}$-${}_*{}^*$-quasideterminant as an entry of the matrix*

$$\det({}_*{}^*)\, a = \mathcal{H}a^{-1_*{}^*}$$

*which is called* **${}_*{}^*$-quasideterminant**. □

THEOREM 3.3.4. *Expression for ${}_*{}^*$-inverse matrix has form*

$$a^{-1_*{}^*} = \mathcal{H}\det({}_*{}^*)a \tag{3.3.11}$$

PROOF. (3.3.11) follows from (3.3.10). □

[3.5] The notation in the equality (3.3.8) means that we exchange rows and columns in inverse of matrix.

[3.6] Professor Kyrchei uses double determinant (see the definition [10]-5.28) to solve system of linear equations in quaternion algebra. Professor Kyrchei also solved eigenvalues problem in quaternion algebra using double determinant (see the section [16]-2.5). I confine myself by consideration of quasideterminant, because I am interested in a wider set of algebras.

THEOREM 3.3.5. *Expression for* $\binom{j}{i}$*-*${}_*{}^*$*-quasideterminant can be evaluated by either form*[3.7]

$$\det({}_*{}^*)^j_i\ a = a^j_i - a^j_{[i]}\ {}_*{}^*\left(a^{[j]}_{[i]}\right)^{-1_*{}^*}\ {}_*{}^* a^{[j]}_i \tag{3.3.12}$$

$$\det({}_*{}^*)^j_i\ a = a^j_i - a^j_{[i]}\ {}_*{}^*\mathcal{H}\det({}_*{}^*)a^{[j]}_{[i]}\ {}_*{}^* a^{[j]}_i \tag{3.3.13}$$

PROOF. Statement follows from (3.3.9) and (3.3.10). □

THEOREM 3.3.6. *Consider matrix*

$$\begin{pmatrix} a^1_1 & a^1_2 \\ a^2_1 & a^2_2 \end{pmatrix}$$

*Then*

$$\det({}_*{}^*)a = \begin{pmatrix} a^1_1 - a^1_2(a^2_2)^{-1}a^2_1 & a^1_2 - a^1_1(a^2_1)^{-1}a^2_2 \\ a^2_1 - a^2_2(a^1_2)^{-1}a^1_1 & a^2_2 - a^2_1(a^1_1)^{-1}a^1_2 \end{pmatrix} \tag{3.3.14}$$

$$\det({}^*{}_*)a = \begin{pmatrix} a^1_1 - a^2_1(a^2_2)^{-1}a^1_2 & a^1_2 - a^2_2(a^2_1)^{-1}a^1_1 \\ a^2_1 - a^1_1(a^1_2)^{-1}a^2_2 & a^2_2 - a^1_2(a^1_1)^{-1}a^2_1 \end{pmatrix} \tag{3.3.15}$$

$$a^{-1_*{}^*} = \begin{pmatrix} (a^1_1 - a^2_1(a^2_2)^{-1}a^1_2)^{-1} & (a^2_1 - a^1_1(a^1_2)^{-1}a^2_2)^{-1} \\ (a^1_2 - a^2_2(a^2_1)^{-1}a^1_1)^{-1} & (a^2_2 - a^1_2(a^1_1)^{-1}a^2_1)^{-1} \end{pmatrix} \tag{3.3.16}$$

PROOF. According to the equality (3.3.12)

$$\det({}_*{}^*)^1_1\ a = a^1_1 - a^1_2(a^2_2)^{-1}a^2_1 \tag{3.3.17}$$

$$\det({}_*{}^*)^2_1\ a = a^1_2 - a^1_1(a^2_1)^{-1}a^2_2 \tag{3.3.18}$$

$$\det({}_*{}^*)^1_2\ a = a^2_1 - a^2_2(a^1_2)^{-1}a^1_1 \tag{3.3.19}$$

$$\det({}_*{}^*)^2_2\ a = a^2_2 - a^2_1(a^1_1)^{-1}a^1_2 \tag{3.3.20}$$

(3.3.17), (3.3.18), (3.3.19), (3.3.20). □

THEOREM 3.3.7.

$$\det({}_*{}^*)^i_j\ a^T = \det({}^*{}_*)^j_i\ a \tag{3.3.21}$$

[3.7] We can provide similar proof for $\binom{j}{i}$**-**${}^*{}_*$**-quasideterminant**. However we can write corresponding statement using the duality principle. Thus, if we read equality (3.3.12) from right to left, we get equality

$$\det({}^*{}_*)a^j_i = a^j_i - a^{[j]}_i\ {}^*{}_*\left(a^{[j]}_{[i]}\right)^{-1^*{}_*}\ {}^*{}_* a^j_{[i]}$$

$$\det({}^*{}_*)a^j_i = a^j_i - a^{[j]}_i\ {}^*{}_*\mathcal{H}\det({}^*{}_*)^{[j]}A_{[i]}\ {}^*{}_* a^j_{[i]}$$

Proof. According to (3.3.10) and (3.2.6)

$$\det({}_*{}^*)^i_j\, a^T = (((a^T)^{-1_*{}^*})^{\cdot}{}_i{}^{-}\,{}^j{}_{\cdot -})^{-1}$$

Using theorem 3.2.16 we get

$$\det({}_*{}^*)^i_j\, a^T = (((a^{-1^*{}_*})^T)^{\cdot}{}_i{}^{-}\,{}^j{}_{\cdot -})^{-1}$$

Using (3.2.5) we get

$$\det({}_*{}^*)^i_j\, a^T = ((a^{-1^*{}_*})_{\cdot -}{}^i\,{}^{\cdot -}{}_j)^{-1} \tag{3.3.22}$$

Using (3.3.22), (3.2.6), (3.3.10) we get (3.3.21). □

The theorem 3.3.7 extends the duality principle stated in the theorem 3.2.18 to statements on quasideterminants and tells us that the same expression is ${}_*{}^*$-quasideterminant of matrix $a$ and ${}^*{}_*$-quasideterminant of matrix $a^T$. Using this theorem, we can write any statement for ${}^*{}_*$-matrix on the basis of similar statement for ${}_*{}^*$-matrix.

Theorem 3.3.8 (duality principle). *Let $\mathfrak{A}$ be true statement about matrix biring. If we exchange the same time*

- *row and column*
- *${}_*{}^*$-quasideterminant and ${}^*{}_*$-quasideterminant*

*then we soon get true statement.*

Theorem 3.3.9.

$$(ma)^{-1_*{}^*} = a^{-1_*{}^*} m^{-1} \tag{3.3.23}$$

$$(am)^{-1_*{}^*} = m^{-1} a^{-1_*{}^*} \tag{3.3.24}$$

Proof. To prove the equality (3.3.23) we proceed by induction on size of the matrix.

Since

$$(ma)^{-1_*{}^*} = \left((ma)^{-1}\right) = \left(a^{-1}m^{-1}\right) = \left(a^{-1}\right) m^{-1} = a^{-1_*{}^*} m^{-1}$$

the statement is evident for $1 \times 1$ matrix.

Let the statement holds for $(n-1)\times(n-1)$ matrix. Then from the equality (3.3.1) it follows that

$$\begin{aligned}(((ma)^{-1_*{}^*})^I_J)^{-1_*{}^*} &= (ma)^J_I - (ma)^J_{[I]}{}_*{}^*\left((ma)^{[J]}_{[I]}\right)^{-1_*{}^*}{}_*{}^*(ma)^{[J]}_I\\ &= ma^J_I - m\, a^J_{[I]}{}_*{}^*\left(a^{[J]}_{[I]}\right)^{-1_*{}^*} m^{-1}{}_*{}^* m\, a^{[J]}_I\\ &= m\, a^J_I - m\, a^J_{[I]}{}_*{}^*\left(a^{[J]}_{[I]}\right)^{-1_*{}^*}{}_*{}^* a^{[J]}_I\end{aligned}$$

$$(((ma)^{-1_*{}^*})^I_J)^{-1_*{}^*} = m\,(a^{-1_*{}^*})^I_J \tag{3.3.25}$$

The equality (3.3.23) follows from the equality (3.3.25). In the same manner we prove the equality (3.3.24). □

CHAPTER 4

# Polynomial

## 4.1. Polynomial over Associative $D$-Algebra

Let $D$ be commutative ring and $A$ be associative $D$-algebra with unit.

THEOREM 4.1.1. *Let $p_k(x)$ be* **monomial of power** *$k$ over $D$-algebra $A$. Then*

4.1.1.1: *Monomial of power 0 has form* $p_0(x) = a_0$, $a_0 \in A$.

4.1.1.2: *If $k > 0$, then*

$$p_k(x) = p_{k-1}(x)xa_k$$

*where $a_k \in A$.*

PROOF. We prove the theorem by induction over power $n$ of monomial.

Let $n = 0$. We get the statement 4.1.1.1 since monomial $p_0(x)$ is constant.

Let $n = k$. Last factor of monomial $p_k(x)$ is either $a_k \in A$, or has form $x^l$, $l \geq 1$. In the later case we assume $a_k = 1$. Factor preceding $a_k$ has form $x^l$, $l \geq 1$. We can represent this factor as $x^{l-1}x$. Therefore, we proved the statement. □

DEFINITION 4.1.2. *We denote* $A_k[x]$ *Abelian group generated by the set of monomials of power $k$. Element $p_k(x)$ of Abelian group $A_k[x]$ is called* **homogeneous polynomial of power** *$k$.* □

THEOREM 4.1.3. *The map*

$$f : A^{k+1} \to A_k[x]$$

*defined by the equality*

$$f(a_0, a_1, ..., a_k) = (a_0, a_1, ..., a_k) \circ x^k$$

*is polylinear map.*

PROOF. Let $d \in D$. From equalities

$$\begin{aligned}
f(a_0, ..., da_i, ..., a_k) &= a_0x...(da_i)...xa_k = d(a_0x...a_i...xa_k)\\
&= df(a_0, ..., a_i, ..., a_k)\\
f(a_0, ..., a_i + b_i, ..., a_k) &= a_0x...(a_i + b_i)...xa_k\\
&= (a_0x...a_i...xa_k) + (a_0x...b_i...xa_k)\\
&= f(a_0, ..., a_i, ..., a_k) + f(a_0, ..., b_i, ..., a_k)
\end{aligned}$$

it follows that the map $f$ is linear with respect to $a_i$. Therefore, the map $f$ is polylinear map. □

THEOREM 4.1.4. *There exists linear map*

$$(a_0 \otimes a_1 \otimes ... \otimes a_k) \circ x^k = (a_0, a_1, ..., a_k) \circ x^k$$

PROOF. The theorem follows from theorems [9]-3.6.4, 4.1.3. □

Corollary 4.1.5. *We can present homogeneous polynomial $p(x)$ in the following form*

$$p(x) = a \circ x^k \quad a \in A^{(k+1)\otimes}$$

□

Definition 4.1.6. *We denote*

$$A[x] = \bigoplus_{n=0}^{\infty} A_n[x]$$

*direct sum*[4.1] *of $A$-modules $A_n[x]$. An element $p(x)$ of $A$-module $A[x]$ is called* **polynomial** *over $D$-algebra $A$.* □

Therefore, we can present polynomial of power $n$ in the following form

$$p(x) = a_0 + a_1 \circ x + ... + a_n \circ x^n \quad a_i \in A^{i+1\otimes} \quad i = 0, ..., n$$

Definition 4.1.7. *Bilinear map*

$$* : A^{n\otimes} \times A^{m\otimes} \to A^{n+m-1\otimes}$$

*is defined by the equality*

$$(a_1 \otimes ... \otimes a_n) * (b_1 \otimes ... \otimes b_n) = a_1 \otimes ... \otimes a_{n-1} \otimes a_n b_1 \otimes b_2 \otimes ... \otimes b_n$$

□

Theorem 4.1.8. *For any tensors $a \in A^{n\otimes}$, $b \in A^{m\otimes}$, $c \in A^{k\otimes}$*

$$(a * b) * c = a * (b * c) \tag{4.1.1}$$

Proof. According to the definition 4.1.7, it is enough to prove the theorem for tensors $a = a_0 \otimes ... \otimes a_n$, $b = b_0 \otimes ... \otimes b_m$, $c = c_0 \otimes ... \otimes c_k$. The equality (4.1.1) follows from the equality

$$\begin{aligned}
&((a_0 \otimes ... \otimes a_n) * (b_0 \otimes b_m)) * (c_0 \otimes ... \otimes c_k)\\
=&(a_0 \otimes ... \otimes a_n b_0 \otimes ... \otimes b_m) * (c_0 \otimes ... \otimes c_k)\\
=&a_0 \otimes ... \otimes a_n b_0 \otimes ... \otimes b_m c_0 \otimes ... \otimes c_k\\
=&(a_0 \otimes ... \otimes a_n) * (b_0 \otimes ... \otimes b_m c_0 \otimes ... \otimes c_k)\\
=&(a_0 \otimes ... \otimes a_n) * ((b_0 \otimes ... \otimes b_m) * (c_0 \otimes ... \otimes c_k)
\end{aligned}$$

□

Theorem 4.1.9. *For any tensors $a \in A^{n\otimes}$, $b \in A^{m\otimes}$ , Product of homogeneous polynomials $a \circ x^n$, $b \circ x^m$ is defined by the equality*

$$(a * b) \circ x^{n+m} = (a \circ x^n)(b \circ x^m) \tag{4.1.2}$$

Proof. According to the definition 4.1.7, it is enough to prove the theorem for tensors $a = a_0 \otimes ... \otimes a_n$, $b = b_0 \otimes ... \otimes b_m$. We proceed by induction on $m$.

From the equality

$$\begin{aligned}
((a_0 \otimes a_1 \otimes ... \otimes a_n) * b_0) \circ x^n &= (a_0 \otimes a_1 \otimes ... \otimes a_n b_0) \circ x^n\\
&= a_0 x a_1 x ... x a_n b_0 = (a_0 x a_1 x ... x a_n) b_0\\
&= ((a_0 \otimes a_1 \otimes ... \otimes a_n) \circ x^n)(b_0 \circ x^0)
\end{aligned}$$

[4.1] See the definition 2.2.5.

it follows that the theorem is true for $m=0$.

From the equality

$$\begin{aligned}
&((a_0\otimes a_1\otimes ...\otimes a_n)*(b_0\otimes b_1))\circ x^{n+1}\\
&=(a_0\otimes a_1\otimes ...\otimes a_nb_0\otimes b_1)\circ x^{n+1}\\
&=a_0xa_1x...xa_nb_0xb_1=(a_0xa_1x...xa_n)(b_0xb_1)\\
&=((a_0\otimes a_1\otimes ...\otimes a_n)\circ x^n)((b_0\otimes b_1)\circ x^1)
\end{aligned} \tag{4.1.3}$$

it follows that the theorem is true for $m=1$.

Let the theorem be true for $m=k$

$$(a*(b_0\otimes ...\otimes b_k))\circ x^{n+k}=(a\circ x^n)((b_0\otimes ...\otimes b_k)\circ x^k) \tag{4.1.4}$$

The equality

$$\begin{aligned}
&(a*(b_0\otimes ...\otimes b_k\otimes b_{k+1}))\circ x^{n+k+1}\\
&=(a*((b_0\otimes ...\otimes b_k)*(1\otimes b_{k+1}))\circ x^{n+k+1}\\
&=((a*(b_0\otimes ...\otimes b_k))*(1\otimes b_{k+1}))\circ x^{n+k+1}\\
&=((a*(b_0\otimes ...\otimes b_k))\circ x^{n+k})((1\otimes b_{k+1})\circ x)\\
&=((a\circ x^n)((b_0\otimes ...\otimes b_k)\circ x^k))((1\otimes b_{k+1})\circ x)\\
&=(a\circ x^n)(((b_0\otimes ...\otimes b_k)\circ x^k)((1\otimes b_{k+1})\circ x))\\
&=(a\circ x^n)(((b_0\otimes ...\otimes b_k)*(1\otimes b_{k+1}))\circ x^{k+1})\\
&=(a\circ x^n)((b_0\otimes ...\otimes b_k\otimes b_{k+1})\circ x^{k+1})
\end{aligned} \tag{4.1.5}$$

follows from the equality (4.1.1). (4.1.3), (4.1.4). From the equality (4.1.5), it follows that theorem is true for $m=k+1$. □

## 4.2. Derivative of Polynomial

CONVENTION 4.2.1. *For $n\ge 0$, let $SO(k,n)$ be set of permutations*

$$\sigma=\begin{pmatrix}y_1&...&y_k&x_{k+1}&...&x_n\\ \sigma(y_1)&...&\sigma(y_k)&\sigma(x_{k+1})&...&\sigma(x_n)\end{pmatrix}$$

*such that each permutation $\sigma$ preserves the order of variables $x_i$: since $i<j$, then in the tuple*

$$(\sigma(y_1),...,\sigma(y_k),\sigma(x_{k+1}),...,\sigma(x_n))$$

*$x_i$ precedes $x_j$.* □

LEMMA 4.2.2. *We can enumerate the set of permutations $SO(1,n)$ by index $i$, $1\le i\le n$, such way that*

4.2.2.1: *$\sigma_1(y)=y$.*

4.2.2.2: *Since $i>1$, then $\sigma_i(x_i)=y$.*

PROOF. Since the order of variables $x_2$, ..., $x_n$ in permutation $\sigma\in SO(1,n)$ does not depend on permutation, permutations $\sigma\in SO(1,n)$ are different by position which variable $y$ has. Accordingly, we can enumerate the set of permutations $SO(1,n)$ by index whose value corresponds to the number of position of variable $y$. □

The lemma 4.2.2 has simple interpretation. Let $n-1$ white balls and 1 black ball be in narrow box. The black ball is the most left ball; white balls are numbered from 2 to $n$ in the order as we put them into the box. The essence of the permutation $\sigma_k$ is that we take out the black ball from the box and then we put it into cell with number $k$. At the same time, white balls with the number not exceeding $k$ shift to the left.

LEMMA 4.2.3. *For $n>0$, let*

$$SO^+(1,n)=\{\sigma:\sigma=(\tau(y),\tau(x_2),...,\tau(x_n),x_{n+1}),\tau\in SO(1,n)\} \tag{4.2.1}$$

*Then*

$$SO(1,n+1)=SO^+(1,n)\cup\{(x_2,...,x_{n+1},y)\} \tag{4.2.2}$$

PROOF. Let $\sigma\in SO^+(1,n)$. According to the definition (4.2.1), there exists permutation $\tau\in SO(1,n)$ such that

$$(\sigma(y),\sigma(x_1),...,\sigma(x_{n+1}))=(\tau(y),\tau(x_2),...,\tau(x_n),x_{n+1})$$

According to the convention 4.2.1, the statement $i<j<n+1$ implyes that in the tuple

$$(\sigma(y),\sigma(x_1),...,\sigma(x_{n+1}))=(\tau(y),\tau(x_2),...,\tau(x_n),x_{n+1})$$

the variable $x_j$ is located between variables $x_i$ and $x_{n+1}$. According to the convention 4.2.1, $\sigma\in SO(1,n+1)$. Therefore

$$SO^+(1,n)\subseteq SO(1,n+1) \tag{4.2.3}$$

According to the lemma 4.2.2, the set $SO^+(1,n)$ has $n$ permutations.

Let $\sigma=(x_2,...,x_{n+1},y)$. According to the convention 4.2.1,

$$(x_2,...,x_{n+1},y)\in SO(1,n+1) \tag{4.2.4}$$

According to the definition (4.2.1), $\sigma\notin SO^+(1,n)$.

Therefore, we have listed $n+1$ elements of the set $SO(1,n+1)$. According to the lemma 4.2.2, the statement (4.2.2) follows from statements (4.2.3), (4.2.4). □

THEOREM 4.2.4. *For any monomial*

$$p_n(x)=(a_0\otimes...\otimes a_n)\circ x^n$$

*derivative has form*

$$\begin{matrix}\displaystyle\frac{dp_n(x)}{dx}\circ dx=\sum_{\sigma\in SO(1,n)}(a_0\otimes...\otimes a_n)\circ(\sigma(dx),\sigma(x_2),...,\sigma(x_n))\\ x_2=...=x_n=x\end{matrix} \tag{4.2.5}$$

PROOF. For $n=1$, the map

$$p_1(x)=(a_0\otimes a_1)\circ x=a_0xa_1$$

is linear map. According to the theorem A.1.5 and convention 4.2.1

$$\frac{dp_1(x)}{dx}\circ dx=a_0dxa_1=(a_0\otimes a_1)\circ dx=\sum_{\sigma\in SO(1,1)}(a_0\otimes a_1)\circ\sigma(dx) \tag{4.2.6}$$

Let the statement be true for $n-1$

$$\frac{dp_{n-1}(x)}{dx}\circ dx=\sum_{\sigma\in SO(1,n-1)}(a_0\otimes...\otimes a_{n-1})\circ\sigma(dx,x_2,...,x_{n-1}) \tag{4.2.7}$$
$$x_2=...=x_{n-1}=x$$

Since

$$p_n(x)=p_{n-1}(x)xa_n \tag{4.2.8}$$

then according to the theorem A.1.10 and the definition (4.2.8)

$$\frac{dp_n(x)}{dx}\circ dx=\left(\frac{dp_{n-1}(x)}{dx}\circ dx\right)\ xa_n+p_{n-1}(x)\ \left(\frac{dxa_n}{dx}\circ dx\right) \tag{4.2.9}$$

The equality

$$\begin{aligned}&\frac{dp_n(x)}{dx}\circ dx\\ &=\sum_{\sigma\in SO(1,n-1)}(a_0\otimes...\otimes a_{n-1})\circ(\sigma(dx),\sigma(x_2),...,\sigma(x_{n-1}))xa_n\\ &+((a_0\otimes...\otimes a_{n-1})\circ(x_2,...,x_n))dxa_n\\ &x_2=...=x_n=x\end{aligned} \tag{4.2.10}$$

follows from (A.1.1), (A.1.2), (4.2.7), (4.2.9). The equality

$$\begin{aligned}&\frac{dp_n(x)}{dx}\circ dx\\ &=\sum_{\sigma\in SO(1,n-1)}(a_0\otimes...\otimes a_n)\circ(\sigma(dx),\sigma(x_2),...,\sigma(x_{n-1}),x_n))\\ &+(a_0\otimes...\otimes a_n)\circ(x_2,...,x_n,dx)\\ &x_2=...=x_n=x\end{aligned} \tag{4.2.11}$$

follows from (4.2.10) and multiplication rule of monomials (the definition 4.1.7 and the theorem 4.1.9). According to the lemma 4.2.3, the equality (4.2.5) follows from the equality (4.2.11). □

LEMMA 4.2.5. *Let* $1<k<n$. *For any pair of permutations* $\mu\in SO(k,n)$, $\nu\in SO(1,n-k)$, *there exists unique permutation* $\sigma\in SO(k+1,n)$ *such that*

$$\sigma=(\mu(y_1),...,\mu(y_k),\nu(\mu(y_{k+1})),\nu(\mu(x_{k+2})),...,\nu(\mu(x_n)))$$

PROOF. According to the statement of the lemma, we first perform the permutation $\mu$

$$\begin{pmatrix}\mu(y_1) & ... & \mu(y_k) & \mu(x_{k+1}=y_{k+1}) & \mu(x_{k+2}) & ... & \mu(x_n)\end{pmatrix}$$

and than perform the permutation $\nu$

$$\begin{pmatrix}y_1 & ... & y_k & \nu(y_{k+1}) & \nu(x_{k+2}) & ... & \nu(x_n)\end{pmatrix}$$

which we apply only to variables $y_{k+1}$, $x_{k+2}$, ..., $x_n$. According to the convention 4.2.1, the permutation $\mu$ preserves the order of variables $x_{k+1}=y_{k+1}$, $x_{k+2}$, ..., $x_n$. Therefore, the permutation $\mu$ preserves the order of variables $x_{k+2}$, ..., $x_n$. According to the convention 4.2.1, the permutation $\nu$ preserves the order of variables $x_{k+2}$, ..., $x_n$. Therefore, the permutation $\sigma$ preserves the order of variables $x_{k+2}$, ..., $x_n$. According to the convention 4.2.1, $\sigma\in SO(k+1,n)$ . □

Lemma 4.2.6. *Let* $1 < k < n$. *For any permutation* $\sigma \in SO(k+1,n)$, *there exists unique pair of permutations* $\mu \in SO(k,n)$, $\nu \in SO(1,n-k)$ *such that*

$$\sigma = (\mu(y_1), ..., \mu(y_k), \nu(\mu(y_{k+1})), \nu(\mu(x_{k+2})), ..., \nu(\mu(x_n)))$$

Proof. If in the tuple

$$(\sigma(y_1), ..., \sigma(y_{k+1}), \sigma(x_{k+2}), ..., \sigma(x_n))$$

$y_{k+1}$ precedes $x_{k+2}$, then we set

$$\nu = \begin{pmatrix} y_{k+1} & x_{k+2} & ... & x_n \\ y_{k+1} & x_{k+2} & ... & x_n \end{pmatrix}$$

$$\mu = \begin{pmatrix} y_1 & ... & y_k & x_{k+1} = y_{k+1} & x_{k+2} & ... & x_n \\ \sigma(y_1) & ... & \sigma(y_k) & \sigma(y_{k+1}) & \sigma(x_{k+2}) & ... & \sigma(x_n) \end{pmatrix}$$

Since permutation $\nu$ preserves the order of variables $x_i$, then, according to the convention 4.2.1, $\nu \in SO(1, n-k)$. In the tuple

$$(\mu(y_1), ..., \mu(y_{k+1}), \mu(x_{k+2}), ..., \mu(x_n))$$

$x_{k+1} = y_{k+1}$ precedes $x_{k+2}$, According to the convention 4.2.1, $x_{k+1} = y_{k+1}$ precedes $x_{k+2}$, ..., $x_n$. Therefore, $\mu \in SO(k,n)$.

Let in the tuple

$$(\sigma(y_1), ..., \sigma(y_{k+1}), \sigma(x_{k+2}), ..., \sigma(x_n))$$

$x_{k+2}$, ..., $x_j$ precede $y_{k+1}$. Then we set

$$\nu = \begin{pmatrix} x_{k+2} & x_{k+3} & ... & x_j & y_{k+1} & ... & x_n \\ y_{k+1} & x_{k+2} & ... & x_{j-1} & x_j & ... & x_n \end{pmatrix}$$

$$\mu = \begin{pmatrix} y_1 & ... & y_k & x_{k+1} = y_{k+1} & x_{k+2} & ... & x_n \\ \sigma(y_1) & ... & \sigma(y_k) & \sigma(y_{k+1}) & \sigma(x_{k+2}) & ... & \sigma(x_n) \end{pmatrix}$$

Since permutations $\mu$, $\nu$ preserve the order of variables $x_i$, then $\mu \in SO(k,n)$, $\nu \in SO(1, n-k)$. $\square$

Theorem 4.2.7. *For any monomial*

$$p_n(x) = (a_0 \otimes ... \otimes a_n) \circ x^n$$

*derivative of order* $k$ *has form*

$$\begin{aligned} &\frac{d^k p_n(x)}{dx^k} \circ (dx^1; ...; dx^k) \\ &= \sum_{\sigma \in SO(k,n)} (a_0 \otimes ... \otimes a_n, \sigma) \circ (dx_1; ...; dx_k; x_{k+1}; ...; x_n) \\ &x_{k+1} = ... = x_n = x \end{aligned} \tag{4.2.12}$$

Proof. For $k = 1$, the statement of the theorem is the statement of the theorem 4.2.4. Let the theorem be true for $k-1$. Then

$$(4.2.13)\quad \begin{aligned}&\frac{d^{k-1}p_n(x)}{dx^{k-1}}\circ(dx_1;...;dx_{k-1})\\ =&\sum_{\mu\in SO(k-1,n)}(a_0\otimes...\otimes a_n,\mu)\circ(dx_1;...;dx_{k-1};x_k;...;x_n)\\ &x_k=...=x_n=x\end{aligned}$$

According to the definition [15]-4.1.4, the equality

$$(4.2.14)\quad \begin{aligned}&\frac{d^kp_n(x)}{dx^k}\circ(dx_1;...;dx_{k-1};dx_k)\\ =&\frac{d}{dx}\left(\frac{d^{k-1}p_n(x)}{dx^{k-1}}\circ(dx_1;...;dx_{k-1})\right)\circ dx^k\\ =&\frac{d}{dx}\left(\sum_{\mu\in SO(k-1,n)}(a_0\otimes...\otimes a_n,\mu)\circ(dx_1;...;dx_{k-1};x_k;...;x_n)\right)\circ dx^k\\ &x_k=...=x_n=x\end{aligned}$$

follows from the equality (4.2.13). According to the theorem 4.2.4 the equality

$$(4.2.15)\quad \begin{aligned}&\frac{d^kp_n(x)}{dx^k}\circ(dx_1;...;dx_{k-1};dx_k)\\ =&\sum_{\substack{\mu\in SO(k-1,n)\\ \nu\in SO(1,n-k+1)}}\begin{aligned}&(a_0\otimes...\otimes a_n)\circ(\mu(dx_1);...;\mu(dx_{k-1});\nu(\mu(dx_k));\\ &\nu(\mu(x_{k+1}));...;\nu(\mu(x_n))\end{aligned}\\ &x_{k+1}=...=x_n=x\end{aligned}$$

follows from the equality (4.2.14). According to the lemmas 4.2.5, 4.2.6, the equality

$$(4.2.16)\quad \begin{aligned}&\frac{d^kp_n(x)}{dx^k}\circ(dx_1;...;dx_{k-1};dx_k)\\ =&\sum_{\sigma\in SO(k,n)}(a_0\otimes...\otimes a_n,\sigma)\circ(dx_1;...;dx_{k-1};dx_k;x_{k+1};...;x_n)\\ &x_{k+1}=...=x_n=x\end{aligned}$$

follows from the equality (4.2.15). Therefore, the theorem is true for $k$. □

Theorem 4.2.8. *The derivative* $\frac{d^mp_n(x)}{dx^m}\circ(h_1;...;h_m)$ *is symmetric polynomial with respect to variables* $h_1$, ..., $h_m$.

Proof. without loss of generality, we can assume that $p_n$ is monomial. We will prove the theorem by induction over $m$.

The theorem is evident for $m = 0$ and $m = 1$.

Lemma 4.2.9. *The derivative* $\frac{d^2p_n(x)}{dx^2}$ *is symmetric bilinear map.*

Proof. We consider monomial $p_n(x)$ in the following form

$$p_n(x)=a_0xa_1...a_{n_1}xa_n \tag{4.2.17}$$

According to the theorem 4.2.4 and the lemma 4.2.2, the equality

$$\begin{aligned}\frac{dp_n(x)}{dx}\circ c_1=&a_0c_1a_1...a_{n-1}xa_n+...+a_0xa_1...a_{i-1}c_1a_i...a_{n-1}xa_n+...\\&+a_0xa_1...a_{n-1}c_1a_n\end{aligned} \tag{4.2.18}$$

follows from the equality (4.2.17). We denote $p_{n1}(c_1,x)$, ..., $p_{ni}(c_1,x)$, ..., $p_{nn}(c_1,x)$ terms on the right side of the equality (4.2.18)

$$\frac{dp_n(x)}{dx}\circ c_1=p_{n1}(c_1,x)+...+p_{ni}(c_1,x)+...+p_{nn}(c_1,x) \tag{4.2.19}$$

$$p_{ni}(c_1,x)=a_0xa_1...a_{i-1}c_1a_i...a_{n-1}xa_n \tag{4.2.20}$$

According to the definition [15]-4.1.1 and the theorem A.1.4, the equality

$$\begin{aligned}\frac{d^2p_n(x)}{dx^2}\circ(c_1,c_2)=&\frac{d}{dx}\left(\frac{dp_n(x)}{dx}\circ c_1\right)\circ c_2\\=&\frac{d}{dx}(p_{n1}(c_1,x)+...+p_{ni}(c_1,x)+...+p_{nn}(c_1,x))\circ c_2\\=&\frac{dp_{n1}(c_1,x)}{dx}\circ c_2+...+\frac{dp_{ni}(c_1,x)}{dx}\circ c_2+...+\frac{dp_{nn}(c_1,x)}{dx}\circ c_2\end{aligned} \tag{4.2.21}$$

follows from the equality (4.2.19). In the right side of the equality (4.2.21), consider term with number $i$, $1\le i\le n$. According to the theorem 4.2.4 and the lemma 4.2.2, the equality[4.2]

$$\begin{aligned}\frac{dp_{ni}(c_1,x)}{dx}\circ c_2=&a_0c_2a_1...a_{i-1}c_1a_i...a_{n-1}xa_n+...\\&+a_0xa_1...a_{i-1}c_1a_i......a_{j-1}c_2a_j...a_{n-1}xa_n+...\\&+a_0xa_1...a_{i-1}c_1a_i...a_{n-1}c_2a_n\end{aligned} \tag{4.2.22}$$

follows from the equality (4.2.20). We denote $p_{ni1}(c_1,c_2,x)$, ..., $p_{ni\,i-1}(c_1,c_2,x)$, $p_{ni\,i+1}(c_1,c_2,x)$, ..., $p_{nij}(c_1,c_2,x)$, ..., $p_{nn}(c_1,c_2,x)$ terms on the right side of the equality (4.2.22)

$$\frac{dp_{ni}(c_1,x)}{dx}\circ c_2=p_{ni1}(c_1,c_2,x)+...+p_{nij}(c_1,c_2,x)+...+p_{nin}(c_1,c_2,x) \tag{4.2.23}$$

$$p_{nij}(c_1,c_2,x)=a_0xa_1...a_{i-1}c_1a_i......a_{j-1}c_2a_j...a_{n-1}xa_n \tag{4.2.24}$$

From the equality

$$p_{nji}(c_1,c_2,x)=a_0xa_1...a_{i-1}c_2a_i......a_{j-1}c_1a_j...a_{n-1}xa_n \tag{4.2.25}$$

and the equality (4.2.24), it follows that

$$\frac{d^2p_n(x)}{dx^2}\circ(c_1,c_2)=\frac{d^2p_n(x)}{dx^2}\circ(c_2,c_1) \tag{4.2.26}$$

[4.2] We accept some conventions in the equality (4.2.22).

- If $i=1$, then there is no first term.
- If $i=n$, then there is no last term.
- We assumed that $i<j$. In case $j<i$, the order of factors will be different.

The lemma follows from the equality (4.2.26). ⊙

LEMMA 4.2.10. *For any monomial $p_n(x)$ and for any $m > 3$*

$$\begin{aligned}&\frac{d^m p_n(x)}{dx^m}\circ(c_1,...,c_{i-1},c_i,c_{i+1},...,c_m)\\=&\frac{d^m p_n(x)}{dx^m}\circ(c_1,...,c_{i-1}c_{i+1},c_i,,...,c_m)\end{aligned}\tag{4.2.27}$$

PROOF. According to the theorem A.1.13,

$$\begin{aligned}&\frac{d^m p_n(x)}{dx^m}\circ(c_1,...,c_m)\\=&\frac{d^{m-i-1}p_n(x)}{dx^{m-i-1}}\left(\frac{d^2}{dx^2}\left(\frac{d^{i-1}p_n(x)}{dx^{i-1}}\circ(c_1,...,c_{i-1})\right)\circ(c_i,c_{i+1})\right)\\&\circ(c_{i+2},...,c_m)\end{aligned}\tag{4.2.28}$$

The expression

$$p(x)=\frac{d^{i-1}p_n(x)}{dx^{i-1}}\circ(c_1,...,c_{i-1})\tag{4.2.29}$$

is polynomial and according to the lemma 4.2.9

$$\frac{d^2p(x)}{dx^2}\circ(c_i,c_{i+1})=\frac{d^2p(x)}{dx^2}\circ(c_{i+1},c_i)\tag{4.2.30}$$

The lemma follows from equalities (4.2.28), (4.2.29), (4.2.30). ⊙

The theorem follows from the lemma 4.2.10.

We proved the statement of theorem. □

THEOREM 4.2.11. *Let $D$ be the complete commutative ring of characteristic 0. Let $A$ be associative Banach $D$-algebra. Then*

$$\frac{dx^{n+1}}{dx}=\sum_{i=0}^n x^i\otimes x^{n-i}\tag{4.2.31}$$

$$dx^{n+1}=\sum_{i=0}^n x^i dx x^{n-i}\tag{4.2.32}$$

PROOF. Since $x^{n+1}=(1\otimes...\otimes 1)\circ x$, then the theorem follows from the theorem 4.2.4. □

## 4.3. Matrix of tensors

Let $D$ be commutative ring and $A$ be associative $D$-algebra with unit.

DEFINITION 4.3.1. *Let $a$ be a matrix and $a^i_j\in A^{n\otimes}$. The matrix $a$ is called matrix of tensors $A^{n\otimes}$.* □

DEFINITION 4.3.2. *Following operations are defined on the set of matrices of tensors.*

- *Let $b\in A^{m\otimes}$.*

$$b*\begin{pmatrix}a^1_1&...&a^1_n\\...&...&...\\a^m_1&...&a^m_n\end{pmatrix}=\begin{pmatrix}b*a^1_1&...&b*a^1_n\\...&...&...\\b*a^m_1&...&b*a^m_n\end{pmatrix}$$

$$\begin{pmatrix} a_1^1 & ... & a_n^1 \\ ... & ... & ... \\ a_1^m & ... & a_n^m \end{pmatrix} * b = \begin{pmatrix} a_1^1 * b & ... & a_n^1 * b \\ ... & ... & ... \\ a_1^m * b & ... & a_n^m * b \end{pmatrix}$$

•
$$\begin{pmatrix} a_1^1 & ... & a_n^1 \\ ... & ... & ... \\ a_1^m & ... & a_n^m \end{pmatrix} {}_*{}^* \begin{pmatrix} b_1^1 & ... & b_k^1 \\ ... & ... & ... \\ b_1^n & ... & b_k^n \end{pmatrix} = \begin{pmatrix} a_i^1 * b_1^i & ... & a_i^1 * b_k^i \\ ... & ... & ... \\ a_i^m * b_1^i & ... & a_i^m * b_k^i \end{pmatrix}$$

•
$$\begin{pmatrix} a_1^1 & ... & a_n^1 \\ ... & ... & ... \\ a_1^m & ... & a_n^m \end{pmatrix} {}^*{}_* \begin{pmatrix} b_1^1 & ... & b_m^1 \\ ... & ... & ... \\ b_1^k & ... & b_m^k \end{pmatrix} = \begin{pmatrix} a_1^i * b_i^1 & ... & a_n^i * b_i^k \\ ... & ... & ... \\ a_1^i * b_i^k & ... & a_n^i * b_i^k \end{pmatrix}$$

• *Let* $a_j^i \in A^{p\otimes}, \quad x \in A.$

$$\begin{pmatrix} a_1^1 & ... & a_n^1 \\ ... & ... & ... \\ a_1^m & ... & a_n^m \end{pmatrix} \circ x^p = \begin{pmatrix} a_1^1 \circ x^p & ... & a_n^1 \circ x^p \\ ... & ... & ... \\ a_1^m \circ x^p & ... & a_n^m \circ x^p \end{pmatrix} \tag{4.3.1}$$

• *Let* $a_j^i \in A^{p\otimes}, \quad x_1, ...., x_p \in A.$

$$\begin{pmatrix} a_1^1 & ... & a_n^1 \\ ... & ... & ... \\ a_1^m & ... & a_n^m \end{pmatrix} \circ (x_1, ..., x_p) = \begin{pmatrix} a_1^1 \circ (x_1, ..., x_p) & ... & a_n^1 \circ (x_1, ..., x_p) \\ ... & ... & ... \\ a_1^m \circ (x_1, ..., x_p) & ... & a_n^m \circ (x_1, ..., x_p) \end{pmatrix}$$

• *Let* $a_j^i \in A^{2\otimes}, \quad b_j^i \in A^{2\otimes},$

$$\begin{pmatrix} a_1^1 & ... & a_n^1 \\ ... & ... & ... \\ a_1^m & ... & a_n^m \end{pmatrix} {}_\circ{}^\circ \begin{pmatrix} b_1^1 & ... & b_k^1 \\ ... & ... & ... \\ b_1^n & ... & b_k^n \end{pmatrix} = \begin{pmatrix} a_i^1 \circ b_1^i & ... & a_i^1 \circ b_k^i \\ ... & ... & ... \\ a_i^m \circ b_1^i & ... & a_i^m \circ b_k^i \end{pmatrix}$$

• *Let* $a_j^i \in A^{2\otimes}, \quad b_j^i \in A^{2\otimes},$

$$\begin{pmatrix} a_1^1 & ... & a_n^1 \\ ... & ... & ... \\ a_1^m & ... & a_n^m \end{pmatrix} {}^\circ{}_\circ \begin{pmatrix} b_1^1 & ... & b_m^1 \\ ... & ... & ... \\ b_1^k & ... & b_m^k \end{pmatrix} = \begin{pmatrix} a_1^i \circ b_i^1 & ... & a_n^i \circ b_i^k \\ ... & ... & ... \\ a_1^i \circ b_i^k & ... & a_n^i \circ b_i^k \end{pmatrix}$$

□

It is important to pay attention to the difference between ${}_*{}^*$-product and ${}_\circ{}^\circ$-product

$$\begin{pmatrix}{a_1}^1_1\otimes {a_2}^1_1&...&{a_1}^1_n\otimes {a_2}^1_n\\ ...&...&...\\ {a_1}^m_1\otimes {a_2}^m_1&...&{a_1}^m_n\otimes {a_2}^m_n\end{pmatrix}{}_*{}^*\begin{pmatrix}{b_1}^1_1\otimes {b_2}^1_1&...&{b_1}^1_k\otimes {b_2}^1_k\\ ...&...&...\\ {b_1}^n_1\otimes {b_2}^n_1&...&{b_1}^n_k\otimes {b_2}^n_k\end{pmatrix}$$
$$=\begin{pmatrix}{a_1}^1_i\otimes {a_2}^1_i{b_1}^i_1\otimes {b_2}^i_1&...&{a_1}^1_i\otimes {a_2}^1_i{b_1}^i_k\otimes {b_2}^i_k\\ ...&...&...\\ {a_1}^m_i\otimes {a_2}^m_i{b_1}^i_1\otimes {b_2}^i_1&...&{a_1}^m_i\otimes {a_2}^m_i{b_1}^i_k\otimes {b_2}^i_k\end{pmatrix}$$

$$\begin{pmatrix}{a_1}^1_1\otimes {a_2}^1_1&...&{a_1}^1_n\otimes {a_2}^1_n\\ ...&...&...\\ {a_1}^m_1\otimes {a_2}^m_1&...&{a_1}^m_n\otimes {a_2}^m_n\end{pmatrix}{}_\circ{}^\circ\begin{pmatrix}{b_1}^1_1\otimes {b_2}^1_1&...&{b_1}^1_k\otimes {b_2}^1_k\\ ...&...&...\\ {b_1}^n_1\otimes {b_2}^n_1&...&{b_1}^n_k\otimes {b_2}^n_k\end{pmatrix}$$
$$=\begin{pmatrix}{a_1}^1_i{b_1}^i_1\otimes {b_2}^i_1{a_2}^1_i&...&{a_1}^1_i{b_1}^i_k\otimes {b_2}^i_k{a_2}^1_i\\ ...&...&...\\ {a_1}^m_i{b_1}^i_1\otimes {b_2}^i_1{a_2}^m_i&...&{a_1}^m_i{b_1}^i_k\otimes {b_2}^i_k{a_2}^m_i\end{pmatrix}$$

THEOREM 4.3.3. *Let $A$ be $D$-algebra. Let $a$ be a matrix of tensors $A^{p\otimes}$. Let $b$ be a matrix of tensors $A^{q\otimes}$. Then*

$$(a{}_*{}^*b)\circ x^{p+q}=(a\circ x^p){}_*{}^*(b\circ x^q)\tag{4.3.2}$$

PROOF. The equality

$$((a{}_*{}^*b)\circ x^{p+q})^i_j=(a{}_*{}^*b)^i_j\circ x^{p+q}=(a^i_kb^k_j)\circ x^{p+q}\tag{4.3.3}$$

follows from equalities (3.2.3), (4.3.1). The equality

$$((a{}_*{}^*b)\circ x^{p+q})^i_j=(a^i_k\circ x^p)(b^k_j\circ x^q)=(a\circ x^p)^i_k(b\circ x^q)^k_j\tag{4.3.4}$$

follows from equalities (4.1.2), (4.3.1), (4.3.3). The equality (4.3.2) follows from equalities (3.2.3), (4.3.4). □

THEOREM 4.3.4. *Let $A$ be $D$-algebra and $a$ be a matrix of tensors $A^{2\otimes}$. Then*

$$a^{(k+1){}_*{}^*}\circ x^{k+1}=(a^{k{}_*{}^*}\circ x^k){}_*{}^*(a\circ x)\tag{4.3.5}$$

PROOF. We will prove the theorem by induction over $k$.

The equality (4.3.5) is true, if $k=1$, because

$$a^{0{}_*{}^*}\circ x^0=E$$

Let the theorem be true for $k=l$

$$a^{(l+1){}_*{}^*}\circ x^{l+1}=(a^{l{}_*{}^*}\circ x^l){}_*{}^*(a\circ x)\tag{4.3.6}$$

The equality

$$a^{(l+2){}_*{}^*}\circ x^{l+2}=(a^{(l+1){}_*{}^*}{}_*{}^*a)\circ x^{l+2}=(a^{(l+1){}_*{}^*}\circ x^{l+1}){}_*{}^*(a\circ x)$$

follows from the equality (4.3.2). Therefore, the theorem be true for $k=l+1$. □

Theorem 4.3.5. *Let $A$ be $D$-algebra and $a$ be a matrix of tensors $A^{2\otimes}$. Then*

$$a^{n_*{}^*} \circ x^n = (a \circ x)^{n_*{}^*} \tag{4.3.7}$$

Proof. We will prove the theorem by induction over $n$.

The equality (4.3.7) is true, if $n = 1$ or $n = 1$.

Let the theorem be true for $n = k$

$$a^{k_*{}^*} \circ x^k = (a \circ x)^{k_*{}^*} \tag{4.3.8}$$

According to the theorem 4.3.4, the equality

$$a^{(k+1)_*{}^*} \circ x^{k+1} = (a^k \circ x^k){}_*{}^*(a \circ x) = (a \circ x)^{k_*{}^*}{}_*{}^*(a \circ x) \tag{4.3.9}$$

follows from the equality (4.3.8). The equality

$$a^{(k+1)_*{}^*} \circ x^{k+1} = (a \circ x)^{(k+1)_*{}^*} \tag{4.3.10}$$

follows from equalities (3.2.16), (4.3.9). Therefore, the theorem be true for $n = k + 1$. □

CHAPTER 5

# Vector Space over Division Algebra

## 5.1. Vector Space

Let $D$ be commutative associative ring with unit. Let $A$ be associative $D$-algebra.

THEOREM 5.1.1. *Let $D$-algebra $A$ be division algebra. $D$-algebra $A$ has unit.*

PROOF. The theorem follows from the statement that the equation

$$ax = a$$

has solution for any $a \in A$. □

THEOREM 5.1.2. *Let $D$-algebra $A$ be division algebra. The ring $D$ is the field and subset of the center of $D$-algebra $A$.*

PROOF. Let $e$ be unit of $D$-algebra $A$. Since the map

$$d \in D \to d\,e \in A$$

is embedding of the ring $D$ into $D$-algebra $A$, then the ring $D$ is subset of the center of $D$-algebra $A$. □

DEFINITION 5.1.3. *Let $A$ be division algebra. Effective left-side representation*

$$f : A \xrightarrow{*} V \qquad f(a) : v \in V \to av \in V \quad a \in A$$

*of Abelian group $A$ in $D$-module $V$ is called* **left vector space** *over $D$-algebra $A$. We will also say that $D$-module $V$ is* **left $A$-vector space** *or* **$A*$-vector space***. $V$-number is called* **vector***.* □

DEFINITION 5.1.4. *Let $A$ be division algebra. Effective right-side representation*

$$f : A \xrightarrow{*} V \qquad f(a) : v \in V \to va \in V \quad a \in A$$

*of Abelian group $A$ in $D$-module $V$ is called* **right vector space** *over $D$-algebra $A$. We will also say that $D$-module $V$ is* **left $A$-vector space** *or* **$*A$-vector space***. $V$-number is called* **vector***.* □

Since any statement about left $A$-vector space is equivalent to a statement about right $A$-vector space up to the order of multipliers, in this chapter, we will consider left $A$-vector space. If $D$-algebra $A$ is commutative, then left $A$-vector space coincides with right $A$-vector space.

Theorem 5.1.5. *The following diagram of representations describes left $A$-vector space $V$*

(5.1.1)

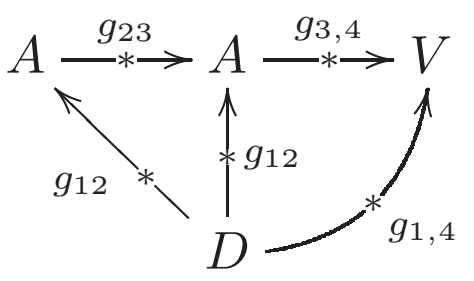

$$g_{12}(d) : a \to d\,a$$
$$g_{23}(v) : w \to C(w,v)$$
$$C \in \mathcal{L}(A^2 \to A)$$
$$g_{3,4}(a) : v \to v\,a$$
$$g_{1,4}(d) : v \to d\,v$$

*The diagram of representations (5.1.1) holds* **commutativity of representations** *of commutative ring $D$ and $D$-algebra $A$ in Abelian group $V$*

$$a(dv) = d(av) \tag{5.1.2}$$

Proof. The diagram of representations (5.1.1) follows from the definition 5.1.3, from the theorem [13]-5.1.11 and the corollary [13]-5.1.12. Since left-side transformation $g_{3,4}(a)$ is endomorphism of $D$-vector space $V$, we obtain the equality (5.1.2). □

Theorem 5.1.6. *Let $V$ be left $A$-vector space. Following conditions hold for left $A$-vector space $V$:*

- *if $D$-algebra $A$ is associative, then* **associative law** *holds*

$$(pq)v = p(qv) \tag{5.1.3}$$

- **distributive law**

$$p(v+w) = pv + pw \tag{5.1.4}$$
$$(p+q)v = pv + qv \tag{5.1.5}$$

- **unitarity law**

$$1v = v \tag{5.1.6}$$

*for any $p, q \in A$, $v, w \in V$.*

Proof. Since left-side transformation $p$ is endomorphism of the Abelian group $V$, we obtain the equation (5.1.4). Since representation $g_{3,4}$ is homomorphism of the additive group of $D$-algebra $A$, we obtain the equation (5.1.5). Since representation $g_{3,4}$ is left-side representation of the multiplicative group of $D$-algebra $A$, we obtain the equality (5.1.6). Since representation $g_{3,4}$ is left-side representation of the multiplicative group of $D$-algebra $A$, we obtain the equality (5.1.3) when $D$-algebra $A$ is associative. □

## 5.2. Basis of Left $A$-Vector Space

Theorem 5.2.1. *Let $V$ be left $A$-vector space. The set of vectors generated by the set of vectors $v = (v_i \in V, i \in I)$ has form*[5.1]

$$J(v) = \left\{ w : w = \sum_{i \in I} c^i v_i, c^i \in A, |\{i : c^i \neq 0\}| < \infty \right\} \tag{5.2.1}$$

[5.1] For a set $A$, we denote by $|A|$ the cardinal number of the set $A$. The notation $|A| < \infty$ means that the set $A$ is finite.

PROOF. We prove the theorem by induction based on the theorem [12]-3.4, Acording to the theorem [12]-3.4, we need to prove following statements:

5.2.1.1: $v_k \in X_0 \subseteq J(v)$
5.2.1.2: $c^k v_k \in J(v)$, $c^k \in A$, $k \in I$
5.2.1.3: $\sum_{k\in I} c^k v_k \in J(v)$, $c^k \in A$, $|\{i : c^i \neq 0\}| < \infty$
5.2.1.4: $v, w \in J(v) \Rightarrow v + w \in J(v)$
5.2.1.5: $a \in A$, $v \in J(v) \Rightarrow av \in J(v)$

- For any $v_k \in v$, let $c^i = \delta^i_k \in A$. Then

$$v_k = \sum_{i\in I} c^i v_i \tag{5.2.2}$$

The statement 5.2.1.1 follows from equalities (5.2.1), (5.2.2).

- The statement 5.2.1.2 follow from the theorems [12]-3.4, 5.1.6 and from the statement 5.2.1.1.
- Since $V$ is Abelian group, then the statement 5.2.1.3 follows from the statement 5.2.1.2 and from the theorem [12]-3.4.
- Let $w_1, w_2 \in X_k \subseteq J(v)$. Since $V$ is Abelian group, then, according to the statement [12]-3.4.3,

$$w_1 + w_2 \in X_{k+1} \tag{5.2.3}$$

According to the equality (5.2.1), there exist $A$-numbers $w^i_1, w^i_2$, $i \in I$, such that

$$w_1 = \sum_{i\in I} w^i_1 v_i \qquad w_2 = \sum_{i\in I} w^i_2 v_i \tag{5.2.4}$$

where sets

$$H_1 = \{i \in I : w^i_1 \neq 0\} \qquad H_2 = \{i \in I : w^i_2 \neq 0\} \tag{5.2.5}$$

are finite. Since $V$ is Abelian group, then from the equality (5.2.4) it follows that

$$w_1 + w_2 = \sum_{i\in I} w^i_1 v_i + \sum_{i\in I} w^i_2 v_i = \sum_{i\in I} (w^i_1 v_i + w^i_2 v_i) \tag{5.2.6}$$

The equality

$$w_1 + w_2 = \sum_{i\in I} (w^i_1 + w^i_2) v_i \tag{5.2.7}$$

follows from equalities (5.1.5), (5.2.6). From the equality (5.2.5), it follows that the set

$$\{i \in I : w^i_1 + w^i_2 \neq 0\} \subseteq H_1 \cup H_2$$

is finite.

- Let $w \in X_k \subseteq J(v)$. According to the statement [12]-3.4.4, for any $A$-number $a$,

$$aw \in X_{k+1} \tag{5.2.8}$$

According to the equality (5.2.1), there exist $A$-numbers $w^i$, $i \in I$, such that

$$w = \sum_{i\in I} w^i v_i \tag{5.2.9}$$

where

$$|\{i \in I : w^i \neq 0\}| < \infty \tag{5.2.10}$$

From the equality (5.2.9) it follows that

$$aw = a\sum_{i\in I} w^i v_i = \sum_{i\in I} a(w^i v_i) = \sum_{i\in I}(aw^i)v_i \tag{5.2.11}$$

From the statement (5.2.10), it follows that the set $\{i \in I : aw^i \neq 0\}$ is finite.

From equalities (5.2.3), (5.2.7), (5.2.8), (5.2.11), it follows that $X_{k+1} \subseteq J(v)$. □

Definition 5.2.2. *$J(v)$ is called vector subspace generated by set $v$, and $v$ is a* **generating set** *of vector subspace $J(v)$. In particular, a* **generating set** *of left $A$-vector space $V$ is a subset $X \subset V$ such that $J(X) = V$.* □

Definition 5.2.3. *If the set $X \subset V$ is generating set of left $A$-vector space $V$, then any set $Y$, $X \subset Y \subset V$ also is generating set of left $A$-vector space $V$. If there exists minimal set $X$ generating the left $A$-vector space $V$, then the set $X$ is called* **basis** *of left $A$-vector space $V$.* □

Definition 5.2.4. *Let $v = (v_i \in V, i \in I)$ be set of vectors. The expression $w^i v_i$ is called* **linear combination** *of vectors $v_i$. A vector $\overline{w} = w^i v_i$ is called* **linearly dependent** *on vectors $v_i$.* □

Theorem 5.2.5. *Let $A$ be associative division $D$-algebra. Since the equation*

$$w^i v_i = 0$$

*implies existence of index $i = j$ such that $w^j \neq 0$, then the vector $v_j$ linearly depends on rest of vectors $v$.*

Proof. The theorem follows from the equality

$$v_j = \sum_{i\in I\setminus\{j\}} (w^j)^{-1} w^i v_i$$

and from the definition 5.2.4. □

It is evident that for any set of vectors $v_i$

$$w^i = 0 \Rightarrow w^*{}_* v = 0$$

Definition 5.2.6. *The set of vectors*[5.2] *$v_i$, $i \in I$, of left $A$-vector space $V$ is* **linearly independent** *if $w = 0$ follows from the equation*

$$w^i v_i = 0$$

*Otherwise the set of vectors $v_i$, $i \in I$, is* **linearly dependent**. □

Theorem 5.2.7. *Let $A$ be associative division $D$-algebra. The set of vectors $\overline{\overline{e}} = (e_i, i \in I)$ is a* **basis of left $A$-vector space** *$V$ if vectors $e_i$ are linearly independent and any vector $v \in V$ linearly depends on vectors $e_i$.*

[5.2] I follow to the definition in [1], page 130.

Proof. Let the set of vectors $e_i$, $i \in I$, be linear dependent. Then the equation

$$w^i e_i = 0$$

implies existence of index $i = j$ such that $w^j \neq 0$. According to the theorem 5.2.5, the vector $e_j$ linearly depends on rest of vectors of the set $\overline{\overline{e}}$. According to the definition 5.2.3, the set of vectors $e_i$, $i \in I$, is not a basis for left $A$-vector space $V$.

Therefore, if the set of vectors $e_i$, $i \in I$, is a basis, then these vectors are linearly independent. Since an arbitrary vector $v \in V$ is linear combination of vectors $e_i$, $i \in I$, , then the set of vectors $v$, $e_i$, $i \in I$, is not linearly independent. □

Definition 5.2.8. *Let $\overline{\overline{e}}$ be the basis of left $A$-vector space $V$ and vector $\overline{v} \in V$ has expansion*

$$\overline{v} = v^*{}_* e = v^i e_i$$

*with respect to the basis $\overline{\overline{e}}$. $A$-numbers $v^i$ are called* **coordinates** *of vector $\overline{v}$ with respect to the basis $\overline{\overline{e}}$. Matrix of $A$-numbers $v = (v^i, i \in I)$ is called* **coordinate matrix of vector** $\overline{v}$ *in basis $\overline{\overline{e}}$.* □

Theorem 5.2.9. *Coordinates of vector $v \in V$ relative to basis $\overline{\overline{e}}$ of left $A$-vector space $V$ are uniquely defined.*

Proof. According to the theorem 5.2.7, the system of vectors $\overline{v}$, $e_i$, $i \in I$, is linearly dependent and in equality

$$b\overline{v} + c^*{}_* e = 0 \tag{5.2.12}$$

at least $b$ is different from 0. Then equality

$$\overline{v} = (-cb^{-1})^*{}_* e \tag{5.2.13}$$

follows from (5.2.12). The equality

$$\overline{v} = v^*{}_* e \tag{5.2.14}$$

follows from the equality (5.2.13).

Assume we get another expansion

$$\overline{v} = v'{}_*{}^* e \tag{5.2.15}$$

We subtract (5.2.14) from (5.2.15) and get

$$0 = (v' - v)^*{}_* e \tag{5.2.16}$$

Since vectors $e_i$ are linearly independent, then the equality

$$v' - v = 0 \tag{5.2.17}$$

follows from the equality (5.2.16). Therefore, the theorem follows from the equality (5.2.17). □

Theorem 5.2.10. *We can consider division $D$-algebra $A$ as left $A$-vector space with basis $\overline{\overline{e}} = \{1\}$.*

Proof. The theorem follows from the equality (5.1.6). □

## 5.3. Vector Space Type

Since a basis of $A$-vector space $V$ is finite, then the power of the basis is called dimention of $A$-vector space either we say that $A$-vector space $V$ finite dimensional. According to the definition 5.2.8, the theory of finite dimensional $A$-vector space $V$ is associated with the theory of matrices.

Example 5.3.1. *In the section 5.2, we represented the set of vectors* $v_1$, ..., $v_n$ *as row of matrix*

$$v = \begin{pmatrix} v_1 & ... & v_n \end{pmatrix}$$

*and the set of $A$-nummbers* $c^1$, ..., $c^n$ *as column of matrix*

$$c = \begin{pmatrix} c^1 \\ ... \\ c^n \end{pmatrix}$$

*Thus we can represent linear combination of vectors* $v_1$, ..., $v_n$ *as* ${}^*{}_*$*-product of matrices*

$$c^*{}_*v = c^i v_i = \begin{pmatrix} c^1 \\ ... \\ c^n \end{pmatrix} {}^*{}_* \begin{pmatrix} v_1 & ... & v_n \end{pmatrix} \tag{5.3.1}$$

*The expression* (5.3.1) *for linear combination of vectors is called* **left linear combination of columns**. *In particular, if we write vectors of basis* $\overline{\overline{e}}$ *as row of matrix*

$$e = \begin{pmatrix} e_1 & ... & e_n \end{pmatrix}$$

*and coordinates of vector* $\overline{w}$ *with respect to basis* $\overline{\overline{e}}$ *as column of matrix*

$$w = \begin{pmatrix} w^1 \\ ... \\ w^n \end{pmatrix}$$

*then we can represent the vector* $\overline{w}$ *as* ${}^*{}_*$*-product of matrices*

$$\overline{w} = w^*{}_*e = w^i e_i = \begin{pmatrix} w^1 \\ ... \\ w^n \end{pmatrix} {}^*{}_* \begin{pmatrix} e_1 & ... & e_n \end{pmatrix}$$

*Corresponding representation of vector space $V$ is called* **left $A^*$-vector space** *or* **left $A$-column space**, *and $V$-number is called* **column vector**. □

The problem which we consider in particular case determines vector representation format.

Example 5.3.2. *We represent the set of vectors* $v_1, ..., v_n$ *as column of matrix*

$$v=\begin{pmatrix}v^1\\ ...\\ v^n\end{pmatrix}$$

*and the set of A-nummbers* $c^1, ..., c^n$ *as row of matrix*

$$c=\begin{pmatrix}c^1\\ ...\\ c^n\end{pmatrix}$$

*Thus we can represent linear combination of vectors* $v_1, ..., v_n$ *as* ${}_*{}^*$*-product of matrices*

$$c_*{}^*v=c_iv^i=\begin{pmatrix}c_1 & ... & c_n\end{pmatrix}{}_*{}^*\begin{pmatrix}v^1\\ ...\\ v^n\end{pmatrix}\tag{5.3.2}$$

*The expression* (5.3.2) *for linear combination of vectors is called* **left linear combination of rows**. *In particular, if we write vectors of basis* $\overline{\overline{e}}$ *as column of matrix*

$$e=\begin{pmatrix}e^1\\ ...\\ e^n\end{pmatrix}$$

*and coordinates of vector* $\overline{w}$ *with respect to basis* $\overline{\overline{e}}$ *as row of matrix*

$$w=\begin{pmatrix}w_1 & ... & w_n\end{pmatrix}$$

*then we can represent the vector* $\overline{w}$ *as* ${}_*{}^*$*-product of matrices*

$$\overline{w}=w_*{}^*e=w_ie^i=\begin{pmatrix}w_1 & ... & w_n\end{pmatrix}{}_*{}^*\begin{pmatrix}e^1\\ ...\\ e^n\end{pmatrix}$$

*Corresponding representation of vector space* $V$ *is called* **left $A_*$-vector space** *or* **left $A$-row space**, *and $V$-number is called* **row vector**. □

Similar types of vector space we can consider for right vector space.

Example 5.3.3. *We represent the set of vectors* $v_1, ..., v_n$ *as row of matrix*

$$v=\begin{pmatrix}v_1 & ... & v_n\end{pmatrix}$$

*and the set of A-nummbers* $c^1, ..., c^n$ *as column of matrix*

$$c=\begin{pmatrix}c^1\\ ...\\ c^n\end{pmatrix}$$

*Thus we can represent linear combination of vectors* $v_1, ..., v_n$ *as* ${}_*{}^*$*-product of matrices*

$$v_*{}^*c = v_i c^i = \begin{pmatrix} v_1 & ... & v_n \end{pmatrix} {}_*{}^* \begin{pmatrix} c^1 \\ ... \\ c^n \end{pmatrix} \tag{5.3.3}$$

*The expression* (5.3.3) *for linear combination of vectors is called* **right linear combination of columns**. *In particular, if we write vectors of basis* $\overline{\overline{e}}$ *as row of matrix*

$$e = \begin{pmatrix} e_1 & ... & e_n \end{pmatrix}$$

*and coordinates of vector* $\overline{w}$ *with respect to basis* $\overline{\overline{e}}$ *as column of matrix*

$$w = \begin{pmatrix} w^1 \\ ... \\ w^n \end{pmatrix}$$

*then we can represent the vector* $\overline{w}$ *as* ${}_*{}^*$*-product of matrices*

$$\overline{w} = e_*{}^*w = e_i w^i = \begin{pmatrix} e_1 & ... & e_n \end{pmatrix} {}_*{}^* \begin{pmatrix} w^1 \\ ... \\ w^n \end{pmatrix}$$

*Corresponding representation of vector space* $V$ *is called* **right** $A^*$**-vector space** *or* **right** $A$**-column space**, *and* $V$*-number is called* **column vector**. □

Example 5.3.4. *We represent the set of vectors* $v_1, ..., v_n$ *as column of matrix*

$$v = \begin{pmatrix} v^1 \\ ... \\ v^n \end{pmatrix}$$

*and the set of* $A$*-nummbers* $c^1, ..., c^n$ *as row of matrix*

$$c = \begin{pmatrix} c^1 \\ ... \\ c^n \end{pmatrix}$$

*Thus we can represent linear combination of vectors* $v_1, ..., v_n$ *as* ${}^*{}_*$*-product of matrices*

$$v^*{}_*c = v^i c_i = \begin{pmatrix} v^1 \\ ... \\ v^n \end{pmatrix} {}^*{}_* \begin{pmatrix} c_1 & ... & c_n \end{pmatrix} \tag{5.3.4}$$

*The expression* (5.3.4) *for linear combination of vectors is called* **right linear combination of rows**. *In particular, if we write vectors of basis* $\overline{\overline{e}}$ *as column of matrix*

$$e = \begin{pmatrix} e^{1} \\ ... \\ e^{n} \end{pmatrix}$$

*and coordinates of vector* $\overline{w}$ *with respect to basis* $\overline{\overline{e}}$ *as row of matrix*

$$w = \begin{pmatrix} w_{1} & ... & w_{n} \end{pmatrix}$$

*then we can represent the vector* $\overline{w}$ *as* ${}^*{}_*$*-product of matrices*

$$\overline{w} = e^*{}_*w = e^{i}w_{i} = \begin{pmatrix} e^{1} \\ ... \\ e^{n} \end{pmatrix} {}^*{}_* \begin{pmatrix} w_{1} & ... & w_{n} \end{pmatrix}$$

*Corresponding representation of vector space* $V$ *is called* **right** $A_*$**-vector space** *or* **right** $A$**-row space**, *and* $V$*-number is called* **row vector**. □

Choice of type of vector space is agreement required to solve a specific problem.

## 5.4. System of Linear Equations

Definition 5.4.1. *Let* $V$ *be a right* $A^*$*-vector space and* $\{a_{i} \in V, \boldsymbol{i} \in \boldsymbol{I}\}$ *be set of vectors.* **Linear span** *in right* $A^*$*-vector space is set* $span(\overline{a}_{i}, \boldsymbol{i} \in \boldsymbol{I})$ *of vectors linearly dependent on vectors* $a_i$. □

Theorem 5.4.2. *Let* $span(\overline{a}_{i}, \boldsymbol{i} \in \boldsymbol{I})$ *be linear span in right* $A_*$*-vector space* $V$*. Then* $span(\overline{a}_{i}, \boldsymbol{i} \in \boldsymbol{I})$ *is subspace of right* $A_*$*-vector space* $V$*.*

Proof. Suppose

$$\overline{b} \in \mathrm{span}(\overline{a}_{i}, \boldsymbol{i} \in \boldsymbol{I})$$
$$\overline{c} \in \mathrm{span}(\overline{a}_{i}, \boldsymbol{i} \in \boldsymbol{I})$$

According to definition 5.4.1

$$\overline{b} = a_*{}^*b$$
$$\overline{c} = a_*{}^*c$$

Then

$$\overline{b} + \overline{c} = a_*{}^*b + a_*{}^*c = a_*{}^*(b + c) \in \mathrm{span}(\overline{a}_{i}, \boldsymbol{i} \in \boldsymbol{I})$$
$$\overline{b}k = (a_*{}^*b)k = a_*{}^*(bk) \in \mathrm{span}(\overline{a}_{i}, \boldsymbol{i} \in \boldsymbol{I})$$

This proves the statement. □

Example 5.4.3. *Let* $V$ *be a right* $A_*$*-vector space and row*

$$\overline{a} = \begin{pmatrix} a_{1} & ... & a_{n} \end{pmatrix} = (a_{i}, \boldsymbol{i} \in \boldsymbol{I})$$

*be set of vectors. To answer the question of whether vector* $\overline{b} \in span(a_{i}, \boldsymbol{i} \in \boldsymbol{I})$ *we write linear equation*

$$\overline{b} = \overline{a}_*{}^*x \tag{5.4.1}$$

*where*

$$x=\begin{pmatrix}x^{1}\\ ...\\ x^{n}\end{pmatrix}$$

*is column of unknown coefficients of expansion.* $\overline{b}\in span(a_{i},i\in I)$ *if equation* (5.4.1) *has a solution. Suppose* $\overline{\overline{f}}=(f_{j},j\in J)$ *is a basis. Then vectors* $\overline{b}$, $\overline{a}_i$ *have expansion*

$$\overline{b}=f_{*}{}^{*}b \tag{5.4.2}$$

$$\overline{a}_i=f_{*}{}^{*}a_i \tag{5.4.3}$$

*If we substitute* (5.4.2) *and* (5.4.3) *into* (5.4.1) *we get*

$$f_{*}{}^{*}b=f_{*}{}^{*}a_{*}{}^{*}x \tag{5.4.4}$$

*Applying theorem 5.2.9 to* (5.4.4) *we get* **system of linear equations**

$$a_{*}{}^{*}x=b \tag{5.4.5}$$

*We can write system of linear equations* (5.4.5) *in one of the next forms*

$$\begin{pmatrix}a_{1}^{1}&...&a_{n}^{1}\\ ...&...&...\\ a_{1}^{m}&...&a_{n}^{m}\end{pmatrix}{}_{*}{}^{*}\begin{pmatrix}x^{1}\\ ...\\ x^{n}\end{pmatrix}=\begin{pmatrix}b^{1}\\ ...\\ b^{m}\end{pmatrix}$$

$$a_{i}^{j}x^{i}=b^{j} \tag{5.4.6}$$

$$\begin{matrix}a_{1}^{1}x^{1}&+...&+a_{n}^{1}x^{n}&=b^{1}\\ ...&...&...&...\\ a_{m}^{1}x^{1}&+...&+a_{n}^{m}x^{n}&=b^{m}\end{matrix}$$

□

Example 5.4.4. *Let* $V$ *be a left* $A_{*}$*-vector space and column*

$$\overline{a}=\begin{pmatrix}\overline{a}^{1}\\ ...\\ \overline{a}^{m}\end{pmatrix}=(\overline{a}^{j},j\in J)$$

*be set of vectors. To answer the question of whether vector* $\overline{b}\in span(\overline{a}^{j},j\in J)$ *we write linear equation*

$$\overline{b}=x_{*}{}^{*}\overline{a} \tag{5.4.7}$$

*where*

$$x=\begin{pmatrix}x_{1}&...&x_{m}\end{pmatrix}$$

*is row of unknown coefficients of expansion.* $\overline{b} \in span(\overline{a}^j, j \in J)$, *if equation* (5.4.7) *has a solution. Suppose* $\overline{\overline{f}} = (f^i, i \in I)$ *is a basis. Then vectors* $\overline{b}$, $a^j$ *have expansion*

$$\overline{b} = b_*{}^* f \tag{5.4.8}$$

$$\overline{a}^j = a^j{}_*{}^* f \tag{5.4.9}$$

*If we substitute* (5.4.8) *and* (5.4.9) *into* (5.4.7) *we get*

$$b_*{}^* f = x_*{}^* a_*{}^* f \tag{5.4.10}$$

*Applying theorem 5.2.9 to* (5.4.10) *we get* **system of linear equations**[5.3]

$$x_*{}^* a = b \tag{5.4.11}$$

*We can write system of linear equations* (5.4.11) *in one of the next forms*

$$\begin{pmatrix} x_1 & ... & x_m \end{pmatrix} {}_*{}^* \begin{pmatrix} a_1^1 & ... & a_n^1 \\ ... & ... & ... \\ a_1^m & ... & a_n^m \end{pmatrix} = \begin{pmatrix} b_1 & ... & b_n \end{pmatrix}$$

$$x_j a_i^j = b_i \tag{5.4.12}$$

$$\begin{matrix} x_1 a_1^1 & ... & x_1 a_n^1 \\ + & ... & + \\ ... & ... & ... \\ + & ... & + \\ x_m a_1^m & ... & x_m a_n^m \\ = & ... & = \\ b_1 & ... & b_n \end{matrix}$$

□

To find a solution of system of linear equations, we need a matrix of this system. From examples 5.3.1, 5.3.3, we see that we may consider a column of matrix as vector of left or right $A$-vector space. To make statements more clear, we will distinguish the format of linear dependence For instance, the statement

- Columns of matrix are $D^*{}_*$-linear dependent.

means that

- Columns of matrix are vectors of $D^*{}_*$-vector space, and corresponding vectors are linear dependent.

In particular, system of linear equations (5.4.5) in ${}_*{}^*D$-vector space is called **system of** $_*{}^*D$**-linear equations**. and system of linear equations (5.4.11) in $D_*{}^*$-vector space is called **system of** $D_*{}^*$**-linear equations**.

Definition 5.4.5. *If* $n \times n$ *matrix* $a$ *has* ${}_*{}^*$*-inverse matrix we call such matrix* ${}_*{}^*$**-nonsingular matrix**. *Otherwise, we call such matrix* ${}_*{}^*$**-singular matrix**. □

[5.3] Reading system of linear equations (5.4.5) in right $A_*$-vector space from bottom up and from left to right we get system of linear equations (5.4.11) in left $A_*$-vector space.

Definition 5.4.6. *Suppose $A$ is ${}_*{}^*$-nonsingular matrix. We call appropriate system of ${}_*{}^*D$-linear equations*

$$a_*{}^*x = b \tag{5.4.13}$$

**nonsingular system of ${}_*{}^*D$-linear equations**. □

Theorem 5.4.7. *Solution of nonsingular system of ${}_*{}^*D$-linear equations (5.4.13) is determined uniquely and can be presented in either form*[5.4]

$$x = A^{-1_*{}^*}{}_*{}^*b \tag{5.4.14}$$

$$x = \mathcal{H}\det({}_*{}^*)a_*{}^*b \tag{5.4.15}$$

Proof. Multiplying both sides of the equation (5.4.13) from left by $a^{-1_*{}^*}$, we get (5.4.14). Using definition (3.3.10) we get (5.4.15). Based on the theorem 3.2.20, the solution is unique. □

## 5.5. Rank of Matrix

In this section, we will make the following assumption.

- $i \in M$, $|M| = m$, $j \in N$, $|N| = n$.
- $a = (a^i_j)$ is an arbitrary matrix.
- $k$, $s \in S \supseteq M$, $l$, $t \in T \supseteq N$, $k = |S| = |T|$.
- $p \in M \setminus S$, $r \in N \setminus T$.

Definition 5.5.1. *Matrix $a^S_T$ is a minor matrix of an order $k$.* □

Definition 5.5.2. *If minor matrix $a^S_T$ is ${}_*{}^*$-nonsingular matrix then we say that ${}_*{}^*$-rank of matrix $a$ is not less then $k$.* **${}_*{}^*$-rank of matrix** $a$

$$\mathrm{rank}_{{}_*{}^*}\, a$$

*is the maximal value of $k$. We call an appropriate minor matrix the* **${}_*{}^*$-major minor matrix**. □

Theorem 5.5.3. *Let matrix $a$ be ${}_*{}^*$-singular matrix and minor matrix $a^S_T$ be ${}_*{}^*$-major minor matrix. Then*

$$\det({}_*{}^*)^p_r\, a^{S\cup\{p\}}_{T\cup\{r\}} = 0$$

Proof. To understand why minor matrix

$$b = a^{S\cup\{p\}}_{T\cup\{r\}} \tag{5.5.1}$$

does not have ${}_*{}^*$-inverse matrix,[5.5] we assume that there exists ${}_*{}^*$-inverse matrix $b^{-1_*{}^*}$. We write down the system of linear equations (3.3.4), (3.3.5) in case $I = \{r\}$, $J = \{p\}$ (then $[I] = T$, $[J] = S$ )

$$b^S_T{}_*{}^*b^{-1_*{}^*}{}^T_p + b^S_r\, b^{-1_*{}^*}{}^r_p = 0 \tag{5.5.2}$$

[5.4]We can see a solution of system (5.4.13) in theorem [4]-1.6.1. I repeat this statement because I slightly changed the notation.

[5.5]It is natural to expect relationship between ${}_*{}^*$-singularity of the matrix and its ${}_*{}^*$-quasideterminant similar to relationship which is known in commutative case. However ${}_*{}^*$-quasideterminant is defined not always. For instance, it is not defined when ${}_*{}^*$-inverse matrix has too much elements equal 0. As it follows from this theorem, the ${}_*{}^*$-quasideterminant is undefined also in case when ${}_*{}^*$-rank of the matrix is less then $n-1$.

$$b^{p}_{T}{}_*{}^*b^{-1_*{}^*}{}^{T}_{p}+b^{p}_{r}\ b^{-1_*{}^*}{}^{r}_{p}=1 \tag{5.5.3}$$

and will try to solve this system. We multiply (5.5.2) by $(b^S_T)^{-1_*{}^*}$

$$b^{-1_*{}^*}{}^{T}_{p}+(b^{S}_{T})^{-1_*{}^*}{}_*{}^*b^{S}_{r}\ b^{-1_*{}^*}{}^{r}_{p}=0 \tag{5.5.4}$$

Now we can substitute (5.5.4) into (5.5.3)

$$-b^{p}_{T}{}_*{}^*(b^{S}_{T})^{-1_*{}^*}{}_*{}^*b^{S}_{r}\ b^{-1_*{}^*}{}^{r}_{p}+b^{p}_{r}\ b^{-1_*{}^*}{}^{r}_{p}=1 \tag{5.5.5}$$

From (5.5.5) it follows that

$$(b^{p}_{r}-b^{p}_{T}{}_*{}^*(b^{S}_{T})^{-1_*{}^*}{}_*{}^*b^{S}_{r})\ b^{-1_*{}^*}{}^{r}_{p}=1 \tag{5.5.6}$$

Expression in brackets is quasideterminant $\det({}_*{}^*)^p_r b$. Substituting this expression into (5.5.6), we get

$$\det({}_*{}^*)^{p}_{r}b\ b^{-1_*{}^*}{}^{r}_{p}=1 \tag{5.5.7}$$

Thus we proved that quasideterminant $\det({}_*{}^*)^p_r b$ is defined and its equation to 0 is necessary and sufficient condition that the matrix $b$ is singular. Since (5.5.1) the statement of theorem is proved. □

THEOREM 5.5.4. *Let $a$ be a matrix,*

$$\mathrm{rank}_{{}_*{}^*}\,a=k<m$$

*and $a^S_T$ is ${}_*{}^*$-major minor matrix. Then row $a^p$ is a left linear composition of rows $a^S$*

$$a^{M\setminus S}=r_*{}^*a^{S} \tag{5.5.8}$$

$$a^{p}=r^{p}{}_*{}^*a^{S} \tag{5.5.9}$$

$$a^{p}_{b}=r^{p}_{s}a^{s}_{b} \tag{5.5.10}$$

PROOF. If the matrix $a$ has $k$ columns, then assuming that row $a^p$ is a left linear combination of rows (5.5.9) of rows $a_s$ with coefficients $R^p_s$ we get system of linear equations (5.5.10). According to theorem 5.4.7 the system of linear equations (5.5.10) has a unique solution[5.6] and this solution is nontrivial because all ${}_*{}^*$-quasideterminants are different from 0.

It remains to prove this statement in case when a number of columns of the matrix $a$ is more then $k$. I get row $a^p$ and column $a_r$. According to assumption, minor matrix $a^{S\cup\{p\}}_{T\cup\{r\}}$ is a ${}_*{}^*$-singular matrix and its ${}_*{}^*$-quasideterminant

$$\det({}_*{}^*)^{p}_{r}\ a^{S\cup\{p\}}_{T\cup\{r\}}=0 \tag{5.5.11}$$

According to (3.3.12) the equation (5.5.11) has form

$$a^{p}_{r}-a^{p}_{T}{}_*{}^*((a^{S}_{T})^{-1_*{}^*})_*{}^*a^{S}_{r}=0$$

Matrix

$$r^{p}=a^{p}_{T}{}_*{}^*((a^{S}_{T})^{-1_*{}^*}) \tag{5.5.12}$$

does not depend on $r$, Therefore, for any $r\in N\setminus T$

$$a^{p}_{r}=r^{p}{}_*{}^*a^{S}_{r} \tag{5.5.13}$$

From equation

$$((a^{S}_{T})^{-1_*{}^*})_*{}^*a^{S}_{l}=\delta^{T}_{l}$$

[5.6] We assume that unknown variables are $x_s=R^p_s$

it follows that

(5.5.14) $$a_l^p = a_T^p{}_*{}^*\delta_l^T = a_T^p{}_*{}^*((a_T^S)^{-1_*{}^*})_*{}^*a_l^S$$

Substituting (5.5.12) into (5.5.14) we get

(5.5.15) $$a_l^p = r^p{}_*{}^*a_l^S$$

(5.5.13) and (5.5.15) finish the proof. □

Corollary 5.5.5. *Let $a$ be a matrix,*

$$\mathrm{rank}_{*}{}^{*}\, a = k < m$$

*Then rows of the matrix are $D_*{}^*$-linearly dependent.*

$$\lambda_*{}^*a = 0$$

Proof. Suppose row $a^p$ is a left linear composition of rows (5.5.9). We assume $\lambda_p = -1$, $\lambda_s = R_s^p$ and the rest $\lambda_c = 0$. □

Theorem 5.5.6. *Let $(a^i, i \in M, |M| = m)$ be set of left linearly independent rows. Then ${}_*{}^*$-rank of their coordinate matrix equal $m$.*

Proof. Let $\overline{\overline{e}}$ be the basis of left vector space of rows. According to model built in example 5.3.2, the coordinate matrix of set of vectors $(\overline{a}^i)$ relative basis $\overline{\overline{e}}$ consists from rows which are coordinate matrices of vectors $\overline{a}^i$ relative the basis $\overline{\overline{e}}$. Therefore ${}_*{}^*$-rank of this matrix cannot be more then $m$.

Let ${}_*{}^*$-rank of the coordinate matrix be less then $m$. According to corollary 5.5.5, rows of matrix are left linear dependent

(5.5.16) $$\lambda_*{}^*a = 0$$

Suppose $c = \lambda_*{}^*A$. From the equation (5.5.16) it follows that left linear composition of rows

$$c_*{}^*\overline{e} = 0$$

of vectors of basis equal 0. This contradicts to statement that vectors $\overline{e}$ form basis. We proved statement of theorem. □

Theorem 5.5.7. *Let $a$ be a matrix,*

$$\mathrm{rank}_{*}{}^{*}\, a = k < n$$

*and $a_T^S$ is ${}_*{}^*$-major minor matrix. Then column $a_r$ is a ${}_*{}^*D$-linear composition of columns $a_t$*

(5.5.17) $$a_{N\setminus T} = a^T{}_*{}^*r$$

(5.5.18) $$a_r = a_T{}_*{}^*r_r$$

(5.5.19) $$a_r^a = a_t^a\ R_r^t$$

Proof. If the matrix $a$ has $k$ rows, then assuming that column $a_r$ is a ${}_*{}^*D$-linear combination (5.5.18) of columns $a_t$ with coefficients $R_r^t$ we get system of ${}_*{}^*D$-linear equations (5.5.19). According to theorem 5.4.7 the system of ${}_*{}^*D$-linear equations (5.5.19) has a unique solution[5.7] and this solution is nontrivial because all ${}_*{}^*$-quasideterminants are different from 0.

[5.7] We assume that unknown variables are ${}_tx = R_r^t$

It remains to prove this statement in case when a number of rows of the matrix $A$ is more then $k$. I get column $a_r$ and row $a^p$. According to assumption, minor matrix $a^{S\cup\{p\}}_{T\cup\{r\}}$ is a ${}_*{}^*$-singular matrix and its $\overline{RC}$-quasideterminant

$$\det({}_*{}^*)^p_r\, a^{S\cup\{p\}}_{T\cup\{r\}} = 0 \tag{5.5.20}$$

According to (3.3.12) (5.5.20) has form

$$a^p_r - a^p_T{}_*{}^*((a^S_T)^{-1_*{}^*})_*{}^*a^S_r = 0$$

Matrix

$$r_r = ((a^S_T)^{-1_*{}^*})_*{}^*a^S_r \tag{5.5.21}$$

does not depend on $p$, Therefore, for any $p \in M \setminus S$

$$a^p_r = a^p_T{}_*{}^*r_r \tag{5.5.22}$$

From equation

$$a^k_T{}_*{}^*((a^S_T)^{-1_*{}^*})_s = \delta^k_s$$

it follows that

$$a^k_r = \delta^k_S{}_*{}^*a^S_r = a^k_T{}_*{}^*((a^S_T)^{-1_*{}^*})^s{}_*{}^*a^S_r \tag{5.5.23}$$

Substituting (5.5.21) into (5.5.23) we get

$$a^k_r = a^k_T{}_*{}^*R_r \tag{5.5.24}$$

(5.5.22) and (5.5.24) finish the proof. □

Corollary 5.5.8. *Let $a$ be a matrix,*

$$\operatorname{rank}_{{}_*{}^*} a = k < m$$

*Then columns of the matrix are ${}_*{}^*D$-linearly dependent.*

$$a_*{}^*\lambda = 0$$

Proof. Suppose column $a_r$ is a right linear composition (5.5.18). We assume $\lambda^r = -1$, $\lambda^t = R^t_r$ and the rest $\lambda^c = 0$. □

Base on theorem 3.3.8 we can write similar statements for ${}^*{}_*$-rank of matrix.

Theorem 5.5.9. *Suppose $a$ is a matrix,*

$$\operatorname{rank}_{{}^*{}_*} a = k < m$$

*and $a^T_S$ is ${}^*{}_*$-major minor matrix. Then row $a^p$ is a right linear composition of rows $a^s$.*

$$a^{M\setminus S} = a^S{}^*{}_*r \tag{5.5.25}$$

$$a^p = a^S{}^*{}_*r^p \tag{5.5.26}$$

$$a^b_p = a^b_s r^s_p \tag{5.5.27}$$

Corollary 5.5.10. *Let $a$ be a matrix,*

$$\operatorname{rank}_{{}^*{}_*} a = k < m$$

*Then rows of matrix are right linearly dependent.*

$$A^*{}_*\lambda = 0$$

Theorem 5.5.11. *Let $a$ be a matrix,*

$$\mathrm{rank}_{*_*} a = \boldsymbol{k} < \boldsymbol{n}$$

*and $a_{\boldsymbol{S}}^{\boldsymbol{T}}$ is ${}^*{}_*$-major minor matrix. Then column $a_{\boldsymbol{r}}$ is a left linear composition of columns $a_{\boldsymbol{t}}$*

$$a_{\boldsymbol{N}\setminus\boldsymbol{T}} = R^*{}_*a_{\boldsymbol{T}} \tag{5.5.28}$$

$$a_{\boldsymbol{r}} = R_{\boldsymbol{r}}{}^*{}_*a_{\boldsymbol{T}} \tag{5.5.29}$$

$$a_{\boldsymbol{r}}^{\boldsymbol{i}} = R_{\boldsymbol{r}}^{\boldsymbol{t}}a_{\boldsymbol{t}}^{\boldsymbol{i}} \tag{5.5.30}$$

Corollary 5.5.12. *Let $a$ be a matrix,*

$$\mathrm{rank}_{*_*} a = \boldsymbol{k} < \boldsymbol{m}$$

*Then columns of matrix are left linearly dependent*

$$\lambda^*{}_*a = 0$$

## 5.6. Eigenvalue of Matrix

Let $a$ be $\boldsymbol{n} \times \boldsymbol{n}$ matrix of $A$-numbers and $E$ be $\boldsymbol{n} \times \boldsymbol{n}$ unit matrix.

Definition 5.6.1. *$A$-number $b$ is called ${}_*{}^*$**-eigenvalue** of the matrix $a$ if the matrix $a - bE$ is ${}_*{}^*$-singular matrix.*

*Column vector $v$ is called ${}_*{}^*$**-eigencolumn** of matrix $a$ corresponding to ${}_*{}^*$-eigenvalue $b$, if the following equality is true*

$$a_*{}^*v = bv$$

*Row vector $v$ is called ${}_*{}^*$**-eigenrow** of matrix $a$ corresponding to ${}_*{}^*$-eigenvalue $b$, if the following equality is true*

$$v_*{}^*a = vb$$

□

Definition 5.6.2. *$A$-number $b$ is called ${}^*{}_*$**-eigenvalue** of the matrix $a$ if the matrix $a - bE$ is ${}^*{}_*$-singular matrix.*

*Column vector $v$ is called ${}^*{}_*$**-eigencolumn** of matrix $a$ corresponding to ${}^*{}_*$-eigenvalue $b$, if the following equality is true*

$$v^*{}_*a = vb$$

*Row vector $v$ is called ${}^*{}_*$**-eigenrow** of matrix $a$ corresponding to ${}^*{}_*$-eigenvalue $b$, if the following equality is true*

$$a^*{}_*v = bv$$

□

Theorem 5.6.3. *Let $a$ be $\boldsymbol{n} \times \boldsymbol{n}$ matrix of $A$-numbers and*

$${}_*{}^*\,\mathrm{rank}\, a = \boldsymbol{k} < \boldsymbol{n}$$

*Then 0 is ${}_*{}^*$-eigenvalue of multiplicity $\boldsymbol{n} - \boldsymbol{k}$.*

Proof. The theorem follows from the equality

$$a = a - 0E$$

□

Theorem 5.6.4. *Let $a$ be $\boldsymbol{n}\times\boldsymbol{n}$ matrix of $A$-numbers and $a^{\boldsymbol{I}}_{\boldsymbol{J}}$ be ${}_*{}^*$-major minor matrix of the matrix $a$. Then non zero ${}_*{}^*$-eigenvalues of the matrix $a$ are roots of any equation*

$$\det({}_*{}^*)^{\boldsymbol{i}}_{\boldsymbol{j}}\,(a-bE)^{\boldsymbol{I}}_{\boldsymbol{J}}=0\qquad \boldsymbol{i}\in\boldsymbol{I},\ \boldsymbol{j}\in\boldsymbol{J}$$

Proof. The theorem follows from the theorem 5.5.3. □

CHAPTER 6

# Differential Equation Solved with Respect to the Derivative

## 6.1. Differential Equation $\frac{dy}{dx} = 1 \otimes x + x \otimes 1$

Definition 6.1.1. *Let $A$, $B$ be Banach $D$-modules. The map*

$$g : A \to \mathcal{L}(D; A \to B)$$

*is called* **integrable***, if there exists a map*

$$f : A \to B$$

*such that*

$$\frac{df(x)}{dx} = g(x)$$

*Then we use notation*

$$f(x) = \int g(x) \circ dx$$

*and the map $f$ is called* **indefinite integral** *of the map $g$.* □

It is convenient to use table of derivatives[6.1] and table of integrals[6.2] to solve differential equation

$$\frac{dy}{dx} = g(x)$$

if the map $g$ is simple.

Definition 6.1.2. *If there exist indefinite integral*

$$\int g(x) \circ dx$$

*then differential equation*

$$\frac{dy}{dx} = g(x)$$

*is called* **integrable**. □

Example 6.1.3. *Let $A$ be Banach algebra. According to the theorem A.2.7, the map*

$$y = x^2 + C$$

*is solution of the differential equation*

$$\frac{dy}{dx} = x \otimes 1 + 1 \otimes x \tag{6.1.1}$$

[6.1] In case of maps of real field, see for instance sections [20]-95, [2]-2.3. In case of maps of Banach algebra, see also section A.1.

[6.2] In case of maps of real field, see for instance sections [21]-265, [2]-4.4. In case of maps of Banach algebra, see also section A.2.

*with initial condition*

$$x_0=0\quad y_0=C$$

□

If $D$-algebra $B$ has finite basis, then we can write the differential equation

$$\frac{dy}{dx}=g(x)$$

as system of differential equations

$$\frac{\partial y^i}{\partial x^j}=g^i_j\quad y=y^ie_{B\cdot i}\quad x=x^ie_{A\cdot i}$$

In order that the calculations do not hide idea how to make transformation, we consider relatively simple example 6.1.3.

We start with the differential equation (6.1.1) in complex field. The theorems 6.1.4, 6.1.6 are simple. But they are necessary to prepare the reader for the proof of the theorem 6.1.10.

Theorem 6.1.4. *The differential equation* (6.1.1) *in complex field has form of system of differential equations*

$$\left\{\begin{array}{ll}\frac{\partial y^0}{\partial x^0}=2x^0 & \frac{\partial y^0}{\partial x^1}=-2x^1\\ \frac{\partial y^1}{\partial x^0}=2x^1 & \frac{\partial y^1}{\partial x^1}=2x^0\end{array}\right. \tag{6.1.2}$$

*with respect to the basis*

$$e_0=1\quad e_1=i$$

Proof. Since product in complex field is commutative, we can write differential equation (6.1.1) as follows

$$dy=2x\,dx \tag{6.1.3}$$

If we represent differentials $dx$, $dy$ as vector-column, then, according to the theorem [9]-2.5.1, the equation (6.1.3) gets following form

$$\begin{pmatrix}dy^0\\ dy^1\end{pmatrix}=\begin{pmatrix}2x^0 & -2x^1\\ 2x^1 & 2x^0\end{pmatrix}\begin{pmatrix}dx^0\\ dx^1\end{pmatrix} \tag{6.1.4}$$

Since the matrix of derivative has following form

$$\frac{dy}{dx}=\begin{pmatrix}\frac{\partial y^0}{\partial x^0} & \frac{\partial y^0}{\partial x^1}\\ \frac{\partial y^1}{\partial x^0} & \frac{\partial y^1}{\partial x^1}\end{pmatrix}$$

then the equality

$$\begin{pmatrix}\frac{\partial y^0}{\partial x^0} & \frac{\partial y^0}{\partial x^1}\\ \frac{\partial y^1}{\partial x^0} & \frac{\partial y^1}{\partial x^1}\end{pmatrix}=\begin{pmatrix}2x^0 & -2x^1\\ 2x^1 & 2x^0\end{pmatrix} \tag{6.1.5}$$

follows from (6.1.4). The system of differential equations (6.1.15) follows from the equality (6.1.5). □

LEMMA 6.1.5. *Let*

$$x = x^0 + x^1 i$$

*be complex number. Then*

$$x^2 = (x^0)^2 - (x^1)^2 + 2x^0 x^1 i \tag{6.1.6}$$

PROOF. The equality (6.1.6) follows from the definition of product in complex field. □

THEOREM 6.1.6. *The map*

$$y = x^2 + C$$

*is solution of the system of differential equations* (6.1.2) *in complex field with initial condition*

$$x = 0 \quad y^0 = C^0 \quad y^1 = C^1 \quad C = C^0 + C^1 i$$

PROOF. The map

$$y^0 = (x^0)^2 + C_1^0(x^1) \tag{6.1.7}$$

is solution of the differential equation

$$\frac{\partial y^0}{\partial x^0} = 2x^0$$

From (6.1.2), (6.1.7), it follows that

$$\frac{\partial y^0}{\partial x^1} = \frac{dC_1^0}{dx^1} = -2x^1 \tag{6.1.8}$$

The map

$$C_1^0(x^1) = -(x^1)^2 + C^0 \tag{6.1.9}$$

is solution of the differential equation (6.1.8). The equality

$$y^0 = (x^0)^2 - (x^1)^2 + C^0 \tag{6.1.10}$$

follows from equalities (6.1.7), (6.1.9).

The map

$$y^1 = 2x^0 x^1 + C_1^1(x^1) \tag{6.1.11}$$

is solution of the differential equation

$$\frac{\partial y^1}{\partial x^0} = 2x^1$$

From (6.1.2), (6.1.11), it follows that

$$\frac{\partial y^1}{\partial x^1} = 2x^0 + \frac{dC_1^1}{dx^1} = 2x^0 \tag{6.1.12}$$

The map

$$C_1^1(x^1) = C^1 \tag{6.1.13}$$

is solution of the differential equation (6.1.12). The equality

$$y^1 = 2x^0 x^1 + C^1 \tag{6.1.14}$$

follows from equalities (6.1.11), (6.1.13).

The theorem follows from the comparison of equality (6.1.6) and equalities (6.1.10), (6.1.14). □

We consider the differential equation (6.1.1) in quaternion algebra in the same way as we considered this differential equation in complex field.

Theorem 6.1.7. *The differential equation (6.1.1) in quaternion algebra has form of system of differential equations*

$$\left\{\begin{array}{llll}\frac{\partial y^0}{\partial x^0}=2x^0 & \frac{\partial y^0}{\partial x^1}=-2x^1 & \frac{\partial y^0}{\partial x^2}=-2x^2 & \frac{\partial y^0}{\partial x^3}=-2x^3\\ \frac{\partial y^1}{\partial x^0}=2x^1 & \frac{\partial y^1}{\partial x^1}=2x^0 & \frac{\partial y^1}{\partial x^2}=0 & \frac{\partial y^1}{\partial x^3}=0\\ \frac{\partial y^2}{\partial x^0}=2x^2 & \frac{\partial y^2}{\partial x^1}=0 & \frac{\partial y^2}{\partial x^2}=2x^0 & \frac{\partial y^2}{\partial x^3}=0\\ \frac{\partial y^3}{\partial x^0}=2x^3 & \frac{\partial y^3}{\partial x^1}=0 & \frac{\partial y^3}{\partial x^2}=0 & \frac{\partial y^3}{\partial x^3}=2x^0\end{array}\right. \tag{6.1.15}$$

*with respect to the basis*

$$e_0 = 1 \quad e_1 = i \quad e_2 = j \quad e_3 = k$$

Proof. We can write differential equation (6.1.1) as follows

$$dy = x\,dx + dx\,x \tag{6.1.16}$$

If we represent differentials $dx$, $dy$ as vector-column, then, according to theorems [11]-5.1, [11]-5.2, the equation (6.1.16) gets following form

$$\begin{aligned}\begin{pmatrix}dy^0\\dy^1\\dy^2\\dy^3\end{pmatrix} &= E_l(x)\begin{pmatrix}dx^0\\dx^1\\dx^2\\dx^3\end{pmatrix} + E_r(x)\begin{pmatrix}dx^0\\dx^1\\dx^2\\dx^3\end{pmatrix}\\ &= \begin{pmatrix}x^0 & -x^1 & -x^2 & -x^3\\ x^1 & x^0 & -x^3 & x^2\\ x^2 & x^3 & x^0 & -x^1\\ x^3 & -x^2 & x^1 & x^0\end{pmatrix}\begin{pmatrix}dx^0\\dx^1\\dx^2\\dx^3\end{pmatrix}\\ &+ \begin{pmatrix}x^0 & -x^1 & -x^2 & -x^3\\ x^1 & x^0 & x^3 & -x^2\\ x^2 & -x^3 & x^0 & x^1\\ x^3 & x^2 & -x^1 & x^0\end{pmatrix}\begin{pmatrix}dx^0\\dx^1\\dx^2\\dx^3\end{pmatrix}\\ &= \begin{pmatrix}2x^0 & -2x^1 & -2x^2 & -2x^3\\ 2x^1 & 2x^0 & 0 & 0\\ 2x^2 & 0 & 2x^0 & 0\\ 2x^3 & 0 & 0 & 2x^0\end{pmatrix}\begin{pmatrix}dx^0\\dx^1\\dx^2\\dx^3\end{pmatrix}\end{aligned} \tag{6.1.17}$$

Since the matrix of derivative has following form

$$\frac{dy}{dx}=\begin{pmatrix}\frac{\partial y^0}{\partial x^0}&\frac{\partial y^0}{\partial x^1}&\frac{\partial y^0}{\partial x^2}&\frac{\partial y^0}{\partial x^3}\\ \frac{\partial y^1}{\partial x^0}&\frac{\partial y^1}{\partial x^1}&\frac{\partial y^1}{\partial x^2}&\frac{\partial y^1}{\partial x^3}\\ \frac{\partial y^2}{\partial x^0}&\frac{\partial y^2}{\partial x^1}&\frac{\partial y^2}{\partial x^2}&\frac{\partial y^2}{\partial x^3}\\ \frac{\partial y^3}{\partial x^0}&\frac{\partial y^3}{\partial x^1}&\frac{\partial y^3}{\partial x^2}&\frac{\partial y^3}{\partial x^3}\end{pmatrix}$$

then the equality

$$\begin{pmatrix}\frac{\partial y^0}{\partial x^0}&\frac{\partial y^0}{\partial x^1}&\frac{\partial y^0}{\partial x^2}&\frac{\partial y^0}{\partial x^3}\\ \frac{\partial y^1}{\partial x^0}&\frac{\partial y^1}{\partial x^1}&\frac{\partial y^1}{\partial x^2}&\frac{\partial y^1}{\partial x^3}\\ \frac{\partial y^2}{\partial x^0}&\frac{\partial y^2}{\partial x^1}&\frac{\partial y^2}{\partial x^2}&\frac{\partial y^2}{\partial x^3}\\ \frac{\partial y^3}{\partial x^0}&\frac{\partial y^3}{\partial x^1}&\frac{\partial y^3}{\partial x^2}&\frac{\partial y^3}{\partial x^3}\end{pmatrix}=\begin{pmatrix}2x^0&-2x^1&-2x^2&-2x^3\\ 2x^1&2x^0&0&0\\ 2x^2&0&2x^0&0\\ 2x^3&0&0&2x^0\end{pmatrix}\tag{6.1.18}$$

follows from (6.1.17). The system of differential equations (6.1.15) follows from the equality (6.1.18). □

Lemma 6.1.8. *Let*

$$\begin{aligned}x&=x^0+x^1i+x^2j+x^3k\\ y&=y^0+y^1i+y^2j+y^3k\end{aligned}\tag{6.1.19}$$

*be quaternions. Then*

$$\begin{aligned}xy=&x^0y^0-x^1y^1-x^2y^2-x^3y^3+(x^0y^1+x^1y^0+x^2y^3-x^3y^2)i\\ &+(x^0y^2+x^2y^0+x^3y^1-x^1y^3)j+(x^0y^3+x^3y^0+x^1y^2-x^2y^1)k\end{aligned}\tag{6.1.20}$$

Proof. The lemma follows from the definition [7]-11.2.1. □

Lemma 6.1.9. *Let*

$$x=x^0+x^1i+x^2j+x^3k$$

*be quaternion. Then*

$$x^2=(x^0)^2-(x^1)^2-(x^2)^2-(x^3)^2+2x^0x^1i+2x^0x^2j+2x^0x^3k\tag{6.1.21}$$

Proof. The equality

$$\begin{aligned}x^2=&(x^0)^2-(x^1)^2-(x^2)^2-(x^3)^2+(x^0x^1+x^1x^0+x^2x^3-x^3x^2)i\\ &+(x^0x^2+x^2x^0+x^3x^1-x^1x^3)j+(x^0x^3+x^3x^0+x^1x^2-x^2x^1)k\end{aligned}\tag{6.1.22}$$

follows from the equality (6.1.20). The equality (6.1.21) follows from the equality (6.1.22). □

THEOREM 6.1.10. *The map*

$$y = x^2 + C$$

*is solution of the system of differential equations* (6.1.15) *with initial condition*

$$x = 0 \quad y^0 = C^0 \quad y^1 = C^1 \quad y^2 = C^2 \quad y^3 = C^3 \quad C = C^0 + C^1 i + C^2 j + C^3 k$$

PROOF. The map

$$y^0 = (x^0)^2 + C_1^0(x^1, x^2, x^3) \tag{6.1.23}$$

is solution of the differential equation

$$\frac{\partial y^0}{\partial x^0} = 2x^0$$

From (6.1.15), (6.1.23), it follows that

$$\frac{\partial y^0}{\partial x^1} = \frac{\partial C_1^0}{\partial x^1} = -2x^1 \tag{6.1.24}$$

The map

$$C_1^0(x^1, x^2, x^3) = -(x^1)^2 + C_2^0(x^2, x^3) \tag{6.1.25}$$

is solution of the differential equation (6.1.24). The equality

$$y^0 = (x^0)^2 - (x^1)^2 + C_2^0(x^2, x^3) \tag{6.1.26}$$

follows from equalities (6.1.23), (6.1.25). From (6.1.15), (6.1.26), it follows that

$$\frac{\partial y^0}{\partial x^2} = \frac{\partial C_2^0}{\partial x^2} = -2x^2 \tag{6.1.27}$$

The map

$$C_2^0(x^2, x^3) = -(x^2)^2 + C_3^0(x^3) \tag{6.1.28}$$

is solution of the differential equation (6.1.27). The equality

$$y^0 = (x^0)^2 - (x^1)^2 - (x^2)^2 + C_3^0(x^3) \tag{6.1.29}$$

follows from equalities (6.1.26), (6.1.28). From (6.1.15), (6.1.29), it follows that

$$\frac{\partial y^0}{\partial x^3} = \frac{dC_3^0}{dx^3} = -2x^3 \tag{6.1.30}$$

The map

$$C_3^0(x^3) = -(x^3)^2 + C^0 \tag{6.1.31}$$

is solution of the differential equation (6.1.30). The equality

$$y^0 = (x^0)^2 - (x^1)^2 - (x^2)^2 - (x^3)^2 + C^0 \tag{6.1.32}$$

follows from equalities (6.1.29), (6.1.31).

The map

$$y^1 = 2x^0 x^1 + C_1^1(x^1, x^2, x^3) \tag{6.1.33}$$

is solution of the differential equation

$$\frac{\partial y^1}{\partial x^0} = 2x^1$$

From (6.1.15), (6.1.33), it follows that

$$\frac{\partial y^1}{\partial x^1} = 2x^0 + \frac{\partial C_1^1}{\partial x^1} = 2x^0 \tag{6.1.34}$$

The map

$$C_1^1(x^1, x^2, x^3) = C_2^1(x^2, x^3) \tag{6.1.35}$$

is solution of the differential equation (6.1.34). The equality

$$y^1 = 2x^0x^1 + C_2^1(x^2, x^3) \tag{6.1.36}$$

follows from equalities (6.1.33), (6.1.35). From (6.1.15), (6.1.36), it follows that

$$\frac{\partial y^1}{\partial x^2} = \frac{\partial C_2^1}{\partial x^2} = 0 \quad \frac{\partial y^1}{\partial x^3} = \frac{\partial C_2^1}{\partial x^3} = 0 \tag{6.1.37}$$

The map

$$C_2^1(x^2, x^3) = C^1 \tag{6.1.38}$$

is solution of the system of differential equations (6.1.37). The equality

$$y^1 = 2x^0x^1 + C^1 \tag{6.1.39}$$

follows from equalities (6.1.36), (6.1.38).

The similar way we see that

$$y^2 = 2x^0x^2 + C^2 \tag{6.1.40}$$

$$y^3 = 2x^0x^3 + C^3 \tag{6.1.41}$$

The theorem follows from the comparison of equality (6.1.21) and equalities (6.1.32), (6.1.39), (6.1.40), (6.1.41). □

## 6.2. Differential Equation $\frac{dy}{dx} = 3x \otimes x$

In this section I consider the differential equation in which existence of a solution depends on choice of algebra.

Theorem 6.2.1. *The differential equation*

$$\frac{dy}{dx} = 3x \otimes x \tag{6.2.1}$$

*in complex field has form of system of differential equations*

$$\left\{\begin{array}{ll} \frac{\partial y^0}{\partial x^0} = 3(x^0)^2 - 3(x^1)^2 & \frac{\partial y^0}{\partial x^1} = -6x^0x^1 \\ \frac{\partial y^1}{\partial x^0} = 6x^0x^1 & \frac{\partial y^1}{\partial x^1} = 3(x^0)^2 - 3(x^1)^2 \end{array}\right. \tag{6.2.2}$$

*with respect to the basis*

$$e_0 = 1 \quad e_1 = i$$

Proof. Since product in complex field is commutative, we can write differential equation (6.2.1) as follows

$$dy = 3x^2\,dx \tag{6.2.3}$$

If we represent differentials $dx$, $dy$ as vector-column, then, according to the theorem [9]-2.5.1 and lemma 6.1.5, the equation (6.2.3) gets following form

$$\begin{pmatrix} dy^0 \\ dy^1 \end{pmatrix} = 3 \begin{pmatrix} (x^0)^2 - (x^1)^2 & -2x^0x^1 \\ 2x^0x^1 & (x^0)^2 - (x^1)^2 \end{pmatrix} \begin{pmatrix} dx^0 \\ dx^1 \end{pmatrix} \tag{6.2.4}$$

Since the matrix of derivative has following form

$$\frac{dy}{dx} = \begin{pmatrix} \frac{\partial y^0}{\partial x^0} & \frac{\partial y^0}{\partial x^1} \\ \frac{\partial y^1}{\partial x^0} & \frac{\partial y^1}{\partial x^1} \end{pmatrix}$$

then the equality

$$\begin{pmatrix} \frac{\partial y^0}{\partial x^0} & \frac{\partial y^0}{\partial x^1} \\ \frac{\partial y^1}{\partial x^0} & \frac{\partial y^1}{\partial x^1} \end{pmatrix} = \begin{pmatrix} 3(x^0)^2 - 3(x^1)^2 & -6x^0x^1 \\ 6x^0x^1 & 3(x^0)^2 - 3(x^1)^2 \end{pmatrix} \tag{6.2.5}$$

follows from (6.2.4). The system of differential equations (6.2.2) follows from the equality (6.2.5). □

Lemma 6.2.2. *Let*

$$x = x^0 + x^1 i$$

*be complex number. Then*

$$x^3 = (x^0)^3 - 3x^0(x^1)^2 + (3(x^0)^2x^1 - (x^1)^3)i \tag{6.2.6}$$

Proof. According to the lemma 6.1.5,

$$\begin{aligned} x^3 = x^2x &= ((x^0)^2 - (x^1)^2 + 2x^0x^1 i)(x^0 + x^1 i) \\ &= ((x^0)^2 - (x^1)^2)x^0 + 2(x^0)^2x^1 i + ((x^0)^2 - (x^1)^2)x^1 i - 2x^0(x^1)^2 \\ &= (x^0)^3 - x^0(x^1)^2 + 2(x^0)^2x^1 i + ((x^0)^2x^1 - (x^1)^3)i - 2x^0(x^1)^2 \end{aligned} \tag{6.2.7}$$

The equality (6.2.6) follows from the equality (6.2.7) □

Theorem 6.2.3. *The map*

$$y = x^3 + C$$

*is solution of the system of differential equations* (6.2.2) *in complex field with initial condition*

$$x = 0 \quad y^0 = C^0 \quad y^1 = C^1 \quad C = C^0 + C^1 i$$

Proof. The map

$$y^0 = (x^0)^3 - 3x^0(x^1)^2 + C_1^0(x^1) \tag{6.2.8}$$

is solution of the differential equation

$$\frac{\partial y^0}{\partial x^0} = 3(x^0)^2 - 3(x^1)^2$$

From (6.2.2), (6.2.8), it follows that

$$\frac{\partial y^0}{\partial x^1} = -6x^0x^1 + \frac{dC_1^0}{dx^1} = -6x^0x^1 \tag{6.2.9}$$

The map

$$C_1^0(x^1) = C^0 \tag{6.2.10}$$

is solution of the differential equation (6.2.9). The equality

$$y^0 = (x^0)^3 - 3x^0(x^1)^2 + C^0 \tag{6.2.11}$$

follows from equalities (6.2.8), (6.2.10).

The map

$$y^{1} = 3(x^{0})^2x^{1} + C_1^{1}(x^{1}) \tag{6.2.12}$$

is solution of the differential equation

$$\frac{\partial y^{1}}{\partial x^{0}} = 6x^{0}x^{1}$$

From (6.2.2), (6.2.12), it follows that

$$\frac{\partial y^{1}}{\partial x^{1}} = 3(x^{0})^2 + \frac{dC_1^{1}}{dx^{1}} = 3(x^{0})^2 - 3(x^{1})^2 \tag{6.2.13}$$

The map

$$C_1^{1}(x^{1}) = -(x^{1})^3 + C^{1} \tag{6.2.14}$$

is solution of the differential equation (6.2.13). The equality

$$y^{1} = 3(x^{0})^2x^{1} - (x^{1})^3 + C^{1} \tag{6.2.15}$$

follows from equalities (6.2.12), (6.2.14).

The theorem follows from the comparison of equality (6.2.6) and equalities (6.2.11), (6.2.15). □

THEOREM 6.2.4. *The differential equation*

$$\frac{dy}{dx} = 3x \otimes x \tag{6.2.16}$$

*in non commutative algebra does not possess a solution.*

PROOF. We consider **method of successive differentiation** to solve differential equation (6.2.16). Differentiating one after another equation (6.2.16), we get the chain of equations

$$\frac{d^2y}{dx^2} = 3 \otimes 1 \otimes x + 3x \otimes 1 \otimes 1 \tag{6.2.17}$$

$$\frac{d^3y}{dx^3} = 6 \otimes 1 \otimes 1 \otimes 1 \tag{6.2.18}$$

The expansion into Taylor series

$$y = x^3 + C$$

follows from equations (6.2.16), (6.2.17), (6.2.18) and from initial condition. According to the theorem A.1.9

$$\frac{dx^3}{dx} \neq 3x \otimes x \tag{6.2.19}$$

The theorem follows from the statement (6.2.19). □

LEMMA 6.2.5. *In quaternion algebra*

$$
\begin{array}{rl}
& J_l(a)J_r(b) = J_r(b)J_l(a) \\
= & \left(\begin{array}{c:c:c:c}
a^0b^0 - a^1b^1 & -a^0b^1 - a^1b^0 & -a^0b^2 - a^1b^3 & -a^0b^3 + a^1b^2 \\
-a^2b^2 - a^3b^3 & +a^2b^3 - a^3b^2 & -a^2b^0 + a^3b^1 & -a^2b^1 - a^3b^0 \\ \hdashline
a^1b^0 + a^0b^1 & -a^1b^1 + a^0b^0 & -a^1b^2 + a^0b^3 & -a^1b^3 - a^0b^2 \\
-a^3b^2 + a^2b^3 & +a^3b^3 + a^2b^2 & -a^3b^0 - a^2b^1 & -a^3b^1 + a^2b^0 \\ \hdashline
a^2b^0 + a^3b^1 & -a^2b^1 + a^3b^0 & -a^2b^2 + a^3b^3 & -a^2b^3 - a^3b^2 \\
+a^0b^2 - a^1b^3 & -a^0b^3 - a^1b^2 & +a^0b^0 + a^1b^1 & +a^0b^1 - a^1b^0 \\ \hdashline
a^3b^0 - a^2b^1 & -a^3b^1 - a^2b^0 & -a^3b^2 - a^2b^3 & -a^3b^3 + a^2b^2 \\
+a^1b^2 + a^0b^3 & -a^1b^3 + a^0b^2 & +a^1b^0 - a^0b^1 & +a^1b^1 + a^0b^0
\end{array}\right)
\end{array}
\tag{6.2.20}
$$

PROOF. The product of matrices $J_l(a)$, $J_r(b)$ has form

$$
\begin{array}{rl}
& \begin{pmatrix}
a^0 & -a^1 & -a^2 & -a^3 \\
a^1 & a^0 & -a^3 & a^2 \\
a^2 & a^3 & a^0 & -a^1 \\
a^3 & -a^2 & a^1 & a^0
\end{pmatrix}
\begin{pmatrix}
b^0 & -b^1 & -b^2 & -b^3 \\
b^1 & b^0 & b^3 & -b^2 \\
b^2 & -b^3 & b^0 & b^1 \\
b^3 & b^2 & -b^1 & b^0
\end{pmatrix} \\
= & \left(\begin{array}{c:c:c:c}
a^0b^0 - a^1b^1 & -a^0b^1 - a^1b^0 & -a^0b^2 - a^1b^3 & -a^0b^3 + a^1b^2 \\
-a^2b^2 - a^3b^3 & +a^2b^3 - a^3b^2 & -a^2b^0 + a^3b^1 & -a^2b^1 - a^3b^0 \\ \hdashline
a^1b^0 + a^0b^1 & -a^1b^1 + a^0b^0 & -a^1b^2 + a^0b^3 & -a^1b^3 - a^0b^2 \\
-a^3b^2 + a^2b^3 & +a^3b^3 + a^2b^2 & -a^3b^0 - a^2b^1 & -a^3b^1 + a^2b^0 \\ \hdashline
a^2b^0 + a^3b^1 & -a^2b^1 + a^3b^0 & -a^2b^2 + a^3b^3 & -a^2b^3 - a^3b^2 \\
+a^0b^2 - a^1b^3 & -a^0b^3 - a^1b^2 & +a^0b^0 + a^1b^1 & +a^0b^1 - a^1b^0 \\ \hdashline
a^3b^0 - a^2b^1 & -a^3b^1 - a^2b^0 & -a^3b^2 - a^2b^3 & -a^3b^3 + a^2b^2 \\
+a^1b^2 + a^0b^3 & -a^1b^3 + a^0b^2 & +a^1b^0 - a^0b^1 & +a^1b^1 + a^0b^0
\end{array}\right)
\end{array}
\tag{6.2.21}
$$

The product of matrices $J_r(b)$, $J_l(a)$ has form

$$
\begin{aligned}
&\begin{pmatrix}
b^0 & -b^1 & -b^2 & -b^3\\
b^1 & b^0 & b^3 & -b^2\\
b^2 & -b^3 & b^0 & b^1\\
b^3 & b^2 & -b^1 & b^0
\end{pmatrix}
\begin{pmatrix}
a^0 & -a^1 & -a^2 & -a^3\\
a^1 & a^0 & -a^3 & a^2\\
a^2 & a^3 & a^0 & -a^1\\
a^3 & -a^2 & a^1 & a^0
\end{pmatrix}\\
&=\left(\begin{array}{c|c|c|c}
\begin{array}{r} a^0b^0-a^1b^1\\ -a^2b^2-a^3b^3\end{array} &
\begin{array}{r} -a^0b^1-a^1b^0\\ +a^2b^3-a^3b^2\end{array} &
\begin{array}{r} -a^0b^2-a^1b^3\\ -a^2b^0+a^3b^1\end{array} &
\begin{array}{r} -a^0b^3+a^1b^2\\ -a^2b^1-a^3b^0\end{array}\\
\hline
\begin{array}{r} a^1b^0+a^0b^1\\ -a^3b^2+a^2b^3\end{array} &
\begin{array}{r} -a^1b^1+a^0b^0\\ +a^3b^3+a^2b^2\end{array} &
\begin{array}{r} -a^1b^2+a^0b^3\\ -a^3b^0-a^2b^1\end{array} &
\begin{array}{r} -a^1b^3-a^0b^2\\ -a^3b^1+a^2b^0\end{array}\\
\hline
\begin{array}{r} a^2b^0+a^3b^1\\ +a^0b^2-a^1b^3\end{array} &
\begin{array}{r} -a^2b^1+a^3b^0\\ -a^0b^3-a^1b^2\end{array} &
\begin{array}{r} -a^2b^2+a^3b^3\\ +a^0b^0+a^1b^1\end{array} &
\begin{array}{r} -a^2b^3-a^3b^2\\ +a^0b^1-a^1b^0\end{array}\\
\hline
\begin{array}{r} a^3b^0-a^2b^1\\ +a^1b^2+a^0b^3\end{array} &
\begin{array}{r} -a^3b^1-a^2b^0\\ -a^1b^3+a^0b^2\end{array} &
\begin{array}{r} -a^3b^2-a^2b^3\\ +a^1b^0-a^0b^1\end{array} &
\begin{array}{r} -a^3b^3+a^2b^2\\ +a^1b^1+a^0b^0\end{array}
\end{array}\right)
\end{aligned}
\tag{6.2.22}
$$

The equality (6.2.20) follows from equalities (6.2.21), (6.2.22). □

Lemma 6.2.6. *The product of matrices* $J_r(a)$, $J_l(a)$ *has form*

$$
\left(\begin{array}{c|c|c|c}
\begin{array}{r} a^0a^0-a^1a^1\\ -a^2a^2-a^3a^3\end{array} & -2a^0a^1 & -2a^0a^2 & -2a^0a^3\\
\hline
2a^0a^1 & \begin{array}{r} -a^1a^1+a^0a^0\\ +a^3a^3+a^2a^2\end{array} & -2a^1a^2 & -2a^1a^3\\
\hline
2a^2a^0 & -2a^2a^1 & \begin{array}{r} -a^2a^2+a^3a^3\\ +a^0a^0+a^1a^1\end{array} & -2a^2a^3\\
\hline
2a^3a^0 & -2a^3a^1 & -2a^2a^3 & \begin{array}{r} -a^3a^3+a^2a^2\\ +a^1a^1+a^0a^0\end{array}
\end{array}\right)
\tag{6.2.23}
$$

Proof. According to the lemma 6.2.5, the product of matrices $J_r(a)$, $J_l(a)$ has form

$$
\left(\begin{array}{c|c|c|c}
\begin{matrix} a^0a^0 - a^1a^1 \\ -a^2a^2 - a^3a^3 \end{matrix} & \begin{matrix} -a^0a^1 - a^1a^0 \\ +a^2a^3 - a^3a^2 \end{matrix} & \begin{matrix} -a^0a^2 - a^1a^3 \\ -a^2a^0 + a^3a^1 \end{matrix} & \begin{matrix} -a^0a^3 + a^1a^2 \\ -a^2a^1 - a^3a^0 \end{matrix} \\ \hline
\begin{matrix} a^1a^0 + a^0a^1 \\ -a^3a^2 + a^2a^3 \end{matrix} & \begin{matrix} -a^1a^1 + a^0a^0 \\ +a^3a^3 + a^2a^2 \end{matrix} & \begin{matrix} -a^1a^2 + a^0a^3 \\ -a^3a^0 - a^2a^1 \end{matrix} & \begin{matrix} -a^1a^3 - a^0a^2 \\ -a^3a^1 + a^2a^0 \end{matrix} \\ \hline
\begin{matrix} a^2a^0 + a^3a^1 \\ +a^0a^2 - a^1a^3 \end{matrix} & \begin{matrix} -a^2a^1 + a^3a^0 \\ -a^0a^3 - a^1a^2 \end{matrix} & \begin{matrix} -a^2a^2 + a^3a^3 \\ +a^0a^0 + a^1a^1 \end{matrix} & \begin{matrix} -a^2a^3 - a^3a^2 \\ +a^0a^1 - a^1a^0 \end{matrix} \\ \hline
\begin{matrix} a^3a^0 - a^2a^1 \\ +a^1a^2 + a^0a^3 \end{matrix} & \begin{matrix} -a^3a^1 - a^2a^0 \\ -a^1a^3 + a^0a^2 \end{matrix} & \begin{matrix} -a^3a^2 - a^2a^3 \\ +a^1a^0 - a^0a^1 \end{matrix} & \begin{matrix} -a^3a^3 + a^2a^2 \\ +a^1a^1 + a^0a^0 \end{matrix}
\end{array}\right) \tag{6.2.24}
$$

The equality (6.2.23) follows from the equality (6.2.24). □

Theorem 6.2.7. *The differential equation (6.2.16) in quaternion algebra has form of system of differential equations*

$$
\left\{\begin{aligned}
\frac{\partial y^0}{\partial x^0} &= 3(x^0)^2 - 3(x^1)^2 - 3(x^2)^2 - 3(x^3)^2 \\
\frac{\partial y^0}{\partial x^1} &= -6x^0x^1 \\
\frac{\partial y^0}{\partial x^2} &= -6x^0x^2 \\
\frac{\partial y^0}{\partial x^3} &= -6x^0x^3
\end{aligned}\right. \tag{6.2.25}
$$

$$
\left\{\begin{aligned}
\frac{\partial y^1}{\partial x^0} &= 6x^0x^1 \\
\frac{\partial y^1}{\partial x^1} &= 3(x^0)^2 - 3(x^1)^2 + 3(x^2)^2 + 3(x^3)^2 \\
\frac{\partial y^1}{\partial x^2} &= -6x^1x^2 \\
\frac{\partial y^1}{\partial x^3} &= -6x^1x^3
\end{aligned}\right. \tag{6.2.26}
$$

$$
\left\{\begin{aligned}
\frac{\partial y^2}{\partial x^0} &= 6x^0x^2 \\
\frac{\partial y^2}{\partial x^1} &= -6x^1x^2 \\
\frac{\partial y^2}{\partial x^2} &= -3(x^2)^2 + 3(x^3)^2 + 3(x^0)^2 + 3(x^1)^2 \\
\frac{\partial y^2}{\partial x^3} &= -6x^2x^3
\end{aligned}\right. \tag{6.2.27}
$$

(6.2.28)
$$\begin{cases}\frac{\partial y^3}{\partial x^0}=6x^0x^3\\ \frac{\partial y^3}{\partial x^1}=-6x^1x^3\\ \frac{\partial y^3}{\partial x^2}=-6x^2x^3\\ \frac{\partial y^3}{\partial x^3}=-3(x^3)^2+3(x^2)^2+3(x^1)^2+3(x^0)^2\end{cases}$$

*with respect to the basis*

$$e_0=1\quad e_1=i\quad e_2=j\quad e_3=k$$

Proof. We can write differential equation (6.2.16) as follows

(6.2.29)
$$dy=3x\,dx\,x$$

If we represent differentials $dx$, $dy$ as vector-column, then, according to the lemma 6.2.5, the equation (6.2.29) gets following form

(6.2.30)
$$\begin{pmatrix}dy^0\\ dy^1\\ dy^2\\ dy^3\end{pmatrix}=3J_l(x)J_r(x)\begin{pmatrix}dx^0\\ dx^1\\ dx^2\\ dx^3\end{pmatrix}$$

Since the matrix of derivative has following form

$$\frac{dy}{dx}=\begin{pmatrix}\frac{\partial y^0}{\partial x^0}&\frac{\partial y^0}{\partial x^1}&\frac{\partial y^0}{\partial x^2}&\frac{\partial y^0}{\partial x^3}\\ \frac{\partial y^1}{\partial x^0}&\frac{\partial y^1}{\partial x^1}&\frac{\partial y^1}{\partial x^2}&\frac{\partial y^1}{\partial x^3}\\ \frac{\partial y^2}{\partial x^0}&\frac{\partial y^2}{\partial x^1}&\frac{\partial y^2}{\partial x^2}&\frac{\partial y^2}{\partial x^3}\\ \frac{\partial y^3}{\partial x^0}&\frac{\partial y^3}{\partial x^1}&\frac{\partial y^3}{\partial x^2}&\frac{\partial y^3}{\partial x^3}\end{pmatrix}$$

then the equality

(6.2.31)

$$\frac{1}{3}\begin{pmatrix}\frac{\partial y^0}{\partial x^0} & \frac{\partial y^0}{\partial x^1} & \frac{\partial y^0}{\partial x^2} & \frac{\partial y^0}{\partial x^3}\\ \frac{\partial y^1}{\partial x^0} & \frac{\partial y^1}{\partial x^1} & \frac{\partial y^1}{\partial x^2} & \frac{\partial y^1}{\partial x^3}\\ \frac{\partial y^2}{\partial x^0} & \frac{\partial y^2}{\partial x^1} & \frac{\partial y^2}{\partial x^2} & \frac{\partial y^2}{\partial x^3}\\ \frac{\partial y^3}{\partial x^0} & \frac{\partial y^3}{\partial x^1} & \frac{\partial y^3}{\partial x^2} & \frac{\partial y^3}{\partial x^3}\end{pmatrix} =$$

$$\begin{pmatrix}x^0x^0 - x^1x^1 - x^2x^2 - x^3x^3 & -2x^0x^1 & -2x^0x^2 & -2x^0x^3\\ 2x^0x^1 & -x^1x^1 + x^0x^0 + x^3x^3 + x^2x^2 & -2x^1x^2 & -2x^1x^3\\ 2x^2x^0 & -2x^2x^1 & -x^2x^2 + x^3x^3 + x^0x^0 + x^1x^1 & -2x^2x^3\\ 2x^3x^0 & -2x^3x^1 & -2x^2x^3 & -x^3x^3 + x^2x^2 + x^1x^1 + x^0x^0\end{pmatrix}$$

follows from (6.2.30). The system of differential equations (6.2.25), (6.2.26), (6.2.27), (6.2.28) follows from the equality (6.2.31). □

THEOREM 6.2.8. *The system of differential equations* (6.2.25), (6.2.26), (6.2.27), (6.2.28) *does not possess a solution.*

PROOF. The map

(6.2.32) $$y^1 = 3(x^0)^2x^1 + C_1^1(x^1, x^2, x^3)$$

is solution of the differential equation

$$\frac{\partial y^1}{\partial x^0} = 6x^0x^1$$

From (6.2.26), (6.2.32), it follows that

(6.2.33) $$\frac{\partial y^1}{\partial x^1} = 3(x^0)^2 + \frac{\partial C_1^1}{\partial x^1} = 3(x^0)^2 - 3(x^1)^2 + 3(x^2)^2 + 3(x^3)^2$$

The map

(6.2.34) $$C_1^1(x^1, x^2, x^3) = -(x^1)^3 + 3(x^2)^2x^1 + 3(x^3)^2x^1 + C_2^1(x^2, x^3)$$

is solution of the differential equation (6.2.33). The equality

(6.2.35) $$y^1 = 3(x^0)^2x^1 - (x^1)^3 + 3(x^2)^2x^1 + 3(x^3)^2x^1 + C_2^1(x^2, x^3)$$

follows from equalities (6.2.32), (6.2.34). From (6.2.26), (6.2.35), it follows that

(6.2.36) $$\frac{\partial y^1}{\partial x^2} = 6x^2x^1 + \frac{\partial C_2^1}{\partial x^2} = -6x^1x^2$$

From the equation (6.2.36), it follows that the map $C_2^1$ depends on $x^1$; this contradicts to the statement that the map $C_2^1$ does not depend on $x^1$. Therefore, the system of differential equations (6.2.25), (6.2.26), (6.2.27), (6.2.28) does not possess a solution. □

## 6.3. Before Going Any Further

Before going any further, we ask questions that need to be answered.

QUESTION 6.3.1. *If $D$-algebra $B$ has finite basis $\overline{\overline{e}}$, then we can write the differential equation*

$$\frac{df(x)}{dx} = g(x) \tag{6.3.1}$$

*as system of differential equations*

$$\frac{\partial y^i}{\partial x^j} = g^i_j \quad y = y^i e_{B\cdot i} \quad x = x^i e_{A\cdot i} \tag{6.3.2}$$

*What is the relationship between the differential equation (6.3.1) and system of differential equations (6.3.2)?* □

QUESTION 6.3.2. *What is the condition for the integrability of the differential equation (6.3.1)?* □

QUESTION 6.3.3. *Can we write down a differential equation (6.3.1) given the system of differential equations (6.3.2)?* □

QUESTION 6.3.4. *Systems of differential equations (6.1.2), (6.1.15), (6.2.2) are completely integrable.*[6.3] *The system of differential equations (6.2.25), (6.2.26), (6.2.27), (6.2.28) is not completely integrable. Is the requirement complete integrability of systems of differential equations (6.3.2) equivalent to the requirement of integrability of corresponding differential equation (6.3.1)?* □

## 6.4. Forms of representation of differential equations

The theorem 6.4.1 answers the question 6.3.1.

THEOREM 6.4.1. *Let $B$ be free finite dimensional associative $D$-algebra. Let $\overline{\overline{e}}_A$ be basis of $D$ module $A$. Let $\overline{\overline{e}}_B$ be basis of $D$ module $B$. Let $\overline{\overline{F}}$ be the basis of left $B\otimes B$-module $\mathcal{L}(D;B\to B)$ and*

$$G : A \to B$$

*be linear map of maximal rank such that $\ker G \subseteq \ker g$. Let $g^{k\cdot ij}$ be standard components of the map $g$. Let $C^p_{kl}$ be structural constants of algebra $B$. Then we can write the differential equation (6.3.1) as system of differential equations*

$$\frac{\partial y^k}{\partial x^l} = g^{k\cdot ij} F_{k\cdot r}{}^{m} G^r_l C^p_{im} C^k_{pj} \tag{6.4.1}$$

PROOF. The theorem follows from the theorem [13]-6.4.5. □

Based on the theorem 6.4.1, we consider the question 6.3.3. However, the answer to this question is not straightforward. We begin with the theorem 6.4.2.

[6.3] See the definition of completely integrable system of differential equations on page 2 in [19].

THEOREM 6.4.2. *Let $B$ be free finite dimensional associative $D$-algebra. Let $\overline{\overline{e}}$ be basis of $D$ module $B$. Let $\overline{\overline{F}}$ be the basis of left $B\otimes B$-module $\mathcal{L}(D;B\rightarrow B)$ and*

$$G:A\rightarrow B$$

*be linear map of maximal rank such that $\ker G\subseteq\ker g$. Let $C^{p}_{kl}$ be structural constants of algebra $B$. Consider matrix*

$$\mathcal{C}=\left(\mathcal{C}^{\cdot}{}^{k}_{m}{}^{\cdot}{}_{ij}\right)=\left(C^{p}_{im}C^{k}_{pj}\right)$$

*whose rows and columns are indexed by ${}^{\cdot}{}^{k}_{m}$ and ${}_{\cdot ij}$, respectively. If matrix $\mathcal{C}$ is nonsingular, then we can write the system of differential equations*

$$\frac{\partial y^{i}}{\partial x^{j}}=g^{i}_{j}\qquad y=y^{i}e_{B\cdot i}\qquad x=x^{i}e_{A\cdot i}$$

*as differential equation*

$$\frac{dy}{dx}=g^{k\cdot ij}(e_{B\cdot i}\otimes e_{B\cdot j})\circ F_{k}\circ G \tag{6.4.2}$$

*where standard components $g^{k\cdot ij}$ of map $g$ are solution of system of linear equations*

$$g^{k}_{l}=g^{k\cdot ij}F_{k\cdot}{}^{m}_{r}G^{r}_{l}C_{B\cdot}{}^{p}_{im}C_{B\cdot}{}^{k}_{pj}$$

*If matrix $\mathcal{C}$ is singular, then we can write the differential equation (6.4.2), if condition*

$$\operatorname{rank}\left(\mathcal{C}^{\cdot}{}^{k}_{m}{}^{\cdot}{}_{ij}\quad g^{k}_{m}\right)=\operatorname{rank}\mathcal{C}$$

*is satisfied.*

PROOF. The theorem follows from the theorem [13]-6.4.5. □

Analysis of the theorem 6.4.2 shows that it is not sufficient to write the equation (6.4.2) to answer the question 6.3.3. The reason for this is that, according to the construction, maps $g^{k\cdot ij}$ depend on coordinates of $A$-number $x$. However, it should be noted that the map

$$x=x^{i}e_{i}\rightarrow x^{j}$$

is the linear map. So, if $D$-module $A$ is $D$-algebra, then there exists effective procedure to answer the question 6.3.3.

REMARK 6.4.3. *The statement of this remark is prelimenary, because I did not finished research in this direction. Standard components of the map $g$ is universal form of representation of the map, but not the only and not the most expressive. Let $\overline{\overline{E}}$ be basis of $B$-module $\mathcal{L}(D;B\rightarrow B)$. Then the map $g$ has representation*

$$g(x)=g^{k}(x)\circ E_{k}\qquad g^{k}:A\rightarrow B$$

□

## 6.5. Condition of Integrability

To answer the question 6.3.2, I recall that the map $g$ is differential form (the definition [15]-7.3.2).

THEOREM 6.5.1. *The differential equation*

$$\frac{df(x)}{dx} = g(x)$$

*is integrable iff*

$$dg = 0$$

PROOF. The theorem follows from the theorem [15]-8.3.1. □

Consider first how the theorem 6.5.1 works in case of differential equations (6.1.1), (6.2.16).

EXAMPLE 6.5.2. *Consider the differential form*

$$g(x) = x \otimes 1 + 1 \otimes x$$

*According to theorems A.1.1, A.1.2, A.1.4, A.1.11,*

$$\frac{dg}{dx} = 1 \otimes_2 1 \otimes_1 1 + 1 \otimes_1 1 \otimes_2 1$$

*Index accompanying symbol* $\otimes$, *shows which argument of polylinear map should be written instead of corresponding symbol* $\otimes$. *For instance*

$$\frac{dg}{dx} \circ (h_1, h_2) = (1 \otimes_2 1 \otimes_1 1 + 1 \otimes_1 1 \otimes_2 1) \circ (h_1, h_2) = h_2 h_1 + h_1 h_2 \tag{6.5.1}$$

*From the equality (6.5.1), it follows that bilinear map* $\frac{dg}{dx}$ *is symmetrical. According to the definition [15]-7.4.1*

$$dg = 0$$

*According to the theorem 6.5.1, the differential equation (6.1.1) is integrable (the theorem 6.1.7).* □

EXAMPLE 6.5.3. *Consider the differential form*

$$g(x) = x \otimes x$$

*According to theorems A.1.2, A.1.11,*

$$\frac{dg}{dx} = 1 \otimes_2 1 \otimes_1 x + x \otimes_1 1 \otimes_2 1$$

*Index accompanying symbol* $\otimes$, *shows which argument of polylinear map should be written instead of corresponding symbol* $\otimes$. *For instance*

$$\frac{dg}{dx} \circ (h_1, h_2) = (1 \otimes_2 1 \otimes_1 x + x \otimes_1 1 \otimes_2 1) \circ (h_1, h_2) = h_2 h_1 x + x h_1 h_2 \tag{6.5.2}$$

*From the equality (6.5.2) and from the definition [15]-7.4.1, it follows that*

$$dg = (h_2 h_1 - h_1 h_2) x + x(h_1 h_2 - h_2 h_1) \tag{6.5.3}$$

*From the equality (6.5.3), it follows that*

$$dg \neq 0$$

*According to the theorem 6.5.1, the differential equation (6.2.16). is not integrable (the theorem 6.2.4).* □

The theorem 6.5.4 answers the question 6.3.2.

THEOREM 6.5.4. *Let $B$ be free finite dimensional associative $D$-algebra. Let $\overline{\overline{e}}$ be basis of $D$ module $B$. Let $\overline{\overline{\overline{F}}}$ be the basis of left $B\otimes B$-module $\mathcal{L}(D;B\to B)$ and $G$ be the linear map*

$$G:A\to B$$

*of maximal rank such that* $\ker G\subseteq\ker g$. *Let* $g^{k\cdot ij}$ *be standard components of the map*

$$g:A\to\mathcal{L}(D;A\to B)$$

*The differential equation*

$$\frac{dy}{dx}=g^{k\cdot ij}(e_{B\cdot i}\otimes e_{B\cdot j})\circ F_k\circ G$$

*is integrable iff*

$$\begin{aligned}&\left(\frac{dg^{k\cdot ij}(x)}{dx}\circ h_2\right)(e_{B\cdot i}(F_k\circ G\circ h_1)e_{B\cdot j})\\ =&\left(\frac{dg^{k\cdot ij}(x)}{dx}\circ h_1\right)(e_{B\cdot i}(F_k\circ G\circ h_2)e_{B\cdot j})\end{aligned}\tag{6.5.4}$$

PROOF. Since

$$g\circ h_1=g^{k\cdot ij}(e_i\otimes e_j)\circ F_k\circ G\circ h_1=g^{k\cdot ij}(e_i(F_k\circ G\circ h_1)e_j)\tag{6.5.5}$$

then, according to the theorem [15]-8.2.3,

$$\frac{dg}{dx}\circ(h_1,h_2)=\left(\frac{dg^{k\cdot ij}}{dx}\circ h_2\right)(e_i(F_k\circ G\circ h_1)e_j)\tag{6.5.6}$$

The equality (6.5.4) follows from the equality (6.5.6), from the theorem 6.5.1 and the definition [15]-7.4.1. □

The theorem 6.5.5 answers the question 6.3.4.

THEOREM 6.5.5. *Let $B$ be free finite dimensional associative $D$-algebra. Let $\overline{\overline{e}}_A$ be basis of $D$ module $A$. Let $\overline{\overline{e}}_B$ be basis of $D$ module $B$. The differential equation*

$$\frac{dy}{dx}=g(x)\tag{6.5.7}$$

*is integrable iff corresponding system of differential equations*

$$\frac{\partial y^i}{\partial x^j}=g^i_j\quad y=y^ie_{B\cdot i}\quad x=x^ie_{A\cdot i}\tag{6.5.8}$$

*is completely integrable.*

PROOF. Let $\overline{\overline{\overline{F}}}$ be the basis of left $B\otimes B$-module $\mathcal{L}(D;B\to B)$ and $G$ be the linear map

$$G:A\to B$$

of maximal rank such that $\ker G\subseteq\ker$ .L et $g^{k\cdot ij}$ be standard components of the map

$$g:A\to\mathcal{L}(D;A\to B)$$

Then the differential equation (6.5.7) has form

$$\frac{dy}{dx}=g^{k\cdot ij}(e_{B\cdot i}\otimes e_{B\cdot j})\circ F_k\circ G\tag{6.5.9}$$

According to the theorem 6.5.4, the differential equation (6.5.7) is integrable iff

(6.5.10)
$$\begin{aligned}&\left(\frac{dg^{k\cdot ij}(x)}{dx}\circ h_2\right)(e_{B\cdot i}(F_k\circ G\circ h_1)e_{B\cdot j})\\=&\left(\frac{dg^{k\cdot ij}(x)}{dx}\circ h_1\right)(e_{B\cdot i}(F_k\circ G\circ h_2)e_{B\cdot j})\end{aligned}$$

Let

$$h_1=h_1^ke_k\quad h_2=h_2^le_l$$

Then we can write the equality (6.5.10) as

(6.5.11)
$$\begin{aligned}&\left(\frac{dg^{k\cdot ij}(x)}{dx}\circ (h_2^le_{1\cdot l})\right)(e_{2\cdot i}(F_k\circ G\circ (h_1^ke_{1\cdot k}))e_{2\cdot j})\\=&\left(\frac{dg^{k\cdot ij}(x)}{dx}\circ (h_1^ke_{1\cdot k})\right)(e_{2\cdot i}(F_k\circ G\circ (h_2^le_{1\cdot l}))e_{2\cdot j})\end{aligned}$$

Since $h_1$, $h_2$ are arbitrary $A$-numbers, then the equality

(6.5.12)
$$\begin{aligned}&\left(\frac{dg^{k\cdot ij}(x)}{dx}\circ e_{1\cdot l}\right)(e_{2\cdot i}(F_k\circ G\circ e_{1\cdot k})e_{2\cdot j})\\=&\left(\frac{dg^{k\cdot ij}(x)}{dx}\circ e_{1\cdot k}\right)(e_{2\cdot i}(F_k\circ G\circ e_{1\cdot l})e_{2\cdot j})\end{aligned}$$

follows from the equality (6.5.11). Since

$$F_k\circ G\circ e_{1\cdot k}=F_{k\ r}^{\ m}G_k^re_{2\cdot m}$$

then the equality

(6.5.13)
$$\begin{aligned}&\left(\frac{dg^{k\cdot ij}(x)}{dx}\circ e_{1\cdot l}\right)F_{k\ r}^{\ m}G_k^re_{2\cdot i}e_{2\cdot m}e_{2\cdot j}\\=&\left(\frac{dg^{k\cdot ij}(x)}{dx}\circ e_{1\cdot k}\right)F_{k\ r}^{\ m}G_l^re_{2\cdot i}e_{2\cdot m}e_{2\cdot j}\end{aligned}$$

follows from the equality (6.5.12). The equality

(6.5.14)
$$\begin{aligned}&\left(\frac{dg^{ij}(x)}{dx}\circ e_{1\cdot l}\right)F_{k\ r}^{\ m}G_k^rC_{2\cdot\ im}^{\ \ p}C_{2\cdot\ pj}^{\ \ r}e_{2\cdot r}\\=&\left(\frac{dg^{k\cdot ij}(x)}{dx}\circ e_{1\cdot k}\right)F_{k\ r}^{\ m}G_l^rC_{2\cdot\ im}^{\ \ p}C_{2\cdot\ pj}^{\ \ r}e_{2\cdot r}\end{aligned}$$

follows from the equality (6.5.13). The equality

(6.5.15)
$$\left(\frac{dg_k^r(x)}{dx}\circ e_{1\cdot l}\right)e_{2\cdot r}=\left(\frac{dg_l^r(x)}{dx}\circ e_{1\cdot k}\right)e_{2\cdot r}$$

follows from the equality (6.5.14). The equality

(6.5.16)
$$\frac{dg_k^r(x)}{dx}\circ e_{1\cdot l}=\frac{dg_l^r(x)}{dx}\circ e_{1\cdot k}$$

follows from the equality (6.5.15). □

CHAPTER 7

# Differential Equation of First Order

## 7.1. Differential Equation with Separated Variables

In the theory of differential equations over field, the first differential equation which we study is separable equation

$$\frac{dy}{dx}=\frac{M(x)}{N(y)} \tag{7.1.1}$$

As soon as we separate variables, we can write differential equation (7.1.1) as

$$N(y)dy=M(x)dx \tag{7.1.2}$$

It is easy to integrate differential equation (7.1.2).

In noncommutative $D$-algebra, $dx$ and $dy$ belong in general to different $D$-modules. At the same time, a solution of the differential equation (7.1.2) is implicit function of variables $x$ and $y$. So I formulate the problem as follows.

Let $X$, $Y$ be Banach $D$-modules. Let $A$ be Banach $D$-algebra. Consider maps

$$f:X\rightarrow\mathcal{L}(D;X\rightarrow A)$$
$$g:Y\rightarrow\mathcal{L}(D;Y\rightarrow A)$$

Differential equation

$$\frac{dy}{dx}=g(y)\circ f(x) \tag{7.1.3}$$

is called **differential separable equation**.

In noncommutative $D$-algebra, even this is division algebra, operation of separation of variables may be impossible. However, if there exists the map

$$h:Y\rightarrow\mathcal{L}(D;Y\rightarrow A)$$

such that

$$h\circ g=1\otimes 1$$

then we can write the differential equation (7.1.3) in the following form

$$h(y)\circ dy=f(x)\circ dx$$

Let $X$, $Y$ be Banach $D$-modules. Let $A$ be Banach $D$-algebra. Consider maps

$$M:X\rightarrow\mathcal{L}(D;X\rightarrow A)$$
$$N:Y\rightarrow\mathcal{L}(D;Y\rightarrow A)$$

Differential equation

$$M(x)\circ dx+N(y)\circ dy=0$$

is called **differential equation with separated variables**.

THEOREM 7.1.1. *Let maps*

$$M : X \to \mathcal{L}(D; X \to A)$$
$$N : Y \to \mathcal{L}(D; Y \to A)$$

*be integrable. The solution of differential equation*

$$M(x) \circ dx + N(y) \circ dy = 0 \tag{7.1.4}$$

*is implicit function*[7.1]

$$\int M(x) \circ dx + \int N(y) \circ dy = C \tag{7.1.5}$$

PROOF. The equality

$$N(y) \circ dx = N(y) \circ \frac{dy}{dx} \circ dx \tag{7.1.6}$$

follows from the definition [15]-3.3.2. Integrability of the differential form in right side of the equality (7.1.6) follows from integrability of the differential form in left side of the equality (7.1.6). According to the theorem A.2.4 the differential form

$$M(x) \circ dx + N(y) \circ dx = \left(M(x) + N(y) \circ \frac{dy}{dx}\right) \circ dx \tag{7.1.7}$$

is integrable. The equality

$$M(x) + N(y) \circ \frac{dy}{dx} = 0 \otimes 0 \tag{7.1.8}$$

follows from equalities (7.1.4), (7.1.7). The equality (7.1.5) follows from the equality (7.1.8) and from the theorem A.2.3. □

EXAMPLE 7.1.2. *Consider the differential equation*

$$(1 \otimes x + x \otimes 1) \circ dx + (1 \otimes y + y \otimes 1) \circ dy = 0 \tag{7.1.9}$$

*According to theorems 7.1.1, A.2.6, implicit function*

$$x^2 + y^2 = C$$

*is solution of differential equation (7.1.9).* □

THEOREM 7.1.3. *The solution of the differential equation*

$$M(x) \circ dx + N(y) \circ dy = 0 \tag{7.1.10}$$

*with initial condition*

$$x_0 = 0 \quad y_0 = C$$

*has form*[7.2]

$$\int_{x_0}^{x} M(x) \circ dx + \int_{y_0}^{y} N(y) \circ dy = 0 \tag{7.1.11}$$

---

[7.1] The proof of the theorem is based on the proofs on the page [3]-44 and on the page [6]-24; however the proof is different since $D$-algebra $A$ in general is noncommutative.

[7.2] See also remarks on the page [3]-44 and on the page [6]-24.

Proof. According to the theorem 7.1.1, implicit function

$$\int M(x)\circ dx+\int N(y)\circ dy=C \tag{7.1.12}$$

is solution of the differential equation (7.1.10). According to the definition 6.1.1, there exist maps

$$P:X\to A$$
$$R:Y\to A$$

such that

$$P(x)=\int M(x)\circ dx \tag{7.1.13}$$

$$R(y)=\int N(y)\circ dy \tag{7.1.14}$$

According to the theorem [15]-6.2.1 and to the definition [15]-8.3.7, equalities

$$\int M(x)\circ dx=\int_{x_0}^{x}M(x)\circ dx+P(x_0) \tag{7.1.15}$$

$$\int N(y)\circ dy=\int_{y_0}^{y}N(y)\circ dy+R(y_0) \tag{7.1.16}$$

follow from equalities (7.1.13), (7.1.14). The equality

$$\int_{x_0}^{x}M(x)\circ dx+P(x_0)+\int_{y_0}^{y}N(y)\circ dy+R(x_0)=C \tag{7.1.17}$$

follows from equalities (7.1.12), (7.1.15), (7.1.16). If $x=x_0$, $y=y_0$, then the equality

$$P(x_0)+R(x_0)=C \tag{7.1.18}$$

follows from the equality (7.1.17). The equality (7.1.11) follows from equalities (7.1.17), (7.1.18). □

## 7.2. Exact Differential Equation

Definition 7.2.1. *The differential equation*

$$M(x,y)\circ dx+N(x,y)\circ dy=0$$

*where*

$$M:(x,y)\in X\times Y\to M(x,y)\in\mathcal{L}(D;X\to A)$$
$$N:(x,y)\in X\times Y\to N(x,y)\in\mathcal{L}(D;Y\to A)$$

*is called* **exact differential equation***, if there exists map*

$$u:X\times Y\to A$$

*such that*

$$\frac{\partial u(x,y)}{\partial x}=M(x,y) \tag{7.2.1}$$

$$\frac{\partial u(x,y)}{\partial y}=N(x,y) \tag{7.2.2}$$

□

THEOREM 7.2.2. *The differential equation*

$$M(x,y)\circ dx+N(x,y)\circ dy=0 \tag{7.2.3}$$

*is integrable iff*

$$\frac{\partial M(x,y)}{\partial x}\circ(dx_1,dx_2)=\frac{\partial M(x,y)}{\partial x}\circ(dx_2,dx_1) \tag{7.2.4}$$

$$\frac{\partial N(x,y)}{\partial y}\circ(dy_1,dy_2)=\frac{\partial N(x,y)}{\partial y}\circ(dy_2,dy_1) \tag{7.2.5}$$

$$\frac{\partial M(x,y)}{\partial y}\circ(dx,dy)=\frac{\partial N(x,y)}{\partial x}\circ(dy,dx) \tag{7.2.6}$$

PROOF. Let $Z=X\oplus Y$, $z=x\oplus y$. According to definitions [15]-7.3.2, 2.2.5, 7.2.1 and the theorem 2.2.6, the expression

$$\omega=M(x,y)\circ dx+N(x,y)\circ dy \tag{7.2.7}$$

is differential form of degree 1. According to the theorem [15]-8.3.1, differential form $\omega$ is integrable iff $d\omega=0$. Equalities

$$\begin{aligned}\frac{d\omega}{dz}\circ(z_1,z_2)=&\frac{\partial M(x,y)}{\partial x}\circ(dx_1,dx_2)+\frac{\partial M(x,y)}{\partial y}\circ(dx_1,dy_2)\\&+\frac{\partial N(x,y)}{\partial x}\circ(dy_1,dx_2)+\frac{\partial N(x,y)}{\partial y}\circ(dy_1,dy_2)\end{aligned} \tag{7.2.8}$$

$$\begin{aligned}\frac{d\omega}{dz}\circ(z_2,z_1)=&\frac{\partial M(x,y)}{\partial x}\circ(dx_2,dx_1)+\frac{\partial M(x,y)}{\partial y}\circ(dx_2,dy_1)\\&+\frac{\partial N(x,y)}{\partial x}\circ(dy_2,dx_1)+\frac{\partial N(x,y)}{\partial y}\circ(dy_2,dy_1)\end{aligned} \tag{7.2.9}$$

follow from the equality (7.2.7). The equality

$$\begin{aligned}&\frac{\partial M(x,y)}{\partial x}\circ(dx_1,dx_2)+\frac{\partial M(x,y)}{\partial y}\circ(dx_1,dy_2)\\&+\frac{\partial N(x,y)}{\partial x}\circ(dy_1,dx_2)+\frac{\partial N(x,y)}{\partial y}\circ(dy_1,dy_2)\\&=\frac{\partial M(x,y)}{\partial x}\circ(dx_2,dx_1)+\frac{\partial M(x,y)}{\partial y}\circ(dx_2,dy_1)\\&+\frac{\partial N(x,y)}{\partial x}\circ(dy_2,dx_1)+\frac{\partial N(x,y)}{\partial y}\circ(dy_2,dy_1)\end{aligned} \tag{7.2.10}$$

follows from equalities (7.2.8), (7.2.9). Equalities (7.2.4), (7.2.5), (7.2.6) follow from the equality (7.2.10). □

Let equalities (7.2.4), (7.2.5), (7.2.6) be true. Our goal is to find the implicit function[7.3] $u(x,y)=C$ which satisfies to equalities (7.2.1), (7.2.2). According to theorems [15]-8.3.1, A.2.1, the equality

$$u(x,y)=\int M(x,y)\circ dx+C_1(y) \tag{7.2.11}$$

[7.3] Solving of the equation (7.2.3) is similar to solving of exact differential equation in commutative algebra. See, for instance, the proof of the theorem [3]-2.6.1 as well solving of exact differential equation on pages 33, 34 in the book [6].

follows from equalities (7.2.1), (7.2.4). The equality

$$\frac{\partial}{\partial y}\int M(x,y)\circ dx+\frac{dC_1(y)}{dy}=N(x,y) \tag{7.2.12}$$

follows from equalities (7.2.2), (7.2.11). The equality

$$\frac{dC_1(y)}{dy}=N(x,y)-\frac{\partial}{\partial y}\int M(x,y)\circ dx \tag{7.2.13}$$

follows from the equality (7.2.12). There is variable $x$ in the expression

$$N(x,y)-\frac{\partial}{\partial y}\int M(x,y)\circ dx \tag{7.2.14}$$

of right side of the equation (7.2.13). In order for the equation (7.2.13) to have a solution, it is necessary that expression (7.2.14) does not depend on $x$. Consider derivative of the expression (7.2.14) with respect to $x$. The equality

$$\begin{aligned}
&\frac{\partial}{\partial x}\left(N(x,y)\circ dy-\left(\frac{\partial}{\partial y}\int M(x,y)\circ dx\right)\circ dy\right)\circ dx\\
&=\frac{\partial N(x,y)}{\partial x}\circ(dy,dx)-\left(\frac{\partial^2}{\partial x\partial y}\int M(x,y)\circ dx\right)\circ(dy,dx)\\
&=\frac{\partial N(x,y)}{\partial x}\circ(dy,dx)-\left(\frac{\partial^2}{\partial y\partial x}\int M(x,y)\circ dx\right)\circ(dx,dy)\\
&=\frac{\partial N(x,y)}{\partial x}\circ(dy,dx)-\left(\frac{\partial}{\partial y}\frac{\partial}{\partial x}\int M(x,y)\circ dx\right)\circ dy\\
&=\frac{\partial N(x,y)}{\partial x}\circ(dy,dx)-\frac{\partial M(x,y)}{\partial y}\circ(dx,dy)=0
\end{aligned} \tag{7.2.15}$$

follows from equalities [8]-(9.1.5), (7.2.6), (A.2.2).

From theorems [15]-8.3.1, A.2.1 and from the equality (7.2.5), it follows that there exist the integral

$$\begin{aligned}
&\int\left(N(x,y)-\frac{\partial}{\partial y}\int M(x,y)\circ dx\right)\circ dy\\
&=\int N(x,y)\circ dy-\int\left(\frac{\partial}{\partial y}\int M(x,y)\circ dx\right)\circ dy
\end{aligned}$$

Therefore, we can find the map $C_1$ as solution of the differetial equation (7.2.13).

Example 7.2.3. *Consider differential equation*

$$(1\otimes 1+1\otimes y)\circ dx+(x\otimes 1+1\otimes 1)\circ dy=0 \tag{7.2.16}$$

*Our goal is to find implicit function $u(x,y)$, which is solution of the differential equation (7.2.16). Differential equations*

$$\frac{\partial u}{\partial x}=1\otimes 1+1\otimes y \tag{7.2.17}$$

$$\frac{\partial u}{\partial y}=x\otimes 1+1\otimes 1 \tag{7.2.18}$$

*follow from the differential equation (7.2.16). Then*

$$u(x,y)=\int(1\otimes 1+1\otimes y)\circ dx=x+xy+C_1(y) \tag{7.2.19}$$

*is the solution of the differential equation (7.2.17). The equality*

$$\frac{\partial u(x,y)}{\partial y} = x\otimes 1 + \frac{dC_1(y)}{dy} = x\otimes 1 + 1\otimes 1 \tag{7.2.20}$$

*follows from the equality (7.2.19) and from the equation (7.2.18). The differential equation*

$$\frac{dC_1(y)}{dy} = 1\otimes 1 \tag{7.2.21}$$

*follows from the equality (7.2.20).*

$$C_1(y) = y \tag{7.2.22}$$

*is the solution of the differential equation (7.2.21). From equalities (7.2.19), (7.2.22), it follows that implicit function*

$$x + xy + y = C$$

*is the solution of the differential equation (7.2.16).* □

The example 7.2.4 shows that equalities (7.2.4), (7.2.5) are essential for integralibility of differential equation (7.2.3).

Example 7.2.4. *Consider the differential equation*

$$(3x^2\otimes 1 + 1\otimes y)\circ dx + (x\otimes 1)\circ dy = 0 \tag{7.2.23}$$

*Here*

$$M(x,y) = 3x^2\otimes 1 + 1\otimes y \qquad N(x,y) = x\otimes 1$$

*Therefore*

$$\frac{\partial M(x,y)}{\partial y} = 1\otimes_x 1\otimes_y 1$$
$$\frac{\partial N(x,y)}{\partial x} = 1\otimes_x 1\otimes_y 1$$

*However, the differential equation (7.2.23) is not integrable, because integral*

$$\int (3x^2\otimes 1 + 1\otimes y)\circ dx$$

*does not exist.* □

The example 7.2.5 shows that order of variables is essential in the equality (7.2.6).

Example 7.2.5. *Consider the differential equation*

$$(1\otimes y)\circ dx + (1\otimes x)\circ dy = 0 \tag{7.2.24}$$

*Here*

$$M(x,y) = 1\otimes y \qquad N(x,y) = 1\otimes x$$

*Therefore*

$$\frac{\partial M(x,y)}{\partial y} = 1\otimes_x 1\otimes_y 1$$
$$\frac{\partial N(x,y)}{\partial x} = 1\otimes_y 1\otimes_x 1$$

*Although the derivatives* $\dfrac{\partial M(x,y)}{\partial y}, \dfrac{\partial N(x,y)}{\partial x}$ *are represented by the same tensor, their action on variables* $x$, $y$ *is different*

$$\frac{\partial M(x,y)}{\partial y}\circ(dx,dy)=(1\otimes_x 1\otimes_y 1)\circ(dx,dy)=dx\,dy$$
$$\neq\frac{\partial N(x,y)}{\partial x}\circ(dy,dx)=(1\otimes_y 1\otimes_x 1)\circ(dy,dx)=dy\,dx$$

*It is easy to see that the differential equation* (7.2.24) *is not integrable.* □

CHAPTER 8

# Exponent

## 8.1. Definition of Exponent

Definition 8.1.1. *The map*

$$y = e^x = \sum_{n=0}^{\infty} \frac{1}{n!} x^n \tag{8.1.1}$$

*is called exponent.* □

Theorem 8.1.2. *The derivative of exponent has the following form*

$$\frac{de^x}{dx} = \sum_{n=0}^{\infty} \frac{1}{(n+1)!} \sum_{i=0}^{n} x^i \otimes x^{n-i} \tag{8.1.2}$$

Proof. According to theorems A.1.3, A.1.4, the equality

$$\frac{de^x}{dx} = \sum_{n=0}^{\infty} \frac{1}{(n+1)!} \frac{dx^{n+1}}{dx} \tag{8.1.3}$$

follows from the equality (8.1.1). According to the theorem 4.2.11, the equality (8.1.2) follows from the equality (8.1.3). □

The equality (8.1.2) is the differential equation for exponent. However this differential equation is solved with respect to the derivative. In a field, exponent is a solution of differential equation

$$y' = y \tag{8.1.4}$$

Our goal in non commutative Banach algebra is to consider a differential equation that is analogous to the differential equation (8.1.4).

Theorem 8.1.3. *The map* $y = ae^x b$, $ab = C$, *is solution of the differential equation*

$$\frac{dy}{dx} \circ 1 = y \tag{8.1.5}$$

*with initial condition*

$$x = 0 \quad y = C$$

Proof. The equality

$$\begin{aligned} \frac{de^x}{dx} \circ 1 &= \sum_{n=0}^{\infty} \frac{1}{(n+1)!} \sum_{i=0}^{n} (x^i \otimes x^{n-i}) \circ 1 = \sum_{n=0}^{\infty} \frac{1}{(n+1)!} \sum_{i=0}^{n} x^i x^{n-i} \\ &= \sum_{n=0}^{\infty} \frac{n+1}{(n+1)!} x^n \end{aligned} \tag{8.1.6}$$

follows from the equality (8.1.2). From equalities (8.1.1), (8.1.6), it follows that

(8.1.7) $$\frac{de^x}{dx}\circ 1=e^x$$

and the map $y=e^x$ is solution of the differential equation (8.1.5).

Let

(8.1.8) $$y=ae^xb$$

According to the theorem A.1.3, the equality

(8.1.9) $$\frac{dy}{dx}\circ 1=\frac{dae^xb}{dx}\circ 1=a\left(\frac{de^x}{dx}\circ 1\right)b=ae^xb=y$$

follows from the equality (8.1.8). Since $e^0=1$, then the theorem follows from the equality (8.1.9). □

According to the theorem 8.1.3, a solution of the differential equation (8.1.5) with initial condition

$$x=0\quad y=C$$

is not unique. This is due in particular to the equality

$$\frac{dy}{dx}\circ 1=\frac{\partial y}{\partial x^0}$$

However, this differential equation is very important.

## 8.2. Properties of Exponent

Theorem 8.2.1. *The equality*

(8.2.1) $$e^{a+b}=e^ae^b$$

*is true iff*

(8.2.2) $$ab=ba$$

Proof. To prove the theorem it is enough to consider Taylor series

(8.2.3) $$e^a=\sum_{n=0}^{\infty}\frac{1}{n!}a^n$$

(8.2.4) $$e^b=\sum_{n=0}^{\infty}\frac{1}{n!}b^n$$

(8.2.5) $$e^{a+b}=\sum_{n=0}^{\infty}\frac{1}{n!}(a+b)^n$$

Let us multiply expressions (8.2.3) and (8.2.4). The sum of monomials of order 3 has form

(8.2.6) $$\frac{1}{6}a^3+\frac{1}{2}a^2b+\frac{1}{2}ab^2+\frac{1}{6}b^3$$

and in general does not equal expression

(8.2.7) $$\frac{1}{6}(a+b)^3=\frac{1}{6}a^3+\frac{1}{6}a^2b+\frac{1}{6}aba+\frac{1}{6}ba^2+\frac{1}{6}ab^2+\frac{1}{6}bab+\frac{1}{6}b^2a+\frac{1}{6}b^3$$

The proof of statement that (8.2.1) follows from (8.2.2) is trivial. □

THEOREM 8.2.2. *Let $A$ be associative $D$-algebra. Then*

$$ae^{xa} = e^{ax}a \tag{8.2.8}$$

PROOF. The equality

$$\begin{aligned} e^{ax}a &= (1 + ax + \frac{1}{2}axax + \frac{1}{3!}axaxax + \frac{1}{4!}axaxaxax + ...)a \\ &= a + axa + \frac{1}{2}axaxa + \frac{1}{3!}axaxaxa + \frac{1}{4!}axaxaxaxa + ... \\ &= a(1 + xa + \frac{1}{2}xaxa + \frac{1}{3!}xaxaxa + \frac{1}{4!}xaxaxaxa + ...) = ae^{xa} \end{aligned} \tag{8.2.9}$$

follows from the equality

$$e^{ax} = 1 + ax + \frac{1}{2}axax + \frac{1}{3!}axaxax + \frac{1}{4!}axaxaxax + ... \tag{8.2.10}$$

The equality (8.2.8) follows from the equality (8.2.9). □

THEOREM 8.2.3. *Let $A$ be associative $D$-algebra. If $A$-number $a$ has a reverse element, then*

$$e^{bxa} = a^{-1}e^{abx}a \tag{8.2.11}$$

*In particular*

$$e^{xa} = a^{-1}e^{ax}a \tag{8.2.12}$$

PROOF. The equality (8.2.12) follows from the equality (8.2.8). The equality (8.2.11) follows from the equality (8.2.12), if we write $bx$ imstead of $x$. □

## 8.3. Quasiexponent

Exponent is not only solution of the differential equation (8.1.5).

DEFINITION 8.3.1. *The map*

$$e[c_1, ..., c_n]^x = \frac{d^n e^x}{dx} \circ (c_1, ..., c_n) \tag{8.3.1}$$

*is called* **quasiexponent**. □

THEOREM 8.3.2. *The quasiexponent $y = e[c]^x$ has following Taylor series expansion*

$$e[c]^x = \sum_{n=0}^{\infty} \frac{1}{(n+1)!} \sum_{m=0}^{n} x^m c x^{n-m} = c + \frac{1}{2}(cx + xc) + \frac{1}{6}(cx^2 + xcx + x^2c) + ...$$

PROOF. The theorem follows from the definition 8.3.1 and the theorem 8.1.2. □

EXAMPLE 8.3.3. *Let $y = e[i]^x$. According to the theorem 8.3.2, $(c = i)$,*

$$e[i]^x = \sum_{n=0}^{\infty} \frac{1}{(n+1)!} \sum_{m=0}^{n} x^m i x^{n-m} = i + \frac{1}{2}(ix + xi) + \frac{1}{6}(ix^2 + xix + x^2i) + ...$$

□

THEOREM 8.3.4. *If $c \in \ker A$, then*

$$e[c]^x = ce^x \tag{8.3.2}$$

PROOF. According to the definition 8.3.1 and the theorem 8.1.2

$$\begin{aligned} e[c]^x &= \sum_{n=0}^{\infty} \frac{1}{(n+1)!} \sum_{m=0}^{n} x^m c x^{m-i} = \sum_{n=0}^{\infty} \frac{1}{(n+1)!} \sum_{m=0}^{n} c x^n \\ &= c \sum_{n=0}^{\infty} \frac{1}{n!} x^n = c e^x \end{aligned} \tag{8.3.3}$$

The theorem follows from the equality (8.3.3). □

THEOREM 8.3.5.

$$\frac{d^m e[c_1, ..., c_n]^x}{dx^m} \circ (C_1, ..., C_m) = e[c_1, ..., c_n, C_1, ..., C_m]^x \tag{8.3.4}$$

PROOF. The equality

$$\begin{aligned} \frac{d^m e[c_1, ..., c_n]^x}{dx^m} \circ (C_1, ..., C_m) &= \frac{d^m}{dx^m} \left( \frac{d^n e^x}{dx^n} \circ (c_1, ..., c_n) \right) \circ (C_1, ..., C_m) \\ &= \frac{d^{n+m} e^x}{dx^{n+m}} \circ (c_1, ..., c_n, C_1, ..., C_m) \\ &= e[c_1, ..., c_n, C_1, ..., C_m]^x \end{aligned} \tag{8.3.5}$$

follows from equalities (8.3.1), (A.1.19). □

THEOREM 8.3.6. *Quasiexponent is a solution of the differential equation*

$$\frac{dy}{dx} \circ 1 = y \tag{8.3.6}$$

PROOF. The equality

$$\frac{de[c_1, ..., c_n]^x}{dx} \circ 1 = \frac{d}{dx} \left( \frac{d^n e^x}{dx^n} \circ (c_1, ..., c_n) \right) \circ 1 = \frac{d^{n+1} e^x}{dx^{n+1}} \circ (c_1, ..., c_n, 1) \tag{8.3.7}$$

follows from the equality (8.3.1) and the definition [15]-4.1.4. According to the theorem 4.2.8 and the definition [15]-4.1.4, the equality

$$\frac{de[c_1, ..., c_n]^x}{dx} \circ 1 = \frac{d^{n+1} e^x}{dx^{n+1}} \circ (1, c_1, ..., c_n) = \frac{d^n}{dx^n} \left( \frac{de^x}{dx} \circ 1 \right) \circ (c_1, ..., c_n) \tag{8.3.8}$$

follows from the equality (8.3.7). The theorem follows from equalities (8.1.7), (8.3.1), (8.3.8). □

## 8.4. Differential Equation $\frac{dy}{dx} \circ 1 = y$

Following questions follow from the theorem 8.3.6.

QUESTION 8.4.1. *Why does the differential equation* (8.3.6) *have infinitely many solutions for given initial conditions?* □

QUESTION 8.4.2. *Why it is interesting for us the differential equation* (8.3.6) *if it has infinitely many solutions for given initial conditions?* □

QUESTION 8.4.3. *Is there an analogue in noncommutative algebra for linear homogenius differential equation?* □

**8.4.1. Answer to the question 8.4.1.** The theorem [15]-7.3.7 answers the question 8.4.1. Let $\overline{\overline{e}}$ be a basis of $D$-algebra $A$. According to this theorem the expression[8.1]

$$\frac{\partial y}{\partial x^{0}}=\frac{dy}{dx}\circ e_{0} \tag{8.4.1}$$

is partial derivative with respect to the variable $x^{0}$.

THEOREM 8.4.4. *For any* $n>0$

$$\frac{d^n y}{dx^n}\circ(e_{0},...,e_{0})=\frac{\partial^n y}{(\partial x^{0})^n}$$

PROOF. We will prove the theorem by induction over $n$.

According to the equality (8.4.1), the theorem is true for $n=1$.

Let the theorem be true for $n=k$

$$\frac{d^k y}{dx^k}\circ(e_{0},...,e_{0})=\frac{\partial^k y}{(\partial x^{0})^k} \tag{8.4.2}$$

According to the definition [15]-4.1.4 and the theorem A.1.5, the equality

$$\begin{aligned}\frac{d^{k+1}y}{dx^{k+1}}\circ(e_{0},...,e_{0},e_{0})&=\frac{d}{dx}\left(\frac{d^k y}{dx^k}\circ(e_{0},...,e_{0})\right)\circ e_{0}\\&=\frac{\partial}{\partial x^{0}}\frac{\partial^k y}{(\partial x^{0})^k}=\frac{\partial^{k+1}y}{(\partial x^{0})^{k+1}}\end{aligned}$$

follows from equalities (8.4.1), (8.4.2). Therefore, the theorem is true for $n=k+1$. □

According to the theorem 8.4.4, solution of the differential equation (8.3.6) for given initial condition is unique for $A$-numbers like $de_{0}$, $d\in D$. However partial derivative with respect to $x^{i}$, $i\neq 0$, is not defined; this is a source of ambiguity of solution.

**8.4.2. Answer to the question 8.4.2.** The question 8.4.2 looks natural in the light of the answer in subsection 8.4.1. Initially, I perceived the quasiexponent as an uninvited guest. However, the analysis of the definition 8.3.1 and theorems 8.1.3, 8.3.6 shows that this map is very important. Namely, differential of quasiexponent is quasiexponent.

**8.4.3. Answer to the question 8.4.3.** According to the theorem 8.1.3, the map $y=e^x$ is a solution of the differential equation

$$\frac{dy}{dx}\circ 1=y \tag{8.4.3}$$

What can i say about the map $y=e^{ax}$?

According to the theorem [15]-3.3.23

$$\frac{de^{ax}}{dx}\circ 1=\frac{de^{ax}}{dax}\circ\frac{dax}{dx}\circ 1=\frac{de^{ax}}{dax}\circ(a\otimes 1)\circ 1=\frac{de^{ax}}{dax}\circ a \tag{8.4.4}$$

The equality

$$\frac{de^{ax}}{dx}\circ 1=e[a]^{ax} \tag{8.4.5}$$

[8.1] According to convention 1.5.3, $e_{0}=1$.

follows from equalities (8.3.1), (8.4.4). Therefore, the map $y=e^{ax}$ is not a solution of the differential equation (8.4.3).

If $D$-algebra $A$ is division algebra, then we can write equality

$$\frac{de^{ax}}{dx}\circ a^{-1}=\frac{de^{ax}}{dax}\circ\frac{dax}{dx}\circ a^{-1}=\frac{de^{ax}}{dax}\circ 1=e^{ax} \tag{8.4.6}$$

instead of equality (8.4.4).

Even if I will consider the differential equation

$$\frac{dy}{dx}\circ a^{-1}=y \tag{8.4.7}$$

it is not clear how to move from this differential equation to a differential equation of higher order or to a system of homogenius differential equations, because there is no coefficient before $y$ in the differential equation (8.4.7).

However I am not ready to give negative answer to the question 8.4.3.

## 8.5. Exponent of matrix

Since there are two products ($a_*{}^*b$ and $a^*{}_*b$) on the set of matrices, we need to consider two types of exponent.[8.2]

Theorem 8.5.1. **$_*{}^*$-exponent** $y=e^{x_*{}^*}$ *has following Taylor series expansion*

$$e^{x_*{}^*}=\sum_{n=0}^{\infty}\frac{1}{n!}x^{n_*{}^*} \tag{8.5.1}$$

*and is the solution of the differential equation*

$$\frac{dy}{dx}{}_\circ{}^\circ E=y$$

*with initial condition $y(0)=E$.*

Proof. The theorem follows from the theorem 8.1.3, if we consider operation $a_*{}^*b$ as product. □

Theorem 8.5.2. **$^*{}_*$-exponent** $y=e^{x^*{}_*}$ *has following Taylor series expansion*

$$e^{x^*{}_*}=\sum_{n=0}^{\infty}\frac{1}{n!}x^{n^*{}_*} \tag{8.5.2}$$

*and is the solution of the differential equation*

$$\frac{dy}{dx}{}^\circ{}_\circ E=y$$

*with initial condition $y(0)=E$.*

Proof. The theorem follows from the theorem 8.1.3, if we consider operation $a^*{}_*b$ as product. □

[8.2] I wrote this section under the influence of the paper [16].

CHAPTER 9

# Differential Equation of Function of Real Variable

## 9.1. Differential Equation $\frac{dx}{dt} = ax$

THEOREM 9.1.1. *Let $A$ be Banach $D$-algebra and $a \in A$. The map*

$$f : t \in R \to e^{at} \in A$$

*has the following Taylor series decomposition*

$$e^{at} = \sum_{n=0}^{\infty} \frac{1}{n!} a^n t^n \tag{9.1.1}$$

PROOF. According to the definition 8.1.1

$$e^{at} = \sum_{n=0}^{\infty} \frac{1}{n!} (at)^n \tag{9.1.2}$$

LEMMA 9.1.2. *For any $n \geq 0$*

$$(at)^n = a^n t^n \tag{9.1.3}$$

PROOF. We will prove the lemma by induction over $n$.
The lemma is true for $n = 0$, because

$$(at)^0 = 1 = a^0 t^0$$

Let the lemma be true for $n = k$

$$(at)^k = a^k t^k \tag{9.1.4}$$

Then the equality

$$(at)^{k+1} = (at)^k (at) = (a^k t^k)(at) \tag{9.1.5}$$

follows from the equality (9.1.4). Since the product in $D$-algebra $A$ is bilinear, the equality

$$(at)^{k+1} = (a^k a)(t^k t) = a^{k+1} t^{k+1} \tag{9.1.6}$$

follows from the equality (9.1.5). From the equality (9.1.6), it follows that the lemma is true for $n = k + 1$. ⊙

The equality (9.1.1) follows from equalities (9.1.2), (9.1.3). □

THEOREM 9.1.3. *Let $A$ be Banach associative $D$-algebra and $a$, $c \in A$. The condition*

$$ca = ac \tag{9.1.7}$$

*implies that*

$$e^{at} c = c e^{at} \tag{9.1.8}$$

Proof. Let the equality (9.1.7) be true. According to the theorem 9.1.1, to prove the equality (9.1.8), it is enough to prove the equality

$$a^n c = c a^n \tag{9.1.9}$$

We will prove the equality (9.1.9) by induction over $n$.

For $n = 0$, the equality (9.1.9) is evident since $a^0 = 1$. For $n = 1$, the equality (9.1.9) is the equality (9.1.7).

Let the equality (9.1.9) be true for $n = k$

$$a^k c = c a^k \tag{9.1.10}$$

The equality

$$a^{k+1} c = a a^k c = a c a^k = c a a^k = c a^{k+1}$$

follows from equalities (9.1.7), (9.1.10). Therefore, the equality (9.1.9) is true for $n = k + 1$. □

Theorem 9.1.4. *Let $A$ be Banach $D$-algebra and $a \in A$. Then*

$$e^{at} a = a e^{at} \tag{9.1.11}$$

Proof. The theorem follows from the theorem 9.1.3 if we set $c = a$. □

Theorem 9.1.5. *Let $A$ be non commutative $D$-algebra. For any $b \in A$, there exists submodule $A[b]$ of $D$-algebra $A$ such that*

$$c \in A[b] \Rightarrow cb = bc$$

Proof. The subset $A[b]$ is not empty because $0 \in A[b]$. Let $d_1$, $d_2 \in D$, $c_1$, $c_2 \in A[b]$. Then

$$\begin{aligned}(d_1 c_1 + d_2 c_2) b &= d_1 c_1 b + d_2 c_2 b = b d_1 c_1 + b d_2 c_2 \\ &= b(d_1 c_1 + d_2 c_2)\end{aligned}$$

And this implies $d_1 c_1 + d_2 c_2 \in A[b]$. □

Theorem 9.1.6. *Let $A$ be Banach $D$-algebra and $a \in A$. The derivative of the map*

$$f : t \in R \to e^{at} \in A$$

*has the following form*

$$\frac{de^{at}}{dt} = e^{at} a = a e^{at} \tag{9.1.12}$$

Proof. According to theorems 8.1.2, [15]-3.3.23,

$$\frac{de^{at}}{dt} = \frac{de^{at}}{dat} \circ \frac{dat}{dt} = \sum_{n=0}^{\infty} \frac{1}{(n+1)!} \sum_{i=0}^{n} ((at)^i \otimes (at)^{n-i}) \circ a \tag{9.1.13}$$

According to the lemma 9.1.2, the equality

$$\begin{aligned}\frac{de^{at}}{dt} &= \sum_{n=0}^{\infty} \frac{1}{(n+1)!} \sum_{i=0}^{n} a^i t^i a a^{n-i} t^{n-i} = a \sum_{n=0}^{\infty} \frac{1}{(n+1)!} \sum_{i=0}^{n} a^n t^n \\ &= a \sum_{n=0}^{\infty} \frac{1}{n!} a^n t^n\end{aligned} \tag{9.1.14}$$

follows from the equality (9.1.13). The equality (9.1.12) follows from equalities (9.1.1), (9.1.14). □

Theorem 9.1.7. *Let $A$ be Banach $D$-algebra and $a \in A$. For any $A$-number $c$, the map*

$$x = e^{at}c \tag{9.1.15}$$

*is solution of the differential equation*

$$\frac{dx}{dt} = ax \tag{9.1.16}$$

Proof. The equality

$$\frac{dx}{dt} = \frac{de^{at}c}{dt} = \frac{de^{at}}{dt}c = ae^{at}c = ax$$

follows from theorems A.1.3, 9.1.6. □

Theorem 9.1.8. *Let $A$ be Banach $D$-algebra and $a \in A$. For any $A$-numbers $c_1$, $c_2$, the map*

$$x = c_1e^{at}c_2 \tag{9.1.17}$$

*is solution of the differential equation (9.1.16) iff $c_1 \in A[a]$ and, therefore, the map (9.1.17) has the form (9.1.15).*

Proof. The equality

$$\frac{dx}{dt} = \frac{dc_1e^{at}c_2}{dt} = c_1\frac{de^{at}}{dt}c_2 = c_1ae^{at}c_2$$

follows from theorems A.1.3, 9.1.6. If $c_1 \notin A[a]$, then the condition

$$c_1a = ac_1 \tag{9.1.18}$$

is not true and the map (9.1.17) is not a solution of the differential equation (9.1.16). However, if the condition (9.1.18) is true, then we the map (9.1.17) gets form

$$x = c_1e^{at}c_2 = e^{at}c_1c_2 \tag{9.1.19}$$

and is the map of the form (9.1.15). □

Theorem 9.1.9. *Let $A$ be Banach $D$-algebra and $a \in A$. For any $A$-number $c$, the map*

$$x = ce^{at} \tag{9.1.20}$$

*is solution of the differential equation*

$$\frac{dx}{dt} = xa \tag{9.1.21}$$

Proof. The equality

$$\frac{dx}{dt} = \frac{dce^{at}}{dt} = c\frac{de^{at}}{dt} = ce^{at}a = xa$$

follows from theorems A.1.3, 9.1.6. □

Theorem 9.1.10. *Let $A$ be Banach $D$-algebra and $a \in A$. For any $A$-numbers $c_1$, $c_2$, the map*

$$x = c_1e^{at}c_2 \tag{9.1.22}$$

*is solution of the differential equation (9.1.21) iff $c_2 \in A[a]$ and, therefore, the map (9.1.22) has the form (9.1.20).*

PROOF. The equality

$$\frac{dx}{dt} = \frac{dc_1e^{at}c_2}{dt} = c_1\frac{de^{at}}{dt}c_2 = c_1e^{at}ac_2$$

follows from theorems A.1.3, 9.1.6. If $c_2 \notin A[a]$, then the condition

$$c_2a = ac_2 \tag{9.1.23}$$

is not true and the map (9.1.22) is not a solution of the differential equation (9.1.21). However, if the condition (9.1.23) is true, then we the map (9.1.22) gets form

$$x = c_1e^{at}c_2 = c_1c_2e^{at} \tag{9.1.24}$$

and is the map of the form (9.1.20). □

## 9.2. Quasiexponent

THEOREM 9.2.1. *The map*

$$x = e[c]^{at}$$

*has following Taylor series expansion*

$$e[c]^{at} = \sum_{n=0}^{\infty} \frac{t^n}{(n+1)!} \sum_{m=0}^{n} a^m ca^{n-m} \tag{9.2.1}$$

PROOF. The theorem follows from the theorem 8.3.2 and the lemma 9.1.2. □

THEOREM 9.2.2. *Let* $c \in A[a]$. *Then*

$$e[c]^{at} = ce^{at} \tag{9.2.2}$$

PROOF. The theorem follows from the theorem 9.2.1. □

THEOREM 9.2.3.

$$\frac{de[c]^{at}}{dt} = \sum_{n=0}^{\infty} \frac{t^n}{(n+2)!} \sum_{m=0}^{n+1} a^m ca^{n+1-m} \tag{9.2.3}$$

PROOF. The theorem follows from the theorem 9.2.1. □

## 9.3. Differential Equation $\frac{dx}{dt} = a_*{}^*x$

Let $A$ be Banach $D$-algebra. The system of differential equations

$$\begin{aligned} \frac{dx^1}{dt} &= a^1_1x^1 + ... + a^1_nx^n \\ &...... \\ \frac{dx^n}{dt} &= a^n_1x^1 + ... + a^n_nx^n \end{aligned} \tag{9.3.1}$$

where $a^i_j \in A$ and $x^i : R \to A$ is $A$-valued function of real variable, is called homogeneous system of linear differential equations.

Let

$$x = \begin{pmatrix} x^1 \\ ... \\ x^n \end{pmatrix} \qquad \frac{dx}{dt} = \begin{pmatrix} \frac{dx^1}{dt} \\ ... \\ \frac{dx^n}{dt} \end{pmatrix}$$

$$a=\begin{pmatrix}a_1^1 & ... & a_n^1\\ ... & ... & ...\\ a_1^n & ... & a_n^n\end{pmatrix}$$

Then we can write system of differential equations5 (9.3.1) in matrix form

$$\frac{dx}{dt}=a_*{}^*x \tag{9.3.2}$$

**9.3.1. Solution as exponent $x=e^{bt}c$ .** We will look for a solution of the system of differential equations (9.3.2) in the form of an exponent

$$x=e^{bt}c=\begin{pmatrix}e^{bt}c^1\\ ...\\ e^{bt}c^n\end{pmatrix}\qquad c=\begin{pmatrix}c^1\\ ...\\ c^n\end{pmatrix} \tag{9.3.3}$$

According to the theorem 9.1.6, the equality

$$\frac{dx^i}{dt}=e^{bt}bc^i=be^{bt}c^i=bx^i \tag{9.3.4}$$

follows from the equalitiy (9.3.3). The equality

$$e^{bt}bc=be^{bt}c=a_*{}^*(e^{bt}c) \tag{9.3.5}$$

follows from equalities (9.3.2), (9.3.4).

THEOREM 9.3.1. *Let $b$ be ${}_*{}^*$-eigenvalue of the matrix $a$. The condition*

$$a_j^i\in A[b]\quad i=1,...,n\quad j=1,...,n \tag{9.3.6}$$

*implies that the matrix of maps (9.3.3) is solution of the system of differential equations (9.3.2) for ${}_*{}^*$-eigencolumn $c$.*

PROOF. According to the theorem 9.1.3, the equality

$$e^{bt}bc=e^{bt}(a_*{}^*c) \tag{9.3.7}$$

follows from the equality (9.3.5). Since the expression $e^{bt}$, in general, is different from 0, the equality

$$bc=a_*{}^*c$$

follows from the equality (9.3.7). According to the definition 5.6.1, $A$-number $b$ is ${}_*{}^*$-eigenvalue of the matrix $a$ and the matrix $c$ is ${}_*{}^*$-eigencolumn of matrix $a$ corresponding to ${}_*{}^*$-eigenvalue $b$. □

THEOREM 9.3.2. *Let $b$ be ${}_*{}^*$-eigenvalue of the matrix $a$. Let entries of matrix $a$ do not satisfy to the condition (9.3.6). If entries of ${}_*{}^*$-eigencolumn $c$ satisfy to the condition*

$$c^i\in A[b]\quad i=1,...,n \tag{9.3.8}$$

*then the matrix of maps (9.3.3) is solution of the system of differential equations (9.3.2).*

PROOF. According to the theorem 9.1.3, the equality

$$bce^{bt} = (a_*{}^*c)e^{bt} \tag{9.3.9}$$

follows from the equality (9.3.5). Since the expression $e^{bt}$, in general, is different from 0, the equality

$$bc = a_*{}^*c$$

follows from the equality (9.3.9). According to the definition 5.6.1, $A$-number $b$ is ${}_*{}^*$-eigenvalue of the matrix $a$ and the matrix $c$ is ${}_*{}^*$-eigencolumn of matrix $a$ corresponding to ${}_*{}^*$-eigenvalue $b$. □

Let entries of matrix $a$ do not satisfy to the condition (9.3.6). Let entries of ${}_*{}^*$-eigencolumn $c$ do not satisfy to the condition (9.3.8). Then the matrix of maps (9.3.3) is not solution of the system of differential equations (9.3.2).

**9.3.2. Solution as exponent $x = ce^{bt}$ .** In the paper [17], on page 35, authors suggest to consider solution as exponent

$$x = ce^{bt} \tag{9.3.10}$$

The equality

$$a_*{}^*(ce^{bt}) = cbe^{bt} \tag{9.3.11}$$

follows from equalities (9.3.2), (9.3.10). Since the expression $e^{bt}$, in general, is different from 0, the equality

$$a_*{}^*c = cb \tag{9.3.12}$$

follows from the equality (9.3.11).

Based on the equality (9.3.12), authors of papers [16], [17] introduce the definition of right eigenvalue versus left eigenvalue which we define by the equality

$$bc = a_*{}^*c$$

According to the theorem 9.1.8, if the matrix of maps (9.3.10) is a solution of the system of differential equations (9.3.2), then vector $c$ satisfies to the condition

$$c^{\boldsymbol{i}} \in A[b] \qquad \boldsymbol{i} = \boldsymbol{1}, ..., \boldsymbol{n}$$

In such case, we can consider matrix of maps

$$x = e^{bt}c$$

instead of the matrix of maps (9.3.10) and we do not need to consider the definition of right eigenvalue.

**9.3.3. Method of successive differentiation.**

THEOREM 9.3.3. *Differentiating one after another system of differential equations (9.3.2) we get the chain of systems of differential equations*

$$\frac{d^nx}{dt^n} = a^{n_*{}^*}{}_*{}^*x \tag{9.3.13}$$

PROOF. We will prove the theorem by induction over $n$.

According to equalities (3.2.17), (9.3.2), the theorem is true for $n = 1$.

Let the theorem be true for $n = k$

$$\frac{d^kx}{dt^k} = a^{k_*{}^*}{}_*{}^*x \tag{9.3.14}$$

According to the definition [15]-4.1.4, the equality

$$\frac{d^{k+1}x}{dt^{k+1}}=\frac{d}{dt}\frac{d^kx}{dt^k}=\frac{d}{dt}(a^{k_*{}^*}{}_*{}^*x)=a^{k_*{}^*}{}_*{}^*\frac{dx}{dt}=a^{k_*{}^*}{}_*{}^*a_*{}^*x=a^{k+1_*{}^*}{}_*{}^*x$$

follows from equalities (3.2.16), (9.3.2), (9.3.14). Therefore, the theorem is true for $n=k+1$. □

THEOREM 9.3.4. *The solution of the system of differential equations* (9.3.2) *with initial condition*

$$t=0\qquad x=\begin{pmatrix}x^1\\ ...\\ x^n\end{pmatrix}=c=\begin{pmatrix}c^1\\ ...\\ c^n\end{pmatrix}$$

*has the following form*

$$x=e^{ta_*{}^*}{}_*{}^*c \tag{9.3.15}$$

PROOF. According to the theorem 9.3.3 and the definition of Taylor series in the section [15]-4.2, Taylor series for the solution has form

$$x=\sum_{n=0}^{\infty}\frac{1}{n!}t^na^{n_*{}^*}{}_*{}^*c \tag{9.3.16}$$

The equality

$$x=\sum_{n=0}^{\infty}\frac{1}{n!}(ta)^{n_*{}^*}{}_*{}^*c \tag{9.3.17}$$

follows from the equality (9.3.16). The equality (9.3.15) follows from equalities (8.5.1), (9.3.17). □

## 9.4. Hyperbolic Trigonometry

Consider the system of differential equations

$$\begin{aligned}\frac{dx^1}{dt}&=x^2\\ \frac{dx^2}{dt}&=x^1\end{aligned} \tag{9.4.1}$$

The matrix $a$ has form

$$a=\begin{pmatrix}0&1\\1&0\end{pmatrix}$$

Since entries of matrix $a$ are real numbers, then the equation to find eigenvalue is

$$\begin{vmatrix}-b&1\\1&-b\end{vmatrix}=b^2-1=0 \tag{9.4.2}$$

It is evident that roots of the equation (9.4.2) depend on choice of $D$-algebra $A$.

In quaternion algebra, according to the statement [14]-5.5.1, the equation (9.5.2) has roots $b_1=-1$, $b_2=1$.

Coefficients $c_1^{\mathit{1}}$, $c_1^{\mathit{2}}$ which correspond to eigenvalue $b_1 = -1$, satisfy to the equation

$$c_1^{\mathit{1}} + c_1^{\mathit{2}} = 0$$

Therefore, corresponding solution of the system of differential equations (9.4.1) has form

$$x^{\mathit{1}} = e^{-t} \qquad x^{\mathit{2}} = -e^{-t} \tag{9.4.3}$$

Coefficients $c_2^{\mathit{1}}$, $c_2^{\mathit{2}}$ which correspond to eigenvalue $b_2 = 1$, satisfy to the equation

$$-c_2^{\mathit{1}} + c_2^{\mathit{2}} = 0$$

Therefore, corresponding solution of the system of differential equations (9.4.1) has form

$$x^{\mathit{1}} = e^{t} \qquad x^{\mathit{2}} = e^{t} \tag{9.4.4}$$

We can write the solution of the system of differential equations (9.4.1) as linear combination of solutions (9.4.3), (9.4.4)

$$\begin{aligned} x^{\mathit{1}} &= C_1 e^t + C_2 e^{-t} \\ x^{\mathit{2}} &= C_1 e^t - C_2 e^{-t} \end{aligned} \tag{9.4.5}$$

The solution of the system of differential equations (9.4.1) which satisfies to initial condition

$$t = 0 \quad x^{\mathit{1}} = 0 \quad x^{\mathit{2}} = 1$$

is unique and has form

$$\begin{aligned} x^{\mathit{1}} &= \frac{1}{2}e^t - \frac{1}{2}e^{-t} \\ x^{\mathit{2}} &= \frac{1}{2}e^t + \frac{1}{2}e^{-t} \end{aligned} \tag{9.4.6}$$

We can solve the system of differential equations (9.4.1) using method of successive differentiation. To this end we consider powers of the matrix $a$

$$a^{1_*{}^*} = \begin{pmatrix} 0 & 1 \\ 1 & 0 \end{pmatrix} \qquad a^{2_*{}^*} = \begin{pmatrix} 1 & 0 \\ 0 & 1 \end{pmatrix} \qquad a^{3_*{}^*} = \begin{pmatrix} 0 & 1 \\ 1 & 0 \end{pmatrix} \tag{9.4.7}$$

The equality

$$\begin{aligned} x = \Bigg( &\begin{pmatrix} 1 & 0 \\ 0 & 1 \end{pmatrix} + t \begin{pmatrix} 0 & 1 \\ 1 & 0 \end{pmatrix} \\ &+ \frac{1}{2!}t^2 \begin{pmatrix} 1 & 0 \\ 0 & 1 \end{pmatrix} + \frac{1}{3!}t^3 \begin{pmatrix} 0 & 1 \\ 1 & 0 \end{pmatrix} + \dots \Bigg) {}_*{}^* \begin{pmatrix} 0 \\ 1 \end{pmatrix} \end{aligned} \tag{9.4.8}$$

follows from equalities (9.3.16), (9.4.7). The equality

$$\begin{aligned} x^{\mathit{0}} &= \sum_{n=0}^{\infty} \frac{1}{(2n+1)!} t^{2n+1} = \sinh t \\ x^{\mathit{1}} &= \sum_{n=0}^{\infty} \frac{1}{(2n)!} t^{2n} = \cosh t \end{aligned} \tag{9.4.9}$$

follows from the equality (9.4.8). From the equalities (9.4.6), (9.4.9), it follows that Euler's formula,

$$\begin{aligned}\sinh t &= \frac{1}{2}e^t - \frac{1}{2}e^{-t}\\ \cosh t &= \frac{1}{2}e^t + \frac{1}{2}e^{-t}\end{aligned} \tag{9.4.10}$$

is true for $t \in R$. The equality

$$\begin{aligned}\frac{d\sinh t}{dt} &= \cosh t\\ \frac{d\cosh t}{dt} &= \sinh t\end{aligned} \tag{9.4.11}$$

follows from equalities (9.4.1), (9.4.9).

To see if this formula is correct for any quaternion $f \neq 0$, consider now the system of differential equations

$$\begin{aligned}\frac{dx^1}{dt} &= fx^2\\ \frac{dx^2}{dt} &= fx^1\end{aligned} \tag{9.4.12}$$

The matrix $a$ has form

$$a = \begin{pmatrix} 0 & f\\ f & 0\end{pmatrix}$$

Now entries of matrix $a$ are quaternions, not real numbers. So we cannot use determinant. To answer the question when matrix $a - Eb$ is singular, we calculate the quasideterminants of first rows[9.1]

$$\det({}_*{}^*)^1_1\, a = -b + fb^{-1}f = 0 \tag{9.4.13}$$

$$\det({}_*{}^*)^1_2\, a = f - bf^{-1}b = 0 \tag{9.4.14}$$

It is not difficult to see that equations (9.4.13) and (9.4.14) are equivalent. The equation (9.4.14) is quadratic equation. In the equation (9.4.14), I will write expression $f + g$ instead of unknown value $b$ and will try to solve the equation

$$f = (f+g)f^{-1}(f+g) = (f+g)(1+f^{-1}g) = f + g + ff^{-1}g + gf^{-1}g \tag{9.4.15}$$

with respect to unknown value $g$. The equation

$$0 = 2g + gf^{-1}g \tag{9.4.16}$$

follows from the equation (9.4.15). The equation

$$g(2 + f^{-1}g) = 0 \tag{9.4.17}$$

follows from the equation (9.4.16). From the equation (9.4.17), it follows that $g = 0$ and $g = -2f$ are roots of the equation (9.4.16). Therefore, $b = f$ and $b = -f$ are roots of the equation (9.4.14).

Coefficients $c^1_1$, $c^2_1$ which correspond to eigenvalue $b_1 = -f$, satisfy to the equation

$$c^1_1 + c^2_1 = 0$$

---

[9.1] See the definition 3.3.3 for quasideterminants. See the theorem 5.5.3 how to find rank of the matrix.

Therefore, corresponding solution of the system of differential equations (9.4.12) has form

$$x^{1} = e^{-ft} \qquad x^{2} = -e^{-ft} \tag{9.4.18}$$

Coefficients $c_2^1$, $c_2^2$ which correspond to eigenvalue $b_2 = f$, satisfy to the equation

$$-c_2^1 + c_2^2 = 0$$

Therefore, corresponding solution of the system of differential equations (9.4.12) has form

$$x^{1} = e^{ft} \qquad x^{2} = e^{ft} \tag{9.4.19}$$

We can write the solution of the system of differential equations (9.4.12) as linear combination of solutions (9.4.18), (9.4.19)

$$\begin{aligned} x^{1} &= C_1 e^{ft} + C_2 e^{-ft} \\ x^{2} &= C_1 e^{ft} - C_2 e^{-ft} \end{aligned} \tag{9.4.20}$$

The solution of the system of differential equations (9.4.12) which satisfies to initial condition

$$t = 0 \quad x^{1} = 0 \quad x^{2} = 1$$

is unique and has form

$$\begin{aligned} x^{1} &= \frac{1}{2}e^{ft} - \frac{1}{2}e^{-ft} \\ x^{2} &= \frac{1}{2}e^{ft} + \frac{1}{2}e^{-ft} \end{aligned} \tag{9.4.21}$$

We can solve the system of differential equations (9.4.12) using method of successive differentiation. To this end we consider powers of the matrix $a$

$$a^{1_*{}^*} = \begin{pmatrix} 0 & f \\ f & 0 \end{pmatrix} \qquad a^{2_*{}^*} = \begin{pmatrix} f^2 & 0 \\ 0 & f^2 \end{pmatrix} \qquad a^{3_*{}^*} = \begin{pmatrix} 0 & f^3 \\ f^3 & 0 \end{pmatrix} \tag{9.4.22}$$

The equality

$$\begin{aligned} x = &\left( \begin{pmatrix} 1 & 0 \\ 0 & 1 \end{pmatrix} + t \begin{pmatrix} 0 & f \\ f & 0 \end{pmatrix} \right. \\ &\left. + \frac{1}{2!}t^2 \begin{pmatrix} f^2 & 0 \\ 0 & f^2 \end{pmatrix} + \frac{1}{3!}t^3 \begin{pmatrix} 0 & f^3 \\ f^3 & 0 \end{pmatrix} + \dots \right) {}_*{}^* \begin{pmatrix} 0 \\ 1 \end{pmatrix} \end{aligned} \tag{9.4.23}$$

follows from equalities (9.3.16), (9.4.22). The equality

$$\begin{aligned} x^{0} &= \sum_{n=0}^{\infty} \frac{1}{(2n+1)!} t^{2n+1} f^{2n+1} = \sum_{n=0}^{\infty} \frac{1}{(2n+1)!} (tf)^{2n+1} = \sinh tf \\ x^{1} &= \sum_{n=0}^{\infty} \frac{1}{(2n)!} t^{2n} f^{2n} = \sum_{n=0}^{\infty} \frac{1}{(2n)!} (tf)^{2n} = \cosh tf \end{aligned} \tag{9.4.24}$$

follows from the equality (9.4.23). From the equalities (9.4.21), (9.4.24), it follows that Euler's formula,

$$\begin{aligned}\sinh ft &= \frac{1}{2}e^{ft} - \frac{1}{2}e^{-ft}\\ \cosh ft &= \frac{1}{2}e^{ft} + \frac{1}{2}e^{-ft}\end{aligned} \tag{9.4.25}$$

is true for $tf \in A$. The equality

$$\begin{aligned}\frac{d\sinh tf}{dt} &= f\cosh tf\\ \frac{d\cosh tf}{dt} &= f\sinh tf\end{aligned} \tag{9.4.26}$$

follows from equalities (9.4.12), (9.4.24). Let $t = 1$. From the equality

$$\begin{aligned}\sinh f &= \frac{1}{2}e^{f} - \frac{1}{2}e^{-f}\\ \cosh f &= \frac{1}{2}e^{f} + \frac{1}{2}e^{-f}\end{aligned} \tag{9.4.27}$$

it follows that Euler's formula is true for $f \in A$.

THEOREM 9.4.1.

$$\sinh tf\, f = f\sinh tf \tag{9.4.28}$$

$$\cosh tf\, f = f\cosh tf \tag{9.4.29}$$

PROOF. Equalities

$$\begin{aligned}\sinh tf\, f &= \sum_{n=0}^{\infty}\frac{1}{(2n+1)!}t^{2n+1}f^{2n+1}f = f\sum_{n=0}^{\infty}\frac{1}{(2n+1)!}t^{2n+1}f^{2n+1}\\ &= f\sinh tf\end{aligned} \tag{9.4.30}$$

$$\begin{aligned}\cosh tf\, f &= \sum_{n=0}^{\infty}\frac{1}{(2n)!}t^{2n}f^{2n}f = f\sum_{n=0}^{\infty}\frac{1}{(2n)!}t^{2n}f^{2n}\\ &= f\cosh tf\end{aligned} \tag{9.4.31}$$

follow from the equality (9.4.24). The equality (9.4.28) follows from the equality (9.4.30). The equality (9.4.29) follows from the equality (9.4.31). □

## 9.5. Elliptical Trigonometry

Consider the system of differential equations

$$\begin{aligned}\frac{dx^1}{dt} &= x^2\\ \frac{dx^2}{dt} &= -x^1\end{aligned} \tag{9.5.1}$$

The matrix $a$ has form

$$a = \begin{pmatrix}0 & 1\\ -1 & 0\end{pmatrix}$$

Since entries of matrix $a$ are real numbers, then the equation to find eigenvalue is

$$\begin{vmatrix} -b & 1 \\ -1 & -b \end{vmatrix} = b^2+1=0 \tag{9.5.2}$$

It is evident that roots of the equation (9.5.2) depend on choice of $D$-algebra $A$.

In quaternion algebra, according to the statement [14]-5.5.3, the equation (9.5.2) has infinitely many roots

$$b=b^1i+b^2j+b^3k$$

such that

$$(b^1)^2+(b^2)^2+(b^3)^2=1$$

Coefficients $c_b^1$, $c_b^2$ which correspond to given eigenvalue $b$, satisfy to the equation

$$-bc_b^1+c_b^2=0$$

Therefore, corresponding solution of the system of differential equations (9.5.1) has form

$$x^1=e^{bt} \qquad x^2=be^{bt} \tag{9.5.3}$$

If we want to find the solution of the system of differential equations (9.5.1) which satisfies to initial condition

$$t=0 \quad x^1=0 \quad x^2=1$$

then we see that we have too many choices.

Linear combination of two solutions of the system of differential equations (9.5.1) is solution of the system of differential equations (9.5.1). We will start from consideration of linear combination of two solutions of the form (9.5.3) because, in such case, constants of linear combination are unique.

Thus, we search solution of the form

$$\begin{aligned} x^1&=C_1e^{b_1t}+C_2e^{b_2t}\\ x^2&=C_1b_1e^{b_1t}+C_2b_2e^{b_2t} \end{aligned} \tag{9.5.4}$$

According to initial condition, the system of equations

$$\begin{aligned} C_1+C_2&=0\\ C_1b_1+C_2b_2&=1 \end{aligned} \tag{9.5.5}$$

follows from the equality (9.5.4). The equality

$$\begin{aligned} C_2&=-C_1\\ C_1(b_1-b_2)&=1 \end{aligned} \tag{9.5.6}$$

follows from the system of equations (9.5.5).

Example 9.5.1. *Let* $b_1=i$, $b_2=j$. *The equality*

$$C_1(i-j)=1 \tag{9.5.7}$$

*follows from the system of equations* (9.5.6). *The equality*

$$C_1=\frac{1}{2}(-i+j) \tag{9.5.8}$$

*follows from the equality* (9.5.7)*. Our goal is to verify whether map*

(9.5.9)
$$x^1 = \frac{1}{2}(-i+j)(e^{it}-e^{jt})$$
$$x^2 = \frac{1}{2}(-i+j)(ie^{it}-je^{jt})$$

*is solution of system of differential equations* (9.3.2)*. The equality*

(9.5.10)
$$\frac{dx^1}{dt} = \frac{1}{2}(-i+j)(ie^{it}-je^{jt}) = x^2$$
$$\frac{dx^2}{dt} = \frac{1}{2}(-i+j)(i^2e^{it}-j^2e^{jt}) = -x^1$$

*follows from the equality* (9.5.9)*. Therefore, the map* (9.5.9) *is solution of system of differential equations* (9.3.2) *which satisfies to initial condition*

$$t=0 \quad x^1=0 \quad x^2=1$$

□

If we consider linear combination of three or more solutions of the form (9.5.3), then coefficients of linear combination will be determined ambiguously because corresponding system of linear equations is overdetermined.

Example 9.5.2. *Let* $b_1=i$, $b_2=j$, $b_3=k$. *Thus, we search solution of the form*

(9.5.11)
$$x^1 = C_1e^{it}+C_2e^{jt}+C_3e^{kt}$$
$$x^2 = C_1ie^{it}+C_2je^{jt}+C_3ke^{kt}$$

*According to initial condition, the system of equations*

(9.5.12)
$$C_1+C_2+C_3=0$$
$$C_1i+C_2j+C_3k=1$$

*follows from the equality* (9.5.11)*. The equality*

(9.5.13)
$$C_1=C$$
$$C_2=-C-C_3$$
$$Ci-(C+C_3)j+C_3k=1$$

*follows from the system of equations* (9.5.12)*. The equality*

(9.5.14)
$$C_3(k-j)=1-Ci+Cj$$

*follows from the equality* (9.5.13)*. The equality*

(9.5.15)
$$C_3=\frac{1}{2}(1+C(j-i))(j-k)=\frac{1}{2}((j-k)+C(j-i)(j-k))$$
$$=\frac{1}{2}((j-k)+C(j(j-k)-i(j-k)))$$

*follows from the equality* (9.5.14)*. The equality*

(9.5.16)
$$C_1=C$$
$$C_2=\frac{1}{2}(-(j-k)+C(-1+i+j+k))$$
$$C_3=\frac{1}{2}((j-k)-C(1+i+j+k))$$

*follows from equalities* (9.5.13), (9.5.15). *The equality*

$$\begin{aligned} x^1 &= Ce^{it} + \frac{1}{2}(-(j-k) + C(-1+i+j+k))e^{jt} \\ &\quad + \frac{1}{2}((j-k) - C(1+i+j+k))e^{kt} \\ x^2 &= Cie^{it} + \frac{1}{2}(-(j-k) + C(-1+i+j+k))je^{jt} \\ &\quad + \frac{1}{2}((j-k) - C(1+i+j+k))ke^{kt} \end{aligned} \tag{9.5.17}$$

*follows from equalities* (9.5.11), (9.5.16). □

QUESTION 9.5.3. *Because the set of solutions of system of differential equations* (9.4.12) *is too large, I do not consider Euler's formula. However we need to understand why the same system of differential equations has infinitely many solutions which satisfy to given initial condition; and at the same time solution of this system of differential equations has unique Taylor series expansion. Is this related to identities that I did not consider?* □

## 9.6. Differential Equation $\frac{dx}{dt} = x^*{}_*a$

Let $A$ be Banach $D$-algebra. The system of differential equations

$$\begin{aligned} \frac{dx^1}{dt} &= x^1 a^1_1 + ... + x^n a^1_n \\ &...... \\ \frac{dx^n}{dt} &= x^1 a^n_1 + ... + x^n a^n_n \end{aligned} \tag{9.6.1}$$

where $a^i_j \in A$ and $x^i : R \to A$ is $A$-valued function of real variable, is called homogeneous system of linear differential equations.

Let

$$x = \begin{pmatrix} x^1 \\ ... \\ x^n \end{pmatrix} \qquad \frac{dx}{dt} = \begin{pmatrix} \frac{dx^1}{dt} \\ ... \\ \frac{dx^n}{dt} \end{pmatrix}$$

$$a = \begin{pmatrix} a^1_1 & ... & a^1_n \\ ... & ... & ... \\ a^n_1 & ... & a^n_n \end{pmatrix}$$

Then we can write system of differential equations5 (9.6.1) in matrix form

$$\frac{dx}{dt} = x^*{}_*a \tag{9.6.2}$$

**9.6.1. Solution as exponent** $x = ce^{bt}$ **.** We will look for a solution of the system of differential equations (9.6.2) in the form of an exponent

$$x = ce^{bt} = \begin{pmatrix} c^1 e^{bt} \\ ... \\ c^n e^{bt} \end{pmatrix} \qquad c = \begin{pmatrix} c^1 \\ ... \\ c^n \end{pmatrix} \tag{9.6.3}$$

According to the theorem 9.1.6, the equality

$$\frac{dx^i}{dt} = c^i e^{bt} b = c^i b e^{bt} = x^i b \tag{9.6.4}$$

follows from the equality (9.6.3). The equality

$$ce^{bt} b = cbe^{bt} = (ce^{bt}){}^*{}_* a \tag{9.6.5}$$

follows from equalities (9.6.2), (9.6.4).

THEOREM 9.6.1. *Let $b$ be ${}^*{}_*$-eigenvalue of the matrix $a$. The condition*

$$a^i_j \in A[b] \quad i = 1, ..., n \quad j = 1, ..., n \tag{9.6.6}$$

*implies that the matrix of maps (9.6.3) is solution of the system of differential equations (9.6.2) for ${}^*{}_*$-eigencolumn $c$.*

PROOF. According to the theorem 9.1.3, the equality

$$cbe^{bt} = (c^*{}_* a)e^{bt} \tag{9.6.7}$$

follows from the equality (9.6.5). Since the expression $e^{bt}$, in general, is different from 0, the equality

$$cb = c^*{}_* a$$

follows from the equality (9.6.7). According to the definition 5.6.2, $A$-number $b$ is ${}^*{}_*$-eigenvalue of the matrix $a$ and the matrix $c$ is ${}^*{}_*$-eigencolumn of matrix $a$ corresponding to ${}^*{}_*$-eigenvalue $b$. □

THEOREM 9.6.2. *Let $b$ be ${}^*{}_*$-eigenvalue of the matrix $a$. Let entries of matrix $a$ do not satisfy to the condition (9.6.6). If entries of ${}^*{}_*$-eigencolumn $c$ satisfy to the condition*

$$c^i \in A[b] \quad i = 1, ..., n \tag{9.6.8}$$

*then the matrix of maps (9.6.3) is solution of the system of differential equations (9.6.2).*

PROOF. According to the theorem 9.1.3, the equality

$$e^{bt} cb = e^{bt}(c^*{}_* a) \tag{9.6.9}$$

follows from the equality (9.6.5). Since the expression $e^{bt}$, in general, is different from 0, the equality

$$cb = c^*{}_* a$$

follows from the equality (9.6.9). According to the definition 5.6.2, $A$-number $b$ is ${}^*{}_*$-eigenvalue of the matrix $a$ and the matrix $c$ is ${}^*{}_*$-eigencolumn of matrix $a$ corresponding to ${}^*{}_*$-eigenvalue $b$. □

Let entries of matrix $a$ do not satisfy to the condition (9.6.6). Let entries of ${}_*{}^*$-eigencolumn $c$ do not satisfy to the condition (9.6.8). Then the matrix of maps (9.6.3) is not solution of the system of differential equations (9.6.2).

### 9.6.2. Method of successive differentiation.

THEOREM 9.6.3. *Differentiating one after another system of differential equations* (9.6.2) *we get the chain of systems of differential equations*

$$\frac{d^n x}{dt^n}={x^*}_*{a^n}^*{}_* \tag{9.6.10}$$

PROOF. We will prove the theorem by induction over $n$.
According to equalities (3.2.21), (9.6.2), the theorem is true for $n=1$.
Let the theorem be true for $n=k$

$$\frac{d^k x}{dt^k}={x^*}_*{a^k}^*{}_* \tag{9.6.11}$$

According to the definition [15]-4.1.4, the equality

$$\frac{d^{k+1}x}{dt^{k+1}}=\frac{d}{dt}\frac{d^k x}{dt^k}=\frac{d}{dt}({x^*}_*{a^k}^*{}_*)=\frac{dx}{dt}{}^*{}_*{a^k}^*{}_*={x^*}_*{a^*}_*{a^k}^*{}_*={x^*}_*{a^{k+1}}^*{}_*$$

follows from equalities (3.2.20), (9.6.2), (9.6.11). Therefore, the theorem is true for $n=k+1$. □

THEOREM 9.6.4. *The solution of the system of differential equations* (9.6.2) *with initial condition*

$$t=0\qquad x=\begin{pmatrix}x^1\\ ...\\ x^n\end{pmatrix}=c=\begin{pmatrix}c^1\\ ...\\ c^n\end{pmatrix}$$

*has the following form*

$$x={c^*}_*{e^{ta}}^*{}_* \tag{9.6.12}$$

PROOF. According to the theorem 9.6.3 and the definition of Taylor series in the section [15]-4.2, Taylor series for the solution has form

$$x={c^*}_*\sum_{n=0}^{\infty}\frac{1}{n!}t^n{a^n}^*{}_* \tag{9.6.13}$$

The equality

$$x={c^*}_*\sum_{n=0}^{\infty}\frac{1}{n!}{(ta)^n}^*{}_* \tag{9.6.14}$$

follows from the equality (9.6.13). The equality (9.6.12) follows from equalities (8.5.1), (9.6.14). □

## 9.7. Differential Equation $\frac{dx}{dt}={a^*}_*x$

Let $A$ be Banach $D$-algebra. The system of differential equations

$$\begin{aligned}\frac{dx_1}{dt}&=a_1^1x_1+...+a_1^nx_n\\ &......\\ \frac{dx_n}{dt}&=a_n^1x_1+...+a_n^nx_n\end{aligned} \tag{9.7.1}$$

where $a^i_j \in A$ and $x_i : R \to A$ is $A$-valued function of real variable, is called homogeneous system of linear differential equations.

Let

$$x = \begin{pmatrix} x_1 & ... & x_n \end{pmatrix} \qquad \frac{dx}{dt} = \begin{pmatrix} \frac{dx_1}{dt} & ... & \frac{dx_n}{dt} \end{pmatrix}$$

$$a = \begin{pmatrix} a^1_1 & ... & a^1_n \\ ... & ... & ... \\ a^n_1 & ... & a^n_n \end{pmatrix}$$

Then we can write system of differential equations5 (9.7.1) in matrix form

$$\frac{dx}{dt} = a^*{}_*x \tag{9.7.2}$$

**9.7.1. Solution as exponent $x = e^{bt}c$ .** We will look for a solution of the system of differential equations (9.7.2) in the form of an exponent

$$x = e^{bt}c = \begin{pmatrix} e^{bt}c_1 & ... & e^{bt}c_n \end{pmatrix} \qquad c = \begin{pmatrix} c_1 & ... & c_n \end{pmatrix} \tag{9.7.3}$$

According to the theorem 9.1.6, the equality

$$\frac{dx_i}{dt} = e^{bt}bc_i = be^{bt}c_i = bx_i \tag{9.7.4}$$

follows from the equality (9.7.3). The equality

$$e^{bt}bc = a^*{}_*(e^{bt}c) \tag{9.7.5}$$

follows from equalities (9.7.2), (9.7.4).

Theorem 9.7.1. *Let $b$ be $^*{}_*$-eigenvalue of the matrix $a$. The condition*

$$a^i_j \in A[b] \qquad i = 1, ..., n \quad j = 1, ..., n \tag{9.7.6}$$

*implies that the matrix of maps (9.7.3) is solution of the system of differential equations (9.7.2) for $^*{}_*$-eigenrow $c$.*

Proof. According to the theorem 9.1.3, the equality

$$e^{bt}bc = e^{bt}(a^*{}_*c) \tag{9.7.7}$$

follows from the equality (9.7.5). Since the expression $e^{bt}$, in general, is different from 0, the equality

$$bc = a^*{}_*c$$

follows from the equality (9.7.7). According to the definition 5.6.2, $A$-number $b$ is $^*{}_*$-eigenvalue of the matrix $a$ and the matrix $c$ is $^*{}_*$-eigenrow of matrix $a$ corresponding to $^*{}_*$-eigenvalue $b$. □

Theorem 9.7.2. *Let $b$ be $^*{}_*$-eigenvalue of the matrix $a$. Let entries of matrix $a$ do not satisfy to the condition (9.7.6). If entries of $^*{}_*$-eigenrow $c$ satisfy to the condition*

$$c^i \in A[b] \qquad i = 1, ..., n \tag{9.7.8}$$

*then the matrix of maps (9.7.3) is solution of the system of differential equations (9.7.2).*

Proof. According to the theorem 9.1.3, the equality

$$bce^{bt}=(a^*{}_*c)e^{bt}\tag{9.7.9}$$

follows from the equality (9.7.5). Since the expression $e^{bt}$, in general, is different from 0, the equality

$$bc=a^*{}_*c$$

follows from the equality (9.7.9). According to the definition 5.6.2, $A$-number $b$ is ${}^*{}_*$-eigenvalue of the matrix $a$ and the matrix $c$ is ${}^*{}_*$-eigenrow of matrix $a$ corresponding to ${}^*{}_*$-eigenvalue $b$. □

Let entries of matrix $a$ do not satisfy to the condition (9.7.6). Let entries of ${}_*{}^*$-eigenrow $c$ do not satisfy to the condition (9.7.8). Then the matrix of maps (9.7.3) is not solution of the system of differential equations (9.7.2).

### 9.7.2. Method of successive differentiation.

Theorem 9.7.3. *Differentiating one after another system of differential equations (9.7.2) we get the chain of systems of differential equations*

$$\frac{d^nx}{dt^n}=a^{n^*{}_*{}^*}{}_*x\tag{9.7.10}$$

Proof. We will prove the theorem by induction over $n$.

According to equalities (3.2.21), (9.7.2), the theorem is true for $n=1$.

Let the theorem be true for $n=k$

$$\frac{d^kx}{dt^k}=a^{k^*{}_*{}^*}{}_*x\tag{9.7.11}$$

According to the definition [15]-4.1.4, the equality

$$\frac{d^{k+1}x}{dt^{k+1}}=\frac{d}{dt}\frac{d^kx}{dt^k}=\frac{d}{dt}(a^{k^*{}_*{}^*}{}_*x)=a^{k^*{}_*{}^*}{}_*\frac{dx}{dt}=a^{k^*{}_*{}^*}{}_*a^*{}_*x=a^{k+1^*{}_*{}^*}{}_*x$$

follows from equalities (3.2.20), (9.7.2), (9.7.11). Therefore, the theorem is true for $n=k+1$. □

Theorem 9.7.4. *The solution of the system of differential equations (9.7.2) with initial condition*

$$t=0\qquad x=\begin{pmatrix}x^1\\ \dots\\ x^n\end{pmatrix}=c=\begin{pmatrix}c^1\\ \dots\\ c^n\end{pmatrix}$$

*has the following form*

$$x=e^{ta^*{}_*{}^*}{}_*c\tag{9.7.12}$$

Proof. According to the theorem 9.7.3 and the definition of Taylor series in the section [15]-4.2, Taylor series for the solution has form

$$x=\sum_{n=0}^{\infty}\frac{1}{n!}t^na^{n^*{}_*{}^*}{}_*c\tag{9.7.13}$$

The equality

$$x=\sum_{n=0}^{\infty}\frac{1}{n!}(ta)^{n^*{}_*{}^*}{}_*c\tag{9.7.14}$$

follows from the equality (9.7.13). The equality (9.7.12) follows from equalities (8.5.1), (9.7.14). □

## 9.8. Differential Equation $\frac{dx}{dt}=x_*{}^*a$

Let $A$ be Banach $D$-algebra. The system of differential equations

(9.8.1)
$$\begin{aligned}\frac{dx_1}{dt}&=x_1a^1_1+...+x_na^n_1\\&......\\\frac{dx_n}{dt}&=x_1a^1_n+...+x_na^n_n\end{aligned}$$

where $a^i_j\in A$ and $x^i:R\to A$ is $A$-valued function of real variable, is called homogeneous system of linear differential equations.

Let

$$x=\begin{pmatrix}x_1&...&x_n\end{pmatrix}\quad\frac{dx}{dt}=\begin{pmatrix}\frac{dx_1}{dt}&...&\frac{dx_n}{dt}\end{pmatrix}$$

$$a=\begin{pmatrix}a^1_1&...&a^1_n\\...&...&...\\a^n_1&...&a^n_n\end{pmatrix}$$

Then we can write system of differential equations5 (9.8.1) in matrix form

(9.8.2)
$$\frac{dx}{dt}=x_*{}^*a$$

**9.8.1. Solution as exponent** $x=ce^{bt}$ **.** We will look for a solution of the system of differential equations (9.8.2) in the form of an exponent

(9.8.3)
$$x=ce^{bt}=\begin{pmatrix}c_1e^{bt}&...&c_ne^{bt}\end{pmatrix}\quad c=\begin{pmatrix}c_1&...&c_n\end{pmatrix}$$

According to the theorem 9.1.6, the equality

(9.8.4)
$$\frac{dx_i}{dt}=c_ie^{bt}b=c_ibe^{bt}=bx_i$$

follows from the equality (9.8.3). The equality

(9.8.5)
$$ce^{bt}b=cbe^{bt}=(ce^{bt})_*{}^*a$$

follows from equalities (9.8.2), (9.8.4).

Theorem 9.8.1. *Let $b$ be ${}_*{}^*$-eigenvalue of the matrix $a$. The condition*

(9.8.6)
$$a^i_j\in A[b]\quad i=1,...,n\quad j=1,...,n$$

*implies that the matrix of maps (9.8.3) is solution of the system of differential equations (9.8.2) for ${}_*{}^*$-eigenrow $c$.*

Proof. According to the theorem 9.1.3, the equality

(9.8.7)
$$cbe^{bt}=(c_*{}^*a)e^{bt}$$

follows from the equality (9.8.5). Since the expression $e^{bt}$, in general, is different from 0, the equality

$$cb=c_*{}^*a$$

follows from the equality (9.8.7). According to the definition 5.6.1, $A$-number $b$ is ${}_*{}^*$-eigenvalue of the matrix $a$ and the matrix $c$ is ${}_*{}^*$-eigenrow of matrix $a$ corresponding to ${}_*{}^*$-eigenvalue $b$. □

Theorem 9.8.2. *Let $b$ be ${}_*{}^*$-eigenvalue of the matrix $a$. Let entries of matrix $a$ do not satisfy to the condition (9.8.6). If entries of ${}_*{}^*$-eigenrow $c$ satisfy to the condition*

$$c^i \in A[b] \qquad i = 1, ..., n \tag{9.8.8}$$

*then the matrix of maps (9.8.3) is solution of the system of differential equations (9.8.2).*

Proof. According to the theorem 9.1.3, the equality

$$e^{bt}cb = e^{bt}(c_*{}^*a) \tag{9.8.9}$$

follows from the equality (9.8.5). Since the expression $e^{bt}$, in general, is different from 0, the equality

$$cb = c_*{}^*a$$

follows from the equality (9.8.9). According to the definition 5.6.1, $A$-number $b$ is ${}_*{}^*$-eigenvalue of the matrix $a$ and the matrix $c$ is ${}_*{}^*$-eigenrow of matrix $a$ corresponding to ${}_*{}^*$-eigenvalue $b$. □

Let entries of matrix $a$ do not satisfy to the condition (9.8.6). Let entries of ${}_*{}^*$-eigenrow $c$ do not satisfy to the condition (9.8.8). Then the matrix of maps (9.3.3) is not solution of the system of differential equations (9.8.2).

### 9.8.2. Method of successive differentiation.

Theorem 9.8.3. *Differentiating one after another system of differential equations (9.8.2) we get the chain of systems of differential equations*

$$\frac{d^n x}{dt^n} = x_*{}^*a^{n_*{}^*} \tag{9.8.10}$$

Proof. We will prove the theorem by induction over $n$.

According to equalities (3.2.17), (9.8.2), the theorem is true for $n = 1$.

Let the theorem be true for $n = k$

$$\frac{d^k x}{dt^k} = x_*{}^*a^{k_*{}^*} \tag{9.8.11}$$

According to the definition [15]-4.1.4, the equality

$$\frac{d^{k+1}x}{dt^{k+1}} = \frac{d}{dt}\frac{d^k x}{dt^k} = \frac{d}{dt}(x_*{}^*a^{k_*{}^*}) = \frac{dx}{dt}{}_*{}^*a^{k_*{}^*} = x_*{}^*a_*{}^*a^{k_*{}^*} = x_*{}^*a^{k+1_*{}^*}$$

follows from equalities (3.2.16), (9.8.2), (9.8.11). Therefore, the theorem is true for $n = k + 1$. □

Theorem 9.8.4. *The solution of the system of differential equations (9.8.2) with initial condition*

$$t = 0 \qquad x = \begin{pmatrix} x^1 \\ ... \\ x^n \end{pmatrix} = c = \begin{pmatrix} c^1 \\ ... \\ c^n \end{pmatrix}$$

*has the following form*

$$x = c_*{}^* e^{ta_*{}^*} \tag{9.8.12}$$

PROOF. According to the theorem 9.8.3 and the definition of Taylor series in the section [15]-4.2, Taylor series for the solution has form

$$x = c_*{}^* \sum_{n=0}^{\infty} \frac{1}{n!} t^n a^{n_*{}^*} \tag{9.8.13}$$

The equality

$$x = c_*{}^* \sum_{n=0}^{\infty} \frac{1}{n!} (ta)^{n_*{}^*} \tag{9.8.14}$$

follows from the equality (9.8.13). The equality (9.8.12) follows from equalities (8.5.1), (9.8.14). □

## 9.9. Short Summary

Method of successive differentiation is an important method of solving of differential equations, especially when it is difficult to find another method of solving. However the case is changing, when we consider differential equation over non commutative Banach $D$-algebra.

For the differential equation

$$\frac{dy}{dx} = 3x \otimes x \tag{9.9.1}$$

method of successive differentiation generates the map which is not a solution of the differential equation (9.9.1).

According to the theorem 9.3.4, the solution of the system of differential equations

$$\begin{aligned} \frac{dx^1}{dt} &= x^2 \\ \frac{dx^2}{dt} &= -x^1 \end{aligned} \tag{9.9.2}$$

which satisfies to initial condition

$$t = 0 \quad x^1 = 0 \quad x^2 = 1$$

has the following form

$$\begin{aligned} x^1 &= \sin t \\ x^2 &= \cos t \end{aligned} \tag{9.9.3}$$

According to the example 9.5.1, the map

$$\begin{aligned} x^1 &= \frac{1}{2}(-i + j)(e^{it} - e^{jt}) \\ x^2 &= \frac{1}{2}(-i + j)(ie^{it} - je^{jt}) \end{aligned} \tag{9.9.4}$$

also is the solution of the system of differential equations (9.9.2) which satisfies to initial condition

$$t = 0 \quad x^1 = 0 \quad x^2 = 1$$

However, maps (9.9.3), (9.9.4) have different Taylor series expansion.

From these two extremely opposite examples, it follows that we must be careful when we use method of successive differentiation to solve differential equation over non commutative Banach $D$-algebra.

From the reasoning above, the question follows. If we have an infinite set of solutions, how do we describe this set and how do we choose a solution? To answer this question it would be nice to have more examples how we use considered theory to solve different problems. In particular, I like examples described in the section [17]-1.1.

APPENDIX A

# Summary of Statements

Let $D$ be the complete commutative ring of characteristic 0.

## A.1. Table of Derivatives

THEOREM A.1.1. *For any $b \in A$*

$$\frac{db}{dx} = 0 \otimes 0$$

PROOF. The theorem follows from the theorem [15]-B.1.1. □

THEOREM A.1.2.

$$\frac{dx}{dx} \circ dx = dx \qquad \frac{dx}{dx} = 1 \otimes 1 \tag{A.1.1}$$

PROOF. The theorem follows from the theorem [15]-B.1.2. □

THEOREM A.1.3. *For any $b,\ c \in A$*

$$\left\{\begin{array}{l} \dfrac{dbf(x)c}{dx} = (b \otimes c) \circ \dfrac{df(x)}{dx} \\ \dfrac{dbf(x)c}{dx} \circ dx = b \left(\dfrac{df(x)}{dx} \circ dx\right) c \\ \dfrac{d_{s \cdot 0}bf(x)c}{dx} = b\dfrac{d_{s \cdot 0}f(x)}{dx} \\ \dfrac{d_{s \cdot 1}bf(x)c}{dx} = \dfrac{d_{s \cdot 1}f(x)}{dx}c \end{array}\right. \tag{A.1.2}$$

*For any $F \in A \otimes A$,*

$$\left\{\begin{array}{l} \dfrac{dF \circ f(x)}{dx} = F \circ \dfrac{df(x)}{dx} \\ \dfrac{dF \circ f(x)}{dx} \circ dx = F \circ \left(\dfrac{df(x)}{dx} \circ dx\right) \end{array}\right. \tag{A.1.3}$$

PROOF. The theorem follows from the theorem [15]-B.1.3. □

THEOREM A.1.4. *Let*

$$f : A \to B$$
$$g : A \to B$$

*be maps of Banach $D$-module $A$ into associative Banach $D$-algebra $A$. Since there exist the derivatives $\dfrac{df(x)}{dx}$, $\dfrac{dg(x)}{dx}$, then there exists the derivative $\dfrac{d(f(x)+g(x))}{dx}$*

$$\frac{d(f(x)+g(x))}{dx} \circ dx = \frac{df(x)}{dx} \circ dx + \frac{dg(x)}{dx} \circ dx \tag{A.1.4}$$

(A.1.5)
$$\frac{d(f(x)+g(x))}{dx}=\frac{df(x)}{dx}+\frac{dg(x)}{dx}$$

*Since*

(A.1.6)
$$\frac{df(x)}{dx}=\frac{d_{s\cdot 0}f(x)}{dx}\otimes\frac{d_{s\cdot 1}f(x)}{dx}$$

(A.1.7)
$$\frac{dg(x)}{dx}=\frac{d_{t\cdot 0}g(x)}{dx}\otimes\frac{d_{t\cdot 1}g(x)}{dx}$$

*then*

(A.1.8)
$$\frac{d(f(x)+g(x))}{dx}=\frac{d_{s\cdot 0}f(x)}{dx}\otimes\frac{d_{s\cdot 1}f(x)}{dx}+\frac{d_{t\cdot 0}g(x)}{dx}\otimes\frac{d_{t\cdot 1}g(x)}{dx}$$

PROOF. The theorem follows from the theorem [15]-B.1.4. □

THEOREM A.1.5. *For any $b, c\in A$*

(A.1.9)
$$\begin{cases}\dfrac{dbxc}{dx}=b\otimes c & \dfrac{dbxc}{dx}\circ dx=b\,dx\,c\\ \dfrac{d_{1\cdot 0}bxc}{dx}=b & \dfrac{d_{1\cdot 1}bxc}{dx}=c\end{cases}$$

PROOF. Corollary of theorems A.1.2, A.1.3, when $f(x)=x$. □

THEOREM A.1.6. *Let $f$ be linear map*

$$f\circ x=(a_{s\cdot 0}\otimes a_{s\cdot 1})\circ x=a_{s\cdot 0}\ x\ a_{s\cdot 1}$$

*Then*

$$\frac{\partial f\circ x}{\partial x}=f$$
$$\frac{\partial f\circ x}{\partial x}\circ dx=f\circ dx$$

PROOF. Corollary of theorems A.1.4, A.1.5, [15]-3.3.13. □

COROLLARY A.1.7. *For any $b\in A$*

$$\begin{cases}\dfrac{d(xb-bx)}{dx}=1\otimes b-b\otimes 1\\ \dfrac{d(xb-bx)}{dx}\circ dx=dx\,b-b\,dx\\ \dfrac{d_{1\cdot 0}(xb-bx)}{dx}=1 \quad \dfrac{d_{1\cdot 1}(xb-bx)}{dx}=b\\ \dfrac{d_{2\cdot 0}(xb-bx)}{dx}=-b \quad \dfrac{d_{2\cdot 1}(xb-bx)}{dx}=1\end{cases}$$

□

Theorem A.1.8. *Let $D$ be the complete commutative ring of characteristic* 0*. Let $A$ be associative Banach $D$-algebra. Then*[A.1]

$$\frac{dx^2}{dx}=x\otimes 1+1\otimes x \tag{A.1.10}$$

$$dx^2=x\,dx+dx\,x \tag{A.1.11}$$

$$\left\{\begin{array}{ll}\dfrac{d_{1\cdot 0}x^2}{dx}=x & \dfrac{d_{1\cdot 1}x^2}{dx}=1\\ \dfrac{d_{2\cdot 0}x^2}{dx}=1 & \dfrac{d_{2\cdot 1}x^2}{dx}=x\end{array}\right. \tag{A.1.12}$$

Proof. The theorem follows from the theorem [15]-B.1.15. □

Theorem A.1.9. *Let $D$ be the complete commutative ring of characteristic* 0*. Let $A$ be associative Banach $D$-algebra. Then*

$$\frac{dx^3}{dx}=x^2\otimes 1+x\otimes x+1\otimes x^2 \tag{A.1.13}$$

$$dx^3=x^2\,dx+x\,dx\,x+dx\,x^2 \tag{A.1.14}$$

Proof. According to the theorem A.1.10,

$$\frac{dx^3}{dx}=\frac{dx^2x}{dx}=\frac{dx^2}{dx}x+x^2\frac{dx}{dx}=(x\otimes 1+1\otimes x)x+x^2(1\otimes 1) \tag{A.1.15}$$

The equality (A.1.13) follows from the equality (A.1.15). The equality (A.1.14) follows from the equality (A.1.13) and the definition [15]-3.3.2. □

Theorem A.1.10. *Let $A$ be Banach $D$-module. Let $B$ be Banach $D$-algebra. Let $f$, $g$ be differentiable maps*

$$f:A\to B\quad g:A\to B$$

*The derivative satisfies to relationship*

$$\frac{df(x)g(x)}{dx}\circ dx=\left(\frac{df(x)}{dx}\circ dx\right)\,g(x)+f(x)\,\left(\frac{dg(x)}{dx}\circ dx\right)$$

$$\frac{df(x)g(x)}{dx}=\frac{df(x)}{dx}\,g(x)+f(x)\,\frac{dg(x)}{dx} \tag{A.1.16}$$

Proof. The theorem follows from the theorem [15]-3.3.17. □

---

[A.1] The statement of the theorem is similar to example VIII, [22], p. 451. If product is commutative, then the equality (A.1.10) gets form

$$dx^2\circ dx=2x\;dx$$

$$\frac{dx^2}{dx}=2x$$

THEOREM A.1.11. *Let $A$ be Banach $D$-module. Let $B$, $C$ be Banach $D$-algebras. Let $f$, $g$ be differentiable maps*

$$f:A\to B \quad g:A\to C$$

*The derivative satisfies to relationship*

$$\frac{df(x)\otimes g(x)}{dx}\circ a=\left(\frac{df(x)}{dx}\circ a\right)\otimes g(x)+f(x)\otimes\left(\frac{dg(x)}{dx}\circ a\right)$$

$$\frac{df(x)\otimes g(x)}{dx}=\frac{df(x)}{dx}\otimes g(x)+f(x)\otimes\frac{dg(x)}{dx}$$

PROOF. The theorem follows from the theorem [15]-3.3.20. □

THEOREM A.1.12. *Let $D$-algebra $A$ have unit $e$. Let*

$$f:A\to A$$

*be differentiable map and*

$$P_0:d\in D\to de_0\in A$$

*be immersion of the ring $D$ into $D$-algebra $A$. Then*

$$\frac{df(te)}{dt}=\frac{df(x)}{dx}\circ e \tag{A.1.17}$$

PROOF. According to the theorem [15]-3.3.23,

$$\frac{df(te)}{dt}=\frac{df(te)}{d\,te}\circ\frac{d\,te}{dt} \tag{A.1.18}$$

The equality (A.1.17) follows from the equality (A.1.18). □

THEOREM A.1.13. *Let $D$ be the complete commutative ring of characteristic* 0. *Let $A$ be associative Banach $D$-algebra. Then*

$$\frac{d^n}{dx^n}\frac{d^my}{dx^m}=\frac{d^{n+m}y}{dx^{n+m}} \tag{A.1.19}$$

PROOF. We will prove the theorem by induction over $n$.

The theorem is evident for $n=0$. For $n=1$, the theorem follows from the definition [15]-4.1.4.

Let the theorem be true for $n=k$

$$\frac{d^k}{dx^k}\frac{d^my}{dx^m}=\frac{d^{k+m}y}{dx^{k+m}} \tag{A.1.20}$$

The equality

$$\frac{d^{k+1}}{dx^{k+1}}\frac{d^my}{dx^m}=\frac{d}{dx}\left(\frac{d^k}{dx^k}\frac{d^my}{dx^m}\right)=\frac{d}{dx}\frac{d^{k+m}y}{dx^{k+m}}=\frac{d^{k+1+m}y}{dx^{k+1+m}} \tag{A.1.21}$$

follows from equalities [15]-(4.1.5), (A.1.20). Therefore, the theorem is true for $n=k+1$. □

## A.2. Table of Integrals

THEOREM A.2.1. *Let the map*

$$f : A \to B$$

*be differentiable map. Then*

$$\int \frac{df(x)}{dx} \circ dx = f(x) + C \tag{A.2.1}$$

PROOF. The theorem follows from the definition 6.1.1. □

THEOREM A.2.2. *Let the map*

$$f : A \to \mathcal{L}(D; A; B)$$

*be integrable map. Then*

$$\frac{d}{dx} \int f(x) \circ dx = f(x) \tag{A.2.2}$$

PROOF. The theorem follows from the definition 6.1.1. □

THEOREM A.2.3.

$$\int (0 \otimes 0) \circ dx = C \tag{A.2.3}$$

PROOF. The theorem follows from the theorem A.1.1. □

THEOREM A.2.4.

$$\int (f(x) + g(x)) \circ dx = \int f(x) \circ dx + \int g(x) \circ dx \tag{A.2.4}$$

PROOF. The equlity (A.2.4) follows from theorems A.1.4, A.2.1. □

THEOREM A.2.5.

$$\int (f_{s \cdot 0} \otimes f_{s \cdot 1}) \circ dx = (f_{s \cdot 0} \otimes f_{s \cdot 1}) \circ x + C \tag{A.2.5}$$

$$\int f_{s \cdot 0} \, dx \, f_{s \cdot 1} = f_{s \cdot 0} \, x \, f_{s \cdot 1} + C \tag{A.2.6}$$

$$f_{s \cdot 0} \in A \quad f_{s \cdot 1} \in A$$

PROOF. The theorem follows from the theorem [15]-B.2.2. □

THEOREM A.2.6.

$$\int (1 \otimes x + x \otimes 1) \circ dx = x^2 + C \tag{A.2.7}$$

$$\int dx \, x + x \, dx = x^2 + C \tag{A.2.8}$$

THEOREM A.2.7.

$$\int (1 \otimes x^2 + x \otimes x + x^2 \otimes 1) \circ dx = x^3 + C \tag{A.2.9}$$

$$\int dx \, x^2 + x \, dx \, x + x^2 \, dx = x^3 + C \tag{A.2.10}$$

PROOF. The theorem follows from the theorem [15]-5.1.3. □

THEOREM A.2.8.

$$\int (e^x \otimes 1 + 1 \otimes e^x) \circ dx = 2e^x + C \tag{A.2.11}$$

$$\int e^x\, dx + dx\, e^x = 2e^x + C \tag{A.2.12}$$

PROOF. The theorem follows from the theorem [15]-B.2.3. □

THEOREM A.2.9.

$$\int (\sinh x \otimes 1 + 1 \otimes \sinh x) \circ dx = 2\cosh x + C \tag{A.2.13}$$

$$\int \sinh x\, dx + dx\ \sinh x = 2\cosh x + C \tag{A.2.14}$$

$$\int (\cosh x \otimes 1 + 1 \otimes \cosh x) \circ dx = 2\sinh x + C \tag{A.2.15}$$

$$\int \cosh x\, dx + dx\ \cosh x = 2\sinh x + C \tag{A.2.16}$$

PROOF. The theorem follows from the theorem [15]-B.2.4. □

THEOREM A.2.10.

$$\int (\sin x \otimes 1 + 1 \otimes \sin x) \circ dx = -2\cos x + C \tag{A.2.17}$$

$$\int \sin x\, dx + dx\ \sin x = -2\cos x + C \tag{A.2.18}$$

$$\int (\cos x \otimes 1 + 1 \otimes \cos x) \circ dx = 2\sin x + C \tag{A.2.19}$$

$$\int \cos x\, dx + dx\ \cos x = 2\sin x + C \tag{A.2.20}$$

PROOF. The theorem follows from the theorem [15]-B.2.5. □

# Index

# Special Symbols and Notations

# Дифференциальное уравнение над банаховой алгеброй

Александр Клейн

Aleks_Kleyn@MailAPS.org
http://AleksKleyn.dyndns-home.com:4080/
http://sites.google.com/site/AleksKleyn/
http://arxiv.org/a/kleyn_a_1
http://AleksKleyn.blogspot.com/

Аннотация. В книге рассмотрены дифференциальные уравнения первого порядка над банаховой $D$-алгеброй: дифференциальное уравнение, разрешённое относительно производной; уравнение в полных дифференциалах; линейное однородное уравнение.

Рассмотрел примеры дифференциальных уравнений в алгебре кватернионов.

С целью изучения однородной системы линейных дифференциальных уравнений рассмотрел векторное пространство над $D$-алгеброй с делением, решение линейных уравнений над $D$-алгеброй с делением и теорию собственных значений в некоммутативной $D$-алгебре с делением.

Рассмотрел пример однородной системы линейных дифференциальных уравнений в алгебре кватернионов, для которой задача с начальными значениями имеет бесконечно много решений.

# Оглавление

Глава 1

# Предисловие

## 1.1. Замысел книги

Когда я начал изучать некоммутативную алгебру, я мог записывать линейное отображение только вместе с его аргументом. Тензорное представление линейного отображения позволило мне упростить запись и расширило поле моего исследования. Сейчас я готов к изучению дифференциальных уравнений над банаховой алгеброй.

Поскольку произведение некоммутативно, то множество уравнений, которое я могу рассматривать более ограничено, чем в коммутативном случае. Тем не менее множество дифференциальных уравнений, рассмотренных в этой книге, является началом крайне интересной теории.

## 1.2. Предисловие к изданию 1

Производная отображения в банаховых пространствах может иметь различное представление: если мы рассматриваем отображение $f$ $D$-алгебры $A$, записанное в координатах относительно базиса $D$-модуля $A$, то производная представлена матрицей Якоби отображения $f$; если мы рассматриваем бескоординатную запись отображения $f$, то производная представлена $A\otimes A$- числом. Поэтому я решил написать эту статью таким образом, чтобы ею можно было пользоваться независимо от представления, и рассмотреть в примерах решение дифференциальных уравнений в различных представлениях.

Глава 6 посвящена дифференциальному уравнению, разрешённое относительно производной. Для того, чтобы картина была полной, я рассматриваю разные формы представления дифференциального уравнения и сопоставляю полученные результаты. Примеры, приведенные в тексте, дают возможность читателю увидеть как работают вычисления в некоммутативной алгебре. В тоже время эти примеры являются источником новых результатов и являются неотъемлемой частью книги.

Наиболее интересным и несколько неожиданным утверждением оказалось, что задача с начальным значением для линейного однородного уравнения имеет бесконечно много решений.

Январь, 2018

## 1.3. Предисловие к изданию 2

Меня крайне заинтересовали статьи [16], [18]. Авторы этих статей, так же как и я, рассматривают систему дифференциальных уравнений для отображений в алгебру кватернионов. Однако, отображения в этих статьях зависят от действительной переменной. Это облегчает решение задачи и в то же время позволяет увидеть что-то новое.

Поэтому, не откладывая дело в долгий ящик, я решил рассмотреть систему дифференциальных уравнений

$$\frac{dx^1}{dt} = a^1_1 x^1 + ... + a^1_n x^n$$
$$......$$
$$\frac{dx^n}{dt} = a^n_1 x^1 + ... + a^n_n x^n$$

в общем виде, как это сделали авторы статей [16], [18], а потом попробовать решить конкретные системы дифференциальных уравнений.

У меня не было проблем с системой дифференциальных уравнений

$$\frac{dx^1}{dt} = x^2$$
$$\frac{dx^2}{dt} = x^1$$

Эта система дифференциальных уравнений имеет единственное решение, которое удовлетворяет начальному условию

$$t = 0 \quad x^1 = 0 \quad x^2 = 1$$

и я немедленно получаю формулу Эйлера. Однако система дифференциальных уравнений

$$\frac{dx^1}{dt} = x^2$$
$$\frac{dx^2}{dt} = -x^1$$

имеет бесконечно много решений, которое удовлетворяют начальному условию

$$t = 0 \quad x^1 = 0 \quad x^2 = 1$$

Это утверждение насторожило меня, и я решил вернуться к дифференциальному уравнению

$$\frac{dy}{dx} = \frac{1}{2}(y \otimes 1 + 1 \otimes y) \tag{1.3.1}$$

при начальном условии

$$y(0) = 1$$

Хотя в прошлом году я пришёл к выводу, что эта задача так же имеет бесконечно много решений, мне было важно понять, почему разложение в ряд Тейлора единственно, а решений бесконечно много.

Когда я начал изучать дифференциальное уравнение (1.3.1), я был под сильным влиянием классической традиции: экспонента - это единственное отображение, которое совпадает со своей производной. Однако в некоммутативной банаховой алгебре это неверно. Когда я это понял, я обнаружил большое семейство отображений, отличных от экспоненты и сохраняющих при дифференцировании своё разложение в ряд Тейлора. Так как поведение этих отображений похоже на поведение экспоненты, я назвал эти отображения квазиэкспонентой.

Как всегда, у меня есть не мало не отвеченных вопросов.

Один из этих вопросов связан с тем, что дифференциальное уравнение для экспоненты не является вполне интегрируемым в не коммутативной алгебре.

Меня заинтересовало, каким образом возможно, что система дифференциальных уравнений

$$\frac{dx^1}{dt} = x^2$$
$$\frac{dx^2}{dt} = -x^1$$

имеет не единственное решение, в то время как разложение в ряд Тейлора определенно однозначно. Этот пример так же не согласуется с теоремой [17]-3.1 о единственности решения задачи с начальным значением.

Теория собственных значений тесно связана с теорией многочленов над некоммутативной алгеброй. Это новая и быстро развивающаяся теория. Согласно теореме [14]-7.2, в алгебре кватернионов существует квадратное уравнение, которое имеет 1 корень. Согласно теореме [14]-7.3, в алгебре кватернионов существует квадратное уравнение, которое не имеет корней. Следовательно, $2 \times 2$ матрица может иметь единственное собственное значение или не иметь собственных значений вообще.

Ещё один вопрос, который меня интересует - это вопрос ковариантности. Мы решаем систему дифференциальных уравнений относительно выбранного базиса в модуле над $D$-алгеброй. Как изменится система дифференциальных уравнений и как изменится её решение, если мы выберем другой базис?

Май, 2019

## 1.4. Предисловие к изданию 3

Математика – дедуктивная наука. Это значит, что из утверждения $A$ должно следовать утверждение $B$. Вроде всё просто. Но когда решаешь задачу, дедукции не достаточно. Индукция – это шаг в неизвестное, туда, где ты никогда раньше не был. А может это вообще неизведанная территория, где никогда не ступала нога человека.

В 2008 я начал изучать математический анализ над некоммутативной алгеброй. Но мне потребовался ещё год, чтобы понять структуру линейного отображения. Практически сразу я понял, что производная высшего порядка - это полилинейное отображение, симметричное относительно порядка аргументов. Это позволило мне решать простейшие дифференциальные уравнения. После этого я попытался понять как выглядит дифференциальное уравнение для экспоненты.

Было ясно, что классическое уравнение $y' = y$ не годится, так как слева и справа разный тип выражений. После непродолжительного поиска я обнаружил, что уравнение

$$\frac{dy}{dx} \circ h = \frac{1}{2}(yh + hy) \tag{1.4.1}$$

имеет подходящую симметрию и разложение в ряд Тейлора его решения совпадает с экспонентой. Я был уверен, что я нашёл верное решение.

В прошлом году я начал писать книгу, посвящённую дифференциальным уравнениям над некоммутативной алгеброй. В частности, я записал условие, когда дифференциальное уравнение имеет решение. Оказалось, что дифференциальное уравнение, которое я написал для экспоненты, неверно. В этот раз поиск дифференциального уравнения занял больше времени. Я нашёл дифференциальное уравнение вида

$$\frac{dy}{dx} \circ 1 = y$$

Но оказалось, что это уравнение имеет много других решений кроме экспоненты. Вначале эти отображения мне казались чудовищами , потом я смирился с их существованием. И вдруг я понял, что эти отображения очень важны. Поэтому я назвал новое отображение квазиэкспонентой.

Квазиэкспонента возникает в процессе дифференцирования экспоненты. В коммутативной алгебре дифференциальное уравнение $y' = y$ говорит, что скорость изменения отображения $y = e^x$ равна $y$. Аналогичную роль играет квазиэкспонента.

Поиск не завершён. Впереди ещё немало вопросов, на которые надо дать ответ.

Июнь, 2019

## 1.5. Соглашения

Соглашение 1.5.1. *Мы будем пользоваться соглашением Эйнштейна о сумме, в котором повторяющийся индекс (один вверху и один внизу) подразумевает сумму по повторяющемуся индексу. В этом случае предполагается известным множество индекса суммирования и знак суммы опускается*

$$c^i v_i = \sum_{i \in I} c^i v_i$$

*Я буду явно указывать множество индексов, если это необходимо.* □

Соглашение 1.5.2. *Пусть $A$ - свободная алгебра с конечным или счётным базисом. При разложении элемента алгебры $A$ относительно базиса $\overline{\overline{e}}$ мы пользуемся одной и той же корневой буквой для обозначения этого элемента и его координат. В выражении $a^2$ не ясно - это компонента разложения элемента $a$ относительно базиса или это операция возведения в степень. Для облегчения чтения текста мы будем индекс элемента алгебры выделять цветом. Например,*

$$a = a^i e_i$$

□

Соглашение 1.5.3. *Если свободная конечномерная алгебра имеет единицу, то мы будем отождествлять вектор базиса $e_0$ с единицей алгебры.* □

Глава 2

# Предварительные утверждения

## 2.1. Модуль линейных отображений в $D$-алгебру

Пусть $A$ - $D$-модуль. Пусть $B$ - свободная конечно мерная ассоциативная $D$-алгебра. Согласно теореме [13]-6.4.8, левый $B\otimes B$-модуль $\mathcal{L}(D;A\to B)$ имеет конечный базис $\overline{\overline{F}}_{AB}$. Однако очевидно, что в общем случае выбор этого базиса произволен. В тоже время, выбор базиса $\overline{\overline{F}}$ левого $B\otimes B$-модуля $\mathcal{L}(D;B\to B)$ обычно связан со структурой $D$-алгебры $B$. Поэтому возникает вопрос, есть ли связь между базисами $\overline{\overline{F}}_{AB}$ и $\overline{\overline{F}}$.

Теорема 2.1.1. *Пусть $A$ - $D$-модуль, $n=\dim A$. Пусть $B$ - ассоциативная $D$-алгебра, $m=\dim B$. Пусть $\overline{\overline{F}}$ - базис левого $B\otimes B$-модуля $\mathcal{L}(D;B\to B)$. Пусть $n\le m$. Пусть*

$$G:A\to B$$

*линейное отображение максимального ранга. Множество*

$$\overline{\overline{F}}\circ G=\{F_k\circ G:F_k\in\overline{\overline{F}}\}\tag{2.1.1}$$

*порождает левый $B\otimes B$-модуль*[2.1] $\mathcal{L}(D;A\to B)$.

Доказательство. Пусть

$$g:A\to B$$

линейное отображение. Пусть $\overline{\overline{e}}_A$ - базис $D$-модуля $A$. Пусть $\overline{\overline{e}}_B$ - базис $D$-модуля $B$. Согласно теореме [13]-4.2.2, линейное отображение $G$ имеет координаты

$$G=\begin{pmatrix}G_1^1&...&G_n^1\\...&...&...\\G_1^m&...&G_n^m\end{pmatrix}$$

относительно базисов $\overline{\overline{e}}_A$, $\overline{\overline{e}}_B$ и линейное отображение $g$ имеет координаты

$$g=\begin{pmatrix}g_1^1&...&g_n^1\\...&...&...\\g_1^m&...&g_n^m\end{pmatrix}$$

относительно базисов $\overline{\overline{e}}_A$, $\overline{\overline{e}}_B$. Строка $G_i$ матрицы $G$, также как строка $g_i$ матрицы $g$, является координатами линейной формы $A\to D$. Поскольку матрица $G$ имеет максимальный ранг, то строки матрицы $G$ порождают $D$-модуль

[2.1] Я не утверждаю, что это множество является базисом, так как отображения $F_i\circ G$, $i\in I$, могут быть линейно зависимыми.

$\mathcal{L}(D;A\to D)$ и строки матрицы $g$ являются линейной комбинацией строк матрицы $G$

(2.1.2) $$g^i_k=C^i_jG^j_k$$

Поскольку мы можем рассматривать матрицу $C$ как координаты линейного отображения

$$C:B\to B$$

то равенство

(2.1.3) $$g^i_k=(c_{1.k.l}F_{k\cdot}{}^i_jc_{2.k.l})G^j_k$$

является следствием равенства (2.1.2) и равенства

$$C=c_{1.k.l}F_kc_{2.k.l}$$

Поскольку $G^j_k\in D$, то равенство

(2.1.4) $$g^i_k=c_{1.k.l}(F_{k\cdot}{}^i_jG^j_k)c_{2.k.l}$$

является следствием равенства (2.1.3). Следовательно, отображение $g$ принадлежит линейной оболочке множества отображений (2.1.1). □

Теорема 2.1.2. *Пусть $n=\dim A$, $m=\dim B$. Пусть $\overline{\overline{F}}$ - базис левого $B\otimes B$-модуля $\mathcal{L}(D;B\to B)$. Пусть $n>m$. Пусть*

$$G:A\to B$$

*линейное отображение максимального ранга. Множество*

$$\overline{\overline{F}}\circ G=\{F_k\circ G:F_k\in\overline{\overline{F}}\}$$

*порождает множество отображений*

(2.1.5) $$\{g\in\mathcal{L}(D;A\to B):\ker G\subseteq\ker g\}$$

Доказательство. Доказательство теоремы похоже на доказательство теоремы 2.1.1. Однако, так как число строк матрицы $G$ меньше размерности $D$-модуля $A$, то строки матрицы $G$ не порождают $D$-модуль $\mathcal{L}(D;A\to D)$ и отображение $G$ имеет нетривиальное ядро. В частности, строки матрицы $g$ линейно зависят от строк матрицы $G$ тогда и только тогда, когда $\ker G\subseteq\ker g$. □

Из теоремы 2.1.2, следует, что выбор отображения $G$ зависит от отображения $g$. Нетрудно убедиться, что теорема 2.1.1 является частным случаем теоремы 2.1.2, так как в теореме 2.1.2 $\ker G=\emptyset$.

Теорема 2.1.3. *Пусть $\overline{\overline{e}}$ - базис $D$-модуля $A$. Тогда множество отображений*

$$P_i:d\in D\to de_i\in A$$

*является базисом $D$-модуля $\mathcal{L}(D;D\to A)$.*

Доказательство. Теорема является следствием равенства

$$f\circ t=(f^it)e_i=P_i\circ(f^it)$$

□

## 2.2. Прямая сумма $D$-модулей

Определение 2.2.1. *Пусть $\mathcal{A}$ - категория. Пусть $\{B_i, i \in I\}$ - множество объектов из $\mathcal{A}$. Объект*

$$P = \coprod_{i \in I} B_i$$

*и множество морфизмов*

$$\{f_i : B_i \to P, i \in I\}$$

*называется* **копроизведением множества объектов $\{B_i, i \in I\}$ в категории $\mathcal{A}$**[2.2]*, если для любого объекта $R$ и множество морфизмов*

$$\{g_i : B_i \to R, i \in I\}$$

*существует единственный морфизм*

$$h : P \to R$$

*такой, что диаграмма*

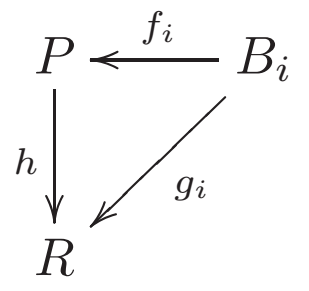

$$h \circ f_i = g_i$$

*коммутативна для всех $i \in I$.*

*Если $|I| = n$, то для копроизведения множества объектов $\{B_i, i \in I\}$ в $\mathcal{A}$ мы так же будем пользоваться записью*

$$P = \coprod_{i=1}^{n} B_i = B_1 \coprod ... \coprod B_n$$

□

Определение 2.2.2. *Копроизведение в категории абелевых групп $Ab$ называется* **прямой суммой**.[2.3] *Мы будем пользоваться записью $A \oplus B$ для прямой суммы абелевых групп $A$ и $B$.* □

Теорема 2.2.3. *Пусть $\{A_i, i \in I\}$ - семейство абелевых групп. Пусть*

$$A \subseteq \prod_{i \in I} A_i$$

*такое множество, что $(x_i, i \in I) \in A$, если $x_i \neq 0$ для конечного числа индексов $i$. Тогда*[2.4]

$$A = \bigoplus_{i \in I} A_i \tag{2.2.1}$$

Доказательство. Согласно построению $A$ является подгруппой абелевой группы $\prod A_i$. Отображение

$$\lambda_j : A_j \to A$$

определённое равенством

$$\lambda_j(x) = (\delta^i_j x, i \in I) \tag{2.2.2}$$

---

[2.2] Определение дано согласно [1], страница 46.
[2.3] Определение дано согласно [1], страница 55.
[2.4] Смотри также предложение [1]-10, страница 55.

является инъективным гомоморфизмом.

Пусть

$$\{f_i : A_i \to B, i \in I\}$$

семейство гомоморфизмов в абелевую группу $B$. Мы определим отображение

$$f : A \to B$$

пользуясь равенством

$$f\left(\bigoplus_{i\in I} x_i\right) = \sum_{i\in I} f_i(x_i) \tag{2.2.3}$$

Сумма в правой части равенства (2.2.3) конечна, так как все слагаемые, кроме конечного числа, равна 0. Из равенства

$$f_i(x_i + y_i) = f_i(x_i) + f_i(y_i)$$

и равенства (2.2.3) следует, что

$$\begin{aligned} f(\bigoplus_{i\in I}(x_i + y_i)) &= \sum_{i\in I} f_i(x_i + y_i) = \sum_{i\in I}(f_i(x_i) + f_i(y_i)) \\ &= \sum_{i\in I} f_i(x_i) + \sum_{i\in I} f_i(y_i) \\ &= f(\bigoplus_{i\in I} x_i) + f(\bigoplus_{i\in I} y_i) \end{aligned}$$

Следовательно, отображение $f$ является гомоморфизмом абелевой группы. Равенство

$$f \circ \lambda_j(x) = \sum_{i\in I} f_i(\delta^i_j x) = f_j(x)$$

является следствием равенств (2.2.2), (2.2.3). Так как отображение $\lambda_i$ инъективно, то отображение $f$ определено однозначно. Следовательно, теорема является следствием определений 2.2.1, 2.2.2. □

Теорема 2.2.4. *Прямая сумма абелевых групп* $A_1$, ..., $A_n$ *совпадает с их прямым произведением*

$$A_1 \oplus ... \oplus A_n = A_1 \times ... \times A_n$$

Доказательство. Теорема является следствием теоремы 2.2.3. □

Пусть

$$A = A_1 \oplus ... \oplus A_n$$

прямая сумма абелевых групп $A_1$, ..., $A_n$. Согласно доказательству теоремы 2.2.3, произвольное $A$-число $a$ имеет вид $(a_1, ..., a_n)$ где $a_i \in A_i$. Мы также будем пользоваться записью

$$a = a_1 \oplus ... \oplus a_n$$

Определение 2.2.5. *Копроизведение в категории $D$-модулей называется* **прямой суммой**.[2.5] *Мы будем пользоваться записью* $A \oplus B$ *для прямой $D$-модулей $A$ и $B$.* □

[2.5] Смотри так же определение прямой суммы модулей в [1], страница 98. Согласно теореме 1 на той же странице, прямая сумма модулей существует.

Теорема 2.2.6. *Пусть* $\{A_i, i \in I\}$ *- семейство $D$-модулей. Тогда представление*

$$D \xrightarrow{*} \bigoplus_{i\in I} A_i \qquad d\left(\bigoplus_{i\in I} a_i\right) = \bigoplus_{i\in I} da_i$$

*кольца $D$ в прямой сумме абелевых групп*

$$A = \bigoplus_{i\in I} A_i$$

*является прямой суммой $D$-модулей*

$$A = \bigoplus_{i\in I} A_i$$

Доказательство.

Пусть

$$\{f_i : A_i \to B, i \in I\}$$

семейство линейных отображений в $D$-модуль $B$. Мы определим отображение

$$f : A \to B$$

пользуясь равенством

$$f\left(\bigoplus_{i\in I} x_i\right) = \sum_{i\in I} f_i(x_i) \tag{2.2.4}$$

Сумма в правой части равенства (2.2.4) конечна, так как все слагаемые, кроме конечного числа, равна 0. Из равенства

$$f_i(x_i + y_i) = f_i(x_i) + f_i(y_i)$$

и равенства (2.2.4) следует, что

$$\begin{aligned} f(\bigoplus_{i\in I}(x_i + y_i)) &= \sum_{i\in I} f_i(x_i + y_i) = \sum_{i\in I}(f_i(x_i) + f_i(y_i)) \\ &= \sum_{i\in I} f_i(x_i) + \sum_{i\in I} f_i(y_i) \\ &= f(\bigoplus_{i\in I} x_i) + f(\bigoplus_{i\in I} y_i) \end{aligned}$$

Из равенства

$$f_i(dx_i) = df_i(x_i)$$

и равенства (2.2.4) следует, что

$$\begin{aligned} f((dx_i, i \in I)) &= \sum_{i\in I} f_i(dx_i) = \sum_{i\in I} df_i(x_i) = d\sum_{i\in I} f_i(x_i) \\ &= df((x_i, i \in I)) \end{aligned}$$

Следовательно, отображение $f$ является линейным отображением. Равенство

$$f \circ \lambda_j(x) = \sum_{i\in I} f_i(\delta^i_j x) = f_j(x)$$

является следствием равенств (2.2.2), (2.2.4). Так как отображение $\lambda_i$ инъективно, то отображение $f$ определено однозначно. Следовательно, теорема является следствием определений 2.2.1, 2.2.2. □

Теорема 2.2.7. *Прямая сумма $D$-модулей $A_1$, ..., $A_n$ совпадает с их прямым произведением*

$$A_1 \oplus ... \oplus A_n = A_1 \times ... \times A_n$$

Доказательство. Теорема является следствием теоремы 2.2.6. □

Теорема 2.2.8. *Пусть $A^1$, ..., $A^n$ - $D$-модули и*

$$A = A^1 \oplus ... \oplus A^n$$

*Представим $A$-число*

$$a = a^1 \oplus ... \oplus a^n$$

*как вектор столбец*

$$a = \begin{pmatrix} a^1 \\ ... \\ a^n \end{pmatrix}$$

*Представим линейное отображение*

$$f : A \to B$$

*как вектор строку*

$$f = \begin{pmatrix} f_1 & ... & f_n \end{pmatrix}$$

$$f_i : A^i \to B$$

*Тогда значение отображения $f$ в $A$-числе $a$ можно представить как произведение матриц*

$$f \circ a = \begin{pmatrix} f_1 & ... & f_n \end{pmatrix} {}_\circ{}^\circ \begin{pmatrix} a^1 \\ ... \\ a^n \end{pmatrix} = f_i \circ a^i \tag{2.2.5}$$

Доказательство. Теорема является следствием определения (2.2.4). □

Теорема 2.2.9. *Пусть $B^1$, ..., $B^m$ - $D$-модули и*

$$B = B^1 \oplus ... \oplus B^m$$

*Представим $B$-число*

$$b = b^1 \oplus ... \oplus b^m$$

*как вектор столбец*

$$b = \begin{pmatrix} b^1 \\ ... \\ b^m \end{pmatrix}$$

*Тогда линейное отображение*

$$f : A \to B$$

*имеет представление как вектор столбец отображений*

$$f = \begin{pmatrix} f^1 \\ ... \\ f^m \end{pmatrix}$$

*таким образом, что если $b = f \circ a$, то*

$$\begin{pmatrix} b^1 \\ ... \\ b^m \end{pmatrix} = \begin{pmatrix} f^1 \\ ... \\ f^m \end{pmatrix} \circ a = \begin{pmatrix} f^1 \circ a \\ ... \\ f^m \circ a \end{pmatrix}$$

Доказательство. Теорема является следствием теоремы [13]-2.1.5. □

Теорема 2.2.10. *Пусть* $A^1$, ..., $A^n$, $B^1$, ..., $B^m$ *- $D$-модули и*

$$A = A^1 \oplus ... \oplus A^n$$
$$B = B^1 \oplus ... \oplus B^m$$

*Представим $A$-число*

$$a = a^1 \oplus ... \oplus a^n$$

*как вектор столбец*

$$a = \begin{pmatrix} a^1 \\ ... \\ a^n \end{pmatrix}$$

*Представим $B$-число*

$$b = b^1 \oplus ... \oplus b^m$$

*как вектор столбец*

$$b = \begin{pmatrix} b^1 \\ ... \\ b^m \end{pmatrix}$$

*Тогда линейное отображение $f$ имеет представление как матрица отображений*

$$f = \begin{pmatrix} f_1^1 & ... & f_n^1 \\ ... & ... & ... \\ f_1^m & ... & f_n^m \end{pmatrix}$$

*таким образом, что если $b = f \circ a$, то*

$$\begin{pmatrix} b^1 \\ ... \\ b^m \end{pmatrix} = \begin{pmatrix} f_1^1 & ... & f_n^1 \\ ... & ... & ... \\ f_1^m & ... & f_n^m \end{pmatrix} {}_\circ{}^\circ \begin{pmatrix} a^1 \\ ... \\ a^n \end{pmatrix} = \begin{pmatrix} f_i^1 \circ a^i \\ ... \\ f_i^m \circ a^i \end{pmatrix} \tag{2.2.6}$$

*Отображение*

$$f_j^i : A^j \to B^i$$

*является линейным отображением и называется* **частным линейным отображением**

Доказательство. Согласно теореме 2.2.9, существует множество линейных отображений

$$f^i : A \to B^i$$

таких, что

$$\begin{pmatrix} b^1 \\ ... \\ b^m \end{pmatrix} = \begin{pmatrix} f^1 \\ ... \\ f^m \end{pmatrix} \circ a = \begin{pmatrix} f^1 \circ a \\ ... \\ f^m \circ a \end{pmatrix} \tag{2.2.7}$$

Согласно теореме 2.2.8, для каждого $i$, существует множество линейных отображений

$$f^i_j : A^j \to B^i$$

таких, что

$$f^i \circ a = \begin{pmatrix} f^i_1 & ... & f^i_n \end{pmatrix} {}_\circ{}^\circ \begin{pmatrix} a^1 \\ ... \\ a^n \end{pmatrix} = f^i_j \circ a^j \tag{2.2.8}$$

Если мы отождествим матрицы

$$\begin{pmatrix} \begin{pmatrix} f^1_1 & ... & f^1_n \end{pmatrix} \\ ... \\ \begin{pmatrix} f^m_1 & ... & f^m_n \end{pmatrix} \end{pmatrix} = \begin{pmatrix} f^1_1 & ... & f^1_n \\ ... & ... & ... \\ f^m_1 & ... & f^m_n \end{pmatrix}$$

то равенство (2.2.6) является следствием равенств (2.2.7), (2.2.8). □

## 2.3. Прямая сумма банаховых $D$-модулей

Теорема 2.3.1. *Пусть* $A^1$, ..., $A^n$ *- банаховые $D$-модули и*

$$A = A^1 \oplus ... \oplus A^n$$

*Тогда мы можем определить норму в $D$-модуле $A$ такую, что $D$-модуль $A$ становится банаховым $D$-модулем.*

Доказательство. Пусть $\|a^i\|_i$ - норма в $D$-модуле $A^i$.

2.3.1.1: Мы определим норму в $D$-модуле $A$ равенством

$$\|b\| = \max(\|b^i\|_i, i = 1, ..., n)$$

где

$$b = b^1 \oplus ... \oplus b^n$$

2.3.1.2: Пусть $\{a_p\}$, $p = 1$, ..., - фундаментальная последовательность, где

$$a_p = a^1_p \oplus ... \oplus a^n_p$$

2.3.1.3: Следовательно, для любого $\epsilon \in R$, $\epsilon > 0$, существует $N$ такое, что для любых $p, q > N$

$$\|a_p - a_q\| < \epsilon$$

2.3.1.4: Согласно утверждениям 2.3.1.1, 2.3.1.2, 2.3.1.3,

$$\|a_p^i - a_q^i\|_i < \epsilon$$

для любых $p, q > N$ и $i = 1, ..., n$.

2.3.1.5: Следовательно, последовательность $\{a_p^i\}$, $i = 1, ..., n$, $p = 1, ...,$ фундаментальна в $D$-модуле $A^i$ и существует предел

$$a^i = \lim_{p\to\infty} a_p^i$$

2.3.1.6: Пусть

$$a = a^1 \oplus ... \oplus a^n$$

2.3.1.7: Согласно утверждению 2.3.1.5, для любого $\epsilon \in R$, $\epsilon > 0$, существует $N$ такое, что для любого $p > N_i$

$$\|a^i - a_p^i\|_i < \epsilon$$

2.3.1.8: Пусть

$$N = \max(N_1, ..., N_n)$$

2.3.1.9: Согласно утверждениям 2.3.1.6, 2.3.1.7, 2.3.1.8, для любого $\epsilon \in R$, $\epsilon > 0$, существует $N$ такое, что для любого $p > N$

$$\|a - a_p\|_i < \epsilon$$

2.3.1.10: Следовательно,

$$a = \lim_{p\to\infty} a_p$$

Теорема является следствием утверждений 2.3.1.1, 2.3.1.2, 2.3.1.10. □

Опираясь на теорему 2.3.1, мы можем рассмотреть производную отображения

$$f : A^1 \oplus ... \oplus A^n \to B^1 \oplus ... \oplus B^m$$

Теорема 2.3.2. *Пусть* $A^1, ..., A^n$, $B^1, ..., B^m$ *- банаховые $D$-модули и*

$$A = A^1 \oplus ... \oplus A^n$$

$$B = B^1 \oplus ... \oplus B^m$$

*Представим дифференциал*

$$dx = dx^1 \oplus ... \oplus dx^n$$

*как вектор столбец*

$$dx = \begin{pmatrix} dx^1 \\ ... \\ dx^n \end{pmatrix}$$

*Представим дифференциал*

$$dy = dy^1 \oplus ... \oplus dy^m$$

*как вектор столбец*

$$dy = \begin{pmatrix} dy^1 \\ ... \\ dy^m \end{pmatrix}$$

*Тогда производная отображения*

$$f : A \to B$$
$$f = f^1 \oplus ... \oplus f^m$$

*имеет представление*

$$\frac{df}{dx} = \begin{pmatrix} \frac{\partial f^1}{\partial x^1} & ... & \frac{\partial f^1}{\partial x^n} \\ ... & ... & ... \\ \frac{\partial f^m}{\partial x^1} & ... & \frac{\partial f^m}{\partial x^n} \end{pmatrix}$$

*таким образом, что*

$$\begin{pmatrix} dy^1 \\ ... \\ dy^m \end{pmatrix} = \begin{pmatrix} \frac{\partial f^1}{\partial x^1} & ... & \frac{\partial f^1}{\partial x^n} \\ ... & ... & ... \\ \frac{\partial f^m}{\partial x^1} & ... & \frac{\partial f^m}{\partial x^n} \end{pmatrix} {}_\circ{}^\circ \begin{pmatrix} dx^1 \\ ... \\ dx^n \end{pmatrix} = \begin{pmatrix} \frac{\partial f^1}{\partial x^i} \circ dx^i \\ ... \\ \frac{\partial f^m}{\partial x^i} \circ dx^i \end{pmatrix} \tag{2.3.1}$$

Утверждение 2.3.3. *Линейное отображение* $\frac{\partial f^i}{\partial x^j}$ *называется* **частной производной** *и является производной отображения* $f^i$ *по переменной* $x^j$ *при условии, что остальные координаты A-числа* $x$ *постоянны.* ⊙

Доказательство. Равенство (2.3.1) является следствием равенства (2.2.6). Мы можем записать отображение

$$f^i : A \to B^i$$

в виде

$$f^i(x) = f^i(x^1, ..., x^n)$$

Равенство

$$\frac{df^i(x)}{dx} \circ dx = \frac{\partial f^i(x^1, ..., x^n)}{\partial x^j} \circ dx^j \tag{2.3.2}$$

является следствием равенства (2.3.1). Согласно теореме [15]-3.3.4,

$$\begin{aligned} \frac{df^i(x)}{dx} \circ dx = & \lim_{t \to 0,\ t \in R} (t^{-1}(f^(x + tdx) - f(x))) \\ = & \lim_{t \to 0,\ t \in R} (t^{-1}(f^i(x^1 + tdx^1, x^2 + tdx^2, ..., x^n + tdx^n) \\ & - f^i(x^1, x^2 + tdx^2, ..., x^n + tdx^n) \\ & + f^i(x^1, x^2 + tdx^2, ..., x^n + tdx^n) - ... \\ & - f^i(x^1, x^2, ..., x^n))) \\ = & \lim_{t \to 0,\ t \in R} (t^{-1}(f^i(x^1 + tdx^1, x^2 + tdx^2, ..., x^n + tdx^n) \\ & - f^i(x^1, x^2 + tdx^2, ..., x^n + tdx^n))) + ... \\ & + \lim_{t \to 0,\ t \in R} (t^{-1}(f^i(x^1, ..., x^n + tdx^n) - f^i(x^1, ..., x^n))) \\ = & f^i_1 \circ dx^1 + ... + f^i_n \circ dx^n \end{aligned} \tag{2.3.3}$$

где $f^i_j$ - производная отображения $f^i$ по переменной $x^j$ при условии, что остальные координаты $A$-числа $x$ постоянны. Равенство

$$f^i_j = \frac{\partial f^i(x^1, ..., x^n)}{\partial x^j} \tag{2.3.4}$$

является следствием равенств (2.3.2), (2.3.3). Утверждение 2.3.3 является следствием равенства (2.3.4). □

Пример 2.3.4. *Рассмотрим отображение*

$$\begin{aligned} y^1 &= f^1(x^1, x^2, x^3) = (x^1)^2 + x^2x^3 \\ y^2 &= f^2(x^1, x^2, x^3) = x^1x^2 + (x^3)^2 \end{aligned} \tag{2.3.5}$$

*Следовательно*

$$\begin{array}{lll} \dfrac{\partial y^1}{\partial x^1} = x^1 \otimes 1 + 1 \otimes x^1 & \dfrac{\partial y^1}{\partial x^2} = 1 \otimes x^3 & \dfrac{\partial y^1}{\partial x^3} = x^2 \otimes 1 \\ \dfrac{\partial y^2}{\partial x^1} = 1 \otimes x^2 & \dfrac{\partial y^2}{\partial x^2} = x^1 \otimes 1 & \dfrac{\partial y^2}{\partial x^3} = x^3 \otimes 1 + 1 \otimes x^3 \end{array}$$

*и производная отображения* (2.3.5) *имеет вид*

$$\frac{df}{dx} = \begin{pmatrix} x^1 \otimes 1 + 1 \otimes x^1 & 1 \otimes x^3 & x^2 \otimes 1 \\ 1 \otimes x^2 & x^1 \otimes 1 & x^3 \otimes 1 + 1 \otimes x^3 \end{pmatrix} \tag{2.3.6}$$

*Равенство*

$$\begin{aligned} dy^1 &= (x^1 \otimes 1 + 1 \otimes x^1) \circ dx^1 + (1 \otimes x^3) \circ dx^2 + (x^2 \otimes 1) \circ dx^3 \\ &= x^1dx^1 + dx^1x^1 + dx^2x^3 + x^2dx^3 \\ dy^2 &= (1 \otimes x^2) \circ dx^1 + (x^1 \otimes 1) \circ dx^2 + (x^3 \otimes 1 + 1 \otimes x^3) \circ dx^3 \\ &= dx^1x^2 + x^1dx^2 + x^3dx^3 + dx^3x^3 \end{aligned} \tag{2.3.7}$$

*является следствием равенств* (2.3.6). *Мы можем также получить выражение* (2.3.7) *непосредственным вычислением*

$$\begin{aligned} dy^1 &= f^1(x + dx) - f^1(x) \\ &= (x^1 + dx^1)^2 + (x^2 + dx^2)(x^3 + dx^3) - (x^1)^2 - x^2x^3 \\ &= (x^1)^2 + x^1dx^1 + dx^1x^1 + x^2x^3 + dx^2x^3 + x^2dx^3 - (x^1)^2 - x^2x^3 \\ &= x^1dx^1 + dx^1x^1 + dx^2x^3 + x^2dx^3 \\ dy^2 &= f^2(x + dx) - f^2(x) \\ &= (x^1 + dx^1)(x^2 + dx^2) + (x^3 + dx^3)^2 - x^1x^2 - (x^3)^2 \\ &= x^1x^2 + dx^1x^2 + x^1dx^2 + (x^3)^2 + x^3dx^3 + dx^3x^3 - x^1x^2 - (x^3)^2 \\ &= dx^1x^2 + x^1dx^2 + x^3dx^3 + dx^3x^3 \end{aligned} \tag{2.3.8}$$

□

Глава 3

# Бикольцо матриц

## 3.1. Концепция обобщённого индекса

Изучая тензорное исчисление, мы начинаем с изучения одновалентных ковариантного и контравариантного тензоров. Несмотря на различие свойств, оба эти объекта являются элементами соответствующих векторных пространств. Если мы введём обобщённый индекс по правилу $a^{i} = a^i$, $b^{i} = b^{\cdot -}_{\ \, i}$, то мы видим, что эти тензоры ведут себя одинаково. Например, преобразование ковариантного тензора принимает форму

$$b'^{i} = b'^{\cdot -}_{\ \, i} = f^{\cdot - \ \, j}_{\ \, i \, \cdot -} b^{\cdot -}_{\ \, j} = f^{i}_{j} b^{j}$$

Это сходство идёт сколь угодно далеко, так как тензоры также порождают векторное пространство.

Эти наблюдения сходства свойств ковариантного и контравариантного тензоров приводят нас к концепции обобщённого индекса. Я пользуюсь символом $\cdot$ перед обобщённым индексом, когда мне необходимо описать его структуру. Я помещаю символ $'-'$ на месте индекса, позиция которого изменилась. Например, если исходное выражение было $a_{ij}$, я пользуюсь записью $a_{i \cdot -}^{\ \ \, j}$ вместо записи $a^{j}_{i}$.

Хотя структура обобщённого индекса произвольна, мы будем предполагать, что существует взаимно однозначное отображение отрезка натуральных чисел $\mathbf{1}$, ..., $\mathbf{n}$ на множество значений индекса. Пусть $\boldsymbol{I}$ - множество значений индекса $\boldsymbol{i}$. Мы будем обозначать мощность этого множества символом $|\boldsymbol{I}|$ и будем полагать $|\boldsymbol{I}| = \boldsymbol{n}$. Если нам надо перечислить элементы $a_{i}$, мы будем пользоваться обозначением $a_{\mathbf{1}}$, ..., $a_{\boldsymbol{n}}$.

Представление координат вектора в форме матрицы позволяет сделать запись более компактной. Вопрос о представлении вектора как строка или столбец матрицы является вопросом соглашения. Мы можем распространить концепцию обобщённого индекса на элементы матрицы. Матрица - это двумерная таблица, строки и столбцы которой занумерованы обобщёнными индексами. Так как мы пользуемся обобщёнными индексами, мы не можем сказать, нумерует ли индекс $a$ строки матрицы или столбцы, до тех пор, пока мы не знаем структуры индекса. Однако как мы увидим ниже для нас несущественна форма представления матрицы. Для того, чтобы обозначения, предлагаемые ниже, были согласованы с традиционными, мы будем предполагать, что матрица представлена в виде

$$a = \begin{pmatrix} a^{\mathbf{1}}_{\mathbf{1}} & ... & a^{\mathbf{1}}_{\boldsymbol{n}} \\ ... & ... & ... \\ a^{\boldsymbol{m}}_{\mathbf{1}} & ... & a^{\boldsymbol{m}}_{\boldsymbol{n}} \end{pmatrix}$$

Верхний индекс перечисляет строки, нижний индекс перечисляет столбцы.

Определение 3.1.1. *Я использую следующие имена и обозначения различных* **минорных матриц** *матрицы* $a$

- **столбец матрицы** *с индексом* $i$

$$a_i = \begin{pmatrix} a_i^1 \\ \dots \\ a_i^m \end{pmatrix}$$

*Верхний индекс перечисляет элементы столбца, нижний индекс перечисляет столбцы.*

$a_T$ : *минорная матрица, полученная из матрицы* $a$ *выбором столбцов с индексом из множества* $T$

$a_{[i]}$ : *минорная матрица, полученный из матрицы* $a$ *удалением столбца* $a_i$

$a_{[T]}$ : *минорная матрица, полученный из матрицы* $a$ *удалением столбцов с индексом из множества* $T$

- **строка матрицы** *с индексом* $j$

$$a^j = \begin{pmatrix} a_1^j & \dots & a_n^j \end{pmatrix}$$

*Нижний индекс перечисляет элементы строки, верхний индекс перечисляет строки.*

$a^S$ : *минорная матрица, полученный из матрицы* $a$ *выбором строк с индексом из множества* $S$

$a^{[j]}$ : *минорная матрица, полученный из матрицы* $a$ *удалением строки* $a^j$

$a^{[S]}$ : *минорная матрица, полученный из матрицы* $a$ *удалением строк с индексом из множества* $S$

□

Замечание 3.1.2. *Мы будем комбинировать запись индексов. Так* $a_i^j$ *является* $1 \times 1$ *минорной матрицей. Одновременно, это обозначение элемента матрицы. Это позволяет отождествить* $1 \times 1$ *матрицу и её элемент. Индекс* $a$ *является номером столбца матрицы и индекс* $b$ *является номером строки матрицы.* □

Пусть $I$, $|I| = n$, - множество индексов. **Символ Кронекера** определён равенством

$$\delta_j^i = \begin{cases} 1 & i = j \\ 0 & i \neq j \end{cases} \qquad i, j \in I$$

## 3.2. Бикольцо

Мы будем рассматривать матрицы, элементы которых принадлежат ассоциативной $D$-алгебре с делением $A$.

Произведение матриц связано с произведением гомоморфизмов векторных пространств над полем. Согласно традиции произведение матриц $a$ и $b$ определено как произведение строк матрицы $a$ и столбцов матрицы $b$.

Пример 3.2.1. *Пусть* $\overline{\overline{e}}$ *- базис правого векторного пространства* $V$ *над* $D$*-алгеброй* $A$ *(смотри определение 5.1.4 и теорему 5.2.7). Мы представим базис* $\overline{\overline{e}}$ *как строку матрицы*

$$e = \begin{pmatrix} e_1 & \dots & e_n \end{pmatrix}$$

*Мы можем представить координаты вектора* $v$ *как вектор столбец*

$$v = \begin{pmatrix} v^1 \\ \dots \\ v^n \end{pmatrix}$$

*Поэтому мы можем представить вектор* $v$ *как традиционное произведение матриц*

$$v = \begin{pmatrix} e_1 & \dots & e_n \end{pmatrix} \begin{pmatrix} v^1 \\ \dots \\ v^n \end{pmatrix} = e_i v^i$$

*Линейный гомоморфизм правого векторного пространства* $V$ *может быть представлен с помощью матрицы*

$$v'^i = f^i_j v^j \tag{3.2.1}$$

*Равенство* (3.2.1) *выражает традиционное произведение матриц* $f$ *и* $v$. □

Пример 3.2.2. *Пусть* $\overline{\overline{e}}$ *- базис левого векторного пространства* $V$ *над* $D$*-алгеброй* $A$ *(смотри определение 5.1.3 и теорему 5.2.7). Мы представим базис* $\overline{\overline{e}}$ *как строку матрицы*

$$e = \begin{pmatrix} e_1 & \dots & e_n \end{pmatrix}$$

*Мы можем представить координаты вектора* $v$ *как вектор столбец*

$$v = \begin{pmatrix} v^1 \\ \dots \\ v^n \end{pmatrix}$$

*Однако мы не можем представить вектор*

$$v = v^i e_i$$

*как традиционное произведение матриц*

$$v = \begin{pmatrix} v^1 \\ \dots \\ v^n \end{pmatrix} \qquad e = \begin{pmatrix} e_1 & \dots & e_n \end{pmatrix}$$

*так как это произведение не определено. Линейный гомоморфизм левого векторного пространства* $V$ *может быть представлен с помощью матрицы*

$$v'^i = v^j f^i_j \tag{3.2.2}$$

*Равенство* (3.2.2) *не может быть выражено как традиционное произведение матриц* $v$ *и* $f$. □

Из примеров 3.2.1, 3.2.2 следует, что мы не можем ограничиться традиционным произведением матриц и нам нужно определить два вида произведения матриц. Чтобы различать эти произведения, мы вводим новые обозначения. Для совместимости обозначений с существующими мы будем иметь в виду ${}_*{}^*$-произведение, когда нет явных обозначений.

Определение 3.2.3. *Пусть число столбцов матрицы $a$ равно числу строк матрицы $b$.* ${}_*{}^*$**-произведение** *матриц $a$ и $b$ имеет вид*

$$\begin{cases} a{}_*{}^*b = \left(a_k^i b_j^k\right) \\ (a{}_*{}^*b)_j^i = a_k^i b_j^k \end{cases} \tag{3.2.3}$$

*и может быть выражено как произведение строк матрицы $a$ и столбцов матрицы $b$.*[3.1] □

Определение 3.2.4. *Пусть число строк матрицы $a$ равно числу столбцов матрицы $b$.* ${}^*{}_*$**-произведение** *матриц $a$ и $b$ имеет вид*

$$\begin{cases} a{}^*{}_*b = \left(a_i^k b_k^j\right) \\ (a{}^*{}_*b)_j^i = a_i^k b_k^j \end{cases} \tag{3.2.4}$$

*и может быть выражено как произведение столбцов матрицы $a$ и строк матрицы $b$.*[3.2] □

Мы так же определим следующие операции на множестве матриц.

Определение 3.2.5. *Транспонирование $a^T$ матрицы $a$ меняет местами строки и столбцы*

$$(a^T)_j^i = a_i^j \tag{3.2.5}$$

□

Определение 3.2.6. *Сумма матриц $a$ и $b$ определена равенством*

$$(a+b)_j^i = a_j^i + b_j^i$$

□

Пусть

$$a = \begin{pmatrix} a_1^1 & \dots & a_n^1 \\ \dots & \dots & \dots \\ a_1^m & \dots & a_n^m \end{pmatrix}$$

[3.1] Мы будем пользоваться символом ${}_*{}^*$- в последующей терминологии и обозначениях. Мы будем читать символ ${}_*{}^*$ как $rc$-произведение или произведение строки на столбец. Символ произведения строки на столбец сформирован из двух символов операции произведения, которые записываются на месте индекса суммирования. Например, если произведение $A$-чисел имеет вид $a\circ b$, то ${}_*{}^*$-произведение матриц $a$ и $b$ имеет вид $a{}_\circ{}^\circ b$.

[3.2] Мы будем пользоваться символом ${}^*{}_*$- в последующей терминологии и обозначениях. Мы будем читать символ ${}^*{}_*$ как $cr$-произведение или произведение столбца на строку. Символ произведение столбца на строку сформирован из двух символов операции произведения, которые записываются на месте индекса суммирования. Например, если произведение $A$-чисел имеет вид $a\circ b$, то ${}^*{}_*$-произведение матриц $a$ и $b$ имеет вид $a{}^\circ{}_\circ b$.

матрица $A$-чисел. Мы называем матрицу[3.3]

$$\mathcal{H}a = \begin{pmatrix} (a_1^1)^{-1} & ... & (a_1^m)^{-1} \\ ... & ... & ... \\ (a_n^1)^{-1} & ... & (a_n^m)^{-1} \end{pmatrix} \tag{3.2.6}$$

$$(\mathcal{H}a)_j^i = ((a^T)_j^i)^{-1} = (a_i^j)^{-1} \tag{3.2.7}$$

**обращением Адамара матрицы** $a$ ([5]-page 4).

Множество $n \times n$ матриц замкнуто относительно ${}_*{}^*$-произведения и ${}^*{}_*$-произведения, а также относительно суммы.

Теорема 3.2.7.

$$(a_*{}^*b)^T = a^T{}^*{}_*b^T \tag{3.2.8}$$

Доказательство. Цепочка равенств

$$((a_*{}^*b)^T)_i^j = (a_*{}^*b)_j^i = a_k^i b_j^k = (a^T)_i^k (b^T)_k^j = ((a^T)^*{}_*(b^T))_i^j \tag{3.2.9}$$

следует из (3.2.5), (3.2.3) и (3.2.4). Равенство (3.2.8) следует из (3.2.9). □

Матрица $\delta = (\delta_j^i)$ является единицей для обоих произведений.

Определение 3.2.8. **Бикольцо** $\mathcal{A}$ *- это множество, на котором мы определили унарную операцию, называемую транспозицией, и три бинарных операции, называемые* ${}_*{}^*$*-произведение,* ${}^*{}_*$*-произведение и сумма, такие что*

- ${}_*{}^*$*-произведение и сумма определяют структуру кольца на* $\mathcal{A}$
- ${}^*{}_*$*-произведение и сумма определяют структуру кольца на* $\mathcal{A}$
- *оба произведения имеют общую единицу* $\delta$
- *произведения удовлетворяют равенству*
$$(a_*{}^*b)^T = a^T{}^*{}_*b^T$$
- *транспозиция единицы есть единица*
$$\delta^T = \delta \tag{3.2.10}$$
- *двойная транспозиция есть исходный элемент*
$$(a^T)^T = a \tag{3.2.11}$$

□

Теорема 3.2.9.

$$(a^*{}_*b)^T = (a^T)_*{}^*(b^T) \tag{3.2.12}$$

Доказательство. Мы можем доказать (3.2.12) в случае матриц тем же образом, что мы доказали (3.2.8). Тем не менее для нас более важно показать, что (3.2.12) следует непосредственно из (3.2.8).

[3.3] Запись в равенстве (3.2.7) означает, что при обращении Адамара столбцы и строки меняются местами. Мы можем формально записать правую часть равенства (3.2.7) следующим образом

$$(a_i^j)^{-1} = \frac{1}{a_j^i}$$

Применяя (3.2.11) к каждому слагаемому в левой части (3.2.12), мы получим

$$(a^*{}_*b)^T = ((a^T)^{T*}{}_*(b^T)^T)^T \tag{3.2.13}$$

Из (3.2.13) и (3.2.8) следует, что

$$(a^*{}_*b)^T = ((a^T{}_*{}^*b^T)^T)^T \tag{3.2.14}$$

(3.2.12) следует из (3.2.14) и (3.2.11). □

Определение 3.2.10. *Мы определим* ${}_*{}^*$**-степень** $\mathcal{A}$*-числа* $a$*, пользуясь рекурсивным правилом*

$$a^{0_*{}^*} = \delta \tag{3.2.15}$$

$$a^{n_*{}^*} = a^{n-1_*{}^*}{}_*{}^*a \tag{3.2.16}$$

□

Теорема 3.2.11.

$$a^{1_*{}^*} = a \tag{3.2.17}$$

Доказательство. Равенство

$$a^{1_*{}^*} = a^{0_*{}^*}{}_*{}^*a = \delta_*{}^*a = a \tag{3.2.18}$$

является следствием равенств (3.2.15), (3.2.16). Равенство (3.2.17) является следствием равенства (3.2.18). □

Определение 3.2.12. *Мы определим* ${}^*{}_*$**-степень** $\mathcal{A}$*-числа* $a$*, пользуясь рекурсивным правилом*

$$a^{0^*{}_*} = \delta \tag{3.2.19}$$

$$a^{n^*{}_*} = a^{n-1^*{}_*}{}^*{}_*a \tag{3.2.20}$$

□

Теорема 3.2.13.

$$a^{1^*{}_*} = a \tag{3.2.21}$$

Доказательство. Равенство

$$a^{1^*{}_*} = a^{0^*{}_*}{}^*{}_*a = \delta^*{}_*a = a \tag{3.2.22}$$

является следствием равенств (3.2.19), (3.2.20). Равенство (3.2.21) является следствием равенства (3.2.22). □

Теорема 3.2.14.

$$(a^T)^{n_*{}^*} = (a^{n^*{}_*})^T \tag{3.2.23}$$

$$(a^T)^{n^*{}_*} = (a^{n_*{}^*})^T \tag{3.2.24}$$

Доказательство. Мы проведём доказательство индукцией по $n$.

При $n = 0$ утверждение непосредственно следует из равенств (3.2.15), (3.2.19) и (3.2.10).

Допустим утверждение справедливо при $n = k - 1$

$$(a^T)^{n-1_*{}^*} = (a^{n-1^*{}_*})^T \tag{3.2.25}$$

Из (3.2.16) следует

$$(a^T)^{k_*{}^*} = (a^T)^{k-1_*{}^*}{}_*{}^*a^T \tag{3.2.26}$$

Из (3.2.26) и (3.2.25) следует

$$(a^T)^{k_*{}^*} = (a^{k-1^*{}_*})^T{}_*{}^*a^T \tag{3.2.27}$$

Из (3.2.27) и (3.2.12) следует

$$(a^T)^{k_*{}^*} = (a^{k-1^*{}_*}{}^*{}_*a)^T \tag{3.2.28}$$

Из (3.2.26) и (3.2.20) следует (3.2.23).

Мы можем доказать (3.2.24) подобным образом. □

Определение 3.2.15. *$\mathcal{A}$-число $a^{-1_*{}^*}$ - это* ${}_*{}^*$**-обратный элемент** *$\mathcal{A}$-числа $a$, если*

$$a_*{}^*a^{-1_*{}^*} = \delta \tag{3.2.29}$$

*$\mathcal{A}$-число $a^{-1^*{}_*}$ - это* ${}^*{}_*$**-обратный элемент** *$\mathcal{A}$-числа $a$, если*

$$a^*{}_*a^{-1^*{}_*} = \delta \tag{3.2.30}$$

□

Теорема 3.2.16. *Предположим, что $\mathcal{A}$-число $a$ имеет ${}_*{}^*$-обратное $\mathcal{A}$-число. Тогда транспонированное $\mathcal{A}$-число $a^T$ имеет ${}^*{}_*$-обратное $\mathcal{A}$-число и эти $\mathcal{A}$-числа удовлетворяют равенству*

$$(a^T)^{-1^*{}_*} = (a^{-1_*{}^*})^T \tag{3.2.31}$$

*Предположим, что $\mathcal{A}$-число $a$ имеет ${}^*{}_*$-обратное $\mathcal{A}$-число. Тогда транспонированное $\mathcal{A}$-число $a^T$ имеет ${}_*{}^*$-обратное $\mathcal{A}$-число и эти $\mathcal{A}$-числа удовлетворяет равенству*

$$(a^T)^{-1_*{}^*} = (a^{-1^*{}_*})^T \tag{3.2.32}$$

Доказательство. Если мы возьмём транспонирование обеих частей (3.2.29) и применим (3.2.10), мы получим

$$(a_*{}^*a^{-1_*{}^*})^T = \delta^T = \delta$$

Применяя (3.2.8), мы получим

$$\delta = a^T{}^*{}_*(a^{-1_*{}^*})^T \tag{3.2.33}$$

(3.2.31) следует из сравнения (3.2.30) и (3.2.33).

Мы можем доказать (3.2.32) подобным образом. □

Теоремы 3.2.7, 3.2.9, 3.2.14 и 3.2.16 показывают, что существует двойственность между ${}_*{}^*$-произведением и ${}^*{}_*$-произведением. Мы можем объединить эти утверждения.

Теорема 3.2.17 (**принцип двойственности для бикольца**). *Пусть $\mathfrak{A}$ - истинное утверждение о бикольце $\mathcal{A}$. Если мы заменим одновременно*

- *$a \in \mathcal{A}$ и $a^T$*
- *${}_*{}^*$-произведение и ${}^*{}_*$-произведение*

*то мы снова получим истинное утверждение.*

Теорема 3.2.18 (**принцип двойственности для бикольца матриц**). *Пусть $\mathcal{A}$ является бикольцом матриц. Пусть $\mathfrak{A}$ - истинное утверждение о матрицах. Если мы заменим одновременно*

- *столбцы и строки всех матриц*
- *${}_*{}^*$-произведение и ${}^*{}_*$-произведение*

*то мы снова получим истинное утверждение.*

Доказательство. Непосредственное следствие теоремы 3.2.17. □

Замечание 3.2.19. *В выражении*

$$a{}_*{}^*b{}_*{}^*c \tag{3.2.34}$$

*мы выполняем операцию умножения слева направо. Однако мы можем выполнять операцию умножения справа налево. Тогда выражение (3.2.34) примет вид*

$$c{}^*{}_*b{}^*{}_*a$$

*Мы сохраним правило, что показатель степени и индексы записываются справа от выражения. Например, если исходное выражение имеет вид*

$$a^{-1_*{}^*}{}_*{}^*b^i$$

*то выражение, читаемое справа налево, примет вид*

$$b^i{}^*{}_*a^{-1^*{}_*}$$

*Если задать порядок, в котором мы записываем индексы, то мы можем утверждать, что мы читаем выражение сверху вниз, читая сперва верхние индексы, потом нижние. Мы можем прочесть выражение снизу вверх. При этом мы дополним правило, что символы операции также читаются в том же направлении, что и индексы. Например, выражение*

$$a^i{}_*{}^*b^{-1_*{}^*}=c^i$$

*прочтённое снизу вверх, имеет вид*

$$a_i{}^*{}_*b^{-1^*{}_*}=c_i$$

*Согласно принципу двойственности, если верно одно утверждение, то верно и другое.* □

Теорема 3.2.20. *Если матрица $a$ имеет ${}_*{}^*$-обратную матрицу, то для любых матриц $b$ и $c$ из равенства*

$$b{}_*{}^*a=c{}_*{}^*a \tag{3.2.35}$$

*следует равенство*

$$b=c \tag{3.2.36}$$

Доказательство. Равенство (3.2.36) следует из (3.2.35), если обе части равенства (3.2.35) умножить на $a^{-1_*{}^*}$. □

## 3.3. Квазидетерминант

Теорема 3.3.1. *Предположим, что* $n\times n$ *матрица* $a$ *имеет* ${}_*{}^*$*-обратную матрицу.*[3.4] *Тогда* $k\times k$ *минорная матрица* ${}_*{}^*$*-обратной матрицы удовлетворяет равенству*

$$\left((a^{-1_*{}^*})^I_J\right)^{-1_*{}^*}=a^J_I-a^J_{[I]}{}_*{}^*\left(a^{[J]}_{[I]}\right)^{-1_*{}^*}{}_*{}^*a^{[J]}_I \tag{3.3.1}$$

Доказательство. Рассмотрим минорную матрицу $(a^{-1_*{}^*})_J$, полученную из матрицы $a^{-1_*{}^*}$ выбором столбцов с индексом из множества $J$. Определение (3.2.29) ${}_*{}^*$-обратной матрицы приводит к системе линейных уравнений

$$a^{[J]}{}_*{}^*(a^{-1_*{}^*})_J=0 \tag{3.3.2}$$

$$a^J{}_*{}^*(a^{-1_*{}^*})_J=\delta \tag{3.3.3}$$

где $a^J$ - минорная матрица, полученную из матрицы $a$ выбором строк с индексом из множества $J$, и $a^{[J]}$ - минорная матрица, полученный из матрицы $a$ удалением строк с индексом из множества $J$. Равенство

$$a^{[J]}_{[I]}{}_*{}^*(a^{-1_*{}^*})^{[I]}_J+a^{[J]}_I{}_*{}^*(a^{-1_*{}^*})^I_J=0 \tag{3.3.4}$$

является следствием равенства (3.3.2), если мы рассмотрим минорную матрицу $a^{[J]}$ как объединение минорных матриц $a^{[J]}_I$ и $a^{[J]}_{[I]}$. Равенство

$$a^J_{[I]}{}_*{}^*(a^{-1_*{}^*})^{[I]}_J+a^J_I{}_*{}^*(a^{-1_*{}^*})^I_J=\delta \tag{3.3.5}$$

является следствием равенства (3.3.3), если мы рассмотрим минорную матрицу $a^J$ как объединение минорных матриц $a^J_I$ и $a^J_{[I]}$. Мы умножим равенство (3.3.4) на $\left(a^{[J]}_{[I]}\right)^{-1_*{}^*}$

$$(a^{-1_*{}^*})^{[I]}_J+\left(a^{[J]}_{[I]}\right)^{-1_*{}^*}{}_*{}^*a^{[J]}_I{}_*{}^*(a^{-1_*{}^*})^I_J=0 \tag{3.3.6}$$

Равенство

$$-a^J_{[I]}{}_*{}^*\left(a^{[J]}_{[I]}\right)^{-1_*{}^*}{}_*{}^*a^{[J]}_I{}_*{}^*(a^{-1_*{}^*})^I_J+a^J_I{}_*{}^*(a^{-1_*{}^*})^I_J=\delta \tag{3.3.7}$$

является следствием равенств (3.3.6), (3.3.5). Равенство (3.3.1) следует из (3.3.7), если мы обе части равенства (3.3.7) умножим на $\left((a^{-1_*{}^*})^I_J\right)^{-1_*{}^*}$. □

Теорема 3.3.2. *Предположим, что* $n\times n$ *матрица* $a$ *имеет* ${}_*{}^*$*-обратную матрицу. Тогда элементы* ${}_*{}^*$*-обратной матрицы удовлетворяют равенству*[3.5]

$$(a^{-1_*{}^*})^i_j=\left(a^j_i-a^j_{[i]}{}_*{}^*\left(a^{[j]}_{[i]}\right)^{-1_*{}^*}{}_*{}^*a^{[j]}_i\right)^{-1} \tag{3.3.8}$$

$$(\mathcal{H}a^{-1_*{}^*})^j_i=a^j_i-a^j_{[i]}{}_*{}^*\left(a^{[j]}_{[i]}\right)^{-1_*{}^*}{}_*{}^*a^{[j]}_i \tag{3.3.9}$$

---

[3.4] Это утверждение и его доказательство основаны на утверждении 1.2.1 из [4] (page 8) для матриц над свободным кольцом с делением.

[3.5] Запись в равенстве (3.3.8) означает, что при обращении матрицы столбцы и строки меняются местами.

Доказательство. Равенство (3.3.8) является следствием равенства (3.3.1), если мы положим $I = i$, $J = j$. Равенство (3.3.9) является следствием равенств (3.2.7), (3.3.8). □

Согласно [4], page 3 у нас нет определения детерминанта в случае алгебры с делением.[3.6] Тем не менее, мы можем определить квазидетерминант, который в конечном итоге даёт похожую картину. В определении, данном ниже, мы следуем определению [4]-1.2.2.

Определение 3.3.3. **$\binom{j}{i}$-${}_*{}^*$-квазидетерминант** *$n \times n$ матрицы $a$ - это формальное выражение*

$$\det({}_*{}^*){}_i^j\, a = (\mathcal{H}a^{-1_*{}^*})_i^j \tag{3.3.10}$$

*Согласно замечанию 3.1.2 мы можем рассматривать $\binom{j}{i}$-${}_*{}^*$-квазидетерминант как элемент матрицы*

$$\det({}_*{}^*)\, a = \mathcal{H}a^{-1_*{}^*}$$

*которую мы будем называть* **${}_*{}^*$-квазидетерминантом**. □

Теорема 3.3.4. *Выражение для ${}_*{}^*$-обратной матрицы имеет вид*

$$a^{-1_*{}^*} = \mathcal{H}\det({}_*{}^*)a \tag{3.3.11}$$

Доказательство. (3.3.11) следует из (3.3.10). □

Теорема 3.3.5. *Выражение для $\binom{j}{i}$-${}_*{}^*$-квазидетерминанта имеет любую из следующих форм*[3.7]

$$\det({}_*{}^*){}_i^j\, a = a_i^j - a_{[i]}^j{}_*{}^* \left(a_{[i]}^{[j]}\right)^{-1_*{}^*} {}_*{}^* a_i^{[j]} \tag{3.3.12}$$

$$\det({}_*{}^*){}_i^j\, a = a_i^j - a_{[i]}^j{}_*{}^* \mathcal{H}\det({}_*{}^*) a_{[i]}^{[j]}{}_*{}^* a_i^{[j]} \tag{3.3.13}$$

Доказательство. Утверждение следует из (3.3.9) и (3.3.10). □

Теорема 3.3.6. *Рассмотрим матрицу*

$$\begin{pmatrix} a_1^1 & a_2^1 \\ a_1^2 & a_2^2 \end{pmatrix}$$

*Тогда*

$$\det({}_*{}^*)a = \begin{pmatrix} a_1^1 - a_2^1(a_2^2)^{-1}a_1^2 & a_2^1 - a_1^1(a_1^2)^{-1}a_2^2 \\ a_1^2 - a_2^2(a_2^1)^{-1}a_1^1 & a_2^2 - a_1^2(a_1^1)^{-1}a_2^1 \end{pmatrix} \tag{3.3.14}$$

---

[3.6] Профессор Кирчей пользуется двойным определителем (смотри определение [10]-5.28) для решения системы линейных уравнений в алгебре кватернионов. Профессор Кирчей решил также задачу о собственных значениях в алгебре кватернионов пользуясь двойным определителем (смотри раздел [16]-2.5). Я ограничился рассмотрением квазидетерминантов, так как меня интересует более широкий класс алгебр.

[3.7] Мы можем дать подобное доказательство для **$\binom{j}{i}$-${}^*{}_*$-квазидетерминанта**. Однако мы можем записать соответствующие утверждения, опираясь на принцип двойственности. Так, если прочесть равенство (3.3.12) справа налево, то мы получим равенство

$$\det({}^*{}_*)a_i^j = a_i^j - a_i^{[j]}{}^*{}_* \left(a_{[i]}^{[j]}\right)^{-1^*{}_*} {}^*{}_* a_{[i]}^j$$

$$\det({}^*{}_*)a_i^j = a_i^j - a_i^{[j]}{}^*{}_* \mathcal{H}\det({}^*{}_*)^{[j]} A_{[i]}{}^*{}_* a_{[i]}^j$$

$$\det({}^*{}_*)a = \begin{pmatrix} a_1^1 - a_1^2(a_2^2)^{-1}a_2^1 & a_2^1 - a_2^2(a_1^2)^{-1}a_1^1 \\ a_1^2 - a_1^1(a_2^1)^{-1}a_2^2 & a_2^2 - a_2^1(a_1^1)^{-1}a_1^2 \end{pmatrix} \tag{3.3.15}$$

$$a^{-1_*{}^*} = \begin{pmatrix} (a_1^1 - a_1^2(a_2^2)^{-1}a_2^1)^{-1} & (a_1^2 - a_1^1(a_2^1)^{-1}a_2^2)^{-1} \\ (a_2^1 - a_2^2(a_1^2)^{-1}a_1^1)^{-1} & (a_2^2 - a_2^1(a_1^1)^{-1}a_1^2)^{-1} \end{pmatrix} \tag{3.3.16}$$

Доказательство. Согласно равенству (3.3.12)

$$\det({}_*{}^*)_1^1\, a = a_1^1 - a_2^1(a_2^2)^{-1}a_1^2 \tag{3.3.17}$$

$$\det({}_*{}^*)_1^2\, a = a_2^1 - a_1^1(a_1^2)^{-1}a_2^2 \tag{3.3.18}$$

$$\det({}_*{}^*)_2^1\, a = a_1^2 - a_2^2(a_2^1)^{-1}a_1^1 \tag{3.3.19}$$

$$\det({}_*{}^*)_2^2\, a = a_2^2 - a_1^2(a_1^1)^{-1}a_2^1 \tag{3.3.20}$$

(3.3.17), (3.3.18), (3.3.19), (3.3.20). □

Теорема 3.3.7.

$$\det({}_*{}^*)_j^i\, a^T = \det({}^*{}_*)_i^j\, a \tag{3.3.21}$$

Доказательство. Согласно (3.3.10) и (3.2.6)

$$\det({}_*{}^*)_j^i\, a^T = (((a^T)^{-1_*{}^*})^{\cdot -}_{i}{}^{j}_{\cdot -})^{-1}$$

Пользуясь теоремой 3.2.16, мы получим

$$\det({}_*{}^*)_j^i\, a^T = (((a^{-1^*{}_*})^T)^{\cdot -}_{i}{}^{j}_{\cdot -})^{-1}$$

Пользуясь (3.2.5), мы имеем

$$\det({}_*{}^*)_j^i\, a^T = ((a^{-1^*{}_*})^{i}_{\cdot -}{}^{\cdot -}_{j})^{-1} \tag{3.3.22}$$

Пользуясь (3.3.22), (3.2.6), (3.3.10), мы получим (3.3.21). □

Теорема 3.3.7 расширяет принцип двойственности, теорема 3.2.18, на утверждения о квазидетерминантах и утверждает, что одно и тоже выражение является ${}_*{}^*$-квазидетерминантом матрицы $a$ и ${}^*{}_*$-квазидетерминантом матрицы $a^T$. Пользуясь этой теоремой, мы можем записать любое утверждение о ${}^*{}_*$-матрице, опираясь на подобное утверждение о ${}_*{}^*$-матрице.

Теорема 3.3.8 (принцип двойственности). *Пусть $\mathfrak{A}$ - истинное утверждение о бикольце матриц. Если мы одновременно заменим*

- *строку и столбец*
- *${}_*{}^*$-квазидетерминант и ${}^*{}_*$-квазидетерминант*

*то мы снова получим истинное утверждение.*

Теорема 3.3.9.

$$(ma)^{-1_*{}^*} = a^{-1_*{}^*}m^{-1} \tag{3.3.23}$$

$$(am)^{-1_*{}^*} = m^{-1}a^{-1_*{}^*} \tag{3.3.24}$$

Доказательство. Мы докажем равенство (3.3.23) индукцией по размеру матрицы.

Для $1 \times 1$ матрицы утверждение очевидно, так как

$$(ma)^{-1_*{}^*} = \left((ma)^{-1}\right) = \left(a^{-1}m^{-1}\right) = \left(a^{-1}\right) m^{-1} = a^{-1_*{}^*} m^{-1}$$

Допустим утверждение справедливо для $(n-1) \times (n-1)$ матрицы. Тогда из равенства (3.3.1) следует

$$\begin{aligned}
(((ma)^{-1_*{}^*})^{I}_{J})^{-1_*{}^*} =& (ma)^{J}_{I} - (ma)^{J}_{[I]}{}_*{}^* \left((ma)^{[J]}_{[I]}\right)^{-1_*{}^*} {}_*{}^*(ma)^{[J]}_{I} \\
=& ma^{J}_{I} - m\ a^{J}_{[I]}{}_*{}^* \left(a^{[J]}_{[I]}\right)^{-1_*{}^*}\ m^{-1}{}_*{}^* m\ a^{[J]}_{I} \\
=& m\ a^{J}_{I} - m\ a^{J}_{[I]}{}_*{}^* \left(a^{[J]}_{[I]}\right)^{-1_*{}^*} {}_*{}^* a^{[J]}_{I}
\end{aligned}$$

$$(((ma)^{-1_*{}^*})^{I}_{J})^{-1_*{}^*} = m\ (a^{-1_*{}^*})^{I}_{J} \tag{3.3.25}$$

Из равенства (3.3.25) следует равенство (3.3.23). Аналогично доказывается равенство (3.3.24). □

Глава 4

# Многочлен

## 4.1. Многочлен над ассоциативной $D$-алгеброй

Пусть $D$ - коммутативное кольцо и $A$ - ассоциативная $D$-алгебра с единицей.

Теорема 4.1.1. *Пусть $p_k(x)$ -* **одночлен степени** *$k$ над $D$-алгеброй $A$. Тогда*

4.1.1.1: *Одночлен степени 0 имеет вид $p_0(x) = a_0$, $a_0 \in A$.*
4.1.1.2: *Если $k > 0$, то*

$$p_k(x) = p_{k-1}(x)xa_k$$

*где $a_k \in A$.*

Доказательство. Мы докажем утверждение теоремы индукцией по степени $n$ одночлена.

Пусть $n = 0$. Так как одночлен $p_0(x)$ является константой, то мы получаем утверждение 4.1.1.1.

Пусть $n = k$. Последний множитель одночлена $p_k(x)$ является либо $a_k \in A$, либо имеет вид $x^l$, $l \geq 1$. В последнем случае мы положим $a_k = 1$. Множитель, предшествующий $a_k$, имеет вид $x^l$, $l \geq 1$. Мы можем представить этот множитель в виде $x^{l-1}x$. Следовательно, утверждение доказано. □

Определение 4.1.2. *Обозначим $A_k[x]$ абелевую группу, порождённую множеством одночленов степени $k$. Элемент $p_k(x)$ абелевой группы $A_k[x]$ называется* **однородным многочленом степени** *$k$.* □

Теорема 4.1.3. *Отображение*

$$f : A^{k+1} \to A_k[x]$$

*определённое равенством*

$$f(a_0, a_1, ..., a_k) = (a_0, a_1, ..., a_k) \circ x^k$$

*является полилинейным отображением.*

Доказательство. Пусть $d \in D$. Из равенств

$$\begin{aligned} f(a_0, ..., da_i, ..., a_k) &= a_0x...(da_i)...xa_k = d(a_0x...a_i...xa_k) \\ &= df(a_0, ..., a_i, ..., a_k) \\ f(a_0, ..., a_i + b_i, ..., a_k) &= a_0x...(a_i + b_i)...xa_k \\ &= (a_0x...a_i...xa_k) + (a_0x...b_i...xa_k) \\ &= f(a_0, ..., a_i, ..., a_k) + f(a_0, ..., b_i, ..., a_k) \end{aligned}$$

следует что отображение $f$ линейно по $a_i$. Следовательно, отображение $f$ - полилинейное отображение. □

Теорема 4.1.4. *Определено линейное отображение*

$$(a_0 \otimes a_1 \otimes ... \otimes a_k) \circ x^k = (a_0, a_1, ..., a_k) \circ x^k$$

Доказательство. Утверждение теоремы является следствием теорем [9]-3.6.4, 4.1.3. □

Следствие 4.1.5. *Однородный многочлен $p(x)$ может быть записан в виде*

$$p(x) = a \circ x^k \quad a \in A^{(k+1)\otimes}$$

□

Определение 4.1.6. *Обозначим*

$$A[x] = \bigoplus_{n=0}^{\infty} A_n[x]$$

*прямую сумму*[4.1] *$A$-модулей $A_n[x]$. Элемент $p(x)$ $A$-модуля $A[x]$ называется* **многочленом** *над $D$-алгеброй $A$.* □

Следовательно, многочлен степени $n$ может быть записан в виде

$$p(x) = a_0 + a_1 \circ x + ... + a_n \circ x^n \quad a_i \in A^{i+1\otimes} \quad i = 0, ..., n$$

Определение 4.1.7. *Билинейное отображение*

$$* : A^{n\otimes} \times A^{m\otimes} \to A^{n+m-1\otimes}$$

*определено равенством*

$$(a_1 \otimes ... \otimes a_n) * (b_1 \otimes ... \otimes b_n) = a_1 \otimes ... \otimes a_{n-1} \otimes a_n b_1 \otimes b_2 \otimes ... \otimes b_n$$

□

Теорема 4.1.8. *Для любых тензоров $a \in A^{n\otimes}$, $b \in A^{m\otimes}$, $c \in A^{k\otimes}$*

$$(a * b) * c = a * (b * c) \tag{4.1.1}$$

Доказательство. Согласно определению 4.1.7, теорему достаточно доказать для тензоров $a = a_0 \otimes ... \otimes a_n$, $b = b_0 \otimes ... \otimes b_m$, $c = c_0 \otimes ... \otimes c_k$. Равенство (4.1.1) является следствием равенства

$$\begin{aligned}
&((a_0 \otimes ... \otimes a_n) * (b_0 \otimes b_m)) * (c_0 \otimes ... \otimes c_k) \\
=&(a_0 \otimes ... \otimes a_n b_0 \otimes ... \otimes b_m) * (c_0 \otimes ... \otimes c_k) \\
=&a_0 \otimes ... \otimes a_n b_0 \otimes ... \otimes b_m c_0 \otimes ... \otimes c_k \\
=&(a_0 \otimes ... \otimes a_n) * (b_0 \otimes ... \otimes b_m c_0 \otimes ... \otimes c_k) \\
=&(a_0 \otimes ... \otimes a_n) * ((b_0 \otimes ... \otimes b_m) * (c_0 \otimes ... \otimes c_k)
\end{aligned}$$

□

Теорема 4.1.9. *Для любых тензоров $a \in A^{n\otimes}$, $b \in A^{m\otimes}$ произведение однородных многочленов $a \circ x^n$, $b \circ x^m$ определено равенством*

$$(a * b) \circ x^{n+m} = (a \circ x^n)(b \circ x^m) \tag{4.1.2}$$

[4.1] Смотри определение 2.2.5.

Доказательство. Согласно определению 4.1.7, теорему достаточно доказать для тензоров $a = a_0 \otimes ... \otimes a_n$, $b = b_0 \otimes ... \otimes b_m$. Мы проведём доказательство индукцией по $m$.

Из равенства

$$\begin{aligned}((a_0 \otimes a_1 \otimes ... \otimes a_n) * b_0) \circ x^n &= (a_0 \otimes a_1 \otimes ... \otimes a_n b_0) \circ x^n\\ &= a_0 x a_1 x ... x a_n b_0 = (a_0 x a_1 x ... x a_n) b_0\\ &= ((a_0 \otimes a_1 \otimes ... \otimes a_n) \circ x^n)(b_0 \circ x^0)\end{aligned}$$

следует, что теорема верна при $m = 0$.

Из равенства

$$\begin{aligned}&((a_0 \otimes a_1 \otimes ... \otimes a_n) * (b_0 \otimes b_1)) \circ x^{n+1}\\ =&(a_0 \otimes a_1 \otimes ... \otimes a_n b_0 \otimes b_1) \circ x^{n+1}\\ =&a_0 x a_1 x ... x a_n b_0 x b_1 = (a_0 x a_1 x ... x a_n)(b_0 x b_1)\\ =&((a_0 \otimes a_1 \otimes ... \otimes a_n) \circ x^n)((b_0 \otimes b_1) \circ x^1)\end{aligned} \tag{4.1.3}$$

следует, что теорема верна при $m = 1$.

Пусть теорема верна при $m = k$

$$(a * (b_0 \otimes ... \otimes b_k)) \circ x^{n+k} = (a \circ x^n)((b_0 \otimes ... \otimes b_k) \circ x^k) \tag{4.1.4}$$

Равенство

$$\begin{aligned}&(a * (b_0 \otimes ... \otimes b_k \otimes b_{k+1})) \circ x^{n+k+1}\\ =&(a * ((b_0 \otimes ... \otimes b_k) * (1 \otimes b_{k+1})) \circ x^{n+k+1}\\ =&((a * (b_0 \otimes ... \otimes b_k)) * (1 \otimes b_{k+1})) \circ x^{n+k+1}\\ =&((a * (b_0 \otimes ... \otimes b_k)) \circ x^{n+k})((1 \otimes b_{k+1}) \circ x)\\ =&((a \circ x^n)((b_0 \otimes ... \otimes b_k) \circ x^k))((1 \otimes b_{k+1}) \circ x)\\ =&(a \circ x^n)(((b_0 \otimes ... \otimes b_k) \circ x^k)((1 \otimes b_{k+1}) \circ x))\\ =&(a \circ x^n)(((b_0 \otimes ... \otimes b_k) * (1 \otimes b_{k+1})) \circ x^{k+1})\\ =&(a \circ x^n)((b_0 \otimes ... \otimes b_k \otimes b_{k+1}) \circ x^{k+1})\end{aligned} \tag{4.1.5}$$

является следствием равенств (4.1.1). (4.1.3), (4.1.4). Из равенства (4.1.5) следует, что теорема верна при $m = k + 1$. □

## 4.2. Производная многочлена

Соглашение 4.2.1. *Для $n \ge 0$, пусть $SO(k, n)$ - множество перестановок*

$$\sigma = \begin{pmatrix} y_1 & ... & y_k & x_{k+1} & ... & x_n\\ \sigma(y_1) & ... & \sigma(y_k) & \sigma(x_{k+1}) & ... & \sigma(x_n)\end{pmatrix}$$

*таких, что каждая перестановка $\sigma$ сохраняет порядок переменных $x_i$: если $i < j$, то в кортеже*

$$(\sigma(y_1), ..., \sigma(y_k), \sigma(x_{k+1}), ..., \sigma(x_n))$$

*$x_i$ предшествует $x_j$.* □

Лемма 4.2.2. *Мы можем пронумеровать множество перестановок $SO(1, n)$ индексом $i$, $1 \le i \le n$, таким образом, что*

4.2.2.1: $\sigma_1(y) = y$.
4.2.2.2: *Если* $i > 1$, *то* $\sigma_i(x_i) = y$.

Доказательство. Так как порядок переменных $x_2, ..., x_n$ в перестановке $\sigma \in SO(1,n)$ не зависит от перестановки, перестановки $\sigma \in SO(1,n)$ отличаются местом, которое занимает переменная $y$. Соответственно, мы можем пронумеровать множество перестановок $SO(1,n)$ индексом, значение которого соответствует номеру позиции переменной $y$. □

Лемма 4.2.2 имеет простую интерпретацию. Предположим, что $n-1$ белых шаров и один чёрный лежат в пенале. Чёрный шар - самый левый шар, и белые шары пронумерованы от 2 до $n$ в том порядке, как они лежат в пенале. Суть перестановки $\sigma_k$ состоит в том, что мы вынимаем чёрный шар из пенала и затем мы ложим его в ячейку с номером $k$. При этом, белые шары с номером, не превышающим $k$, сдвигаются влево.

Лемма 4.2.3. *Для* $n > 0$, *пусть*

$$SO^+(1,n) = \{\sigma : \sigma = (\tau(y), \tau(x_2), ..., \tau(x_n), x_{n+1}), \tau \in SO(1,n)\} \tag{4.2.1}$$

*Тогда*

$$SO(1,n+1) = SO^+(1,n) \cup \{(x_2, ..., x_{n+1}, y)\} \tag{4.2.2}$$

Доказательство. Пусть $\sigma \in SO^+(1,n)$. Согласно определению (4.2.1), существует перестановка $\tau \in SO(1,n)$ такая, что

$$(\sigma(y), \sigma(x_1), ..., \sigma(x_{n+1})) = (\tau(y), \tau(x_2), ..., \tau(x_n), x_{n+1})$$

Согласно соглашению 4.2.1, из утверждения $i < j < n+1$ следует, что в кортеже

$$(\sigma(y), \sigma(x_1), ..., \sigma(x_{n+1})) = (\tau(y), \tau(x_2), ..., \tau(x_n), x_{n+1})$$

переменная $x_j$ расположена между переменными $x_i$ и $x_{n+1}$. Согласно соглашению 4.2.1, $\sigma \in SO(1,n+1)$. Следовательно

$$SO^+(1,n) \subseteq SO(1,n+1) \tag{4.2.3}$$

Согласно лемме 4.2.2, множество $SO^+(1,n)$ имеет $n$ перестановок.

Пусть $\sigma = (x_2, ..., x_{n+1}, y)$. Согласно соглашению 4.2.1,

$$(x_2, ..., x_{n+1}, y) \in SO(1,n+1) \tag{4.2.4}$$

Согласно определению (4.2.1), $\sigma \notin SO^+(1,n)$.

Следовательно, мы перечислили $n+1$ элемент множества $SO(1,n+1)$. Согласно лемме 4.2.2, утверждение (4.2.2) является следствием утверждений (4.2.3), (4.2.4). □

Теорема 4.2.4. *Для произвольного одночлена*

$$p_n(x) = (a_0 \otimes ... \otimes a_n) \circ x^n$$

*производная имеет вид*

$$\begin{array}{c}\displaystyle\frac{dp_n(x)}{dx} \circ dx = \sum_{\sigma \in SO(1,n)} (a_0 \otimes ... \otimes a_n) \circ (\sigma(dx), \sigma(x_2), ..., \sigma(x_n))\\ x_2 = ... = x_n = x\end{array} \tag{4.2.5}$$

Доказательство. Для $n=1$, отображение

$$p_1(x)=(a_0\otimes a_1)\circ x=a_0xa_1$$

является линейным отображением. Согласно теореме A.1.5 и соглашению 4.2.1

$$\frac{dp_1(x)}{dx}\circ dx=a_0dxa_1=(a_0\otimes a_1)\circ dx=\sum_{\sigma\in SO(1,1)}(a_0\otimes a_1)\circ\sigma(dx) \tag{4.2.6}$$

Пусть утверждение справедливо для $n-1$

$$\begin{array}{c}\displaystyle\frac{dp_{n-1}(x)}{dx}\circ dx=\sum_{\sigma\in SO(1,n-1)}(a_0\otimes...\otimes a_{n-1})\circ\sigma(dx,x_2,...,x_{n-1})\\ x_2=...=x_{n-1}=x\end{array} \tag{4.2.7}$$

Если

$$p_n(x)=p_{n-1}(x)xa_n \tag{4.2.8}$$

то согласно теореме A.1.10 и определению (4.2.8)

$$\frac{dp_n(x)}{dx}\circ dx=\left(\frac{dp_{n-1}(x)}{dx}\circ dx\right)xa_n+p_{n-1}(x)\left(\frac{dxa_n}{dx}\circ dx\right) \tag{4.2.9}$$

Равенство

$$\begin{array}{l}\displaystyle\frac{dp_n(x)}{dx}\circ dx\\ \displaystyle=\sum_{\sigma\in SO(1,n-1)}(a_0\otimes...\otimes a_{n-1})\circ(\sigma(dx),\sigma(x_2),...,\sigma(x_{n-1}))xa_n\\ +((a_0\otimes...\otimes a_{n-1})\circ(x_2,...,x_n))dxa_n\\ x_2=...=x_n=x\end{array} \tag{4.2.10}$$

следует из (A.1.1), (A.1.2), (4.2.7), (4.2.9). Равенство

$$\begin{array}{l}\displaystyle\frac{dp_n(x)}{dx}\circ dx\\ \displaystyle=\sum_{\sigma\in SO(1,n-1)}(a_0\otimes...\otimes a_n)\circ(\sigma(dx),\sigma(x_2),...,\sigma(x_{n-1}),x_n))\\ +(a_0\otimes...\otimes a_n)\circ(x_2,...,x_n,dx)\\ x_2=...=x_n=x\end{array} \tag{4.2.11}$$

следует из (4.2.10) и правила перемножения одночленов (определение 4.1.7 и теорема 4.1.9). Согласно лемме 4.2.3, равенство (4.2.5) является следствием равенства (4.2.11). □

Лемма 4.2.5. *Пусть* $1<k<n$. *Для любой пары перестановок* $\mu\in SO(k,n)$, $\nu\in SO(1,n-k)$ *существует единственная перестановка* $\sigma\in SO(k+1,n)$ *такая, что*

$$\sigma=(\mu(y_1),...,\mu(y_k),\nu(\mu(y_{k+1})),\nu(\mu(x_{k+2})),...,\nu(\mu(x_n)))$$

Доказательство. Согласно утверждению леммы мы сперва выполняем перестановку $\mu$

$$\begin{pmatrix}\mu(y_1) & ... & \mu(y_k) & \mu(x_{k+1}=y_{k+1}) & \mu(x_{k+2}) & ... & \mu(x_n)\end{pmatrix}$$

а затем выполняем перестановку $\nu$

$$\begin{pmatrix}y_1 & ... & y_k & \nu(y_{k+1}) & \nu(x_{k+2}) & ... & \nu(x_n)\end{pmatrix}$$

которую применяем только к переменным $y_{k+1}$, $x_{k+2}$, ..., $x_n$. Согласно соглашению 4.2.1, перестановка $\mu$ сохраняет порядок переменных $x_{k+1}=y_{k+1}$, $x_{k+2}$, ..., $x_n$. Следовательно, перестановка $\mu$ сохраняет порядок переменных $x_{k+2}$, ..., $x_n$. Согласно соглашению 4.2.1, перестановка $\nu$ сохраняет порядок переменных $x_{k+2}$, ..., $x_n$. Следовательно, перестановка $\sigma$ сохраняет порядок переменных $x_{k+2}$, ..., $x_n$. Согласно соглашению 4.2.1, $\sigma\in SO(k+1,n)$ . □

Лемма 4.2.6. *Пусть $1<k<n$. Для любой перестановки $\sigma\in SO(k+1,n)$ существует единственная пара перестановок $\mu\in SO(k,n)$, $\nu\in SO(1,n-k)$ такая, что*

$$\sigma=(\mu(y_1),...,\mu(y_k),\nu(\mu(y_{k+1})),\nu(\mu(x_{k+2})),...,\nu(\mu(x_n)))$$

Доказательство. Если в кортеже

$$(\sigma(y_1),...,\sigma(y_{k+1}),\sigma(x_{k+2}),...,\sigma(x_n))$$

$y_{k+1}$ предшествует $x_{k+2}$, то мы положим

$$\nu=\begin{pmatrix}y_{k+1} & x_{k+2} & ... & x_n\\ y_{k+1} & x_{k+2} & ... & x_n\end{pmatrix}$$

$$\mu=\begin{pmatrix}y_1 & ... & y_k & x_{k+1}=y_{k+1} & x_{k+2} & ... & x_n\\ \sigma(y_1) & ... & \sigma(y_k) & \sigma(y_{k+1}) & \sigma(x_{k+2}) & ... & \sigma(x_n)\end{pmatrix}$$

Поскольку перестановка $\nu$ сохраняет порядок переменных $x_i$, то, согласно соглашению 4.2.1, $\nu\in SO(1,n-k)$. В кортеже

$$(\mu(y_1),...,\mu(y_{k+1}),\mu(x_{k+2}),...,\mu(x_n))$$

$x_{k+1}=y_{k+1}$ предшествует $x_{k+2}$. Согласно соглашению 4.2.1, $x_{k+1}=y_{k+1}$ предшествует $x_{k+2}$, ..., $x_n$. Следовательно, $\mu\in SO(k,n)$.

Пусть в кортеже

$$(\sigma(y_1),...,\sigma(y_{k+1}),\sigma(x_{k+2}),...,\sigma(x_n))$$

$x_{k+2}$, ..., $x_j$ предшествует $y_{k+1}$. Тогда мы положим

$$\nu=\begin{pmatrix}x_{k+2} & x_{k+3} & ... & x_j & y_{k+1} & ... & x_n\\ y_{k+1} & x_{k+2} & ... & x_{j-1} & x_j & ... & x_n\end{pmatrix}$$

$$\mu=\begin{pmatrix}y_1 & ... & y_k & x_{k+1}=y_{k+1} & x_{k+2} & ... & x_n\\ \sigma(y_1) & ... & \sigma(y_k) & \sigma(y_{k+1}) & \sigma(x_{k+2}) & ... & \sigma(x_n)\end{pmatrix}$$

Так как перестановки $\mu$, $\nu$ сохраняет порядок переменных $x_i$, то $\mu\in SO(k,n)$, $\nu\in SO(1,n-k)$. □

Теорема 4.2.7. *Для произвольного одночлена*

$$p_n(x)=(a_0\otimes...\otimes a_n)\circ x^n$$

*производная порядка $k$ имеет вид*

$$\begin{aligned}&\frac{d^kp_n(x)}{dx^k}\circ(dx^1;...;dx^k)\\ &=\sum_{\sigma\in SO(k,n)}(a_0\otimes...\otimes a_n,\sigma)\circ(dx_1;...;dx_k;x_{k+1};...;x_n)\\ &x_{k+1}=...=x_n=x\end{aligned}\tag{4.2.12}$$

Доказательство. Для $k=1$ утверждение теоремы является утверждением теоремы 4.2.4. Пусть теорема верна для $k-1$. Тогда

$$\begin{aligned}&\frac{d^{k-1}p_n(x)}{dx^{k-1}}\circ(dx_1;...;dx_{k-1})\\ &=\sum_{\mu\in SO(k-1,n)}(a_0\otimes...\otimes a_n,\mu)\circ(dx_1;...;dx_{k-1};x_k;...;x_n)\\ &x_k=...=x_n=x\end{aligned}\tag{4.2.13}$$

Согласно определению [15]-4.1.4, равенство

$$\begin{aligned}&\frac{d^kp_n(x)}{dx^k}\circ(dx_1;...;dx_{k-1};dx_k)\\ &=\frac{d}{dx}\left(\frac{d^{k-1}p_n(x)}{dx^{k-1}}\circ(dx_1;...;dx_{k-1})\right)\circ dx^k\\ &=\frac{d}{dx}\left(\sum_{\mu\in SO(k-1,n)}(a_0\otimes...\otimes a_n,\mu)\circ(dx_1;...;dx_{k-1};x_k;...;x_n)\right)\circ dx^k\\ &x_k=...=x_n=x\end{aligned}\tag{4.2.14}$$

является следствием равенства (4.2.13). Согласно теореме 4.2.4 равенство

$$\begin{aligned}&\frac{d^kp_n(x)}{dx^k}\circ(dx_1;...;dx_{k-1};dx_k)\\ &=\sum_{\substack{\mu\in SO(k-1,n)\\ \nu\in SO(1,n-k+1)}}\begin{aligned}&(a_0\otimes...\otimes a_n)\circ(\mu(dx_1);...;\mu(dx_{k-1});\nu(\mu(dx_k));\\ &\nu(\mu(x_{k+1}));...;\nu(\mu(x_n))\end{aligned}\\ &x_{k+1}=...=x_n=x\end{aligned}\tag{4.2.15}$$

является следствием равенства (4.2.14). Согласно леммам 4.2.5, 4.2.6 равенство

$$\begin{aligned}&\frac{d^kp_n(x)}{dx^k}\circ(dx_1;...;dx_{k-1};dx_k)\\ &=\sum_{\sigma\in SO(k,n)}(a_0\otimes...\otimes a_n,\sigma)\circ(dx_1;...;dx_{k-1};dx_k;x_{k+1};...;x_n)\\ &x_{k+1}=...=x_n=x\end{aligned}\tag{4.2.16}$$

является следствием равенства (4.2.15). Следовательно, теорема верна для $k$. □

Теорема 4.2.8. *Производная* $\frac{d^mp_n(x)}{dx^m}\circ(h_1;...;h_m)$ *является симметричным многочленом по переменным* $h_1$, ..., $h_m$.

Доказательство. Не нарушая общности, мы можем полагать, что $p_n$ является одночленом. Мы докажем теорему индукцией по $m$.

Теорема очевидна для $m=0$ и $m=1$.

Лемма 4.2.9. *Производная* $\dfrac{d^2p_n(x)}{dx^2}$ *является симметричным билинейным отображением.*

Доказательство. Рассмотрим одночлен $p_n(x)$ в виде

$$p_n(x)=a_0xa_1...a_{n_1}xa_n \tag{4.2.17}$$

Согласно теореме 4.2.4 и лемме 4.2.2, равенство

$$\begin{aligned}\frac{dp_n(x)}{dx}\circ c_1=&a_0c_1a_1...a_{n-1}xa_n+...+a_0xa_1...a_{i-1}c_1a_i...a_{n-1}xa_n+...\\&+a_0xa_1...a_{n-1}c_1a_n\end{aligned} \tag{4.2.18}$$

является следствием равенства (4.2.17). Мы обозначим $p_{n1}(c_1,x)$, ..., $p_{ni}(c_1,x)$, ..., $p_{nn}(c_1,x)$ слагаемые в правой части равенства (4.2.18)

$$\frac{dp_n(x)}{dx}\circ c_1=p_{n1}(c_1,x)+...+p_{ni}(c_1,x)+...+p_{nn}(c_1,x) \tag{4.2.19}$$

$$p_{ni}(c_1,x)=a_0xa_1...a_{i-1}c_1a_i...a_{n-1}xa_n \tag{4.2.20}$$

Согласно определению [15]-4.1.1 и теореме A.1.4, равенство

$$\begin{aligned}\frac{d^2p_n(x)}{dx^2}\circ(c_1,c_2)=&\frac{d}{dx}\left(\frac{dp_n(x)}{dx}\circ c_1\right)\circ c_2\\=&\frac{d}{dx}(p_{n1}(c_1,x)+...+p_{ni}(c_1,x)+...+p_{nn}(c_1,x))\circ c_2\\=&\frac{dp_{n1}(c_1,x)}{dx}\circ c_2+...+\frac{dp_{ni}(c_1,x)}{dx}\circ c_2+...+\frac{dp_{nn}(c_1,x)}{dx}\circ c_2\end{aligned} \tag{4.2.21}$$

является следствием равенства (4.2.19). В правой части равенства (4.2.21) рассмотрим слагаемое с номером $i$, $1\le i\le n$. Согласно теореме 4.2.4 и лемме 4.2.2, равенство[4.2]

$$\begin{aligned}\frac{dp_{ni}(c_1,x)}{dx}\circ c_2=&a_0c_2a_1...a_{i-1}c_1a_i...a_{n-1}xa_n+...\\&+a_0xa_1...a_{i-1}c_1a_i......a_{j-1}c_2a_j...a_{n-1}xa_n+...\\&+a_0xa_1...a_{i-1}c_1a_i...a_{n-1}c_2a_n\end{aligned} \tag{4.2.22}$$

является следствием равенства (4.2.20). Мы обозначим $p_{ni1}(c_1,c_2,x)$, ..., $p_{ni\,i-1}(c_1,c_2,x)$, $p_{ni\,i+1}(c_1,c_2,x)$, ..., $p_{nij}(c_1,c_2,x)$, ..., $p_{nn}(c_1,c_2,x)$ слагаемые в правой части равенства (4.2.22)

$$\frac{dp_{ni}(c_1,x)}{dx}\circ c_2=p_{ni1}(c_1,c_2,x)+...+p_{nij}(c_1,c_2,x)+...+p_{nin}(c_1,c_2,x) \tag{4.2.23}$$

$$p_{nij}(c_1,c_2,x)=a_0xa_1...a_{i-1}c_1a_i......a_{j-1}c_2a_j...a_{n-1}xa_n \tag{4.2.24}$$

[4.2] В равенстве (4.2.22) приняты некоторые соглашения.

- Если $i=1$, то нет первого слагаемого.
- Если $i=n$, то нет последнего слагаемого.
- Предполагается, что $i<j$. Если $j<i$, то порядок множителей будет другой.

Из равенства

$$p_{nji}(c_1,c_2,x)=a_0xa_1...a_{i-1}c_2a_i......a_{j-1}c_1a_j...a_{n-1}xa_n \tag{4.2.25}$$

и равенства (4.2.24) следует, что

$$\frac{d^2p_n(x)}{dx^2}\circ(c_1,c_2)=\frac{d^2p_n(x)}{dx^2}\circ(c_2,c_1) \tag{4.2.26}$$

Лемма является следствием равенства (4.2.26). ⊙

Лемма 4.2.10. *Для любого одночлена $p_n(x)$ и любого $m>3$*

$$\begin{aligned}&\frac{d^mp_n(x)}{dx^m}\circ(c_1,...,c_{i-1},c_i,c_{i+1},...,c_m)\\ =&\frac{d^mp_n(x)}{dx^m}\circ(c_1,...,c_{i-1}c_{i+1},c_i,,...,c_m)\end{aligned} \tag{4.2.27}$$

Доказательство. Согласно теореме A.1.13

$$\begin{aligned}&\frac{d^mp_n(x)}{dx^m}\circ(c_1,...,c_m)\\ =&\frac{d^{m-i-1}p_n(x)}{dx^{m-i-1}}\left(\frac{d^2}{dx^2}\left(\frac{d^{i-1}p_n(x)}{dx^{i-1}}\circ(c_1,...,c_{i-1})\right)\circ(c_i,c_{i+1})\right)\\ &\circ(c_{i+2},...,c_m)\end{aligned} \tag{4.2.28}$$

Выражение

$$p(x)=\frac{d^{i-1}p_n(x)}{dx^{i-1}}\circ(c_1,...,c_{i-1}) \tag{4.2.29}$$

является многочленом и согласно лемме 4.2.9

$$\frac{d^2p(x)}{dx^2}\circ(c_i,c_{i+1})=\frac{d^2p(x)}{dx^2}\circ(c_{i+1},c_i) \tag{4.2.30}$$

The lemma follows from equalities (4.2.28), (4.2.29), (4.2.30). ⊙

Теорема является следствием леммы 4.2.10. □

Теорема 4.2.11. *Пусть $D$ - полное коммутативное кольцо характеристики 0. Пусть $A$ - ассоциативная банахова $D$-алгебра. Тогда*

$$\frac{dx^{n+1}}{dx}=\sum_{i=0}^n x^i\otimes x^{n-i} \tag{4.2.31}$$

$$dx^{n+1}=\sum_{i=0}^n x^idxx^{n-i} \tag{4.2.32}$$

Доказательство. Поскольку $x^{n+1}=(1\otimes...\otimes 1)\circ x$, то теорема является следствием теоремы 4.2.4. □

## 4.3. Матрица тензоров

Пусть $D$ - коммутативное кольцо и $A$ - ассоциативная $D$-алгебра с единицей.

Определение 4.3.1. *Пусть $a$ - матрица и $a^i_j\in A^{n\otimes}$. Матрица $a$ называется матрицей тензоров $A^{n\otimes}$.* □

Определение 4.3.2. *На множестве матриц тензоров определены следующие операции.*

- *Пусть* $b\in A^{m\otimes}$.

$$b*\begin{pmatrix}a_1^1&...&a_n^1\\...&...&...\\a_1^m&...&a_n^m\end{pmatrix}=\begin{pmatrix}b*a_1^1&...&b*a_n^1\\...&...&...\\b*a_1^m&...&b*a_n^m\end{pmatrix}$$

$$\begin{pmatrix}a_1^1&...&a_n^1\\...&...&...\\a_1^m&...&a_n^m\end{pmatrix}*b=\begin{pmatrix}a_1^1*b&...&a_n^1*b\\...&...&...\\a_1^m*b&...&a_n^m*b\end{pmatrix}$$

- 

$$\begin{pmatrix}a_1^1&...&a_n^1\\...&...&...\\a_1^m&...&a_n^m\end{pmatrix}{}_*{}^*\begin{pmatrix}b_1^1&...&b_k^1\\...&...&...\\b_1^n&...&b_k^n\end{pmatrix}=\begin{pmatrix}a_i^1*b_1^i&...&a_i^1*b_k^i\\...&...&...\\a_i^m*b_1^i&...&a_i^m*b_k^i\end{pmatrix}$$

- 

$$\begin{pmatrix}a_1^1&...&a_n^1\\...&...&...\\a_1^m&...&a_n^m\end{pmatrix}{}^*{}_*\begin{pmatrix}b_1^1&...&b_m^1\\...&...&...\\b_1^k&...&b_m^k\end{pmatrix}=\begin{pmatrix}a_1^i*b_i^1&...&a_n^i*b_i^k\\...&...&...\\a_1^i*b_i^k&...&a_n^i*b_i^k\end{pmatrix}$$

- *Пусть* $a_j^i\in A^{p\otimes}$, $x\in A$.

$$\begin{pmatrix}a_1^1&...&a_n^1\\...&...&...\\a_1^m&...&a_n^m\end{pmatrix}\circ x^p=\begin{pmatrix}a_1^1\circ x^p&...&a_n^1\circ x^p\\...&...&...\\a_1^m\circ x^p&...&a_n^m\circ x^p\end{pmatrix}\tag{4.3.1}$$

- *Пусть* $a_j^i\in A^{p\otimes}$, $x_1$, ...., $x_p\in A$.

$$\begin{pmatrix}a_1^1&...&a_n^1\\...&...&...\\a_1^m&...&a_n^m\end{pmatrix}\circ(x_1,...,x_p)=\begin{pmatrix}a_1^1\circ(x_1,...,x_p)&...&a_n^1\circ(x_1,...,x_p)\\...&...&...\\a_1^m\circ(x_1,...,x_p)&...&a_n^m\circ(x_1,...,x_p)\end{pmatrix}$$

- *Пусть* $a_j^i\in A^{2\otimes}$, $b_j^i\in A^{2\otimes}$,

$$\begin{pmatrix}a_1^1&...&a_n^1\\...&...&...\\a_1^m&...&a_n^m\end{pmatrix}{}_\circ{}^\circ\begin{pmatrix}b_1^1&...&b_k^1\\...&...&...\\b_1^n&...&b_k^n\end{pmatrix}=\begin{pmatrix}a_i^1\circ b_1^i&...&a_i^1\circ b_k^i\\...&...&...\\a_i^m\circ b_1^i&...&a_i^m\circ b_k^i\end{pmatrix}$$

- *Пусть* $a^{i}_{j} \in A^{2\otimes}$, $b^{i}_{j} \in A^{2\otimes}$,

$$\begin{pmatrix} a^{1}_{1} & ... & a^{1}_{n} \\ ... & ... & ... \\ a^{m}_{1} & ... & a^{m}_{n} \end{pmatrix} {}_{\circ}{}^{\circ} \begin{pmatrix} b^{1}_{1} & ... & b^{1}_{m} \\ ... & ... & ... \\ b^{k}_{1} & ... & b^{k}_{m} \end{pmatrix} = \begin{pmatrix} a^{i}_{1} \circ b^{1}_{i} & ... & a^{i}_{n} \circ b^{k}_{i} \\ ... & ... & ... \\ a^{i}_{1} \circ b^{k}_{i} & ... & a^{i}_{n} \circ b^{k}_{i} \end{pmatrix}$$

□

Важно обратить внимание на различие между ${}_{*}{}^{*}$-произведением и ${}_{\circ}{}^{\circ}$-произведением

$$\begin{pmatrix} a_1{}^{1}_{1} \otimes a_2{}^{1}_{1} & ... & a_1{}^{1}_{n} \otimes a_2{}^{1}_{n} \\ ... & ... & ... \\ a_1{}^{m}_{1} \otimes a_2{}^{m}_{1} & ... & a_1{}^{m}_{n} \otimes a_2{}^{m}_{n} \end{pmatrix} {}_{*}{}^{*} \begin{pmatrix} b_1{}^{1}_{1} \otimes b_2{}^{1}_{1} & ... & b_1{}^{1}_{k} \otimes b_2{}^{1}_{k} \\ ... & ... & ... \\ b_1{}^{n}_{1} \otimes b_2{}^{n}_{1} & ... & b_1{}^{n}_{k} \otimes b_2{}^{n}_{k} \end{pmatrix}$$

$$= \begin{pmatrix} a_1{}^{1}_{i} \otimes a_2{}^{1}_{i} b_1{}^{i}_{1} \otimes b_2{}^{i}_{1} & ... & a_1{}^{1}_{i} \otimes a_2{}^{1}_{i} b_1{}^{i}_{k} \otimes b_2{}^{i}_{k} \\ ... & ... & ... \\ a_1{}^{m}_{i} \otimes a_2{}^{m}_{i} b_1{}^{i}_{1} \otimes b_2{}^{i}_{1} & ... & a_1{}^{m}_{i} \otimes a_2{}^{m}_{i} b_1{}^{i}_{k} \otimes b_2{}^{i}_{k} \end{pmatrix}$$

$$\begin{pmatrix} a_1{}^{1}_{1} \otimes a_2{}^{1}_{1} & ... & a_1{}^{1}_{n} \otimes a_2{}^{1}_{n} \\ ... & ... & ... \\ a_1{}^{m}_{1} \otimes a_2{}^{m}_{1} & ... & a_1{}^{m}_{n} \otimes a_2{}^{m}_{n} \end{pmatrix} {}_{\circ}{}^{\circ} \begin{pmatrix} b_1{}^{1}_{1} \otimes b_2{}^{1}_{1} & ... & b_1{}^{1}_{k} \otimes b_2{}^{1}_{k} \\ ... & ... & ... \\ b_1{}^{n}_{1} \otimes b_2{}^{n}_{1} & ... & b_1{}^{n}_{k} \otimes b_2{}^{n}_{k} \end{pmatrix}$$

$$= \begin{pmatrix} a_1{}^{1}_{i} b_1{}^{i}_{1} \otimes b_2{}^{i}_{1} a_2{}^{1}_{i} & ... & a_1{}^{1}_{i} b_1{}^{i}_{k} \otimes b_2{}^{i}_{k} a_2{}^{1}_{i} \\ ... & ... & ... \\ a_1{}^{m}_{i} b_1{}^{i}_{1} \otimes b_2{}^{i}_{1} a_2{}^{m}_{i} & ... & a_1{}^{m}_{i} b_1{}^{i}_{k} \otimes b_2{}^{i}_{k} a_2{}^{m}_{i} \end{pmatrix}$$

Теорема 4.3.3. *Пусть $A$ - $D$-алгебра. Пусть $a$ - матрица тензоров $A^{p\otimes}$. Пусть $b$ - матрица тензоров $A^{q\otimes}$. Тогда*

(4.3.2) $$(a_{*}{}^{*}b) \circ x^{p+q} = (a \circ x^{p})_{*}{}^{*}(b \circ x^{q})$$

Доказательство. Равенство

(4.3.3) $$((a_{*}{}^{*}b) \circ x^{p+q})^{i}_{j} = (a_{*}{}^{*}b)^{i}_{j} \circ x^{p+q} = (a^{i}_{k} b^{k}_{j}) \circ x^{p+q}$$

является следствием равенств (3.2.3), (4.3.1). Равенство

(4.3.4) $$((a_{*}{}^{*}b) \circ x^{p+q})^{i}_{j} = (a^{i}_{k} \circ x^{p})(b^{k}_{j} \circ x^{q}) = (a \circ x^{p})^{i}_{k}(b \circ x^{q})^{k}_{j}$$

является следствием равенств (4.1.2), (4.3.1), (4.3.3). Равенство (4.3.2) является следствием равенств (3.2.3), (4.3.4). □

Теорема 4.3.4. *Пусть $A$ - $D$-алгебра и $a$ - матрица тензоров $A^{2\otimes}$. Тогда*

(4.3.5) $$a^{(k+1)_{*}{}^{*}} \circ x^{k+1} = (a^{k_{*}{}^{*}} \circ x^{k})_{*}{}^{*}(a \circ x)$$

Доказательство. Мы докажем теорему индукцией по $k$.

Равенство (4.3.5) верно, если $k = 1$, так как

$$a^{0_{*}{}^{*}} \circ x^{0} = E$$

Пусть теорема верна при $k=l$

$$a^{(l+1)_*{}^*}\circ x^{l+1}=(a^{l_*{}^*}\circ x^l)_*{}^*(a\circ x) \tag{4.3.6}$$

Равенство

$$a^{(l+2)_*{}^*}\circ x^{l+2}=(a^{(l+1)_*{}^*}{}_*{}^*a)\circ x^{l+2}=(a^{(l+1)_*{}^*}\circ x^{l+1})_*{}^*(a\circ x)$$

является следствием равенства (4.3.2). Следовательно, теорема верна при $k=l+1$. □

Теорема 4.3.5. *Пусть $A$ - $D$-алгебра и $a$ - матрица тензоров $A^{2\otimes}$. Тогда*

$$a^{n_*{}^*}\circ x^n=(a\circ x)^{n_*{}^*} \tag{4.3.7}$$

Доказательство. Мы докажем теорему индукцией по $n$.

Равенство (4.3.7) верно, если $n=0$ или $n=1$.

Пусть теорема верна при $n=k$

$$a^{k_*{}^*}\circ x^k=(a\circ x)^{k_*{}^*} \tag{4.3.8}$$

Согласно теореме 4.3.4, равенство

$$a^{(k+1)_*{}^*}\circ x^{k+1}=(a^k\circ x^k)_*{}^*(a\circ x)=(a\circ x)^{k_*{}^*}{}_*{}^*(a\circ x) \tag{4.3.9}$$

является следствием равенства (4.3.8). Равенство

$$a^{(k+1)_*{}^*}\circ x^{k+1}=(a\circ x)^{(k+1)_*{}^*} \tag{4.3.10}$$

является следствием равенств (3.2.16), (4.3.9). Следовательно, теорема верна при $n=k+1$. □

Глава 5

# Векторное пространство над алгеброй с делением

## 5.1. Векторное пространство

Пусть $D$ - коммутативное ассоциативное кольцо с единицей. Пусть $A$ - ассоциативная $D$-алгебра.

Теорема 5.1.1. *Пусть $D$-алгебра $A$ является алгеброй с делением. $D$-алгебра $A$ имеет единицу.*

Доказательство. Теорема является следствием утверждения, что уравнение

$$ax = a$$

имеет решение для любого $a \in A$. □

Теорема 5.1.2. *Пусть $D$-алгебра $A$ является алгеброй с делением. Кольцо $D$ является полем и подмножеством центра $D$-алгебры $A$.*

Доказательство. Пусть $e$ - единица $D$-алгебры $A$. Так как отображение

$$d \in D \rightarrow d\,e \in A$$

является вложением кольца $D$ в $D$-алгебру $A$, то кольцо $D$ является подмножеством центра $D$-алгебры $A$. □

Определение 5.1.3. *Пусть $A$ - алгебра с делением. Эффективное левостороннее представление*

$$f : A \overset{*}{\longrightarrow} V \qquad f(a) : v \in V \rightarrow av \in V \quad a \in A$$

*абелевой группы $A$ в $D$-модуле $V$ называется* **левым векторным пространством** *над $D$-алгеброй $A$. Мы также будем говорить, что $D$-модуль $V$ является* **левым $A$-векторным пространством** *или* **$A*$-векторным пространством**. *$V$-число называется* **вектором**. □

Определение 5.1.4. *Пусть $A$ - алгебра с делением. Эффективное правостороннее представление*

$$f : A \overset{*}{\longrightarrow} V \qquad f(a) : v \in V \rightarrow va \in V \quad a \in A$$

*абелевой группы $A$ в $D$-модуле $V$ называется* **правым векторным пространством** *над $D$-алгеброй $A$. Мы также будем говорить, что $D$-модуль $V$ является* **правым $A$-векторным пространством** *или* **$*A$-векторным пространством**. *$V$-число называется* **вектором**. □

Так как любое утверждение о левом $A$-векторном пространстве эквивалентно утверждению о правом $A$-векторном пространстве с точностью до порядка множителей, мы в этой главе будем рассматривать левое $A$-векторное пространство. Если $D$-алгебра $A$ коммутативна, то левое $A$-векторное пространство совпадает с правым $A$-векторным пространством.

Теорема 5.1.5. *Следующая диаграмма представлений описывает левое $A$-векторное пространство $V$*

(5.1.1)

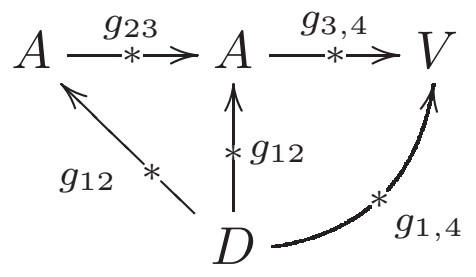

$$
\begin{array}{l}
g_{12}(d) : a \to d\,a \\
g_{23}(v) : w \to C(w,v) \\
\qquad C \in \mathcal{L}(A^2 \to A) \\
g_{3,4}(a) : v \to v\,a \\
g_{1,4}(d) : v \to d\,v
\end{array}
$$

*В диаграмме представлений* (5.1.1) *верна* **коммутативность представлений** *коммутативного кольца $D$ и $D$-алгебры $A$ в абелевой группе $V$*

$$a(dv) = d(av) \tag{5.1.2}$$

Доказательство. Диаграмма представлений (5.1.1) является следствием определения 5.1.3, теоремы [13]-5.1.11 и следствия [13]-5.1.12. Равенство (5.1.2) является следствием утверждения, что левостороннее преобразование $g_{3,4}(a)$ является эндоморфизмом $D$-векторного пространства $V$. □

Теорема 5.1.6. *Пусть $V$ является левым $A$-векторным пространством. Элементы левого $A$-векторного пространства $V$ удовлетворяют соотношениям*

- *если $D$-алгебра $A$ ассоциативна, то верен* **закон ассоциативности**

$$(pq)v = p(qv) \tag{5.1.3}$$

- **закон дистрибутивности**

$$p(v+w) = pv + pw \tag{5.1.4}$$

$$(p+q)v = pv + qv \tag{5.1.5}$$

- **закон унитарности**

$$1v = v \tag{5.1.6}$$

*для любых $p, q \in A$, $v, w \in V$.*

Доказательство. Равенство (5.1.4) следует из утверждения, что левостороннее преобразование $p$ является эндоморфизмом абелевой группы $V$. Равенство (5.1.5) следует из утверждения, что представление $g_{3,4}$ является гомоморфизмом аддитивной группы $D$-алгебры $A$. Равенство (5.1.6) следует из утверждения, что представление $g_{3,4}$ является левосторонним представлением мультипликативной группы $D$-алгебры $A$. Равенство (5.1.3) следует из утверждения, что представление $g_{3,4}$ является левосторонним представлением мультипликативной группы $D$-алгебры $A$, если $D$-алгебра $A$ ассоциативна. □

## 5.2. Базис левого $A$-векторного пространства

Теорема 5.2.1. *Пусть $V$ - левое $A$-векторное пространство. Множество векторов, порождённое множеством векторов $v=(v_i\in V, i\in I)$ имеет вид*[5.1]

$$J(v)=\left\{w:w=\sum_{i\in I}c^iv_i, c^i\in A, |\{i:c^i\neq 0\}|<\infty\right\}\tag{5.2.1}$$

Доказательство. Мы докажем теорему по индукции, опираясь на теорему [12]-3.4, Согласно теореме [12]-3.4, мы должны доказать следующие утверждения:

5.2.1.1: $v_k\in X_0\subseteq J(v)$
5.2.1.2: $c^kv_k\in J(v)$, $c^k\in A$, $k\in I$
5.2.1.3: $\sum_{k\in I}c^kv_k\in J(v)$, $c^k\in A$, $|\{i:c^i\neq 0\}|<\infty$
5.2.1.4: $v,\ w\in J(v)\Rightarrow v+w\in J(v)$
5.2.1.5: $a\in A,\ v\in J(v)\Rightarrow av\in J(v)$

- Для произвольного $v_k\in v$, положим $c^i=\delta^i_k\in A$. Тогда
$$v_k=\sum_{i\in I}c^iv_i\tag{5.2.2}$$
Утверждение 5.2.1.1 является следствием равенств (5.2.1), (5.2.2).
- Утверждение 5.2.1.2 являются следствием теорем [12]-3.4, 5.1.6 и утверждения 5.2.1.1.
- Так как $V$ является абелевой группой, то утверждение 5.2.1.3 следует из утверждения 5.2.1.2 и теоремы [12]-3.4.
- Пусть $w_1,\ w_2\in X_k\subseteq J(v)$. Так как $V$ является абелевой группой, то, согласно утверждению [12]-3.4.3,
$$w_1+w_2\in X_{k+1}\tag{5.2.3}$$
Согласно равенству (5.2.1), существуют $A$-числа $w^i_1,\ w^i_2,\ i\in I$, такие, что
$$w_1=\sum_{i\in I}w^i_1v_i\qquad w_2=\sum_{i\in I}w^i_2v_i\tag{5.2.4}$$
где множества
$$H_1=\{i\in I:w^i_1\neq 0\}\qquad H_2=\{i\in I:w^i_2\neq 0\}\tag{5.2.5}$$
конечны. Так как $V$ является абелевой группой, то из равенства (5.2.4) следует, что
$$w_1+w_2=\sum_{i\in I}w^i_1v_i+\sum_{i\in I}w^i_2v_i=\sum_{i\in I}(w^i_1v_i+w^i_2v_i)\tag{5.2.6}$$
Равенство
$$w_1+w_2=\sum_{i\in I}(w^i_1+w^i_2)v_i\tag{5.2.7}$$

[5.1] Для множества $A$, мы обозначим $|A|$ мощность множества $A$. Запись $|A|<\infty$ означает, что множество $A$ конечно.

является следствием равенств (5.1.5), (5.2.6). Из равенства (5.2.5) следует, что множество

$$\{i \in I : w_1^i + w_2^i \neq 0\} \subseteq H_1 \cup H_2$$

конечно.

- Пусть $w \in X_k \subseteq J(v)$. Согласно утверждению [12]-3.4.4, для любого $A$-числа $a$,

(5.2.8) $$aw \in X_{k+1}$$

Согласно равенству (5.2.1), существуют $A$-числа $w^i$, $i \in I$, такие, что

(5.2.9) $$w = \sum_{i \in I} w^i v_i$$

где

(5.2.10) $$|\{i \in I : w^i \neq 0\}| < \infty$$

Из равенства (5.2.9) следует, что

(5.2.11) $$aw = a \sum_{i \in I} w^i v_i = \sum_{i \in I} a(w^i v_i) = \sum_{i \in I} (aw^i) v_i$$

Из утверждения (5.2.10) следует, что множество $\{i \in I : aw^i \neq 0\}$ конечно.

Из равенств (5.2.3), (5.2.7), (5.2.8), (5.2.11) следует, что $X_{k+1} \subseteq J(v)$. □

Определение 5.2.2. *$J(v)$ называется векторным подпространством, порождённым множеством $v$, а $v$ -* **множеством образующих** *векторного подпространства $J(v)$. В частности,* **множеством образующих** *левого $A$-векторного пространства $V$ будет такое подмножество $X \subset V$, что $J(X) = V$.* □

Определение 5.2.3. *Если множество $X \subset V$ является множеством образующих левого $A$-векторного пространства $V$, то любое множество $Y$, $X \subset Y \subset V$ также является множеством образующих левого $A$-векторного пространства $V$. Если существует минимальное множество $X$, порождающее левый $A$-векторное пространство $V$, то такое множество $X$ называется* **базисом** *левого $A$-векторного пространства $V$.* □

Определение 5.2.4. *Пусть $v = (v_i \in V, i \in I)$ - множество векторов. Выражение $w^i v_i$ называется* **линейной комбинацией** *векторов $v_i$. Вектор $\overline{w} = w^i v_i$ называется* **линейно зависимым** *от векторов $v_i$.* □

Теорема 5.2.5. *Пусть $A$ - ассоциативная $D$-алгебра с делением. Если уравнение*

$$w^i v_i = 0$$

*предполагает существования индекса $i = j$ такого, что $w^j \neq 0$, то вектор $v_j$ линейно зависит от остальных векторов $v$.*

Доказательство. Теорема является следствием равенства

$$v_j = \sum_{i \in I \setminus \{j\}} (w^j)^{-1} w^i v_i$$

и определения 5.2.4. □

Очевидно, что для любого множества векторов $v_i$,

$$w^i = 0 \Rightarrow w^*{}_*v = 0$$

Определение 5.2.6. *Множество векторов*[5.2] $v_i$, $i \in I$, *левого $A$-векторного пространства $V$* **линейно независимо**, *если $w = 0$ следует из уравнения*

$$w^i v_i = 0$$

*В противном случае, множество векторов $v_i$, $i \in I$,* **линейно зависимо**. □

Теорема 5.2.7. *Пусть $A$ - ассоциативная $D$-алгебра с делением. Множество векторов $\overline{\overline{e}} = (e_i, i \in I)$ является* **базисом левого $A$-векторного пространства** $V$, *если векторы $e_i$ линейно независимы и любой вектор $v \in V$ линейно зависит от векторов $e_i$.*

Доказательство. Пусть множество векторов $e_i$, $i \in I$, линейно зависимо. Тогда в равенстве

$$w^i e_i = 0$$

существует индекс $i = j$ такой, что $w^j \neq 0$. Согласно теореме 5.2.5, вектор $e_j$ линейно зависит от остальных векторов множества $\overline{\overline{e}}$. Согласно определению 5.2.3, множество векторов $e_i$, $i \in I$, не является базисом левого $A$-векторного пространства $V$.

Следовательно, если множество векторов $e_i$, $i \in I$, является базисом, то эти векторы линейно независимы. Так как произвольный вектор $v \in V$ является линейной комбинацией векторов $e_i$, $i \in I$, , то множество векторов $v$, $e_i$, $i \in I$, не является линейно независимым. □

Определение 5.2.8. *Пусть $\overline{\overline{e}}$ - базис левого $A$-векторного пространства $V$, и вектор $\overline{v} \in V$ имеет разложение*

$$\overline{v} = v^*{}_*e = v^i e_i$$

*относительно базиса $\overline{\overline{e}}$. $A$-числа $v^i$ называются* **координатами** *вектора $\overline{v}$ относительно базиса $\overline{\overline{e}}$. Матрица $A$-чисел $v = (v^i, i \in I)$ называется* **координатной матрицей вектора** $\overline{v}$ *в базисе $\overline{\overline{e}}$.* □

Теорема 5.2.9. *Координаты вектора $v \in V$ относительно базиса $\overline{\overline{e}}$ левого $A$-векторного пространства $V$ определены однозначно.*

Доказательство. Согласно теореме 5.2.7, система векторов $\overline{v}$, $e_i$, $i \in I$, линейно зависима и в уравнении

$$b\overline{v} + c^*{}_*e = 0 \tag{5.2.12}$$

по крайней мере $b$ отлично от 0. Тогда равенство

$$\overline{v} = (-cb^{-1})^*{}_*e \tag{5.2.13}$$

следует из (5.2.12). Равенство

$$\overline{v} = v^*{}_*e \tag{5.2.14}$$

является следствием из равенства (5.2.13).

[5.2] Я следую определению в [1], страница 100.

Допустим мы имеем другое разложение

$$\overline{v} = v'{}_*{}^*e \tag{5.2.15}$$

Вычтя (5.2.14) из (5.2.15), мы получим

$$0 = (v' - v)^*{}_*e \tag{5.2.16}$$

Так как векторы $e_i$ линейно независимы, то равенство

$$v' - v = 0 \tag{5.2.17}$$

является следствием равенства (5.2.16). Утверждение теоремы является следствием равенства (5.2.17). □

Теорема 5.2.10. *Мы можем рассматривать $D$-алгебру $A$ с делением как левое $A$-векторное пространство с базисом $\overline{\overline{e}} = \{1\}$.*

Доказательство. Теорема является следствием равенства (5.1.6). □

## 5.3. Тип векторного пространства

Если базис $A$-векторного пространства $V$ конечен, то мощность базиса называется размерностью $A$-векторного пространства либо мы говорим что $A$-векторное пространство $V$ конечномерно. Согласно определению 5.2.8 теория конечномерного $A$-векторного пространства $V$ связана с теорией матриц.

Пример 5.3.1. *В разделе 5.2 мы представили множество векторов $v_1$, ..., $v_n$ в виде строки матрицы*

$$v = \begin{pmatrix} v_1 & \dots & v_n \end{pmatrix}$$

*и множество $A$-чисел $c^1$, ..., $c^n$ в виде столбца матрицы*

$$c = \begin{pmatrix} c^1 \\ \dots \\ c^n \end{pmatrix}$$

*Поэтому мы можем записать линейную комбинацию векторов $v_1$, ..., $v_n$ в виде ${}^*{}_*$-произведения матриц*

$$c^*{}_*v = c^i v_i = \begin{pmatrix} c^1 \\ \dots \\ c^n \end{pmatrix} {}^*{}_* \begin{pmatrix} v_1 & \dots & v_n \end{pmatrix} \tag{5.3.1}$$

*Выражение (5.3.1) для линейной комбинации векторов называется* **левой линейной комбинацией столбцов**. *В частности, если векторы базиса $\overline{\overline{e}}$ записать в виде строки матрицы*

$$e = \begin{pmatrix} e_1 & \dots & e_n \end{pmatrix}$$

*и координаты вектора* $\overline{w}$ *относительно базиса* $\overline{\overline{e}}$ *записать в виде столбца матрицы*

$$w = \begin{pmatrix} w^1 \\ ... \\ w^n \end{pmatrix}$$

*то мы можем представить вектор* $\overline{w}$ *в виде* ${}^*{}_*$*-произведения матриц*

$$\overline{w} = w^*{}_*e = w^i e_i = \begin{pmatrix} w^1 \\ ... \\ w^n \end{pmatrix} {}^*{}_* \begin{pmatrix} e_1 & ... & e_n \end{pmatrix}$$

*Соответствующее представление векторного пространства* $V$ *мы будем называть* **левым** $A^*$**-векторным пространством** *или* **левым пространством** $A$**-столбцов**, *а* $V$*-число мы будем называть* **вектор-столбец**. □

Формат представления векторов определяется рассматриваемой задачей.

Пример 5.3.2. *Мы представим множество векторов* $v_1, ..., v_n$ *в виде столбца матрицы*

$$v = \begin{pmatrix} v^1 \\ ... \\ v^n \end{pmatrix}$$

*и множество* $A$*-чисел* $c^1, ..., c^n$ *в виде строки матрицы*

$$c = \begin{pmatrix} c_1 & ... & c_n \end{pmatrix}$$

*Поэтому мы можем записать линейную комбинацию векторов* $v_1, ..., v_n$ *в виде* ${}_*{}^*$*-произведения матриц*

$$c_*{}^*v = c_i v^i = \begin{pmatrix} c_1 & ... & c_n \end{pmatrix} {}_*{}^* \begin{pmatrix} v^1 \\ ... \\ v^n \end{pmatrix} \tag{5.3.2}$$

*Выражение* (5.3.2) *для линейной комбинации векторов называется* **левой линейной комбинацией строк**. *В частности, если векторы базиса* $\overline{\overline{e}}$ *записать в виде столбца матрицы*

$$e = \begin{pmatrix} e^1 \\ ... \\ e^n \end{pmatrix}$$

*и координаты вектора* $\overline{w}$ *относительно базиса* $\overline{\overline{e}}$ *записать в виде строки матрицы*

$$w = \begin{pmatrix} w_1 & ... & w_n \end{pmatrix}$$

*то мы можем представить вектор* $\overline{w}$ *в виде* ${}_*{}^*$*-произведения матриц*

$$\overline{w}=w_*{}^*e=w_ie^i=\begin{pmatrix}w_1 & ... & w_n\end{pmatrix}{}_*{}^*\begin{pmatrix}e^1\\ ...\\ e^n\end{pmatrix}$$

*Соответствующее представление векторного пространства* $V$ *мы будем называть* **левым** $A_*$**-векторным пространством** *или* **левым пространством** $A$**-строк**, *а* $V$*-число мы будем называть* **вектор-строка**. □

Аналогичные типы векторного пространства мы можем рассмотреть для правого векторного пространства.

Пример 5.3.3. *Мы представим множество векторов* $v_1$, ..., $v_n$ *в виде строки матрицы*

$$v=\begin{pmatrix}v_1 & ... & v_n\end{pmatrix}$$

*и множество* $A$*-чисел* $c^1$, ..., $c^n$ *в виде столбца матрицы*

$$c=\begin{pmatrix}c^1\\ ...\\ c^n\end{pmatrix}$$

*Поэтому мы можем записать линейную комбинацию векторов* $v_1$, ..., $v_n$ *в виде* ${}_*{}^*$*-произведения матриц*

$$v_*{}^*c=v_ic^i=\begin{pmatrix}v_1 & ... & v_n\end{pmatrix}{}_*{}^*\begin{pmatrix}c^1\\ ...\\ c^n\end{pmatrix}\tag{5.3.3}$$

*Выражение* (5.3.3) *для линейной комбинации векторов называется* **правой линейной комбинацией столбцов**. *В частности, если векторы базиса* $\overline{\overline{e}}$ *записать в виде строки матрицы*

$$e=\begin{pmatrix}e_1 & ... & e_n\end{pmatrix}$$

*и координаты вектора* $\overline{w}$ *относительно базиса* $\overline{\overline{e}}$ *записать в виде столбца матрицы*

$$w=\begin{pmatrix}w^1\\ ...\\ w^n\end{pmatrix}$$

*то мы можем представить вектор* $\overline{w}$ *в виде* ${}_*{}^*$*-произведения матриц*

$$\overline{w}=e_*{}^*w=e_iw^i=\begin{pmatrix}e_1 & ... & e_n\end{pmatrix}{}_*{}^*\begin{pmatrix}w^1\\ ...\\ w^n\end{pmatrix}$$

*Соответствующее представление векторного пространства $V$ мы будем называть* **правым $A^*$-векторным пространством** *или* **правым пространством $A$-столбцов**, *а $V$-число мы будем называть* **вектор-столбец**. □

Пример 5.3.4. *Мы представим множество векторов $v_1, ..., v_n$ в виде столбца матрицы*

$$v=\begin{pmatrix}v^1\\ ...\\ v^n\end{pmatrix}$$

*и множество $A$-чисел $c^1, ..., c^n$ в виде строки матрицы*

$$c=\begin{pmatrix}c_1 & ... & c_n\end{pmatrix}$$

*Поэтому мы можем записать линейную комбинацию векторов $v_1, ..., v_n$ в виде ${}^*{}_*$-произведения матриц*

$$v^*{}_*c=v^ic_i=\begin{pmatrix}v^1\\ ...\\ v^n\end{pmatrix}{}^*{}_*\begin{pmatrix}c_1 & ... & c_n\end{pmatrix}\tag{5.3.4}$$

*Выражение (5.3.4) для линейной комбинации векторов называется* **правой линейной комбинацией строк**. *В частности, если векторы базиса $\overline{\overline{e}}$ записать в виде столбца матрицы*

$$e=\begin{pmatrix}e^1\\ ...\\ e^n\end{pmatrix}$$

*и координаты вектора $\overline{w}$ относительно базиса $\overline{\overline{e}}$ записать в виде строки матрицы*

$$w=\begin{pmatrix}w_1 & ... & w_n\end{pmatrix}$$

*то мы можем представить вектор $\overline{w}$ в виде ${}^*{}_*$-произведения матриц*

$$\overline{w}=e^*{}_*w=e^iw_i=\begin{pmatrix}e^1\\ ...\\ e^n\end{pmatrix}{}^*{}_*\begin{pmatrix}w_1 & ... & w_n\end{pmatrix}$$

*Соответствующее представление векторного пространства $V$ мы будем называть* **правым $A_*$-векторным пространством** *или* **правым пространством $A$-строк**, *а $V$-число мы будем называть* **вектор-строка**. □

Выбор типа векторного пространства - это соглашение, необходимое для решения конкретной задачи.

## 5.4. Система линейных уравнений

Определение 5.4.1. *Пусть* $V$ *- правое* $A^*$*-векторное пространство и* $\{a_i \in V, i \in I\}$ *- множество векторов.* **Линейная оболочка** *в правом* $A^*$*-векторном пространстве - это множество* $span(\overline{a}_i, i \in I)$ *векторов, линейно зависимых от векторов* $a_i$*.* □

Теорема 5.4.2. *Пусть* $span(\overline{a}_i, i \in I)$ *- линейная оболочка в правом* $A_*$*-векторном пространстве* $V$*. Тогда* $span(\overline{a}_i, i \in I)$ *- подпространство правого* $A_*$*-векторного пространства* $V$*.*

Доказательство. Предположим, что

$$\overline{b} \in \mathrm{span}(\overline{a}_i, i \in I)$$

$$\overline{c} \in \mathrm{span}(\overline{a}_i, i \in I)$$

Согласно определению 5.4.1

$$\overline{b} = a_*{}^*b$$

$$\overline{c} = a_*{}^*c$$

Тогда

$$\overline{b} + \overline{c} = a_*{}^*b + a_*{}^*c = a_*{}^*(b + c) \in \mathrm{span}(\overline{a}_i, i \in I)$$

$$\overline{b}k = (a_*{}^*b)k = a_*{}^*(bk) \in \mathrm{span}(\overline{a}_i, i \in I)$$

Это доказывает утверждение. □

Пример 5.4.3. *Пусть* $V$ *- правое* $A_*$*-векторное пространство и строка*

$$\overline{a} = \begin{pmatrix} a_1 & \dots & a_n \end{pmatrix} = (a_i, i \in I)$$

*задаёт множество векторов. Чтобы ответить на вопрос, или вектор* $\overline{b} \in span(a_i, i \in I)$ *, мы запишем линейное уравнение*

$$\overline{b} = \overline{a}_*{}^*x \tag{5.4.1}$$

*где*

$$x = \begin{pmatrix} x^1 \\ \dots \\ x^n \end{pmatrix}$$

*столбец неизвестных коэффициентов разложения.* $\overline{b} \in span(a_i, i \in I)$ *, если уравнение (5.4.1) имеет решение. Предположим, что* $\overline{\overline{f}} = (f_j, j \in J)$ *- базис. Тогда векторы* $\overline{b}$*,* $\overline{a}_i$ *имеют разложение*

$$\overline{b} = f_*{}^*b \tag{5.4.2}$$

$$\overline{a}_i = f_*{}^*a_i \tag{5.4.3}$$

*Если мы подставим (5.4.2) и (5.4.3) в (5.4.1), мы получим*

$$f_*{}^*b = f_*{}^*a_*{}^*x \tag{5.4.4}$$

*Применяя теорему 5.2.9 к (5.4.4), мы получим* **систему линейных уравнений**

$$a_*{}^*x = b \tag{5.4.5}$$

*Мы можем записать систему линейных уравнений* (5.4.5) *в одной из следующих форм*

$$\begin{pmatrix} a_1^1 & ... & a_n^1 \\ ... & ... & ... \\ a_1^m & ... & a_n^m \end{pmatrix} {}_*{}^* \begin{pmatrix} x^1 \\ ... \\ x^n \end{pmatrix} = \begin{pmatrix} b^1 \\ ... \\ b^m \end{pmatrix}$$

(5.4.6) $$a_i^j x^i = b^j$$

$$\begin{matrix} a_1^1 x^1 & +... & +a_n^1 x^n & = b^1 \\ ... & ... & ... & ... \\ a_m^1 x^1 & +... & +a_n^m x^n & = b^m \end{matrix}$$

□

Пример 5.4.4. *Пусть $V$ - левое $A_*$-векторное пространство и столбец*

$$\overline{a} = \begin{pmatrix} \overline{a}^1 \\ ... \\ \overline{a}^m \end{pmatrix} = (\overline{a}^j, j \in J)$$

*задаёт множество векторов. Чтобы ответить на вопрос, или вектор $\overline{b} \in span(\overline{a}^j, j \in J)$, мы запишем линейное уравнение*

(5.4.7) $$\overline{b} = x_*{}^* \overline{a}$$

*где*

$$x = \begin{pmatrix} x_1 & ... & x_m \end{pmatrix}$$

*строка неизвестных коэффициентов разложения. $\overline{b} \in span(\overline{a}^j, j \in J)$, если уравнение* (5.4.7) *имеет решение. Предположим, что $\overline{\overline{f}} = (f^i, i \in I)$ - базис. Тогда векторы $\overline{b}$, $a^j$ имеют разложение*

(5.4.8) $$\overline{b} = b_*{}^* f$$

(5.4.9) $$\overline{a}^j = a^j{}_*{}^* f$$

*Если мы подставим* (5.4.8) *и* (5.4.9) *в* (5.4.7), *мы получим*

(5.4.10) $$b_*{}^* f = x_*{}^* a_*{}^* f$$

*Применяя теорему 5.2.9 к* (5.4.10), *мы получим* **системy линейных уравнений**[5.3]

(5.4.11) $$x_*{}^* a = b$$

[5.3] Читая систему линейных уравнений (5.4.5) в правом $A_*$-векторном пространстве снизу вверх и слева направо, мы получим систему линейных уравнений (5.4.11) в левом $A_*$-векторном пространстве.

*Мы можем записать систему линейных уравнений* (5.4.11) *в одной из следующих форм*

$$\begin{pmatrix} x_1 & ... & x_m \end{pmatrix} {}_*{}^* \begin{pmatrix} a_1^1 & ... & a_n^1 \\ ... & ... & ... \\ a_1^m & ... & a_n^m \end{pmatrix} = \begin{pmatrix} b_1 & ... & b_n \end{pmatrix}$$

$$x_j a_i^j = b_i \tag{5.4.12}$$

$$\begin{matrix} x_1 a_1^1 & ... & x_1 a_n^1 \\ + & ... & + \\ ... & ... & ... \\ + & ... & + \\ x_m a_1^m & ... & x_m a_n^m \\ = & ... & = \\ b_1 & ... & b_n \end{matrix}$$

□

Чтобы найти решение системы линейных уравнений, мы должны рассмотреть матрицу этой системы. Из примеров 5.3.1, 5.3.3, мы видим, что мы можем рассматривать столбец матрицы как вектор левого или правого $A$-векторного пространства. Чтобы сделать утверждения проще, мы будем различать формат линейной зависимости. Например, утверждение

- Столбцы матрицы $D^*{}_*$-линейно зависимы.

означает, что

- Столбцы матрицы являются векторами $D^*{}_*$-векторного пространства, и соответствующие векторы линейно зависимы.

В частности, система линейных уравнений (5.4.5) в ${}_*{}^*D$-векторном пространстве называется **системой ${}_*{}^*D$-линейных уравнений** и система линейных уравнений (5.4.11) в $D_*{}^*$-векторном пространстве называется **системой $D_*{}^*$-линейных уравнений**.

Определение 5.4.5. *Если $n \times n$ матрица $a$ имеет ${}_*{}^*$-обратную матрицу, мы будем называть такую матрицу* **${}_*{}^*$-невырожденной матрицей**. *В противном случае, мы будем называть такую матрицу* **${}_*{}^*$-вырожденной матрицей**. □

Определение 5.4.6. *Предположим, что $A$ - ${}_*{}^*$-невырожденная матрица. Мы будем называть соответствующую систему ${}_*{}^*D$-линейных уравнений*

$$a_*{}^* x = b \tag{5.4.13}$$

**невырожденной системой ${}_*{}^*D$-линейных уравнений**. □

Теорема 5.4.7. *Решение невырожденной системы ${}_*{}^*D$-линейных уравнений (5.4.13) определено однозначно и может быть записано в любой из следующих форм*[5.4]

$$x = A^{-1_*{}^*}{}_*{}^*b \tag{5.4.14}$$

$$x = \mathcal{H}\det({}_*{}^*)a_*{}^*b \tag{5.4.15}$$

Доказательство. Умножая обе части равенства (5.4.13) слева на $a^{-1_*{}^*}$, мы получим (5.4.14). Пользуясь определением (3.3.10), мы получим (5.4.15). Решение системы единственно в силу теоремы 3.2.20. □

## 5.5. Ранг матрицы

Мы делаем следующие предположения в этом разделе

- $i \in M$, $|M| = m$, $j \in N$, $|N| = n$.
- $a = (a^i_j)$ - произвольная матрица.
- $k, s \in S \supseteq M$, $l, t \in T \supseteq N$, $k = |S| = |T|$.
- $p \in M \setminus S$, $r \in N \setminus T$.

Определение 5.5.1. *Мы будем называть матрицу $a^S_T$ минорной матрицей порядка $k$.* □

Определение 5.5.2. *Если минорная матрица $a^S_T$ - ${}_*{}^*$-невырожденная матрица, то мы будем говорить, что ${}_*{}^*$-ранг матрицы $a$ не меньше, чем $k$.* **${}_*{}^*$-ранг матрицы** $a$

$$\mathrm{rank}_{{}_*{}^*}\, a$$

*это максимальное значение $k$. Мы будем называть соответствующую минорную матрицу* **${}_*{}^*$-главной минорной матрицей**. □

Теорема 5.5.3. *Пусть матрица $a$ - ${}_*{}^*$-вырожденная матрица и минорная матрица $a^S_T$ - ${}_*{}^*$-главная минорная матрица, тогда*

$$\det({}_*{}^*)^p_r\, a^{S\cup\{p\}}_{T\cup\{r\}} = 0$$

Доказательство. Чтобы понять, почему минорная матрица

$$b = a^{S\cup\{p\}}_{T\cup\{r\}} \tag{5.5.1}$$

не имеет ${}_*{}^*$-обратной матрицы,[5.5] мы предположим, что существует ${}_*{}^*$-обратная матрица $b^{-1_*{}^*}$. Запишем систему линейных уравнений (3.3.4), (3.3.5) полагая $I = \{r\}$, $J = \{p\}$ (в этом случае $[I] = T$, $[J] = S$ )

$$b^S_{T*}{}^*b^{-1_*{}^*}{}^T_p + b^S_r\, b^{-1_*{}^*}{}^r_p = 0 \tag{5.5.2}$$

$$b^p_{T*}{}^*b^{-1_*{}^*}{}^T_p + b^p_r\, b^{-1_*{}^*}{}^r_p = 1 \tag{5.5.3}$$

---

[5.4] Мы можем найти решение системы (5.4.13) в теореме [4]-1.6.1. Я повторяю это утверждение, так как я слегка изменил обозначения.

[5.5] Естественно ожидать связь между ${}_*{}^*$-вырожденностью матрицы и её ${}_*{}^*$-квазидетерминантом, подобную связи, известной в коммутативном случае. Однако ${}_*{}^*$-квазидетерминант определён не всегда. Например, если ${}_*{}^*$-обратная матрица имеет слишком много элементов, равных 0. Как следует из этой теоремы, не определён ${}_*{}^*$-квазидетерминант также в случае, когда ${}_*{}^*$-ранг матрицы меньше $n-1$.

и попробуем решить эту систему. Мы умножим (5.5.2) на $(b^S_T)^{-1_*{}^*}$

$$b^{-1_*{}^*}{}^T_p + (b^S_T)^{-1_*{}^*}{}_*{}^*b^S_r\ b^{-1_*{}^*}{}^r_p = 0 \tag{5.5.4}$$

Теперь мы можем подставить (5.5.4) в (5.5.3)

$$-b^p_T{}_*{}^*(b^S_T)^{-1_*{}^*}{}_*{}^*b^S_r\ b^{-1_*{}^*}{}^r_p + b^p_r\ b^{-1_*{}^*}{}^r_p = 1 \tag{5.5.5}$$

Из (5.5.5) следует

$$(b^p_r - b^p_T{}_*{}^*(b^S_T)^{-1_*{}^*}{}_*{}^*b^S_r)\ b^{-1_*{}^*}{}^r_p = 1 \tag{5.5.6}$$

Выражение в скобках является квазидетерминантом $\det({}_*{}^*)^p_r b$. Подставляя это выражение в (5.5.6), мы получим

$$\det({}_*{}^*)^p_r b\ b^{-1_*{}^*}{}^r_p = 1 \tag{5.5.7}$$

Тем самым мы доказали, что квазидетерминант $\det({}_*{}^*)^p_r b$ определён и условие его обращения в 0 необходимое и достаточное условие вырожденности матрицы $b$. Теорема доказана в силу соглашения (5.5.1). □

Теорема 5.5.4. *Предположим, что $a$ - матрица,*

$$\mathrm{rank}_{{}_*{}^*}\, a = \boldsymbol{k} < \boldsymbol{m}$$

*и $a^S_T$ - ${}_*{}^*$-главная минорная матрица. Тогда строка $a^p$ является левой линейной комбинацией строк $a^S$*

$$a^{M\setminus S} = r_*{}^*a^S \tag{5.5.8}$$

$$a^p = r^p{}_*{}^*a^S \tag{5.5.9}$$

$$a^p_b = r^p_s a^s_b \tag{5.5.10}$$

Доказательство. Если матрица $a$ имеет $k$ столбцов, то, полагая, что строка $a^p$ - левая линейная комбинация (5.5.9) строк $a^S$ с коэффициентами $R^p_s$, мы получим систему линейных уравнений (5.5.10). Согласно теореме 5.4.7 система линейных уравнений (5.5.10) имеет единственное решение[5.6] и это решение нетривиально потому, что все ${}_*{}^*$-квазидетерминанты отличны от 0.

Остаётся доказать утверждение в случае, когда число столбцов матрицы $a$ больше чем $k$. Пусть нам даны строка $a^p$ и столбец $a_r$. Согласно предположению, минорная матрица $a^{S\cup\{p\}}_{T\cup\{r\}}$ - ${}_*{}^*$-вырожденная матрица и его ${}_*{}^*$-квазидетерминант

$$\det({}_*{}^*)^p_r\, a^{S\cup\{p\}}_{T\cup\{r\}} = 0 \tag{5.5.11}$$

Согласно (3.3.12) равенство (5.5.11) имеет вид

$$a^p_r - a^p_T{}_*{}^*((a^S_T)^{-1_*{}^*})_*{}^*a^S_r = 0$$

Матрица

$$r^p = a^p_T{}_*{}^*((a^S_T)^{-1_*{}^*}) \tag{5.5.12}$$

не зависит от $r$. Следовательно, для любых $r \in N\setminus T$

$$a^p_r = r^p{}_*{}^*a^S_r \tag{5.5.13}$$

Из равенства

$$((a^S_T)^{-1_*{}^*})_*{}^*a^S_l = \delta^T_l$$

---

[5.6] Мы положим, что неизвестные переменные здесь - это $x_s = R^p_s$

следует, что

$$a_l^p = a_T^p{}_*{}^*\delta_l^T = a_T^p{}_*{}^*((a_T^S)^{-1_*{}^*})_*{}^*a_l^S \tag{5.5.14}$$

Подставляя (5.5.12) в (5.5.14), мы получим

$$a_l^p = r^p{}_*{}^*a_l^S \tag{5.5.15}$$

(5.5.13) и (5.5.15) завершают доказательство. □

Следствие 5.5.5. *Предположим, что $a$ - матрица,*

$$\mathrm{rank}_*{}^*\, a = k < m$$

*Тогда строки матрицы линейно зависимы слева.*

$$\lambda_*{}^*a = 0$$

Доказательство. Предположим, что строка $a^p$ - левая линейная комбинация строк (5.5.9). Мы положим $\lambda_p = -1$, $\lambda_s = R_s^p$ и остальные $\lambda_c = 0$. □

Теорема 5.5.6. *Пусть $(a^i, i \in M, |M| = m)$ семейство строк линейно независимых слева. Тогда ${}_*{}^*$-ранг их координатной матрицы равен $m$.*

Доказательство. Пусть $\overline{\overline{e}}$ - базис левого векторного пространства строк. Согласно конструкции, изложенной в примере 5.3.2, координатная матрица семейства векторов $(\overline{a}^i)$ относительно базиса $\overline{\overline{e}}$ состоит из строк, являющихся координатными матрицами векторов $\overline{a}^i$ относительно базиса $\overline{\overline{e}}$. Поэтому ${}_*{}^*$-ранг этой матрицы не может превышать $m$.

Допустим ${}_*{}^*$-ранг координатной матрицы меньше $m$. Согласно следствию 5.5.5, строки матрицы линейно зависимы слева

$$\lambda_*{}^*a = 0 \tag{5.5.16}$$

Положим $c = \lambda_*{}^*A$. Из равенства (5.5.16) следует, что левая линейная комбинация строк

$$c_*{}^*\overline{e} = 0$$

векторов базиса равна 0. Это противоречит утверждению, что векторы $\overline{e}$ образуют базис. Утверждение теоремы доказано. □

Теорема 5.5.7. *Предположим, что $a$ - матрица,*

$$\mathrm{rank}_*{}^*\, a = k < n$$

*и $a_T^S$ - ${}_*{}^*$-главная минорная матрица. Тогда столбец $a_r$ является ${}_*{}^*D$-линейной композицией столбцов $a_t$*

$$a_{N\setminus T} = a^T{}_*{}^*r \tag{5.5.17}$$

$$a_r = a_T{}_*{}^*r_r \tag{5.5.18}$$

$$a_r^a = a_t^a\ R_r^t \tag{5.5.19}$$

Доказательство. Если матрица $a$ имеет $k$ строк, то, полагая, что столбец $a_r$ - ${}_*{}^*D$-линейная комбинация (5.5.18) столбцов $a_t$ с коэффициентами $R_r^t$, мы получим систему ${}_*{}^*D$-линейных уравнений (5.5.19). Согласно теореме 5.4.7 система ${}_*{}^*D$-линейных уравнений (5.5.19) имеет единственное решение[5.7] и это решение нетривиально потому, что все ${}_*{}^*$-квазидетерминанты отличны от 0.

[5.7] Мы положим, что неизвестные переменные здесь - это ${}^tx = R_r^t$

Остаётся доказать утверждение в случае, когда число строк матрицы $A$ больше чем $k$. Пусть нам даны столбец $a_r$ и строка $a^p$. Согласно предположению, минорная матрица $a^{S\cup\{p\}}_{T\cup\{r\}}$ - ${}_*{}^*$-вырожденная матрица и его ${}_*{}^*$-квазидетерминант

$$\det({}_*{}^*)^p_r\, a^{S\cup\{p\}}_{T\cup\{r\}} = 0 \tag{5.5.20}$$

Согласно (3.3.12) (5.5.20) имеет вид

$$a^p_r - a^p_T{}_*{}^*((a^S_T)^{-1}{}_*{}^*)_*{}^*a^S_r = 0$$

Матрица

$$r_r = ((a^S_T)^{-1}{}_*{}^*)_*{}^*a^S_r \tag{5.5.21}$$

не зависит от $p$, Следовательно, для любых $p \in M \setminus S$

$$a^p_r = a^p_T{}_*{}^*r_r \tag{5.5.22}$$

Из равенства

$$a^k_T{}_*{}^*((a^S_T)^{-1}{}_*{}^*)_s = \delta^k_s$$

следует, что

$$a^k_r = \delta^k_S{}_*{}^*a^S_r = a^k_T{}_*{}^*((a^S_T)^{-1}{}_*{}^*)^s{}_*{}^*a^S_r \tag{5.5.23}$$

Подставляя (5.5.21) в (5.5.23), мы получим

$$a^k_r = a^k_T{}_*{}^*R_r \tag{5.5.24}$$

(5.5.22) и (5.5.24) завершают доказательство. □

Следствие 5.5.8. *Предположим, что $a$ - матрица,*

$$\mathrm{rank}_{{}_*{}^*}\, a = k < m$$

*Тогда столбцы матрицы ${}_*{}^*D$-линейно зависимы.*

$$a_*{}^*\lambda = 0$$

Доказательство. Предположим, что столбец $a_r$ является правой линейной комбинацией (5.5.18). Мы положим $\lambda^r = -1$, $\lambda^t = R^t_r$ и остальные $\lambda^c = 0$. □

Опираясь на теорему 3.3.8, мы можем записать подобные утверждения для ${}^*{}_*$-ранга матрицы.

Теорема 5.5.9. *Предположим, что $a$ - матрица,*

$$\mathrm{rank}_{{}^*{}_*}\, a = k < m$$

*и $a^T_S$ - ${}^*{}_*$-главная минорная матрица. Тогда строка $a^p$ является правой линейной композицией строк $a^s$*

$$a^{M\setminus S} = a^S{}^*{}_*r \tag{5.5.25}$$

$$a^p = a^S{}^*{}_*r^p \tag{5.5.26}$$

$$a^b_p = a^b_s r^s_p \tag{5.5.27}$$

Следствие 5.5.10. *Предположим, что $a$ - матрица,*

$$\mathrm{rank}_{{}^*{}_*}\, a = k < m$$

*Тогда строки матрицы линейно зависимы справа.*

$$A^*{}_*\lambda = 0$$

Теорема 5.5.11. *Предположим, что $a$ - матрица,*

$$\mathrm{rank}_{*_*} a = k < n$$

*и $a_S^T$ - ${}^*{}_*$-главная минорная матрица. Тогда столбец $a_r$ является левой линейной композицией столбцов $a_t$*

$$a_{N\setminus T} = R^*{}_* a_T \tag{5.5.28}$$

$$a_r = R_r{}^*{}_* a_T \tag{5.5.29}$$

$$a_r^i = R_r^t a_t^i \tag{5.5.30}$$

Следствие 5.5.12. *Предположим, что $a$ - матрица,*

$$\mathrm{rank}_{*_*} a = k < m$$

*Тогда столбцы матрицы линейно зависимы слева*

$$\lambda^*{}_* a = 0$$

## 5.6. Собственное значение матрицы

Пусть $a$ - $n \times n$ матрица $A$-чисел и $E$ - $n \times n$ единичная матрица.

Определение 5.6.1. *$A$-число $b$ называется* **${}_*{}^*$-собственным значением** *матрицы $a$, если матрица $a - bE$ - ${}_*{}^*$-вырожденная матрица.*

*Вектор-столбец $v$ называется* **${}_*{}^*$-собственным столбцом** *матрицы $a$, соответствующим ${}_*{}^*$-собственному значению $b$, если верно равенство*

$$a_*{}^* v = bv$$

*Вектор-строка $v$ называется* **${}_*{}^*$-собственной строкой** *матрицы $a$, соответствующей ${}_*{}^*$-собственному значению $b$, если верно равенство*

$$v_*{}^* a = vb$$

□

Определение 5.6.2. *$A$-число $b$ называется* **${}^*{}_*$-собственным значением** *матрицы $a$, если матрица $a - bE$ - ${}^*{}_*$-вырожденная матрица.*

*Вектор-столбец $v$ называется* **${}^*{}_*$-собственным столбцом** *матрицы $a$, соответствующим ${}^*{}_*$-собственному значению $b$, если верно равенство*

$$v^*{}_* a = vb$$

*Вектор-строка $v$ называется* **${}^*{}_*$-собственной строкой** *матрицы $a$, соответствующей ${}^*{}_*$-собственному значению $b$, если верно равенство*

$$a^*{}_* v = bv$$

□

Теорема 5.6.3. *Пусть $a$ - $n \times n$ матрица $A$-чисел и*

$${}_*{}^* \mathrm{rank}\, a = k < n$$

*Тогда $0$ является ${}_*{}^*$-собственным значением кратности $n - k$.*

Доказательство. Теорема является следствием равенства

$$a = a - 0E$$

□

Теорема 5.6.4. *Пусть $a$ - $\boldsymbol{n}\times\boldsymbol{n}$ матрица $A$-чисел и $a^{\boldsymbol{I}}_{\boldsymbol{J}}$ - ${}_*{}^*$-главная минорная матрица матрицы $a$. Тогда ненулевые ${}_*{}^*$-собственные значения матрицы $a$ являются корнями любого уравнения*

$$\det({}_*{}^*)^{\boldsymbol{i}}_{\boldsymbol{j}}\,(a-bE)^{\boldsymbol{I}}_{\boldsymbol{J}}=0 \qquad \boldsymbol{i}\in\boldsymbol{I},\ \boldsymbol{j}\in\boldsymbol{J}$$

Доказательство. Теорема является следствием теоремы 5.5.3. □

Глава 6

# Дифференциальное уравнение, разрешённое относительно производной

## 6.1. Дифференциальное уравнение $\frac{dy}{dx}=1\otimes x+x\otimes 1$

Определение 6.1.1. *Пусть $A$, $B$ - банаховые $D$-модули. Отображение*

$$g:A\to\mathcal{L}(D;A\to B)$$

*называется* **интегрируемым**, *если существует отображение*

$$f:A\to B$$

*такое, что*

$$\frac{df(x)}{dx}=g(x)$$

*Тогда мы пользуемся записью*

$$f(x)=\int g(x)\circ dx$$

*и отображение $f$ называется* **неопределённым интегралом** *отображения $g$.* □

Таблица производных [6.1] и таблица интегралов [6.2] удобный способ интегрирования дифференциального уравнения

$$\frac{dy}{dx}=g(x)$$

если отображение $g$ имеет простой вид.

Определение 6.1.2. *Если существует неопределённый интеграл*

$$\int g(x)\circ dx$$

*то дифференциальное уравнение*

$$\frac{dy}{dx}=g(x)$$

*называется* **интегрируемым**. □

[6.1] Для отображений поля действительных чисел смотри например раздел [20]-95, [2]-2.3. Для отображений банаховой алгебры смотри талже раздел A.1.

[6.2] Для отображений поля действительных чисел смотри например раздел [21]-265, [2]-4.4. Для отображений банаховой алгебры смотри талже раздел A.2.

Пример 6.1.3. *Пусть $A$ - банахова алгебра. Согласно теореме A.2.7, отображение*

$$y=x^2+C$$

*является решением дифференциального уравнения*

$$\frac{dy}{dx}=x\otimes 1+1\otimes x \tag{6.1.1}$$

*при начальном условии*

$$x_0=0\quad y_0=C$$

□

Если $D$-алгебра $B$ имеет конечный базис, то мы можем записать дифференциальное уравнение

$$\frac{dy}{dx}=g(x)$$

в виде системы дифференциальных уравнений

$$\frac{\partial y^i}{\partial x^j}=g^i_j\quad y=y^i e_{B\cdot i}\quad x=x^i e_{A\cdot i}$$

Чтобы вычисления не заслоняли детали преобразования, мы рассмотрим относительно простой пример 6.1.3.

Мы начнём с дифференциального уравнения (6.1.1) в поле комплексных чисел. Теоремы 6.1.4, 6.1.6 просты. Но они необходимы для подготовки читателя к доказательству теоремы 6.1.10.

Теорема 6.1.4. *Дифференциальное уравнение* (6.1.1) *в поле комплексных чисел имеет вид системы дифференциальных уравнений*

$$\left\{\begin{array}{ll}\dfrac{\partial y^0}{\partial x^0}=2x^0 & \dfrac{\partial y^0}{\partial x^1}=-2x^1\\ \dfrac{\partial y^1}{\partial x^0}=2x^1 & \dfrac{\partial y^1}{\partial x^1}=2x^0\end{array}\right. \tag{6.1.2}$$

*относительно базиса*

$$e_0=1\quad e_1=i$$

Доказательство. Поскольку произведение в поле комплексных чисел коммутативно, мы можем записать дифференциальное уравнение (6.1.1) в виде

$$dy=2x\,dx \tag{6.1.3}$$

Если дифференциалы $dx$, $dy$ записать как вектор-столбец, то, согласно теореме [9]-2.5.1, уравнение (6.1.3) принимает вид

$$\begin{pmatrix}dy^0\\ dy^1\end{pmatrix}=\begin{pmatrix}2x^0 & -2x^1\\ 2x^1 & 2x^0\end{pmatrix}\begin{pmatrix}dx^0\\ dx^1\end{pmatrix} \tag{6.1.4}$$

Так как матрица производной имеет вид

$$\frac{dy}{dx}=\begin{pmatrix}\dfrac{\partial y^0}{\partial x^0} & \dfrac{\partial y^0}{\partial x^1}\\ \dfrac{\partial y^1}{\partial x^0} & \dfrac{\partial y^1}{\partial x^1}\end{pmatrix}$$

то равенство

$$\begin{pmatrix}\frac{\partial y^0}{\partial x^0} & \frac{\partial y^0}{\partial x^1}\\ \frac{\partial y^1}{\partial x^0} & \frac{\partial y^1}{\partial x^1}\end{pmatrix}=\begin{pmatrix}2x^0 & -2x^1\\ 2x^1 & 2x^0\end{pmatrix}\tag{6.1.5}$$

является следствием (6.1.4). Система дифференциальных уравнений (6.1.2) является следствием равенства (6.1.5). □

Лемма 6.1.5. *Пусть*

$$x=x^0+x^1i$$

*комплексное число. Тогда*

$$x^2=(x^0)^2-(x^1)^2+2x^0x^1i\tag{6.1.6}$$

Доказательство. Равенство (6.1.6) является следствием определения произведения в поле комплексных чисел. □

Теорема 6.1.6. *Отображение*

$$y=x^2+C$$

*является решением системы дифференциальных уравнений* (6.1.2) *в поле комплексных чисел при начальном условии*

$$x=0\quad y^0=C^0\quad y^1=C^1\quad C=C^0+C^1i$$

Доказательство. Отображение

$$y^0=(x^0)^2+C_1^0(x^1)\tag{6.1.7}$$

является решением дифференциального уравнения

$$\frac{\partial y^0}{\partial x^0}=2x^0$$

Из (6.1.2), (6.1.7) следует, что

$$\frac{\partial y^0}{\partial x^1}=\frac{dC_1^0}{dx^1}=-2x^1\tag{6.1.8}$$

Отображение

$$C_1^0(x^1)=-(x^1)^2+C^0\tag{6.1.9}$$

является решением дифференциального уравнения (6.1.8). Равенство

$$y^0=(x^0)^2-(x^1)^2+C^0\tag{6.1.10}$$

является следствием равенств (6.1.7), (6.1.9).

Отображение

$$y^1=2x^0x^1+C_1^1(x^1)\tag{6.1.11}$$

является решением дифференциального уравнения

$$\frac{\partial y^1}{\partial x^0}=2x^1$$

Из (6.1.2), (6.1.11) следует, что

$$\frac{\partial y^1}{\partial x^1}=2x^0+\frac{dC_1^1}{dx^1}=2x^0\tag{6.1.12}$$

Отображение

(6.1.13) $$C^1_1(x^1) = C^1$$

является решением дифференциального уравнения (6.1.12). Равенство

(6.1.14) $$y^1 = 2x^0x^1 + C^1$$

является следствием равенств (6.1.11), (6.1.13).

Теорема является следствие сравнения равенства (6.1.6) и равенств (6.1.10), (6.1.14). □

Мы рассмотрим дифференциальное уравнение (6.1.1) в алгебре кватернионов аналогично тому как мы рассмотрели это уравнение в поле комплексных чисел.

Теорема 6.1.7. *Дифференциальное уравнение (6.1.1) в алгебре кватернионов имеет вид системы дифференциальных уравнений*

(6.1.15)
$$\left\{\begin{array}{llll}
\frac{\partial y^0}{\partial x^0} = 2x^0 & \frac{\partial y^0}{\partial x^1} = -2x^1 & \frac{\partial y^0}{\partial x^2} = -2x^2 & \frac{\partial y^0}{\partial x^3} = -2x^3\\
\frac{\partial y^1}{\partial x^0} = 2x^1 & \frac{\partial y^1}{\partial x^1} = 2x^0 & \frac{\partial y^1}{\partial x^2} = 0 & \frac{\partial y^1}{\partial x^3} = 0\\
\frac{\partial y^2}{\partial x^0} = 2x^2 & \frac{\partial y^2}{\partial x^1} = 0 & \frac{\partial y^2}{\partial x^2} = 2x^0 & \frac{\partial y^2}{\partial x^3} = 0\\
\frac{\partial y^3}{\partial x^0} = 2x^3 & \frac{\partial y^3}{\partial x^1} = 0 & \frac{\partial y^3}{\partial x^2} = 0 & \frac{\partial y^3}{\partial x^3} = 2x^0
\end{array}\right.$$

*относительно базиса*

$$e_0 = 1 \quad e_1 = i \quad e_2 = j \quad e_3 = k$$

Доказательство. Мы можем записать дифференциальное уравнение (6.1.1) в виде

(6.1.16) $$dy = x\,dx + dx\,x$$

Если дифференциалы $dx$, $dy$ записать как вектор-столбец, то, согласно теоремам [11]-5.1, [11]-5.2, уравнение (6.1.16) принимает вид

$$
\begin{aligned}
\begin{pmatrix} dy^0 \\ dy^1 \\ dy^2 \\ dy^3 \end{pmatrix} &= E_l(x) \begin{pmatrix} dx^0 \\ dx^1 \\ dx^2 \\ dx^3 \end{pmatrix} + E_r(x) \begin{pmatrix} dx^0 \\ dx^1 \\ dx^2 \\ dx^3 \end{pmatrix} \\
&= \begin{pmatrix} x^0 & -x^1 & -x^2 & -x^3 \\ x^1 & x^0 & -x^3 & x^2 \\ x^2 & x^3 & x^0 & -x^1 \\ x^3 & -x^2 & x^1 & x^0 \end{pmatrix} \begin{pmatrix} dx^0 \\ dx^1 \\ dx^2 \\ dx^3 \end{pmatrix} \\
&+ \begin{pmatrix} x^0 & -x^1 & -x^2 & -x^3 \\ x^1 & x^0 & x^3 & -x^2 \\ x^2 & -x^3 & x^0 & x^1 \\ x^3 & x^2 & -x^1 & x^0 \end{pmatrix} \begin{pmatrix} dx^0 \\ dx^1 \\ dx^2 \\ dx^3 \end{pmatrix} \\
&= \begin{pmatrix} 2x^0 & -2x^1 & -2x^2 & -2x^3 \\ 2x^1 & 2x^0 & 0 & 0 \\ 2x^2 & 0 & 2x^0 & 0 \\ 2x^3 & 0 & 0 & 2x^0 \end{pmatrix} \begin{pmatrix} dx^0 \\ dx^1 \\ dx^2 \\ dx^3 \end{pmatrix}
\end{aligned}
\tag{6.1.17}
$$

Так как матрица производной имеет вид

$$
\frac{dy}{dx} = \begin{pmatrix}
\frac{\partial y^0}{\partial x^0} & \frac{\partial y^0}{\partial x^1} & \frac{\partial y^0}{\partial x^2} & \frac{\partial y^0}{\partial x^3} \\
\frac{\partial y^1}{\partial x^0} & \frac{\partial y^1}{\partial x^1} & \frac{\partial y^1}{\partial x^2} & \frac{\partial y^1}{\partial x^3} \\
\frac{\partial y^2}{\partial x^0} & \frac{\partial y^2}{\partial x^1} & \frac{\partial y^2}{\partial x^2} & \frac{\partial y^2}{\partial x^3} \\
\frac{\partial y^3}{\partial x^0} & \frac{\partial y^3}{\partial x^1} & \frac{\partial y^3}{\partial x^2} & \frac{\partial y^3}{\partial x^3}
\end{pmatrix}
$$

то равенство

$$\begin{pmatrix} \frac{\partial y^0}{\partial x^0} & \frac{\partial y^0}{\partial x^1} & \frac{\partial y^0}{\partial x^2} & \frac{\partial y^0}{\partial x^3} \\ \frac{\partial y^1}{\partial x^0} & \frac{\partial y^1}{\partial x^1} & \frac{\partial y^1}{\partial x^2} & \frac{\partial y^1}{\partial x^3} \\ \frac{\partial y^2}{\partial x^0} & \frac{\partial y^2}{\partial x^1} & \frac{\partial y^2}{\partial x^2} & \frac{\partial y^2}{\partial x^3} \\ \frac{\partial y^3}{\partial x^0} & \frac{\partial y^3}{\partial x^1} & \frac{\partial y^3}{\partial x^2} & \frac{\partial y^3}{\partial x^3} \end{pmatrix} = \begin{pmatrix} 2x^0 & -2x^1 & -2x^2 & -2x^3 \\ 2x^1 & 2x^0 & 0 & 0 \\ 2x^2 & 0 & 2x^0 & 0 \\ 2x^3 & 0 & 0 & 2x^0 \end{pmatrix} \tag{6.1.18}$$

является следствием (6.1.17). Система дифференциальных уравнений (6.1.15) является следствием равенства (6.1.18). □

Лемма 6.1.8. *Пусть*

$$\begin{aligned} x &= x^0 + x^1 i + x^2 j + x^3 k \\ y &= y^0 + y^1 i + y^2 j + y^3 k \end{aligned} \tag{6.1.19}$$

*кватернионы. Тогда*

$$\begin{aligned} xy = {} & x^0 y^0 - x^1 y^1 - x^2 y^2 - x^3 y^3 + (x^0 y^1 + x^1 y^0 + x^2 y^3 - x^3 y^2) i \\ & + (x^0 y^2 + x^2 y^0 + x^3 y^1 - x^1 y^3) j + (x^0 y^3 + x^3 y^0 + x^1 y^2 - x^2 y^1) k \end{aligned} \tag{6.1.20}$$

Доказательство. Лемма является следствием определения [7]-11.2.1. □

Лемма 6.1.9. *Пусть*

$$x = x^0 + x^1 i + x^2 j + x^3 k$$

*кватернион. Тогда*

$$x^2 = (x^0)^2 - (x^1)^2 - (x^2)^2 - (x^3)^2 + 2x^0 x^1 i + 2x^0 x^2 j + 2x^0 x^3 k \tag{6.1.21}$$

Доказательство. Равенство

$$\begin{aligned} x^2 = {} & (x^0)^2 - (x^1)^2 - (x^2)^2 - (x^3)^2 + (x^0 x^1 + x^1 x^0 + x^2 x^3 - x^3 x^2) i \\ & + (x^0 x^2 + x^2 x^0 + x^3 x^1 - x^1 x^3) j + (x^0 x^3 + x^3 x^0 + x^1 x^2 - x^2 x^1) k \end{aligned} \tag{6.1.22}$$

является следствием равенства (6.1.20). Равенство (6.1.21) является следствием равенства (6.1.22). □

Теорема 6.1.10. *Отображение*

$$y = x^2 + C$$

*является решением системы дифференциальных уравнений* (6.1.15) *при начальном условии*

$$x = 0 \quad y^0 = C^0 \quad y^1 = C^1 \quad y^2 = C^2 \quad y^3 = C^3 \quad C = C^0 + C^1 i + C^2 j + C^3 k$$

Доказательство. Отображение

$$y^0 = (x^0)^2 + C_1^0(x^1, x^2, x^3) \tag{6.1.23}$$

является решением дифференциального уравнения

$$\frac{\partial y^0}{\partial x^0} = 2x^0$$

Из (6.1.15), (6.1.23) следует, что

(6.1.24) $$\frac{\partial y^0}{\partial x^1} = \frac{\partial C_1^0}{\partial x^1} = -2x^1$$

Отображение

(6.1.25) $$C_1^0(x^1, x^2, x^3) = -(x^1)^2 + C_2^0(x^2, x^3)$$

является решением дифференциального уравнения (6.1.24). Равенство

(6.1.26) $$y^0 = (x^0)^2 - (x^1)^2 + C_2^0(x^2, x^3)$$

является следствием равенств (6.1.23), (6.1.25). Из (6.1.15), (6.1.26) следует, что

(6.1.27) $$\frac{\partial y^0}{\partial x^2} = \frac{\partial C_2^0}{\partial x^2} = -2x^2$$

Отображение

(6.1.28) $$C_2^0(x^2, x^3) = -(x^2)^2 + C_3^0(x^3)$$

является решением дифференциального уравнения (6.1.27). Равенство

(6.1.29) $$y^0 = (x^0)^2 - (x^1)^2 - (x^2)^2 + C_3^0(x^3)$$

является следствием равенств (6.1.26), (6.1.28). Из (6.1.15), (6.1.29) следует, что

(6.1.30) $$\frac{\partial y^0}{\partial x^3} = \frac{dC_3^0}{dx^3} = -2x^3$$

Отображение

(6.1.31) $$C_3^0(x^3) = -(x^3)^2 + C^0$$

является решением дифференциального уравнения (6.1.30). Равенство

(6.1.32) $$y^0 = (x^0)^2 - (x^1)^2 - (x^2)^2 - (x^3)^2 + C^0$$

является следствием равенств (6.1.29), (6.1.31).

Отображение

(6.1.33) $$y^1 = 2x^0 x^1 + C_1^1(x^1, x^2, x^3)$$

является решением дифференциального уравнения

$$\frac{\partial y^1}{\partial x^0} = 2x^1$$

Из (6.1.15), (6.1.33) следует, что

(6.1.34) $$\frac{\partial y^1}{\partial x^1} = 2x^0 + \frac{\partial C_1^1}{\partial x^1} = 2x^0$$

Отображение

(6.1.35) $$C_1^1(x^1, x^2, x^3) = C_2^1(x^2, x^3)$$

является решением дифференциального уравнения (6.1.34). Равенство

(6.1.36) $$y^1 = 2x^0 x^1 + C_2^1(x^2, x^3)$$

является следствием равенств (6.1.33), (6.1.35). Из (6.1.15), (6.1.36) следует, что

(6.1.37) $$\frac{\partial y^1}{\partial x^2} = \frac{\partial C_2^1}{\partial x^2} = 0 \qquad \frac{\partial y^1}{\partial x^3} = \frac{\partial C_2^1}{\partial x^3} = 0$$

Отображение

(6.1.38) $$C_2^1(x^2, x^3) = C^1$$

является решением системы дифференциальных уравнений (6.1.37). Равенство

$$y^1=2x^0x^1+C^1 \tag{6.1.39}$$

является следствием равенств (6.1.36), (6.1.38).

Аналогично мы находим, что

$$y^2=2x^0x^2+C^2 \tag{6.1.40}$$

$$y^3=2x^0x^3+C^3 \tag{6.1.41}$$

Теорема является следствие сравнения равенства (6.1.21) и равенств (6.1.32), (6.1.39), (6.1.40), (6.1.41). □

## 6.2. Дифференциальное уравнение $\frac{dy}{dx}=3x\otimes x$

В этом разделе я рассмотрю дифференциальное уравнение, у которого существование решения зависит от выбора алгебры.

Теорема 6.2.1. *Дифференциальное уравнение*

$$\frac{dy}{dx}=3x\otimes x \tag{6.2.1}$$

*в поле комплексных чисел имеет вид системы дифференциальных уравнений*

$$\left\{\begin{array}{ll}\frac{\partial y^0}{\partial x^0}=3(x^0)^2-3(x^1)^2 & \frac{\partial y^0}{\partial x^1}=-6x^0x^1\\ \frac{\partial y^1}{\partial x^0}=6x^0x^1 & \frac{\partial y^1}{\partial x^1}=3(x^0)^2-3(x^1)^2\end{array}\right. \tag{6.2.2}$$

*относительно базиса*

$$e_0=1\quad e_1=i$$

Доказательство. Поскольку произведение в поле комплексных чисел коммутативно, мы можем записать дифференциальное уравнение (6.2.1) в виде

$$dy=3x^2\,dx \tag{6.2.3}$$

Если дифференциалы $dx$, $dy$ записать как вектор-столбец, то, согласно теореме [9]-2.5.1 и лемме 6.1.5, уравнение (6.2.3) принимает вид

$$\begin{pmatrix}dy^0\\ dy^1\end{pmatrix}=3\begin{pmatrix}(x^0)^2-(x^1)^2 & -2x^0x^1\\ 2x^0x^1 & (x^0)^2-(x^1)^2\end{pmatrix}\begin{pmatrix}dx^0\\ dx^1\end{pmatrix} \tag{6.2.4}$$

Так как матрица производной имеет вид

$$\frac{dy}{dx}=\begin{pmatrix}\frac{\partial y^0}{\partial x^0} & \frac{\partial y^0}{\partial x^1}\\ \frac{\partial y^1}{\partial x^0} & \frac{\partial y^1}{\partial x^1}\end{pmatrix}$$

то равенство

$$\begin{pmatrix}\frac{\partial y^0}{\partial x^0} & \frac{\partial y^0}{\partial x^1}\\ \frac{\partial y^1}{\partial x^0} & \frac{\partial y^1}{\partial x^1}\end{pmatrix}=\begin{pmatrix}3(x^0)^2-3(x^1)^2 & -6x^0x^1\\ 6x^0x^1 & 3(x^0)^2-3(x^1)^2\end{pmatrix} \tag{6.2.5}$$

является следствием (6.2.4). Система дифференциальных уравнений (6.2.2) является следствием равенства (6.2.5). □

Лемма 6.2.2. *Пусть*

$$x = x^0 + x^1 i$$

*комплексное число. Тогда*

$$x^3 = (x^0)^3 - 3x^0(x^1)^2 + (3(x^0)^2 x^1 - (x^1)^3)i \tag{6.2.6}$$

Доказательство. Согласно лемме 6.1.5,

$$\begin{aligned} x^3 &= x^2 x = ((x^0)^2 - (x^1)^2 + 2x^0 x^1 i)(x^0 + x^1 i) \\ &= ((x^0)^2 - (x^1)^2)x^0 + 2(x^0)^2 x^1 i + ((x^0)^2 - (x^1)^2)x^1 i - 2x^0(x^1)^2 \\ &= (x^0)^3 - x^0(x^1)^2 + 2(x^0)^2 x^1 i + ((x^0)^2 x^1 - (x^1)^3)i - 2x^0(x^1)^2 \end{aligned} \tag{6.2.7}$$

Равенство (6.2.6) является следствием равенства (6.2.7). □

Теорема 6.2.3. *Отображение*

$$y = x^3 + C$$

*является решением системы дифференциальных уравнений* (6.2.2) *в поле комплексных чисел при начальном условии*

$$x = 0 \quad y^0 = C^0 \quad y^1 = C^1 \quad C = C^0 + C^1 i$$

Доказательство. Отображение

$$y^0 = (x^0)^3 - 3x^0(x^1)^2 + C_1^0(x^1) \tag{6.2.8}$$

является решением дифференциального уравнения

$$\frac{\partial y^0}{\partial x^0} = 3(x^0)^2 - 3(x^1)^2$$

Из (6.2.2), (6.2.8) следует, что

$$\frac{\partial y^0}{\partial x^1} = -6x^0 x^1 + \frac{dC_1^0}{dx^1} = -6x^0 x^1 \tag{6.2.9}$$

Отображение

$$C_1^0(x^1) = C^0 \tag{6.2.10}$$

является решением дифференциального уравнения (6.2.9). Равенство

$$y^0 = (x^0)^3 - 3x^0(x^1)^2 + C^0 \tag{6.2.11}$$

является следствием равенств (6.2.8), (6.2.10).

Отображение

$$y^1 = 3(x^0)^2 x^1 + C_1^1(x^1) \tag{6.2.12}$$

является решением дифференциального уравнения

$$\frac{\partial y^1}{\partial x^0} = 6x^0 x^1$$

Из (6.2.2), (6.2.12) следует, что

$$\frac{\partial y^1}{\partial x^1} = 3(x^0)^2 + \frac{dC_1^1}{dx^1} = 3(x^0)^2 - 3(x^1)^2 \tag{6.2.13}$$

Отображение

$$C^{1}_{1}(x^{1}) = -(x^{1})^3 + C^{1} \tag{6.2.14}$$

является решением дифференциального уравнения (6.2.13). Равенство

$$y^{1} = 3(x^{0})^2 x^{1} - (x^{1})^3 + C^{1} \tag{6.2.15}$$

является следствием равенств (6.2.12), (6.2.14).

Теорема является следствие сравнения равенства (6.2.6) и равенств (6.2.11), (6.2.15). □

Теорема 6.2.4. *Дифференциальное уравнение*

$$\frac{dy}{dx} = 3x \otimes x \tag{6.2.16}$$

*в некоммутативной алгебре не имеет решений.*

Доказательство. Рассмотрим **метод последовательного дифференцирования** для решения дифференциального уравнения (6.2.16). Последовательно дифференцируя уравнение (6.2.16), мы получаем цепочку уравнений

$$\frac{d^2y}{dx^2} = 3 \otimes 1 \otimes x + 3x \otimes 1 \otimes 1 \tag{6.2.17}$$

$$\frac{d^3y}{dx^3} = 6 \otimes 1 \otimes 1 \otimes 1 \tag{6.2.18}$$

Разложение в ряд Тейлора

$$y = x^3 + C$$

следует из уравнений (6.2.16), (6.2.17), (6.2.18) и начального условия. Согласно теореме A.1.9

$$\frac{dx^3}{dx} \neq 3x \otimes x \tag{6.2.19}$$

Теорема является следствием утверждения (6.2.19). □

Лемма 6.2.5. *В алгебре кватернионов*

$$J_l(a)J_r(b) = J_r(b)J_l(a)$$

$$= \begin{pmatrix} a^0b^0 - a^1b^1 & -a^0b^1 - a^1b^0 & -a^0b^2 - a^1b^3 & -a^0b^3 + a^1b^2 \\ -a^2b^2 - a^3b^3 & +a^2b^3 - a^3b^2 & -a^2b^0 + a^3b^1 & -a^2b^1 - a^3b^0 \\ \hline a^1b^0 + a^0b^1 & -a^1b^1 + a^0b^0 & -a^1b^2 + a^0b^3 & -a^1b^3 - a^0b^2 \\ -a^3b^2 + a^2b^3 & +a^3b^3 + a^2b^2 & -a^3b^0 - a^2b^1 & -a^3b^1 + a^2b^0 \\ \hline a^2b^0 + a^3b^1 & -a^2b^1 + a^3b^0 & -a^2b^2 + a^3b^3 & -a^2b^3 - a^3b^2 \\ +a^0b^2 - a^1b^3 & -a^0b^3 - a^1b^2 & +a^0b^0 + a^1b^1 & +a^0b^1 - a^1b^0 \\ \hline a^3b^0 - a^2b^1 & -a^3b^1 - a^2b^0 & -a^3b^2 - a^2b^3 & -a^3b^3 + a^2b^2 \\ +a^1b^2 + a^0b^3 & -a^1b^3 + a^0b^2 & +a^1b^0 - a^0b^1 & +a^1b^1 + a^0b^0 \end{pmatrix} \tag{6.2.20}$$

Доказательство. Произведение матриц $J_l(a)$, $J_r(b)$ имеет вид

$$
\begin{aligned}
&\begin{pmatrix}
a^0 & -a^1 & -a^2 & -a^3\\
a^1 & a^0 & -a^3 & a^2\\
a^2 & a^3 & a^0 & -a^1\\
a^3 & -a^2 & a^1 & a^0
\end{pmatrix}
\begin{pmatrix}
b^0 & -b^1 & -b^2 & -b^3\\
b^1 & b^0 & b^3 & -b^2\\
b^2 & -b^3 & b^0 & b^1\\
b^3 & b^2 & -b^1 & b^0
\end{pmatrix}\\
&=\left(\begin{array}{c:c:c:c}
\begin{matrix}a^0b^0-a^1b^1\\-a^2b^2-a^3b^3\end{matrix} & \begin{matrix}-a^0b^1-a^1b^0\\+a^2b^3-a^3b^2\end{matrix} & \begin{matrix}-a^0b^2-a^1b^3\\-a^2b^0+a^3b^1\end{matrix} & \begin{matrix}-a^0b^3+a^1b^2\\-a^2b^1-a^3b^0\end{matrix}\\ \hdashline
\begin{matrix}a^1b^0+a^0b^1\\-a^3b^2+a^2b^3\end{matrix} & \begin{matrix}-a^1b^1+a^0b^0\\+a^3b^3+a^2b^2\end{matrix} & \begin{matrix}-a^1b^2+a^0b^3\\-a^3b^0-a^2b^1\end{matrix} & \begin{matrix}-a^1b^3-a^0b^2\\-a^3b^1+a^2b^0\end{matrix}\\ \hdashline
\begin{matrix}a^2b^0+a^3b^1\\+a^0b^2-a^1b^3\end{matrix} & \begin{matrix}-a^2b^1+a^3b^0\\-a^0b^3-a^1b^2\end{matrix} & \begin{matrix}-a^2b^2+a^3b^3\\+a^0b^0+a^1b^1\end{matrix} & \begin{matrix}-a^2b^3-a^3b^2\\+a^0b^1-a^1b^0\end{matrix}\\ \hdashline
\begin{matrix}a^3b^0-a^2b^1\\+a^1b^2+a^0b^3\end{matrix} & \begin{matrix}-a^3b^1-a^2b^0\\-a^1b^3+a^0b^2\end{matrix} & \begin{matrix}-a^3b^2-a^2b^3\\+a^1b^0-a^0b^1\end{matrix} & \begin{matrix}-a^3b^3+a^2b^2\\+a^1b^1+a^0b^0\end{matrix}
\end{array}\right)
\end{aligned}
\tag{6.2.21}
$$

Произведение матриц $J_r(b)$, $J_l(a)$ имеет вид

$$
\begin{aligned}
&\begin{pmatrix}
b^0 & -b^1 & -b^2 & -b^3\\
b^1 & b^0 & b^3 & -b^2\\
b^2 & -b^3 & b^0 & b^1\\
b^3 & b^2 & -b^1 & b^0
\end{pmatrix}
\begin{pmatrix}
a^0 & -a^1 & -a^2 & -a^3\\
a^1 & a^0 & -a^3 & a^2\\
a^2 & a^3 & a^0 & -a^1\\
a^3 & -a^2 & a^1 & a^0
\end{pmatrix}\\
&=\left(\begin{array}{c:c:c:c}
\begin{matrix}a^0b^0-a^1b^1\\-a^2b^2-a^3b^3\end{matrix} & \begin{matrix}-a^0b^1-a^1b^0\\+a^2b^3-a^3b^2\end{matrix} & \begin{matrix}-a^0b^2-a^1b^3\\-a^2b^0+a^3b^1\end{matrix} & \begin{matrix}-a^0b^3+a^1b^2\\-a^2b^1-a^3b^0\end{matrix}\\ \hdashline
\begin{matrix}a^1b^0+a^0b^1\\-a^3b^2+a^2b^3\end{matrix} & \begin{matrix}-a^1b^1+a^0b^0\\+a^3b^3+a^2b^2\end{matrix} & \begin{matrix}-a^1b^2+a^0b^3\\-a^3b^0-a^2b^1\end{matrix} & \begin{matrix}-a^1b^3-a^0b^2\\-a^3b^1+a^2b^0\end{matrix}\\ \hdashline
\begin{matrix}a^2b^0+a^3b^1\\+a^0b^2-a^1b^3\end{matrix} & \begin{matrix}-a^2b^1+a^3b^0\\-a^0b^3-a^1b^2\end{matrix} & \begin{matrix}-a^2b^2+a^3b^3\\+a^0b^0+a^1b^1\end{matrix} & \begin{matrix}-a^2b^3-a^3b^2\\+a^0b^1-a^1b^0\end{matrix}\\ \hdashline
\begin{matrix}a^3b^0-a^2b^1\\+a^1b^2+a^0b^3\end{matrix} & \begin{matrix}-a^3b^1-a^2b^0\\-a^1b^3+a^0b^2\end{matrix} & \begin{matrix}-a^3b^2-a^2b^3\\+a^1b^0-a^0b^1\end{matrix} & \begin{matrix}-a^3b^3+a^2b^2\\+a^1b^1+a^0b^0\end{matrix}
\end{array}\right)
\end{aligned}
\tag{6.2.22}
$$

Равенство (6.2.20) является следствием равенств (6.2.21), (6.2.22). □

Лемма 6.2.6. *Произведение матриц* $J_r(a)$, $J_l(a)$ *имеет вид*

$$
\left(\begin{array}{r:r:r:r}
\begin{matrix} a^0a^0 - a^1a^1 \\ -a^2a^2 - a^3a^3 \end{matrix} & -2a^0a^1 & -2a^0a^2 & -2a^0a^3 \\ \hdashline
2a^0a^1 & \begin{matrix} -a^1a^1 + a^0a^0 \\ +a^3a^3 + a^2a^2 \end{matrix} & -2a^1a^2 & -2a^1a^3 \\ \hdashline
2a^2a^0 & -2a^2a^1 & \begin{matrix} -a^2a^2 + a^3a^3 \\ +a^0a^0 + a^1a^1 \end{matrix} & -2a^2a^3 \\ \hdashline
2a^3a^0 & -2a^3a^1 & -2a^2a^3 & \begin{matrix} -a^3a^3 + a^2a^2 \\ +a^1a^1 + a^0a^0 \end{matrix}
\end{array}\right) \tag{6.2.23}
$$

Доказательство. Согласно лемме 6.2.5, произведение матриц $J_r(a)$, $J_l(a)$ имеет вид

$$
\left(\begin{array}{r:r:r:r}
\begin{matrix} a^0a^0 - a^1a^1 \\ -a^2a^2 - a^3a^3 \end{matrix} & \begin{matrix} -a^0a^1 - a^1a^0 \\ +a^2a^3 - a^3a^2 \end{matrix} & \begin{matrix} -a^0a^2 - a^1a^3 \\ -a^2a^0 + a^3a^1 \end{matrix} & \begin{matrix} -a^0a^3 + a^1a^2 \\ -a^2a^1 - a^3a^0 \end{matrix} \\ \hdashline
\begin{matrix} a^1a^0 + a^0a^1 \\ -a^3a^2 + a^2a^3 \end{matrix} & \begin{matrix} -a^1a^1 + a^0a^0 \\ +a^3a^3 + a^2a^2 \end{matrix} & \begin{matrix} -a^1a^2 + a^0a^3 \\ -a^3a^0 - a^2a^1 \end{matrix} & \begin{matrix} -a^1a^3 - a^0a^2 \\ -a^3a^1 + a^2a^0 \end{matrix} \\ \hdashline
\begin{matrix} a^2a^0 + a^3a^1 \\ +a^0a^2 - a^1a^3 \end{matrix} & \begin{matrix} -a^2a^1 + a^3a^0 \\ -a^0a^3 - a^1a^2 \end{matrix} & \begin{matrix} -a^2a^2 + a^3a^3 \\ +a^0a^0 + a^1a^1 \end{matrix} & \begin{matrix} -a^2a^3 - a^3a^2 \\ +a^0a^1 - a^1a^0 \end{matrix} \\ \hdashline
\begin{matrix} a^3a^0 - a^2a^1 \\ +a^1a^2 + a^0a^3 \end{matrix} & \begin{matrix} -a^3a^1 - a^2a^0 \\ -a^1a^3 + a^0a^2 \end{matrix} & \begin{matrix} -a^3a^2 - a^2a^3 \\ +a^1a^0 - a^0a^1 \end{matrix} & \begin{matrix} -a^3a^3 + a^2a^2 \\ +a^1a^1 + a^0a^0 \end{matrix}
\end{array}\right) \tag{6.2.24}
$$

Равенство (6.2.23) является следствием равенства (6.2.24). □

Теорема 6.2.7. *Дифференциальное уравнение* (6.2.16) *в алгебре кватернионов имеет вид системы дифференциальных уравнений*

$$
\begin{cases}
\dfrac{\partial y^0}{\partial x^0} = 3(x^0)^2 - 3(x^1)^2 - 3(x^2)^2 - 3(x^3)^2 \\
\dfrac{\partial y^0}{\partial x^1} = -6x^0x^1 \\
\dfrac{\partial y^0}{\partial x^2} = -6x^0x^2 \\
\dfrac{\partial y^0}{\partial x^3} = -6x^0x^3
\end{cases} \tag{6.2.25}
$$

(6.2.26)
$$\left\{\begin{array}{l}\frac{\partial y^1}{\partial x^0}=6x^0x^1\\ \frac{\partial y^1}{\partial x^1}=3(x^0)^2-3(x^1)^2+3(x^2)^2+3(x^3)^2\\ \frac{\partial y^1}{\partial x^2}=-6x^1x^2\\ \frac{\partial y^1}{\partial x^3}=-6x^1x^3\end{array}\right.$$

(6.2.27)
$$\left\{\begin{array}{l}\frac{\partial y^2}{\partial x^0}=6x^0x^2\\ \frac{\partial y^2}{\partial x^1}=-6x^1x^2\\ \frac{\partial y^2}{\partial x^2}=-3(x^2)^2+3(x^3)^2+3(x^0)^2+3(x^1)^2\\ \frac{\partial y^2}{\partial x^3}=-6x^2x^3\end{array}\right.$$

(6.2.28)
$$\left\{\begin{array}{l}\frac{\partial y^3}{\partial x^0}=6x^0x^3\\ \frac{\partial y^3}{\partial x^1}=-6x^1x^3\\ \frac{\partial y^3}{\partial x^2}=-6x^2x^3\\ \frac{\partial y^3}{\partial x^3}=-3(x^3)^2+3(x^2)^2+3(x^1)^2+3(x^0)^2\end{array}\right.$$

*относительно базиса*

$$e_0=1\quad e_1=i\quad e_2=j\quad e_3=k$$

Доказательство. Мы можем записать дифференциальное уравнение (6.2.16) в виде

(6.2.29)
$$dy=3x\,dx\,x$$

Если дифференциалы $dx$, $dy$ записать как вектор-столбец, то, согласно лемме 6.2.5, уравнение (6.2.29) принимает вид

(6.2.30)
$$\begin{pmatrix}dy^0\\ dy^1\\ dy^2\\ dy^3\end{pmatrix}=3J_l(x)J_r(x)\begin{pmatrix}dx^0\\ dx^1\\ dx^2\\ dx^3\end{pmatrix}$$

Так как матрица производной имеет вид

$$\frac{dy}{dx} = \begin{pmatrix} \frac{\partial y^0}{\partial x^0} & \frac{\partial y^0}{\partial x^1} & \frac{\partial y^0}{\partial x^2} & \frac{\partial y^0}{\partial x^3} \\ \frac{\partial y^1}{\partial x^0} & \frac{\partial y^1}{\partial x^1} & \frac{\partial y^1}{\partial x^2} & \frac{\partial y^1}{\partial x^3} \\ \frac{\partial y^2}{\partial x^0} & \frac{\partial y^2}{\partial x^1} & \frac{\partial y^2}{\partial x^2} & \frac{\partial y^2}{\partial x^3} \\ \frac{\partial y^3}{\partial x^0} & \frac{\partial y^3}{\partial x^1} & \frac{\partial y^3}{\partial x^2} & \frac{\partial y^3}{\partial x^3} \end{pmatrix}$$

то равенство

(6.2.31)

$$\frac{1}{3}\begin{pmatrix} \frac{\partial y^0}{\partial x^0} & \frac{\partial y^0}{\partial x^1} & \frac{\partial y^0}{\partial x^2} & \frac{\partial y^0}{\partial x^3} \\ \frac{\partial y^1}{\partial x^0} & \frac{\partial y^1}{\partial x^1} & \frac{\partial y^1}{\partial x^2} & \frac{\partial y^1}{\partial x^3} \\ \frac{\partial y^2}{\partial x^0} & \frac{\partial y^2}{\partial x^1} & \frac{\partial y^2}{\partial x^2} & \frac{\partial y^2}{\partial x^3} \\ \frac{\partial y^3}{\partial x^0} & \frac{\partial y^3}{\partial x^1} & \frac{\partial y^3}{\partial x^2} & \frac{\partial y^3}{\partial x^3} \end{pmatrix} =$$

$$\begin{pmatrix} \begin{matrix} x^0x^0 - x^1x^1 \\ -x^2x^2 - x^3x^3 \end{matrix} & -2x^0x^1 & -2x^0x^2 & -2x^0x^3 \\ 2x^0x^1 & \begin{matrix} -x^1x^1 + x^0x^0 \\ +x^3x^3 + x^2x^2 \end{matrix} & -2x^1x^2 & -2x^1x^3 \\ 2x^2x^0 & -2x^2x^1 & \begin{matrix} -x^2x^2 + x^3x^3 \\ +x^0x^0 + x^1x^1 \end{matrix} & -2x^2x^3 \\ 2x^3x^0 & -2x^3x^1 & -2x^2x^3 & \begin{matrix} -x^3x^3 + x^2x^2 \\ +x^1x^1 + x^0x^0 \end{matrix} \end{pmatrix}$$

является следствием (6.2.30). Система дифференциальных уравнений (6.2.25), (6.2.26), (6.2.27), (6.2.28) является следствием равенства (6.2.31). □

Теорема 6.2.8. *Система дифференциальных уравнений* (6.2.25), (6.2.26), (6.2.27), (6.2.28) *не имеет решений.*

Доказательство. Отображение

$$y^1 = 3(x^0)^2x^1 + C_1^1(x^1, x^2, x^3) \tag{6.2.32}$$

является решением дифференциального уравнения

$$\frac{\partial y^1}{\partial x^0} = 6x^0x^1$$

Из (6.2.26), (6.2.32) следует, что

$$\frac{\partial y^1}{\partial x^1}=3(x^0)^2+\frac{\partial C_1^1}{\partial x^1}=3(x^0)^2-3(x^1)^2+3(x^2)^2+3(x^3)^2 \tag{6.2.33}$$

Отображение

$$C_1^1(x^1,x^2,x^3)=-(x^1)^3+3(x^2)^2x^1+3(x^3)^2x^1+C_2^1(x^2,x^3) \tag{6.2.34}$$

является решением дифференциального уравнения (6.2.33). Равенство

$$y^1=3(x^0)^2x^1-(x^1)^3+3(x^2)^2x^1+3(x^3)^2x^1+C_2^1(x^2,x^3) \tag{6.2.35}$$

является следствием равенств (6.2.32), (6.2.34). Из (6.2.26), (6.2.35) следует, что

$$\frac{\partial y^1}{\partial x^2}=6x^2x^1+\frac{\partial C_2^1}{\partial x^2}=-6x^1x^2 \tag{6.2.36}$$

Из уравнения (6.2.36), следует, что отображение $C_2^1$ зависит от $x^1$, что противоречит утверждению, что отображение $C_2^1$ не зависит от $x^1$. Следовательно, система дифференциальных уравнений (6.2.25), (6.2.26), (6.2.27), (6.2.28) не имеет решений. □

## 6.3. Прежде чем двигаться дальше

Прежде чем двигаться дальше, мы сформулируем вопросы, на которые необходимо ответить.

Вопрос 6.3.1. *Если $D$-алгебра $B$ имеет конечный базис $\overline{\overline{e}}$, то мы можем записать дифференциальное уравнение*

$$\frac{df(x)}{dx}=g(x) \tag{6.3.1}$$

*в виде системы дифференциальных уравнений*

$$\frac{\partial y^i}{\partial x^j}=g^i_j \quad y=y^ie_{B\cdot i} \quad x=x^ie_{A\cdot i} \tag{6.3.2}$$

*Какова связь между дифференциальным уравнением (6.3.1) и системой дифференциальных уравнений (6.3.2)?* □

Вопрос 6.3.2. *Каково условие интегрируемости дифференциального уравнения (6.3.1)?* □

Вопрос 6.3.3. *Если дана система дифференциальных уравнений (6.3.2), можем ли мы записать дифференциальное уравнение (6.3.1)?* □

Вопрос 6.3.4. *Системы дифференциальных уравнений (6.1.2), (6.1.15), (6.2.2) являются вполне интегрируемыми.*[6.3] *Система дифференциальных уравнений (6.2.25), (6.2.26), (6.2.27), (6.2.28) не является вполне интегрируемой системой. Равносильно ли требование полной интегрируемости системы дифференциальных уравнений (6.3.2) требованию интегрируемости соответствующего дифференциального уравнения (6.3.1)?* □

[6.3] Смотри определение вполне интегрируемой системы дифференциальных уравнений на страницах 7, 8 в [19].

## 6.4. Формы представления дифференциального уравнения

Теорема 6.4.1 отвечает на вопрос 6.3.1.

Теорема 6.4.1. *Пусть $B$ - свободная конечно мерная ассоциативная $D$-алгебра. Пусть $\overline{\overline{e}}_A$ - базис $D$-модуля $A$. Пусть $\overline{\overline{e}}_B$ - базис $D$-модуля $B$. Пусть $\overline{\overline{F}}$ - базис левого $B\otimes B$-модуля $\mathcal{L}(D;B\to B)$ и*

$$G:A\to B$$

*линейное отображение максимального ранга такое, что $\ker G\subseteq \ker g$. Пусть $g^{k\cdot ij}$ - стандартные компоненты отображения $g$. Пусть $C^p_{kl}$ - структурные константы алгебры $B$. Тогда мы можем записать дифференциальное уравнение (6.3.1) в виде системы дифференциальных уравнений*

$$\frac{\partial y^k}{\partial x^l}=g^{k\cdot ij}F_{k\cdot}{}^m_r G^r_l C^p_{im}C^k_{pj} \tag{6.4.1}$$

Доказательство. Теорема является следствием теоремы [13]-6.4.5. □

Опираясь на теорему 6.4.1, мы можем рассмотреть вопрос 6.3.3. Однако ответ на этот вопрос не однозначен. Мы начнём с теоремы 6.4.2.

Теорема 6.4.2. *Пусть $B$ - свободная конечно мерная ассоциативная $D$-алгебра. Пусть $\overline{\overline{e}}$ - базис $D$-модуля $B$. Пусть $\overline{\overline{F}}$ - базис левого $B\otimes B$-модуля $\mathcal{L}(D;B\to B)$ и*

$$G:A\to B$$

*линейное отображение максимального ранга такое, что $\ker G\subseteq \ker g$. Пусть $C^p_{kl}$ - структурные константы алгебры $B$. Рассмотрим матрицу*

$$\mathcal{C}=\left(\mathcal{C}^{\cdot k}_{m\cdot ij}\right)=\left(C^p_{im}C^k_{pj}\right)$$

*строки и столбцы которой пронумерованы индексами ${}^{\cdot k}_{m}$ и ${}_{\cdot ij}$ соответственно. Если матрица $\mathcal{C}$ невырождена, то мы можем записать систему дифференциальных уравнений*

$$\frac{\partial y^i}{\partial x^j}=g^i_j\qquad y=y^i e_{B\cdot i}\qquad x=x^i e_{A\cdot i}$$

*в виде дифференциального уравнения*

$$\frac{dy}{dx}=g^{k\cdot ij}(e_{B\cdot i}\otimes e_{B\cdot j})\circ F_k\circ G \tag{6.4.2}$$

*где стандартные компоненты $g^{k\cdot ij}$ отображения $g$ являются решением системы линейных уравнений*

$$g^k_l=g^{k\cdot ij}F_{k\cdot}{}^m_r G^r_l C_{B\cdot}{}^p_{im}C_{B\cdot}{}^k_{pj}$$

*Если матрица $\mathcal{C}$ вырождена, то мы можем записать дифференциальное уравнение (6.4.2), если выполнено условие*

$$\operatorname{rank}\begin{pmatrix}\mathcal{C}^{\cdot k}_{m\cdot ij} & g^k_m\end{pmatrix}=\operatorname{rank}\mathcal{C}$$

Доказательство. Теорема является следствием теоремы [13]-6.4.5. □

Анализ теоремы 6.4.2 показывает, что не достаточно записать уравнение (6.4.2), чтобы ответить на вопрос 6.3.3. Причина этого в том, что, согласно построению, отображения $g^{k\cdot ij}$ зависят от координат $A$-числа $x$. Однако необходимо отметить, что отображение

$$x = x^i e_i \to x^j$$

является линейным отображением. Поэтому, если $D$-модуль $A$ является $D$-алгеброй, то существует эффективная процедура ответа на вопрос 6.3.3.

Замечание 6.4.3. *Утверждение этого замечания предварительно, так как исследование в этом направлении не завершено. Стандартные компоненты отображения $g$ - это универсальная форма записи отображения, но не единственная и не самая выразительная. Пусть $\overline{\overline{E}}$ - базис $B$-модуля $\mathcal{L}(D; B \to B)$. Тогда отображение $g$ можно представить в виде*

$$g(x) = g^k(x) \circ E_k \qquad g^k : A \to B$$

□

## 6.5. Условие интегрируемости

Чтобы ответить на вопрос 6.3.2, напомню, что отображение $g$ является дифференциальной формой (определение [15]-7.3.2).

Теорема 6.5.1. *Дифференциальное уравнение*

$$\frac{df(x)}{dx} = g(x)$$

*интегрируемо тогда и только тогда, когда*

$$dg = 0$$

Доказательство. Теорема является следствием теоремы [15]-8.3.1. □

Рассмотрим сперва как работает теорема 6.5.1 в случае дифференциальных уравнений (6.1.1), (6.2.16).

Пример 6.5.2. *Рассмотрим дифференциальную форму*

$$g(x) = x \otimes 1 + 1 \otimes x$$

*Согласно теоремам A.1.1, A.1.2, A.1.4, A.1.11,*

$$\frac{dg}{dx} = 1 \otimes_2 1 \otimes_1 1 + 1 \otimes_1 1 \otimes_2 1$$

*Индекс, сопровождающий символ $\otimes$, показывает, какой аргумент полилинейного отображения должен быть записан вместо соответствующего символа $\otimes$. Например*

$$\frac{dg}{dx} \circ (h_1, h_2) = (1 \otimes_2 1 \otimes_1 1 + 1 \otimes_1 1 \otimes_2 1) \circ (h_1, h_2) = h_2 h_1 + h_1 h_2 \tag{6.5.1}$$

*Из равенства (6.5.1) следует, что билинейное отображение $\frac{dg}{dx}$ симметрично. Согласно определению [15]-7.4.1*

$$dg = 0$$

*Согласно теореме 6.5.1, дифференциальное уравнение (6.1.1) интегрируемо (теорема 6.1.7).* □

Пример 6.5.3. *Рассмотрим дифференциальную форму*

$$g(x) = x \otimes x$$

*Согласно теоремам A.1.2, A.1.11,*

$$\frac{dg}{dx} = 1 \otimes_2 1 \otimes_1 x + x \otimes_1 1 \otimes_2 1$$

*Индекс, сопровождающий символ $\otimes$, показывает, какой аргумент полилинейного отображения должен быть записан вместо соответствующего символа $\otimes$. Например*

$$\frac{dg}{dx} \circ (h_1, h_2) = (1 \otimes_2 1 \otimes_1 x + x \otimes_1 1 \otimes_2 1) \circ (h_1, h_2) = h_2 h_1 x + x h_1 h_2 \tag{6.5.2}$$

*Из равенства (6.5.2) и определения [15]-7.4.1 следует, что*

$$dg = (h_2 h_1 - h_1 h_2) x + x (h_1 h_2 - h_2 h_1) \tag{6.5.3}$$

*Из равенства (6.5.3) следует, что*

$$dg \neq 0$$

*Согласно теореме 6.5.1, дифференциальное уравнение (6.2.16) не интегрируемо (теорема 6.2.4).* □

Теорема 6.5.4 отвечает на вопрос 6.3.2.

Теорема 6.5.4. *Пусть $B$ - свободная конечно мерная ассоциативная $D$-алгебра. Пусть $\overline{\overline{e}}$ - базис $D$-модуля $B$. Пусть $\overline{\overline{\overline{F}}}$ - базис левого $B \otimes B$-модуля $\mathcal{L}(D; B \to B)$ и $G$ - линейное отображение*

$$G : A \to B$$

*максимального ранга такое, что $\ker G \subseteq \ker g$. Пусть $g^{k \cdot ij}$ - стандартные компоненты отображения*

$$g : A \to \mathcal{L}(D; A \to B)$$

*Дифференциальное уравнение*

$$\frac{dy}{dx} = g^{k \cdot ij} (e_{B \cdot i} \otimes e_{B \cdot j}) \circ F_k \circ G$$

*интегрируемо тогда и только тогда, когда*

$$\begin{aligned} &\left( \frac{dg^{k \cdot ij}(x)}{dx} \circ h_2 \right) (e_{B \cdot i} (F_k \circ G \circ h_1) e_{B \cdot j}) \\ = &\left( \frac{dg^{k \cdot ij}(x)}{dx} \circ h_1 \right) (e_{B \cdot i} (F_k \circ G \circ h_2) e_{B \cdot j}) \end{aligned} \tag{6.5.4}$$

Доказательство. Поскольку

$$g \circ h_1 = g^{k \cdot ij} (e_i \otimes e_j) \circ F_k \circ G \circ h_1 = g^{k \cdot ij} (e_i (F_k \circ G \circ h_1) e_j) \tag{6.5.5}$$

то, согласно теореме [15]-8.2.3,

$$\frac{dg}{dx} \circ (h_1, h_2) = \left( \frac{dg^{k \cdot ij}}{dx} \circ h_2 \right) (e_i (F_k \circ G \circ h_1) e_j) \tag{6.5.6}$$

Равенство (6.5.4) является следствием равенства (6.5.6), теоремы 6.5.1 и определения [15]-7.4.1. □

Теорема 6.5.5 отвечает на вопрос 6.3.4.

Теорема 6.5.5. *Пусть $B$ - свободная конечно мерная ассоциативная $D$-алгебра. Пусть $\overline{\overline{e}}_A$ - базис $D$-модуля $A$. Пусть $\overline{\overline{e}}_B$ - базис $D$-модуля $B$. Дифференциальное уравнение*

$$\frac{dy}{dx}=g(x) \tag{6.5.7}$$

*интегрируемо тогда и только тогда, когда вполне интегрируема соответствующая система дифференциальных уравнений*

$$\frac{\partial y^i}{\partial x^j}=g^i_j \quad y=y^i e_{B\cdot i} \quad x=x^i e_{A\cdot i} \tag{6.5.8}$$

Доказательство. Пусть $\overline{\overline{F}}$ - базис левого $B\otimes B$-модуля $\mathcal{L}(D;B\to B)$ и $G$ - линейное отображение

$$G:A\to B$$

максимального ранга такое, что $\ker G\subseteq\ker g$. Пусть $g^{k\cdot ij}$ - стандартные компоненты отображения

$$g:A\to\mathcal{L}(D;A\to B)$$

Тогда дифференциальное уравнение (6.5.7) имеет вид

$$\frac{dy}{dx}=g^{k\cdot ij}(e_{B\cdot i}\otimes e_{B\cdot j})\circ F_k\circ G \tag{6.5.9}$$

Согласно теореме 6.5.4, дифференциальное уравнение (6.5.7) интегрируемо тогда и только тогда, когда

$$\begin{aligned}&\left(\frac{dg^{k\cdot ij}(x)}{dx}\circ h_2\right)(e_{B\cdot i}(F_k\circ G\circ h_1)e_{B\cdot j})\\ =&\left(\frac{dg^{k\cdot ij}(x)}{dx}\circ h_1\right)(e_{B\cdot i}(F_k\circ G\circ h_2)e_{B\cdot j})\end{aligned} \tag{6.5.10}$$

Пусть

$$h_1=h_1^k e_k \quad h_2=h_2^l e_l$$

Тогда мы можем записать равенство (6.5.10) в виде

$$\begin{aligned}&\left(\frac{dg^{k\cdot ij}(x)}{dx}\circ(h_2^l e_{1\cdot l})\right)(e_{2\cdot i}(F_k\circ G\circ(h_1^k e_{1\cdot k}))e_{2\cdot j})\\ =&\left(\frac{dg^{k\cdot ij}(x)}{dx}\circ(h_1^k e_{1\cdot k})\right)(e_{2\cdot i}(F_k\circ G\circ(h_2^l e_{1\cdot l}))e_{2\cdot j})\end{aligned} \tag{6.5.11}$$

Так как $h_1$, $h_2$ - произвольные $A$-числа, то равенство

$$\begin{aligned}&\left(\frac{dg^{k\cdot ij}(x)}{dx}\circ e_{1\cdot l}\right)(e_{2\cdot i}(F_k\circ G\circ e_{1\cdot k})e_{2\cdot j})\\ =&\left(\frac{dg^{k\cdot ij}(x)}{dx}\circ e_{1\cdot k}\right)(e_{2\cdot i}(F_k\circ G\circ e_{1\cdot l})e_{2\cdot j})\end{aligned} \tag{6.5.12}$$

является следствием равенства (6.5.11). Так как

$$F_k\circ G\circ e_{1\cdot k}=F_{k\,r}^{\ m}G^r_k e_{2\cdot m}$$

то равенство

$$\begin{aligned}(6.5.13)\qquad &\left(\frac{dg^{k\cdot ij}(x)}{dx}\circ e_{1\cdot l}\right)F_{k}{}^{m}_{r}G^{r}_{k}e_{2\cdot i}e_{2\cdot m}e_{2\cdot j}\\ &=\left(\frac{dg^{k\cdot ij}(x)}{dx}\circ e_{1\cdot k}\right)F_{k}{}^{m}_{r}G^{r}_{l}e_{2\cdot i}e_{2\cdot m}e_{2\cdot j}\end{aligned}$$

является следствием равенства (6.5.12). Равенство

$$\begin{aligned}(6.5.14)\qquad &\left(\frac{dg^{ij}(x)}{dx}\circ e_{1\cdot l}\right)F_{k}{}^{m}_{r}G^{r}_{k}C_{2\cdot}{}^{p}_{im}C_{2\cdot}{}^{r}_{pj}e_{2\cdot r}\\ &=\left(\frac{dg^{k\cdot ij}(x)}{dx}\circ e_{1\cdot k}\right)F_{k}{}^{m}_{r}G^{r}_{l}C_{2\cdot}{}^{p}_{im}C_{2\cdot}{}^{r}_{pj}e_{2\cdot r}\end{aligned}$$

является следствием равенства (6.5.13). Равенство

$$(6.5.15)\qquad \left(\frac{dg^{r}_{k}(x)}{dx}\circ e_{1\cdot l}\right)e_{2\cdot r}=\left(\frac{dg^{r}_{l}(x)}{dx}\circ e_{1\cdot k}\right)e_{2\cdot r}$$

является следствием равенства (6.5.14). Равенство

$$(6.5.16)\qquad \frac{dg^{r}_{k}(x)}{dx}\circ e_{1\cdot l}=\frac{dg^{r}_{l}(x)}{dx}\circ e_{1\cdot k}$$

является следствием равенства (6.5.15). □

Глава 7

# Дифференциальное уравнение первого порядка

## 7.1. Дифференциальное уравнение с разделёнными переменными

В теории дифференциальных уравнений над полем первое уравнение, которое мы изучаем, - это уравнение с разделяющимися переменными

$$\frac{dy}{dx} = \frac{M(x)}{N(y)} \tag{7.1.1}$$

Разделив переменные, мы можем записать дифференциальное уравнение (7.1.1) в виде

$$N(y)dy = M(x)dx \tag{7.1.2}$$

Дифференциальное уравнение (7.1.2) легко интегрировать.

В некоммутативной $D$-алгебре, $dx$ и $dy$ принадлежат, вообще говоря, различным $D$-модулям. В тоже время, решение дифференциального уравнения (7.1.2) является неявной функцией переменных $x$ и $y$. Поэтому я сформулирую задачу следующим образом.

Пусть $X$, $Y$ - банаховы $D$-модули. Пусть $A$ - банахова $D$-алгебра. Рассмотрим отображения

$$f : X \to \mathcal{L}(D; X \to A)$$
$$g : Y \to \mathcal{L}(D; Y \to A)$$

Дифференциальное уравнение

$$\frac{dy}{dx} = g(y) \circ f(x) \tag{7.1.3}$$

называется **дифференциальным уравнением с разделяющимися переменными**.

В некоммутативной $D$-алгебре, даже если это алгебра с делением, операция разделения переменных может оказаться невыполнимой. Однако, если существует отображение

$$h : Y \to \mathcal{L}(D; Y \to A)$$

такое, что

$$h \circ g = 1 \otimes 1$$

то мы можем записать дифференциальное уравнение (7.1.3) в следующем виде

$$h(y) \circ dy = f(x) \circ dx$$

Пусть $X$, $Y$ - банаховы $D$-модули. Пусть $A$ - банахова $D$-алгебра. Рассмотрим отображения

$$M : X \to \mathcal{L}(D; X \to A)$$
$$N : Y \to \mathcal{L}(D; Y \to A)$$

Дифференциальное уравнение

$$M(x)\circ dx+N(y)\circ dy=0$$

называется **дифференциальным уравнением с разделёнными переменными**.

Теорема 7.1.1. *Пусть отображения*

$$M:X\to\mathcal{L}(D;X\to A)$$
$$N:Y\to\mathcal{L}(D;Y\to A)$$

*интегрируемы. Решением дифференциального уравнения*

$$M(x)\circ dx+N(y)\circ dy=0 \tag{7.1.4}$$

*является неявное отображение*[7.1]

$$\int M(x)\circ dx+\int N(y)\circ dy=C \tag{7.1.5}$$

Доказательство. Равенство

$$N(y)\circ dx=N(y)\circ\frac{dy}{dx}\circ dx \tag{7.1.6}$$

является следствием определения [15]-3.3.2. Из интегрируемости дифференциальной формы в левой части равенства (7.1.6) следует интегрируемость дифференциальной формы в правой части равенства (7.1.6). Согласно теореме A.2.4 дифференциальная форма

$$M(x)\circ dx+N(y)\circ dx=\left(M(x)+N(y)\circ\frac{dy}{dx}\right)\circ dx \tag{7.1.7}$$

интегрируема. Равенство

$$M(x)+N(y)\circ\frac{dy}{dx}=0\otimes 0 \tag{7.1.8}$$

является следствием равенств (7.1.4), (7.1.7). Равенство (7.1.5) является следствием равенства (7.1.8) и теоремы A.2.3. □

Пример 7.1.2. *Рассмотрим дифференциальное уравнение*

$$(1\otimes x+x\otimes 1)\circ dx+(1\otimes y+y\otimes 1)\circ dy=0 \tag{7.1.9}$$

*Согласно теоремам 7.1.1, A.2.6, неявное отображение*

$$x^2+y^2=C$$

*является решением дифференциального уравнения (7.1.9).* □

Теорема 7.1.3. *Решение дифференциального уравнения*

$$M(x)\circ dx+N(y)\circ dy=0 \tag{7.1.10}$$

*при начальном условии*

$$x_0=0\quad y_0=C$$

---

[7.1] Доказательство теоремы опирается на доказательства на странице [3]-44 и на странице [6]-19, однако несколько отличается, так как $D$-алгебра $A$ вообще говоря некоммутативна.

*имеет вид*[7.2]

$$\int_{x_0}^{x} M(x) \circ dx + \int_{y_0}^{y} N(y) \circ dy = 0 \tag{7.1.11}$$

Доказательство. Согласно теореме 7.1.1, неявная функция

$$\int M(x) \circ dx + \int N(y) \circ dy = C \tag{7.1.12}$$

является решением дифференциального уравнения (7.1.10). Согласно определению 6.1.1, существуют отображения

$$P : X \to A$$
$$R : Y \to A$$

такие, что

$$P(x) = \int M(x) \circ dx \tag{7.1.13}$$

$$R(y) = \int N(y) \circ dy \tag{7.1.14}$$

Согласно теореме [15]-6.2.1 и определению [15]-8.3.7, равенства

$$\int M(x) \circ dx = \int_{x_0}^{x} M(x) \circ dx + P(x_0) \tag{7.1.15}$$

$$\int N(y) \circ dy = \int_{y_0}^{y} N(y) \circ dy + R(y_0) \tag{7.1.16}$$

являются следствием равенств (7.1.13), (7.1.14). Равенство

$$\int_{x_0}^{x} M(x) \circ dx + P(x_0) + \int_{y_0}^{y} N(y) \circ dy + R(x_0) = C \tag{7.1.17}$$

являются следствием равенств (7.1.12), (7.1.15), (7.1.16). Если $x = x_0$, $y = y_0$, то равенство

$$P(x_0) + R(x_0) = C \tag{7.1.18}$$

являются следствием равенства (7.1.17). Равенство (7.1.11) являются следствием равенств (7.1.17), (7.1.18). □

## 7.2. Уравнение в полных дифференциалах

Определение 7.2.1. *Дифференциальное уравнение*

$$M(x, y) \circ dx + N(x, y) \circ dy = 0$$

*где*

$$M : (x, y) \in X \times Y \to M(x, y) \in \mathcal{L}(D; X \to A)$$
$$N : (x, y) \in X \times Y \to N(x, y) \in \mathcal{L}(D; Y \to A)$$

*называется* **уравнением в полных дифференциалах***, если существует отображение*

$$u : X \times Y \to A$$

[7.2] Смотри также замечания на странице [3]-44 и на странице [6]-20.

*такое, что*

$$\frac{\partial u(x,y)}{\partial x}=M(x,y)\tag{7.2.1}$$

$$\frac{\partial u(x,y)}{\partial y}=N(x,y)\tag{7.2.2}$$

□

Теорема 7.2.2. *Дифференциальное уравнение*

$$M(x,y)\circ dx+N(x,y)\circ dy=0\tag{7.2.3}$$

*интегрируемо тогда и только тогда, когда*

$$\frac{\partial M(x,y)}{\partial x}\circ(dx_1,dx_2)=\frac{\partial M(x,y)}{\partial x}\circ(dx_2,dx_1)\tag{7.2.4}$$

$$\frac{\partial N(x,y)}{\partial y}\circ(dy_1,dy_2)=\frac{\partial N(x,y)}{\partial y}\circ(dy_2,dy_1)\tag{7.2.5}$$

$$\frac{\partial M(x,y)}{\partial y}\circ(dx,dy)=\frac{\partial N(x,y)}{\partial x}\circ(dy,dx)\tag{7.2.6}$$

Доказательство. Пусть $Z=X\oplus Y$, $z=x\oplus y$. Согласно определениям [15]-7.3.2, 2.2.5, 7.2.1 и теореме 2.2.6, выражение

$$\omega=M(x,y)\circ dx+N(x,y)\circ dy\tag{7.2.7}$$

является дифференциальной формой степени 1. Согласно теореме [15]-8.3.1, дифференциальная форма $\omega$ интегрируема тогда и только тогда, когда $d\omega=0$. Равенства

$$\begin{aligned}\frac{d\omega}{dz}\circ(z_1,z_2)=&\frac{\partial M(x,y)}{\partial x}\circ(dx_1,dx_2)+\frac{\partial M(x,y)}{\partial y}\circ(dx_1,dy_2)\\&+\frac{\partial N(x,y)}{\partial x}\circ(dy_1,dx_2)+\frac{\partial N(x,y)}{\partial y}\circ(dy_1,dy_2)\end{aligned}\tag{7.2.8}$$

$$\begin{aligned}\frac{d\omega}{dz}\circ(z_2,z_1)=&\frac{\partial M(x,y)}{\partial x}\circ(dx_2,dx_1)+\frac{\partial M(x,y)}{\partial y}\circ(dx_2,dy_1)\\&+\frac{\partial N(x,y)}{\partial x}\circ(dy_2,dx_1)+\frac{\partial N(x,y)}{\partial y}\circ(dy_2,dy_1)\end{aligned}\tag{7.2.9}$$

являются следствием равенства (7.2.7). Равенство

$$\begin{aligned}&\frac{\partial M(x,y)}{\partial x}\circ(dx_1,dx_2)+\frac{\partial M(x,y)}{\partial y}\circ(dx_1,dy_2)\\&+\frac{\partial N(x,y)}{\partial x}\circ(dy_1,dx_2)+\frac{\partial N(x,y)}{\partial y}\circ(dy_1,dy_2)\\&=\frac{\partial M(x,y)}{\partial x}\circ(dx_2,dx_1)+\frac{\partial M(x,y)}{\partial y}\circ(dx_2,dy_1)\\&+\frac{\partial N(x,y)}{\partial x}\circ(dy_2,dx_1)+\frac{\partial N(x,y)}{\partial y}\circ(dy_2,dy_1)\end{aligned}\tag{7.2.10}$$

является следствием равенств (7.2.8), (7.2.9). Равенства (7.2.4), (7.2.5), (7.2.6) являются следствием равенства (7.2.10). □

Пусть верны равенства (7.2.4), (7.2.5), (7.2.6). Наша задача построить неявную функцию[7.3] $u(x,y)=C$, которая удовлетворяет равенствам (7.2.1), (7.2.2). Согласно теоремам [15]-8.3.1, A.2.1, равенство

$$u(x,y)=\int M(x,y)\circ dx+C_1(y) \tag{7.2.11}$$

является следствием равенств (7.2.1), (7.2.4). Равенство

$$\frac{\partial}{\partial y}\int M(x,y)\circ dx+\frac{dC_1(y)}{dy}=N(x,y) \tag{7.2.12}$$

является следствием равенств (7.2.2), (7.2.11). Равенство

$$\frac{dC_1(y)}{dy}=N(x,y)-\frac{\partial}{\partial y}\int M(x,y)\circ dx \tag{7.2.13}$$

является следствием равенства (7.2.12). В выражении

$$N(x,y)-\frac{\partial}{\partial y}\int M(x,y)\circ dx \tag{7.2.14}$$

правой части уравнения (7.2.13) присутствует переменная $x$. Для того, чтобы уравнение (7.2.13) имело решение, необходимо чтобы выражение (7.2.14) не зависело от $x$. Рассмотрим производную выражения (7.2.14) по $x$. Равенство

$$\begin{aligned}
&\frac{\partial}{\partial x}\left(N(x,y)\circ dy-\left(\frac{\partial}{\partial y}\int M(x,y)\circ dx\right)\circ dy\right)\circ dx\\
&=\frac{\partial N(x,y)}{\partial x}\circ(dy,dx)-\left(\frac{\partial^2}{\partial x\partial y}\int M(x,y)\circ dx\right)\circ(dy,dx)\\
&=\frac{\partial N(x,y)}{\partial x}\circ(dy,dx)-\left(\frac{\partial^2}{\partial y\partial x}\int M(x,y)\circ dx\right)\circ(dx,dy)\\
&=\frac{\partial N(x,y)}{\partial x}\circ(dy,dx)-\left(\frac{\partial}{\partial y}\frac{\partial}{\partial x}\int M(x,y)\circ dx\right)\circ dy\\
&=\frac{\partial N(x,y)}{\partial x}\circ(dy,dx)-\frac{\partial M(x,y)}{\partial y}\circ(dx,dy)=0
\end{aligned} \tag{7.2.15}$$

является следствием равенств [8]-(9.1.5), (7.2.6), (A.2.2).

Из теорем [15]-8.3.1, A.2.1 и равенства (7.2.5) следует существование интеграла

$$\begin{aligned}
&\int\left(N(x,y)-\frac{\partial}{\partial y}\int M(x,y)\circ dx\right)\circ dy\\
&=\int N(x,y)\circ dy-\int\left(\frac{\partial}{\partial y}\int M(x,y)\circ dx\right)\circ dy
\end{aligned}$$

Следовательно, мы можем определить отображение $C_1$ как решение дифференциального уравнения (7.2.13).

---

[7.3] Процедура решения уравнения (7.2.3) похожа на процедуру решения уравнения в полных дифференциалах в коммутативной алгебре. Смотри, например, доказательство теоремы [3]-2.6.1, а также процедуру решения уравнения в полных дифференциалах на страницах 33, 34 в книге [6]-19.

Пример 7.2.3. *Рассмотрим дифференциальное уравнение*

$$(1\otimes 1+1\otimes y)\circ dx+(x\otimes 1+1\otimes 1)\circ dy=0 \tag{7.2.16}$$

*Наша задача - найти неявную функцию* $u(x,y)$*, которая является решением дифференциального уравнения* (7.2.16)*. Дифференциальные уравнения*

$$\frac{\partial u}{\partial x}=1\otimes 1+1\otimes y \tag{7.2.17}$$

$$\frac{\partial u}{\partial y}=x\otimes 1+1\otimes 1 \tag{7.2.18}$$

*являются следствием дифференциального уравнения* (7.2.16)*. Тогда*

$$u(x,y)=\int(1\otimes 1+1\otimes y)\circ dx=x+xy+C_1(y) \tag{7.2.19}$$

*является решением дифференциального уравнения* (7.2.17)*. Равенство*

$$\frac{\partial u(x,y)}{\partial y}=x\otimes 1+\frac{dC_1(y)}{dy}=x\otimes 1+1\otimes 1 \tag{7.2.20}$$

*является следствием равенства* (7.2.19) *и уравнения* (7.2.18)*. Дифференциальное уравнение*

$$\frac{dC_1(y)}{dy}=1\otimes 1 \tag{7.2.21}$$

*является следствием равенства* (7.2.20)*.*

$$C_1(y)=y \tag{7.2.22}$$

*является решением дифференциального уравнения* (7.2.21)*. Из равенств* (7.2.19)*,* (7.2.22) *следует, что неявная функция*

$$x+xy+y=C$$

*является решением дифференциального уравнения* (7.2.16)*.* □

Пример 7.2.4 иллюстрирует, что равенства (7.2.4), (7.2.5) существенны для интегрируемости дифференциального уравнения (7.2.3).

Пример 7.2.4. *Рассмотрим дифференциальное уравнение*

$$(3x^2\otimes 1+1\otimes y)\circ dx+(x\otimes 1)\circ dy=0 \tag{7.2.23}$$

*Здесь*

$$M(x,y)=3x^2\otimes 1+1\otimes y \qquad N(x,y)=x\otimes 1$$

*Следовательно,*

$$\frac{\partial M(x,y)}{\partial y}=1\otimes_x 1\otimes_y 1$$
$$\frac{\partial N(x,y)}{\partial x}=1\otimes_x 1\otimes_y 1$$

*Однако дифференциальная уравнение* (7.2.23) *не является интегрируемым, так как интеграл*

$$\int(3x^2\otimes 1+1\otimes y)\circ dx$$

*не существует.* □

Пример 7.2.5 иллюстрирует, что порядок переменных существеннен в равенстве (7.2.6).

Пример 7.2.5. *Рассмотрим дифференциальное уравнение*

$$(1\otimes y)\circ dx+(1\otimes x)\circ dy=0 \tag{7.2.24}$$

*Здесь*

$$M(x,y)=1\otimes y\qquad N(x,y)=1\otimes x$$

*Следовательно,*

$$\frac{\partial M(x,y)}{\partial y}=1\otimes_x 1\otimes_y 1$$
$$\frac{\partial N(x,y)}{\partial x}=1\otimes_y 1\otimes_x 1$$

*Хотя производные* $\frac{\partial M(x,y)}{\partial y}$, $\frac{\partial N(x,y)}{\partial x}$ *представлены одним и тем же тензором, их действие на переменные* $x$, $y$ *различно*

$$\frac{\partial M(x,y)}{\partial y}\circ(dx,dy)=(1\otimes_x 1\otimes_y 1)\circ(dx,dy)=dx\,dy$$
$$\neq\frac{\partial N(x,y)}{\partial x}\circ(dy,dx)=(1\otimes_y 1\otimes_x 1)\circ(dy,dx)=dy\,dx$$

*Легко видеть, что дифференциальное уравнение (7.2.24) не является интегрируемым.* □

Глава 8

# Экспонента

## 8.1. Определение экспоненты

Определение 8.1.1. *Мы назовём отображение*

$$y = e^x = \sum_{n=0}^{\infty} \frac{1}{n!} x^n \tag{8.1.1}$$

*экспонентой.* □

Теорема 8.1.2. *Производная экспоненты имеет следующий вид*

$$\frac{de^x}{dx} = \sum_{n=0}^{\infty} \frac{1}{(n+1)!} \sum_{i=0}^{n} x^i \otimes x^{n-i} \tag{8.1.2}$$

Доказательство. Согласно теоремам A.1.3, A.1.4, равенство

$$\frac{de^x}{dx} = \sum_{n=0}^{\infty} \frac{1}{(n+1)!} \frac{dx^{n+1}}{dx} \tag{8.1.3}$$

является следствием равенства (8.1.1). Согласно теореме 4.2.11, равенство (8.1.2) является следствием равенства (8.1.3). □

Равенство (8.1.2) является дифференциальным уравнением для экспоненты. Однако это дифференциальное уравнение разрешено относительно производной. В поле экспонента является решением дифференциального уравнения

$$y' = y \tag{8.1.4}$$

Наша задача - рассмотреть дифференциальное уравнение в некоммутативной банаховой алгебре, которое аналогично дифференциальному уравнению (8.1.4).

Теорема 8.1.3. *Отображение* $y = ae^x b$, $ab = C$, *является решением дифференциального уравнения*

$$\frac{dy}{dx} \circ 1 = y \tag{8.1.5}$$

*при начальном условии*

$$x = 0 \quad y = C$$

Доказательство. Равенство

$$\begin{aligned} \frac{de^x}{dx} \circ 1 &= \sum_{n=0}^{\infty} \frac{1}{(n+1)!} \sum_{i=0}^{n} (x^i \otimes x^{n-i}) \circ 1 = \sum_{n=0}^{\infty} \frac{1}{(n+1)!} \sum_{i=0}^{n} x^i x^{n-i} \\ &= \sum_{n=0}^{\infty} \frac{n+1}{(n+1)!} x^n \end{aligned} \tag{8.1.6}$$

является следствием равенства (8.1.2). Из равенств (8.1.1), (8.1.6) следует, что

$$\frac{de^x}{dx}\circ 1=e^x \tag{8.1.7}$$

и отображение $y=e^x$ является решением дифференциального уравнения (8.1.5).

Пусть

$$y=ae^xb \tag{8.1.8}$$

Согласно теореме A.1.3, равенство

$$\frac{dy}{dx}\circ 1=\frac{dae^xb}{dx}\circ 1=a\left(\frac{de^x}{dx}\circ 1\right)b=ae^xb=y \tag{8.1.9}$$

является следствием равенства (8.1.8). Так как $e^0=1$, то теорема является следствием равенства (8.1.9). □

Согласно теореме 8.1.3 решение дифференциального уравнения (8.1.5) при начальном условии

$$x=0\quad y=C$$

определено неоднозначно. Это связано в частности с равенством

$$\frac{dy}{dx}\circ 1=\frac{\partial y}{\partial x^0}$$

Тем не менее, это дифференциальное уравнение очень важно.

## 8.2. Свойства экспоненты

Теорема 8.2.1. *Равенство*

$$e^{a+b}=e^ae^b \tag{8.2.1}$$

*справедливо тогда и только тогда, когда*

$$ab=ba \tag{8.2.2}$$

Доказательство. Для доказательства теоремы достаточно рассмотреть ряды Тейлора

$$e^a=\sum_{n=0}^{\infty}\frac{1}{n!}a^n \tag{8.2.3}$$

$$e^b=\sum_{n=0}^{\infty}\frac{1}{n!}b^n \tag{8.2.4}$$

$$e^{a+b}=\sum_{n=0}^{\infty}\frac{1}{n!}(a+b)^n \tag{8.2.5}$$

Перемножим выражения (8.2.3) и (8.2.4). Сумма одночленов порядка 3 имеет вид

$$\frac{1}{6}a^3+\frac{1}{2}a^2b+\frac{1}{2}ab^2+\frac{1}{6}b^3 \tag{8.2.6}$$

и не совпадает, вообще говоря, с выражением

$$\frac{1}{6}(a+b)^3=\frac{1}{6}a^3+\frac{1}{6}a^2b+\frac{1}{6}aba+\frac{1}{6}ba^2+\frac{1}{6}ab^2+\frac{1}{6}bab+\frac{1}{6}b^2a+\frac{1}{6}b^3 \tag{8.2.7}$$

Доказательство утверждения, что (8.2.1) следует из (8.2.2) тривиально. □

Теорема 8.2.2. *Пусть $A$ - ассоциативная $D$-алгебра. Тогда*

$$ae^{xa} = e^{ax}a \tag{8.2.8}$$

Доказательство. Равенство

$$\begin{aligned} e^{ax}a &= (1 + ax + \frac{1}{2}axax + \frac{1}{3!}axaxax + \frac{1}{4!}axaxaxax + ...)a \\ &= a + axa + \frac{1}{2}axaxa + \frac{1}{3!}axaxaxa + \frac{1}{4!}axaxaxaxa + ... \\ &= a(1 + xa + \frac{1}{2}xaxa + \frac{1}{3!}xaxaxa + \frac{1}{4!}xaxaxaxa + ...) = ae^{xa} \end{aligned} \tag{8.2.9}$$

является следствием равенства

$$e^{ax} = 1 + ax + \frac{1}{2}axax + \frac{1}{3!}axaxax + \frac{1}{4!}axaxaxax + ... \tag{8.2.10}$$

Равенство (8.2.8) является следствием равенства (8.2.9). □

Теорема 8.2.3. *Пусть $A$ - ассоциативная $D$-алгебра. Если $A$-число $a$ имеет обратный элемент, то*

$$e^{bxa} = a^{-1}e^{abx}a \tag{8.2.11}$$

*В частности*

$$e^{xa} = a^{-1}e^{ax}a \tag{8.2.12}$$

Доказательство. Равенство (8.2.12) является следствием равенства (8.2.8). Равенство (8.2.11) является следствием равенства (8.2.12), если положить $bx$ вместо $x$. □

## 8.3. Квазиэкспонента

Экспонента не является единственным решением дифференциального уравнения (8.1.5).

Определение 8.3.1. *Отображение*

$$e[c_1, ..., c_n]^x = \frac{d^n e^x}{dx} \circ (c_1, ..., c_n) \tag{8.3.1}$$

*называется* **квазиэкспонентой**. □

Теорема 8.3.2. *Квазиэкспонента $y = e[c]^x$ имеет следующее разложение в ряд Тейлора*

$$e[c]^x = \sum_{n=0}^{\infty} \frac{1}{(n+1)!} \sum_{m=0}^{n} x^m c x^{n-m} = c + \frac{1}{2}(cx + xc) + \frac{1}{6}(cx^2 + xcx + x^2c) + ...$$

Доказательство. Теорема является определению 8.3.1 и теоремы 8.1.2. □

Пример 8.3.3. *Положим $y = e[i]^x$. Согласно теореме 8.3.2, $(c = i)$,*

$$e[i]^x = \sum_{n=0}^{\infty} \frac{1}{(n+1)!} \sum_{m=0}^{n} x^m i x^{n-m} = i + \frac{1}{2}(ix + xi) + \frac{1}{6}(ix^2 + xix + x^2i) + ...$$

□

Теорема 8.3.4. *Если $c \in \ker A$, то*

$$e[c]^x = ce^x \tag{8.3.2}$$

Доказательство. Согласно определению 8.3.1 и теореме 8.1.2

$$\begin{aligned} e[c]^x &= \sum_{n=0}^{\infty} \frac{1}{(n+1)!} \sum_{m=0}^{n} x^m c x^{m-i} = \sum_{n=0}^{\infty} \frac{1}{(n+1)!} \sum_{m=0}^{n} c x^n \\ &= c \sum_{n=0}^{\infty} \frac{1}{n!} x^n = ce^x \end{aligned} \tag{8.3.3}$$

Теорема является следствием равенства (8.3.3). □

Теорема 8.3.5.

$$\frac{d^m e[c_1, ..., c_n]^x}{dx^m} \circ (C_1, ..., C_m) = e[c_1, ..., c_n, C_1, ..., C_m]^x \tag{8.3.4}$$

Доказательство. Равенство

$$\begin{aligned} \frac{d^m e[c_1, ..., c_n]^x}{dx^m} \circ (C_1, ..., C_m) &= \frac{d^m}{dx^m} \left( \frac{d^n e^x}{dx^n} \circ (c_1, ..., c_n) \right) \circ (C_1, ..., C_m) \\ &= \frac{d^{n+m} e^x}{dx^{n+m}} \circ (c_1, ..., c_n, C_1, ..., C_m) \\ &= e[c_1, ..., c_n, C_1, ..., C_m]^x \end{aligned} \tag{8.3.5}$$

является следствием равенств (8.3.1), (A.1.19). □

Теорема 8.3.6. *Квазиэкспонента является решением дифференциального уравнения*

$$\frac{dy}{dx} \circ 1 = y \tag{8.3.6}$$

Доказательство. Равенство

$$\frac{de[c_1, ..., c_n]^x}{dx} \circ 1 = \frac{d}{dx} \left( \frac{d^n e^x}{dx^n} \circ (c_1, ..., c_n) \right) \circ 1 = \frac{d^{n+1} e^x}{dx^{n+1}} \circ (c_1, ..., c_n, 1) \tag{8.3.7}$$

является следствием равенства (8.3.1) и определения [15]-4.1.4. Согласно теореме 4.2.8 и определению [15]-4.1.4, равенство

$$\frac{de[c_1, ..., c_n]^x}{dx} \circ 1 = \frac{d^{n+1} e^x}{dx^{n+1}} \circ (1, c_1, ..., c_n) = \frac{d^n}{dx^n} \left( \frac{de^x}{dx} \circ 1 \right) \circ (c_1, ..., c_n) \tag{8.3.8}$$

является следствием равенства (8.3.7). Теорема является следствием равенств (8.1.7), (8.3.1), (8.3.8). □

## 8.4. Дифференциальное уравнение $\frac{dy}{dx} \circ 1 = y$

Следующие вопросы являются следствием теоремы 8.3.6.

Вопрос 8.4.1. *Почему дифференциальное уравнение (8.3.6) имеет бесконечно много решений для заданного начального значения?* □

Вопрос 8.4.2. *Почему для нас интересно дифференциальное уравнение (8.3.6), если оно имеет бесконечно много решений для заданного начального значения?* □

Вопрос 8.4.3. *Существует ли аналог в некоммутативной алгебре для линейного однородного дифференциального уравнения?* □

**8.4.1. Ответ на вопрос 8.4.1.** Теорема [15]-7.3.7 даёт ответ на вопрос 8.4.1. Пусть $\overline{\overline{e}}$ - базис $D$-алгебры $A$. Согласно этой теореме выражение[8.1]

$$\frac{\partial y}{\partial x^0}=\frac{dy}{dx}\circ e_0 \tag{8.4.1}$$

является частной производной по переменной $x^0$.

Теорема 8.4.4. *Для любого* $n>0$

$$\frac{d^ny}{dx^n}\circ(e_0,...,e_0)=\frac{\partial^ny}{(\partial x^0)^n}$$

Доказательство. Мы докажем теорему индукцией по $n$.

Согласно равенству (8.4.1), теорема верна для $n=1$.

Допустим теорема верна для $n=k$

$$\frac{d^ky}{dx^k}\circ(e_0,...,e_0)=\frac{\partial^ky}{(\partial x^0)^k} \tag{8.4.2}$$

Согласно определению [15]-4.1.4 и теореме A.1.5, равенство

$$\begin{aligned}\frac{d^{k+1}y}{dx^{k+1}}\circ(e_0,...,e_0,e_0)&=\frac{d}{dx}\left(\frac{d^ky}{dx^k}\circ(e_0,...,e_0)\right)\circ e_0\\&=\frac{\partial}{\partial x^0}\frac{\partial^ky}{(\partial x^0)^k}=\frac{\partial^{k+1}y}{(\partial x^0)^{k+1}}\end{aligned}$$

является следствием равенств (8.4.1), (8.4.2). Следовательно, теорема верна для $n=k+1$. □

Согласно теореме 8.4.4, решение для дифференциального уравнения (8.3.6) при заданном начальном условии определено однозначно для $A$-чисел $de_0$, $d\in D$. Однако частная производная по $x^i$, $i\neq 0$, не определена, и это порождает неоднозначность решения.

**8.4.2. Ответ на вопрос 8.4.2.** Вопрос 8.4.2 выглядит естественным в свете ответа в разделе 8.4.1. Вначале я воспринял квазиэкспоненту как незваного гостя. Однако анализ определения 8.3.1 и теорем 8.1.3, 8.3.6 показывает, что это отображение очень важно. А именно, дифференциал квазиэкспоненты является квазиэкспонентой.

**8.4.3. Ответ на вопрос 8.4.3.** Согласно теореме 8.1.3, отображение $y=e^x$ является решением дифференциального уравнения

$$\frac{dy}{dx}\circ 1=y \tag{8.4.3}$$

Что я могу сказать об отображении $y=e^{ax}$?

Согласно теореме [15]-3.3.23

$$\frac{de^{ax}}{dx}\circ 1=\frac{de^{ax}}{dax}\circ\frac{dax}{dx}\circ 1=\frac{de^{ax}}{dax}\circ(a\otimes 1)\circ 1=\frac{de^{ax}}{dax}\circ a \tag{8.4.4}$$

Равенство

$$\frac{de^{ax}}{dx}\circ 1=e[a]^{ax} \tag{8.4.5}$$

является следствием равенств (8.3.1), (8.4.4). Следовательно, отображение $y=e^{ax}$ не является решением дифференциального уравнения (8.4.3).

[8.1] Согласно соглашению 1.5.3, $e_0=1$.

Если $D$-алгебра $A$ является алгеброй с делением, то мы можем записать равенство

$$\frac{de^{ax}}{dx}\circ a^{-1}=\frac{de^{ax}}{dax}\circ\frac{dax}{dx}\circ a^{-1}=\frac{de^{ax}}{dax}\circ 1=e^{ax} \tag{8.4.6}$$

вместо равенства (8.4.4).

Даже если я буду рассматривать дифференциальное уравнение

$$\frac{dy}{dx}\circ a^{-1}=y \tag{8.4.7}$$

непонятно, как перейти от этого дифференциального уравнения к дифференциальному уравнению более высокого порядка или системе однородных дифференциальных уравнений, так как в дифференциальном уравнении (8.4.7) нет коэффициента перед $y$.

Однако я не готов дать отрицательный ответ на вопрос 8.4.3.

## 8.5. Экспонента матрицы

Поскольку на множестве матриц определено два закона произведения ($a_*{}^*b$ и $a^*{}_*b$), мы должны рассмотреть два вида экспоненты.[8.2]

Теорема 8.5.1. ${}_*{}^*$**-экспонента** $y=e^{x_*{}^*}$ *имеет следующее разложение в ряд Тейлора*

$$e^{x_*{}^*}=\sum_{n=0}^{\infty}\frac{1}{n!}x^{n_*{}^*} \tag{8.5.1}$$

*и является решением дифференциального уравнения*

$$\frac{dy}{dx}{}_\circ{}^\circ E=y$$

*при начальном условии* $y(0)=E$.

Доказательство. Теорема является следствием теоремы 8.1.3, если в качестве умножения рассмотреть операцию $a_*{}^*b$. □

Теорема 8.5.2. ${}^*{}_*$**-экспонента** $y=e^{x^*{}_*}$ *имеет следующее разложение в ряд Тейлора*

$$e^{x^*{}_*}=\sum_{n=0}^{\infty}\frac{1}{n!}x^{n^*{}_*} \tag{8.5.2}$$

*и является решением дифференциального уравнения*

$$\frac{dy}{dx}{}^\circ{}_\circ E=y$$

*при начальном условии* $y(0)=E$.

Доказательство. Теорема является следствием теоремы 8.1.3, если в качестве умножения рассмотреть операцию $a^*{}_*b$. □

[8.2] Этот раздел я написал под влиянием статьи [16].

Глава 9

# Дифференциальное уравнение функции вещественного переменного

## 9.1. Дифференциальное уравнение $\frac{dx}{dt} = ax$

Теорема 9.1.1. *Пусть $A$ - банахова $D$-алгебра и $a \in A$. Отображение*

$$f : t \in R \to e^{at} \in A$$

*имеет разложение в ряд Тейлора*

$$e^{at} = \sum_{n=0}^{\infty} \frac{1}{n!} a^n t^n \tag{9.1.1}$$

Доказательство. Согласно определению 8.1.1

$$e^{at} = \sum_{n=0}^{\infty} \frac{1}{n!} (at)^n \tag{9.1.2}$$

Лемма 9.1.2. *Для любого $n \geq 0$*

$$(at)^n = a^n t^n \tag{9.1.3}$$

Доказательство. Мы докажем лемму индукцией по $n$.

Лемма верна для $n = 0$, так как

$$(at)^0 = 1 = a^0 t^0$$

Допустим лемма верна для $n = k$

$$(at)^k = a^k t^k \tag{9.1.4}$$

Тогда равенство

$$(at)^{k+1} = (at)^k (at) = (a^k t^k)(at) \tag{9.1.5}$$

является следствием равенства (9.1.4). Поскольку произведение в $D$-алгебре $A$ билинейно, равенство

$$(at)^{k+1} = (a^k a)(t^k t) = a^{k+1} t^{k+1} \tag{9.1.6}$$

является следствием равенства (9.1.5). Из равенства (9.1.6) следует, что лемма верна для $n = k + 1$. ⊙

Равенство (9.1.1) является следствием равенств (9.1.2), (9.1.3). □

Теорема 9.1.3. *Пусть $A$ - банахова ассоциативная $D$-алгебра и $a, c \in A$. Из условия*

$$ca = ac \tag{9.1.7}$$

*следует*

$$e^{at} c = c e^{at} \tag{9.1.8}$$

Доказательство. Пусть равенство (9.1.7) верно. Согласно теореме 9.1.1, для доказательства равенства (9.1.8) достаточно доказать равенство

$$a^n c = c a^n \tag{9.1.9}$$

Мы докажем равенство (9.1.9) индукцией по $n$.

Для $n = 0$ равенство (9.1.9) очевидно так как $a^0 = 1$. Для $n = 1$ равенство (9.1.9) является равенством (9.1.7).

Пусть равенство (9.1.9) верно для $n = k$

$$a^k c = c a^k \tag{9.1.10}$$

Равенство

$$a^{k+1}c = aa^k c = aca^k = caa^k = ca^{k+1}$$

является следствием равенств (9.1.7), (9.1.10). Следовательно, равенство (9.1.9) верно для $n = k + 1$. □

Теорема 9.1.4. *Пусть $A$ - банахова $D$-алгебра и $a \in A$. Тогда*

$$e^{at}a = ae^{at} \tag{9.1.11}$$

Доказательство. Теорема является следствием теоремы 9.1.3, если мы положим $c = a$. □

Теорема 9.1.5. *Пусть $A$ - некоммутативная $D$-алгебра. Для любого $b \in A$ существует подмодуль $A[b]$ $D$-алгебры $A$ такой, что*

$$c \in A[b] \Rightarrow cb = bc$$

Доказательство. Подмножество $A[b]$ не пусто так как $0 \in A[b]$. Пусть $d_1, d_2 \in D$, $c_1, c_2 \in A[b]$. Тогда

$$\begin{aligned}(d_1c_1 + d_2c_2)b &= d_1c_1b + d_2c_2b = bd_1c_1 + bd_2c_2\\ &= b(d_1c_1 + d_2c_2)\end{aligned}$$

Откуда следует $d_1c_1 + d_2c_2 \in A[b]$. □

Теорема 9.1.6. *Пусть $A$ - банахова $D$-алгебра и $a \in A$. Производная отображения*

$$f : t \in R \to e^{at} \in A$$

*имеет вид*

$$\frac{de^{at}}{dt} = e^{at}a = ae^{at} \tag{9.1.12}$$

Доказательство. Согласно теоремам 8.1.2, [15]-3.3.23,

$$\frac{de^{at}}{dt} = \frac{de^{at}}{dat} \circ \frac{dat}{dt} = \sum_{n=0}^{\infty} \frac{1}{(n+1)!} \sum_{i=0}^{n} ((at)^i \otimes (at)^{n-i}) \circ a \tag{9.1.13}$$

Согласно лемме 9.1.2, равенство

$$\begin{aligned}\frac{de^{at}}{dt} &= \sum_{n=0}^{\infty} \frac{1}{(n+1)!} \sum_{i=0}^{n} a^i t^i a a^{n-i} t^{n-i} = a \sum_{n=0}^{\infty} \frac{1}{(n+1)!} \sum_{i=0}^{n} a^n t^n\\ &= a \sum_{n=0}^{\infty} \frac{1}{n!} a^n t^n\end{aligned} \tag{9.1.14}$$

является следствием равенства (9.1.13). Равенство (9.1.12) является следствием равенств (9.1.1), (9.1.14). □

Теорема 9.1.7. *Пусть $A$ - банахова $D$-алгебра и $a\in A$. Для любого $A$-числа $c$ отображение*

$$x=e^{at}c \tag{9.1.15}$$

*является решением дифференциального уравнения*

$$\frac{dx}{dt}=ax \tag{9.1.16}$$

Доказательство. Равенство

$$\frac{dx}{dt}=\frac{de^{at}c}{dt}=\frac{de^{at}}{dt}c=ae^{at}c=ax$$

является следствием теорем A.1.3, 9.1.6. □

Теорема 9.1.8. *Пусть $A$ - банахова $D$-алгебра и $a\in A$. Для любых $A$-чисел $c_1$, $c_2$, отображение*

$$x=c_1e^{at}c_2 \tag{9.1.17}$$

*является решением дифференциального уравнения* (9.1.16) *тогда и только тогда, когда $c_1\in A[a]$ и, следовательно, отображение* (9.1.17) *является отображением вида* (9.1.15).

Доказательство. Равенство

$$\frac{dx}{dt}=\frac{dc_1e^{at}c_2}{dt}=c_1\frac{de^{at}}{dt}c_2=c_1ae^{at}c_2$$

является следствием теорем A.1.3, 9.1.6. Если $c_1\notin A[a]$, то условие

$$c_1a=ac_1 \tag{9.1.18}$$

не верно и отображение (9.1.17) не является решением дифференциального уравнения (9.1.16). Однако если условие (9.1.18) верно, то отображение (9.1.17) принимает вид

$$x=c_1e^{at}c_2=e^{at}c_1c_2 \tag{9.1.19}$$

и является отображением вида (9.1.15). □

Теорема 9.1.9. *Пусть $A$ - банахова $D$-алгебра и $a\in A$. Для любого $A$-числа $c$ отображение*

$$x=ce^{at} \tag{9.1.20}$$

*является решением дифференциального уравнения*

$$\frac{dx}{dt}=xa \tag{9.1.21}$$

Доказательство. Равенство

$$\frac{dx}{dt}=\frac{dce^{at}}{dt}=c\frac{de^{at}}{dt}=ce^{at}a=xa$$

является следствием теорем A.1.3, 9.1.6. □

Теорема 9.1.10. *Пусть $A$ - банахова $D$-алгебра и $a \in A$. Для любых $A$-чисел $c_1$, $c_2$, отображение*

$$x = c_1 e^{at} c_2 \tag{9.1.22}$$

*является решением дифференциального уравнения (9.1.21) тогда и только тогда, когда $c_2 \in A[a]$ и, следовательно, отображение (9.1.22) является отображением вида (9.1.20).*

Доказательство. Равенство

$$\frac{dx}{dt} = \frac{dc_1 e^{at} c_2}{dt} = c_1 \frac{de^{at}}{dt} c_2 = c_1 e^{at} a c_2$$

является следствием теорем A.1.3, 9.1.6. Если $c_2 \notin A[a]$, то условие

$$c_2 a = a c_2 \tag{9.1.23}$$

не верно и отображение (9.1.22) не является решением дифференциального уравнения (9.1.21). Однако если условие (9.1.23) верно, то отображение (9.1.22) принимает вид

$$x = c_1 e^{at} c_2 = c_1 c_2 e^{at} \tag{9.1.24}$$

и является отображением вида (9.1.20). □

## 9.2. Квазиэкспонента

Теорема 9.2.1. *Отображение*

$$x = e[c]^{at}$$

*имеет следующее разложение в ряд Тейлора*

$$e[c]^{at} = \sum_{n=0}^{\infty} \frac{t^n}{(n+1)!} \sum_{m=0}^{n} a^m c a^{n-m} \tag{9.2.1}$$

Доказательство. Теорема является следствием теоремы 8.3.2 и леммы 9.1.2. □

Теорема 9.2.2. *Пусть $c \in A[a]$. Тогда*

$$e[c]^{at} = c e^{at} \tag{9.2.2}$$

Доказательство. Теорема является следствием теоремы 9.2.1. □

Теорема 9.2.3.

$$\frac{de[c]^{at}}{dt} = \sum_{n=0}^{\infty} \frac{t^n}{(n+2)!} \sum_{m=0}^{n+1} a^m c a^{n+1-m} \tag{9.2.3}$$

Доказательство. Теорема является следствием теоремы 9.2.1. □

## 9.3. Дифференциальное уравнение $\frac{dx}{dt} = a_*{}^*x$

Пусть $A$ - банахова $D$-алгебра. Система дифференциальных уравнений

$$\begin{array}{l} \frac{dx^1}{dt} = a^1_1 x^1 + ... + a^1_n x^n \\ ...... \\ \frac{dx^n}{dt} = a^n_1 x^1 + ... + a^n_n x^n \end{array} \tag{9.3.1}$$

где $a^i_j \in A$ и $x^i : R \to A$ - $A$-значная функция действительной переменной, называется однородной системой линейных дифференциальных уравнений.

Положим

$$x = \begin{pmatrix} x^1 \\ ... \\ x^n \end{pmatrix} \qquad \frac{dx}{dt} = \begin{pmatrix} \frac{dx^1}{dt} \\ ... \\ \frac{dx^n}{dt} \end{pmatrix}$$

$$a = \begin{pmatrix} a^1_1 & ... & a^1_n \\ ... & ... & ... \\ a^n_1 & ... & a^n_n \end{pmatrix}$$

Тогда мы можем записать систему дифференциальных уравнений (9.3.1) в матричной форме

$$\frac{dx}{dt} = a_*{}^*x \tag{9.3.2}$$

**9.3.1. Решение в виде экспоненты $x = e^{bt}c$ .** Мы будем искать решение системы дифференциальных уравнений (9.3.2) в виде экспоненты

$$x = e^{bt}c = \begin{pmatrix} e^{bt}c^1 \\ ... \\ e^{bt}c^n \end{pmatrix} \qquad c = \begin{pmatrix} c^1 \\ ... \\ c^n \end{pmatrix} \tag{9.3.3}$$

Согласно теореме 9.1.6, равенство

$$\frac{dx^i}{dt} = e^{bt}bc^i = be^{bt}c^i = bx^i \tag{9.3.4}$$

является следствием равенства (9.3.3). Равенство

$$e^{bt}bc = be^{bt}c = a_*{}^*(e^{bt}c) \tag{9.3.5}$$

является следствием равенств (9.3.2), (9.3.4).

Теорема 9.3.1. *Пусть $b$ является ${}_*{}^*$-собственным значением матрицы $a$. Из условия*

$$a^i_j \in A[b] \quad i = 1, ..., n \quad j = 1, ..., n \tag{9.3.6}$$

*следует, что матрица отображений (9.3.3) является решением системы дифференциальных уравнений (9.3.2) для ${}_*{}^*$-собственного столбца $c$.*

Доказательство. Согласно теореме 9.1.3, равенство

$$e^{bt}bc = e^{bt}(a_*{}^*c) \tag{9.3.7}$$

является следствием равенства (9.3.5). Поскольку выражение $e^{bt}$, вообще говоря, отлично от 0, равенство

$$bc = a_*{}^*c$$

является следствием равенства (9.3.7). Согласно определению 5.6.1, $A$-число $b$ является ${}_*{}^*$-собственным значением матрицы $a$ и матрица $c$ является ${}_*{}^*$-собственным столбцом матрицы $a$, соответствующим ${}^*{}_*$-собственному значению $b$. □

Теорема 9.3.2. *Пусть $b$ является ${}_*{}^*$-собственным значением матрицы $a$. Пусть элементы матрицы $a$ не удовлетворяют условию (9.3.6). Если элементы ${}_*{}^*$-собственного столбца $c$ удовлетворяют условию*

$$c^i \in A[b] \qquad i = 1, ..., n \tag{9.3.8}$$

*следует, что матрица отображений (9.3.3) является решением системы дифференциальных уравнений (9.3.2).*

Доказательство. Согласно теореме 9.1.3, равенство

$$bce^{bt} = (a_*{}^*c)e^{bt} \tag{9.3.9}$$

является следствием равенства (9.3.5). Поскольку выражение $e^{bt}$, вообще говоря, отлично от 0, равенство

$$bc = a_*{}^*c$$

является следствием равенства (9.3.9). Согласно определению 5.6.1, $A$-число $b$ является ${}_*{}^*$-собственным значением матрицы $a$ и матрица $c$ является ${}_*{}^*$-собственным столбцом матрицы $a$, соответствующим ${}^*{}_*$-собственному значению $b$. □

Пусть элементы матрицы $a$ не удовлетворяют условию (9.3.6). Пусть элементы ${}_*{}^*$-собственного столбца $c$ не удовлетворяют условию (9.3.8). Тогда матрица отображений (9.3.3) не является решением системы дифференциальных уравнений (9.3.2).

**9.3.2. Решение в виде экспоненты $x = ce^{bt}$.** В статье [17], на странице 35, авторы предлагают рассмотреть решение в виде экспоненты

$$x = ce^{bt} \tag{9.3.10}$$

Равенство

$$a_*{}^*(ce^{bt}) = cbe^{bt} \tag{9.3.11}$$

является следствием равенств (9.3.2), (9.3.10). Поскольку выражение $e^{bt}$, вообще говоря, отлично от 0, равенство

$$a_*{}^*c = cb \tag{9.3.12}$$

является следствием равенства (9.3.11).

На основании равенства (9.3.12) авторы статей [16], [17] вводят определение правого собственного значения, чтобы отличить от левого собственного значения, которое мы определяем равенством

$$bc = a_*{}^*c$$

Согласно теореме 9.1.8, если матрица отображений (9.3.10) является решением системы дифференциальных уравнений (9.3.2), то вектор $c$ удовлетворяет условию

$$c^i \in A[b] \qquad i = 1, ..., n$$

В этом случае, мы можем рассматривать матрицу отображений

$$x = e^{bt}c$$

вместо матрицы отображений (9.3.10) и нам не надо рассматривать определение правого собственного значения.

### 9.3.3. Метод последовательного дифференцирования.

Теорема 9.3.3. *Последовательно дифференцируя систему дифференциальных уравнений (9.3.2) мы получаем цепочку систем дифференциальных уравнений*

$$\frac{d^n x}{dt^n} = a^{n_*{}^*}{}_*{}^*x \tag{9.3.13}$$

Доказательство. Мы докажем теорему индукцией по $n$.

Согласно равенствам (3.2.17), (9.3.2), теорема верна для $n = 1$.

Допустим теорема верна для $n = k$

$$\frac{d^k x}{dt^k} = a^{k_*{}^*}{}_*{}^*x \tag{9.3.14}$$

Согласно определению [15]-4.1.4, равенство

$$\frac{d^{k+1}x}{dt^{k+1}} = \frac{d}{dt}\frac{d^k x}{dt^k} = \frac{d}{dt}(a^{k_*{}^*}{}_*{}^*x) = a^{k_*{}^*}{}_*{}^*\frac{dx}{dt} = a^{k_*{}^*}{}_*{}^*a_*{}^*x = a^{k+1_*{}^*}{}_*{}^*x$$

является следствием равенств (3.2.16), (9.3.2), (9.3.14). Следовательно, теорема верна для $n = k + 1$. □

Теорема 9.3.4. *Решение системы дифференциальных уравнений (9.3.2) при начальном условии*

$$t = 0 \qquad x = \begin{pmatrix} x^1 \\ ... \\ x^n \end{pmatrix} = c = \begin{pmatrix} c^1 \\ ... \\ c^n \end{pmatrix}$$

*имеет вид*

$$x = e^{ta_*{}^*}{}_*{}^*c \tag{9.3.15}$$

Доказательство. Согласно теореме 9.3.3 и определению ряда Тейлора в разделе [15]-4.2, ряд Тейлора для решения имеет вид

$$x = \sum_{n=0}^{\infty}\frac{1}{n!}t^n a^{n_*{}^*}{}_*{}^*c \tag{9.3.16}$$

Равенство

$$x = \sum_{n=0}^{\infty}\frac{1}{n!}(ta)^{n_*{}^*}{}_*{}^*c \tag{9.3.17}$$

является следствием равенства (9.3.16). Равенство (9.3.15) является следствием равенств (8.5.1), (9.3.17). □

## 9.4. Гиперболическая тригонометрия

Рассмотрим систему дифференциальных уравнений

(9.4.1)
$$\frac{dx^1}{dt} = x^2$$
$$\frac{dx^2}{dt} = x^1$$

Матрица $a$ имеет вид

$$a = \begin{pmatrix} 0 & 1 \\ 1 & 0 \end{pmatrix}$$

Поскольку коэффициенты матрицы $a$ - действительные числа, то уравнение для определения собственных значений имеет вид

(9.4.2)
$$\begin{vmatrix} -b & 1 \\ 1 & -b \end{vmatrix} = b^2 - 1 = 0$$

Очевидно, корни уравнения (9.4.2) зависят от выбора $D$-алгебры $A$.

В алгебре кватернионов, согласно утверждению [14]-5.5.1, уравнение (9.5.2) имеет корни $b_1 = -1$, $b_2 = 1$.

Коэффициенты $c_1^1$, $c_1^2$, которые соответствуют собственному значению $b_1 = -1$, удовлетворяют уравнению

$$c_1^1 + c_1^2 = 0$$

Следовательно, соответствующее решение системы дифференциальных уравнений (9.4.1) имеет вид

(9.4.3)
$$x^1 = e^{-t} \quad x^2 = -e^{-t}$$

Коэффициенты $c_2^1$, $c_2^2$, которые соответствуют собственному значению $b_2 = 1$, удовлетворяют уравнению

$$-c_2^1 + c_2^2 = 0$$

Следовательно, соответствующее решение системы дифференциальных уравнений (9.4.1) имеет вид

(9.4.4)
$$x^1 = e^t \quad x^2 = e^t$$

Мы можем записать решение системы дифференциальных уравнений (9.4.1) как линейную комбинацию решений (9.4.3), (9.4.4)

(9.4.5)
$$x^1 = C_1e^t + C_2e^{-t}$$
$$x^2 = C_1e^t - C_2e^{-t}$$

Решение системы дифференциальных уравнений (9.4.1), которое удовлетворяет начальному условию

$$t = 0 \quad x^1 = 0 \quad x^2 = 1$$

определено однозначно и имеет вид

(9.4.6)
$$x^1 = \frac{1}{2}e^t - \frac{1}{2}e^{-t}$$
$$x^2 = \frac{1}{2}e^t + \frac{1}{2}e^{-t}$$

Мы можем решить систему дифференциальных уравнений (9.4.1), пользуясь методом последовательного дифференцирования. С этой целью мы рассмотрим степени матрицы $a$

$$a^{1_*{}^*} = \begin{pmatrix} 0 & 1 \\ 1 & 0 \end{pmatrix} \quad a^{2_*{}^*} = \begin{pmatrix} 1 & 0 \\ 0 & 1 \end{pmatrix} \quad a^{3_*{}^*} = \begin{pmatrix} 0 & 1 \\ 1 & 0 \end{pmatrix} \tag{9.4.7}$$

Равенство

$$\begin{aligned} x = & \left( \begin{pmatrix} 1 & 0 \\ 0 & 1 \end{pmatrix} + t \begin{pmatrix} 0 & 1 \\ 1 & 0 \end{pmatrix} \right. \\ & \left. + \frac{1}{2!} t^2 \begin{pmatrix} 1 & 0 \\ 0 & 1 \end{pmatrix} + \frac{1}{3!} t^3 \begin{pmatrix} 0 & 1 \\ 1 & 0 \end{pmatrix} + \ldots \right) {}_*{}^* \begin{pmatrix} 0 \\ 1 \end{pmatrix} \end{aligned} \tag{9.4.8}$$

является следствием равенств (9.3.16), (9.4.7). Равенство

$$\begin{aligned} x^{0} &= \sum_{n=0}^{\infty} \frac{1}{(2n+1)!} t^{2n+1} = \sinh t \\ x^{1} &= \sum_{n=0}^{\infty} \frac{1}{(2n)!} t^{2n} = \cosh t \end{aligned} \tag{9.4.9}$$

является следствием равенства (9.4.8). Из равенств (9.4.6), (9.4.9) следует, что формула Эйлера

$$\begin{aligned} \sinh t &= \frac{1}{2} e^t - \frac{1}{2} e^{-t} \\ \cosh t &= \frac{1}{2} e^t + \frac{1}{2} e^{-t} \end{aligned} \tag{9.4.10}$$

верна для $t \in R$. Равенство

$$\begin{aligned} \frac{d \sinh t}{dt} &= \cosh t \\ \frac{d \cosh t}{dt} &= \sinh t \end{aligned} \tag{9.4.11}$$

является следствием равенств (9.4.1), (9.4.9).

Чтобы понять, верна ли эта формула для произвольного кватерниона $f \neq 0$, рассмотрим теперь систему дифференциальных уравнений

$$\begin{aligned} \frac{dx^{1}}{dt} &= f x^{2} \\ \frac{dx^{2}}{dt} &= f x^{1} \end{aligned} \tag{9.4.12}$$

Матрица $a$ имеет вид

$$a = \begin{pmatrix} 0 & f \\ f & 0 \end{pmatrix}$$

Теперь элементы матрицы $a$ - кватернионы, а не действительные числа. Поэтому мы не можем считать определитель. Для того, чтобы ответить на вопрос, когда матрица $a - Eb$ вырождена, мы посчитаем квазидетерминанты первой строки[9.1]

$$\det({}_*{}^*)^1_1\, a = -b + fb^{-1}f = 0 \tag{9.4.13}$$

$$\det({}_*{}^*)^1_2\, a = f - bf^{-1}b = 0 \tag{9.4.14}$$

Нетрудно видеть, что уравнения (9.4.13) и (9.4.14) эквивалентны. Уравнение (9.4.14) - это квадратное уравнение. В уравнении (9.4.14) я запишу выражение $f + g$ вместо неизвестной величины $b$ и буду решать уравнение

$$f = (f+g)f^{-1}(f+g) = (f+g)(1+f^{-1}g) = f + g + ff^{-1}g + gf^{-1}g \tag{9.4.15}$$

относительно неизвестной величины $g$. Уравнение

$$0 = 2g + gf^{-1}g \tag{9.4.16}$$

является следствием уравнения (9.4.15). Уравнение

$$g(2 + f^{-1}g) = 0 \tag{9.4.17}$$

является следствием уравнения (9.4.16). Из уравнения (9.4.17) следует, что $g = 0$ и $g = -2f$ - корни уравнения (9.4.16). Следовательно, $b = f$ и $b = -f$ корни уравнения (9.4.14).

Коэффициенты $c_1^1$, $c_1^2$, которые соответствуют собственному значению $b_1 = -f$, удовлетворяют уравнению

$$c_1^1 + c_1^2 = 0$$

Следовательно, соответствующее решение системы дифференциальных уравнений (9.4.12) имеет вид

$$x^1 = e^{-ft} \qquad x^2 = -e^{-ft} \tag{9.4.18}$$

Коэффициенты $c_2^1$, $c_2^2$, которые соответствуют собственному значению $b_2 = f$, удовлетворяют уравнению

$$-c_2^1 + c_2^2 = 0$$

Следовательно, соответствующее решение системы дифференциальных уравнений (9.4.12) имеет вид

$$x^1 = e^{ft} \qquad x^2 = e^{ft} \tag{9.4.19}$$

Мы можем записать решение системы дифференциальных уравнений (9.4.12) как линейную комбинацию решений (9.4.18), (9.4.19)

$$\begin{aligned} x^1 &= C_1e^{ft} + C_2e^{-ft} \\ x^2 &= C_1e^{ft} - C_2e^{-ft} \end{aligned} \tag{9.4.20}$$

Решение системы дифференциальных уравнений (9.4.12), которое удовлетворяет начальному условию

$$t = 0 \quad x^1 = 0 \quad x^2 = 1$$

---

[9.1] Смотри определение 3.3.3 для квазидетерминанта. Смотри теорему 5.5.3 для определения ранга матрицы.

определено однозначно и имеет вид

$$x^{1} = \frac{1}{2}e^{ft} - \frac{1}{2}e^{-ft} \qquad (9.4.21)$$
$$x^{2} = \frac{1}{2}e^{ft} + \frac{1}{2}e^{-ft}$$

Мы можем решить систему дифференциальных уравнений (9.4.12), пользуясь методом последовательного дифференцирования. С этой целью мы рассмотрим степени матрицы $a$

$$a^{1_*{}^*} = \begin{pmatrix} 0 & f \\ f & 0 \end{pmatrix} \quad a^{2_*{}^*} = \begin{pmatrix} f^2 & 0 \\ 0 & f^2 \end{pmatrix} \quad a^{3_*{}^*} = \begin{pmatrix} 0 & f^3 \\ f^3 & 0 \end{pmatrix} \qquad (9.4.22)$$

Равенство

$$x = \left(\begin{pmatrix} 1 & 0 \\ 0 & 1 \end{pmatrix} + t\begin{pmatrix} 0 & f \\ f & 0 \end{pmatrix} + \frac{1}{2!}t^2\begin{pmatrix} f^2 & 0 \\ 0 & f^2 \end{pmatrix} + \frac{1}{3!}t^3\begin{pmatrix} 0 & f^3 \\ f^3 & 0 \end{pmatrix} + ...\right){}_*{}^*\begin{pmatrix} 0 \\ 1 \end{pmatrix} \qquad (9.4.23)$$

является следствием равенств (9.3.16), (9.4.22). Равенство

$$x^{0} = \sum_{n=0}^{\infty}\frac{1}{(2n+1)!}t^{2n+1}f^{2n+1} = \sum_{n=0}^{\infty}\frac{1}{(2n+1)!}(tf)^{2n+1} = \sinh tf \qquad (9.4.24)$$
$$x^{1} = \sum_{n=0}^{\infty}\frac{1}{(2n)!}t^{2n}f^{2n} = \sum_{n=0}^{\infty}\frac{1}{(2n)!}(tf)^{2n} = \cosh tf$$

является следствием равенства (9.4.23). Из равенств (9.4.21), (9.4.24) следует, что формула Эйлера

$$\sinh ft = \frac{1}{2}e^{ft} - \frac{1}{2}e^{-ft} \qquad (9.4.25)$$
$$\cosh ft = \frac{1}{2}e^{ft} + \frac{1}{2}e^{-ft}$$

верна для $tf \in A$. Равенство

$$\frac{d\sinh tf}{dt} = f\cosh tf \qquad (9.4.26)$$
$$\frac{d\cosh tf}{dt} = f\sinh tf$$

является следствием равенств (9.4.12), (9.4.24). Пусть $t = 1$. Из равенства

$$\sinh f = \frac{1}{2}e^{f} - \frac{1}{2}e^{-f} \qquad (9.4.27)$$
$$\cosh f = \frac{1}{2}e^{f} + \frac{1}{2}e^{-f}$$

следует, что формула Эйлера верна для $f \in A$.

Теорема 9.4.1.

$$\sinh tf\, f = f\sinh tf \qquad (9.4.28)$$

$$\cosh tf\, f = f \cosh tf \tag{9.4.29}$$

Доказательство. Равенства

$$\begin{aligned}\sinh tf\, f &= \sum_{n=0}^{\infty} \frac{1}{(2n+1)!} t^{2n+1} f^{2n+1} f = f \sum_{n=0}^{\infty} \frac{1}{(2n+1)!} t^{2n+1} f^{2n+1}\\ &= f \sinh tf\end{aligned} \tag{9.4.30}$$

$$\begin{aligned}\cosh tf\, f &= \sum_{n=0}^{\infty} \frac{1}{(2n)!} t^{2n} f^{2n} f = f \sum_{n=0}^{\infty} \frac{1}{(2n)!} t^{2n} f^{2n}\\ &= f \cosh tf\end{aligned} \tag{9.4.31}$$

являются следствием равенства (9.4.24). Равенство (9.4.28) является следствием равенства (9.4.30). Равенство (9.4.29) является следствием равенства (9.4.31). □

## 9.5. Эллиптическая тригонометрия

Рассмотрим систему дифференциальных уравнений

$$\begin{aligned}\frac{dx^1}{dt} &= x^2\\ \frac{dx^2}{dt} &= -x^1\end{aligned} \tag{9.5.1}$$

Матрица $a$ имеет вид

$$a = \begin{pmatrix} 0 & 1\\ -1 & 0\end{pmatrix}$$

Поскольку коэффициенты матрицы $a$ - действительные числа, то уравнение для определения собственных значений имеет вид

$$\begin{vmatrix} -b & 1\\ -1 & -b\end{vmatrix} = b^2 + 1 = 0 \tag{9.5.2}$$

Очевидно, корни уравнения (9.5.2) зависят от выбора $D$-алгебры $A$.

В алгебре кватернионов, согласно утверждению [14]-5.5.3, уравнение (9.5.2) имеет бесконечно много корней

$$b = b^1 i + b^2 j + b^3 k$$

таких, что

$$(b^1)^2 + (b^2)^2 + (b^3)^2 = 1$$

Коэффициенты $c_b^1$, $c_b^2$, которые соответствуют данному собственному значению $b$, удовлетворяют уравнению

$$-bc_b^1 + c_b^2 = 0$$

Следовательно, соответствующее решение системы дифференциальных уравнений (9.5.1) имеет вид

$$x^1 = e^{bt} \qquad x^2 = be^{bt} \tag{9.5.3}$$

Если мы хотим найти решение системы дифференциальных уравнений (9.5.1), которое удовлетворяет начальному условию

$$t = 0 \quad x^1 = 0 \quad x^2 = 1$$

то мы видим, что у нас слишком большой выбор.

Линейная комбинация двух решений системы дифференциальных уравнений (9.5.1) является решением системы дифференциальных уравнений (9.5.1). Мы начнём с рассмотрения линейной комбинации двух решений вида (9.5.3) так как в этом случае константы линейной комбинации определены однозначно.

Таким образом, мы ищем решение в форме

$$\begin{aligned} x^1 &= C_1 e^{b_1 t} + C_2 e^{b_2 t} \\ x^2 &= C_1 b_1 e^{b_1 t} + C_2 b_2 e^{b_2 t} \end{aligned} \tag{9.5.4}$$

Согласно начальным условиям, система уравнений

$$\begin{aligned} C_1 + C_2 &= 0 \\ C_1 b_1 + C_2 b_2 &= 1 \end{aligned} \tag{9.5.5}$$

является следствием равенства (9.5.4). Равенство

$$\begin{aligned} C_2 &= -C_1 \\ C_1(b_1 - b_2) &= 1 \end{aligned} \tag{9.5.6}$$

является следствием системы уравнений (9.5.5).

Пример 9.5.1. *Положим* $b_1 = i$, $b_2 = j$. *Равенство*

$$C_1(i - j) = 1 \tag{9.5.7}$$

*является следствием системы уравнений* (9.5.6). *Равенство*

$$C_1 = \frac{1}{2}(-i + j) \tag{9.5.8}$$

*является следствием равенства* (9.5.7). *Наша задача - проверить является ли отображение*

$$\begin{aligned} x^1 &= \frac{1}{2}(-i + j)(e^{it} - e^{jt}) \\ x^2 &= \frac{1}{2}(-i + j)(ie^{it} - je^{jt}) \end{aligned} \tag{9.5.9}$$

*решением системы дифференциальных уравнений* (9.3.2). *Равенство*

$$\begin{aligned} \frac{dx^1}{dt} &= \frac{1}{2}(-i + j)(ie^{it} - je^{jt}) = x^2 \\ \frac{dx^2}{dt} &= \frac{1}{2}(-i + j)(i^2 e^{it} - j^2 e^{jt}) = -x^1 \end{aligned} \tag{9.5.10}$$

*является следствием равенства* (9.5.9). *Следовательно, отображение* (9.5.9) *является решением системы дифференциальных уравнений* (9.3.2), *которое удовлетворяет начальному условию*

$$t = 0 \quad x^1 = 0 \quad x^2 = 1$$

□

Если мы рассмотрим линейную комбинацию трёх или более решений вида (9.5.3), то коэффициенты линейной комбинации будут определены неоднозначно, так как соответствующая система линейных уравнений будет переопределена.

Пример 9.5.2. *Положим* $b_1 = i$, $b_2 = j$, $b_3 = k$. *Таким образом, мы ищем решение в форме*

$$\begin{aligned} x^1 &= C_1 e^{it} + C_2 e^{jt} + C_3 e^{kt} \\ x^2 &= C_1 i e^{it} + C_2 j e^{jt} + C_3 k e^{kt} \end{aligned} \tag{9.5.11}$$

*Согласно начальным условиям, система уравнений*

$$\begin{aligned} C_1 + C_2 + C_3 &= 0 \\ C_1 i + C_2 j + C_3 k &= 1 \end{aligned} \tag{9.5.12}$$

*является следствием равенства* (9.5.11). *Равенство*

$$\begin{aligned} &C_1 = C \\ &C_2 = -C - C_3 \\ &Ci - (C + C_3)j + C_3 k = 1 \end{aligned} \tag{9.5.13}$$

*является следствием системы уравнений* (9.5.12). *Равенство*

$$C_3(k - j) = 1 - Ci + Cj \tag{9.5.14}$$

*является следствием равенства* (9.5.13). *Равенство*

$$\begin{aligned} C_3 &= \frac{1}{2}(1 + C(j - i))(j - k) = \frac{1}{2}((j - k) + C(j - i)(j - k)) \\ &= \frac{1}{2}((j - k) + C(j(j - k) - i(j - k))) \end{aligned} \tag{9.5.15}$$

*является следствием равенства* (9.5.14). *Равенство*

$$\begin{aligned} &C_1 = C \\ &C_2 = \frac{1}{2}(-(j - k) + C(-1 + i + j + k)) \\ &C_3 = \frac{1}{2}((j - k) - C(1 + i + j + k)) \end{aligned} \tag{9.5.16}$$

*является следствием равенств* (9.5.13), (9.5.15). *Равенство*

$$\begin{aligned} x^1 &= Ce^{it} + \frac{1}{2}(-(j - k) + C(-1 + i + j + k))e^{jt} \\ &\quad + \frac{1}{2}((j - k) - C(1 + i + j + k))e^{kt} \\ x^2 &= Cie^{it} + \frac{1}{2}(-(j - k) + C(-1 + i + j + k))je^{jt} \\ &\quad + \frac{1}{2}((j - k) - C(1 + i + j + k))ke^{kt} \end{aligned} \tag{9.5.17}$$

*является следствием равенств* (9.5.11), (9.5.16). □

Вопрос 9.5.3. *Поскольку множество решений системы дифференциальных уравнений* (9.4.12) *слишком велико, я не рассматриваю формулу Эйлера. Однако необходимо понять, почему одна и та же система дифференциальных* *уравнений имеет бесконечно много решений, которые удовлетворяют заданному начальному условию, и в тоже время решение этой системы дифференциальных уравнений имеет единственное разложение в ряд Тейлора. Связано ли это с тождествами, которые я не рассматривал?* □

## 9.6. Дифференциальное уравнение $\frac{dx}{dt} = x^*{}_*a$

Пусть $A$ - банахова $D$-алгебра. Система дифференциальных уравнений

(9.6.1)
$$\begin{aligned}\frac{dx^1}{dt} &= x^1a_1^1 + ... + x^na_n^1\\ &......\\ \frac{dx^n}{dt} &= x^1a_1^n + ... + x^na_n^n\end{aligned}$$

где $a_j^i \in A$ и $x^i : R \to A$ - $A$-значная функция действительной переменной, называется однородной системой линейных дифференциальных уравнений.

Положим

$$x = \begin{pmatrix} x^1 \\ ... \\ x^n \end{pmatrix} \qquad \frac{dx}{dt} = \begin{pmatrix} \frac{dx^1}{dt} \\ ... \\ \frac{dx^n}{dt} \end{pmatrix}$$

$$a = \begin{pmatrix} a_1^1 & ... & a_n^1 \\ ... & ... & ... \\ a_1^n & ... & a_n^n \end{pmatrix}$$

Тогда мы можем записать систему дифференциальных уравнений (9.6.1) в матричной форме

(9.6.2)
$$\frac{dx}{dt} = x^*{}_*a$$

**9.6.1. Решение в виде экспоненты** $x = ce^{bt}$ **.** Мы будем искать решение системы дифференциальных уравнений (9.6.2) в виде экспоненты

(9.6.3)
$$x = ce^{bt} = \begin{pmatrix} c^1e^{bt} \\ ... \\ c^ne^{bt} \end{pmatrix} \qquad c = \begin{pmatrix} c^1 \\ ... \\ c^n \end{pmatrix}$$

Согласно теореме 9.1.6, равенство

(9.6.4)
$$\frac{dx^i}{dt} = c^ie^{bt}b = c^ibe^{bt} = x^ib$$

является следствием равенства (9.6.3). Равенство

(9.6.5)
$$ce^{bt}b = cbe^{bt} = (ce^{bt})^*{}_*a$$

является следствием равенств (9.6.2), (9.6.4).

Теорема 9.6.1. *Пусть $b$ является ${}^*{}_*$-собственным значением матрицы $a$. Из условия*

$$a_j^i \in A[b] \quad i = 1, ..., n \quad j = 1, ..., n \tag{9.6.6}$$

*следует, что матрица отображений (9.6.3) является решением системы дифференциальных уравнений (9.6.2) для ${}^*{}_*$-собственного столбца $c$.*

Доказательство. Согласно теореме 9.1.3, равенство

$$cbe^{bt} = (c^*{}_*a)e^{bt} \tag{9.6.7}$$

является следствием равенства (9.6.5). Поскольку выражение $e^{bt}$, вообще говоря, отлично от 0, равенство

$$cb = c^*{}_*a$$

является следствием равенства (9.6.7). Согласно определению 5.6.2, $A$-число $b$ является ${}^*{}_*$-собственным значением матрицы $a$ и матрица $c$ является ${}^*{}_*$-собственным столбцом матрицы $a$, соответствующим ${}^*{}_*$-собственному значению $b$. □

Теорема 9.6.2. *Пусть $b$ является ${}^*{}_*$-собственным значением матрицы $a$. Пусть элементы матрицы $a$ не удовлетворяют условию (9.6.6). Если элементы ${}^*{}_*$-собственного столбца $c$ удовлетворяют условию*

$$c^i \in A[b] \quad i = 1, ..., n \tag{9.6.8}$$

*следует, что матрица отображений (9.6.3) является решением системы дифференциальных уравнений (9.6.2).*

Доказательство. Согласно теореме 9.1.3, равенство

$$e^{bt}cb = e^{bt}(c^*{}_*a) \tag{9.6.9}$$

является следствием равенства (9.6.5). Поскольку выражение $e^{bt}$, вообще говоря, отлично от 0, равенство

$$cb = c^*{}_*a$$

является следствием равенства (9.6.9). Согласно определению 5.6.2, $A$-число $b$ является ${}^*{}_*$-собственным значением матрицы $a$ и матрица $c$ является ${}^*{}_*$-собственным столбцом матрицы $a$, соответствующим ${}_*{}^*$-собственному значению $b$. □

Пусть элементы матрицы $a$ не удовлетворяют условию (9.6.6). Пусть элементы ${}_*{}^*$-собственного столбца $c$ не удовлетворяют условию (9.6.8). Тогда матрица отображений (9.6.3) не является решением системы дифференциальных уравнений (9.6.2).

### 9.6.2. Метод последовательного дифференцирования.

Теорема 9.6.3. *Последовательно дифференцируя систему дифференциальных уравнений (9.6.2) мы получаем цепочку систем дифференциальных уравнений*

$$\frac{d^n x}{dt^n} = x^*{}_*a^{n^*{}_*} \tag{9.6.10}$$

Доказательство. Мы докажем теорему индукцией по $n$.

Согласно равенствам (3.2.21), (9.6.2), теорема верна для $n=1$.

Допустим теорема верна для $n=k$

$$\frac{d^kx}{dt^k} = x^*{}_*a^{k^*}{}_* \tag{9.6.11}$$

Согласно определению [15]-4.1.4, равенство

$$\frac{d^{k+1}x}{dt^{k+1}} = \frac{d}{dt}\frac{d^kx}{dt^k} = \frac{d}{dt}(x^*{}_*a^{k^*}{}_*) = \frac{dx}{dt}{}^*{}_*a^{k^*}{}_* = x^*{}_*a^*{}_*a^{k^*}{}_* = x^*{}_*a^{k+1^*}{}_*$$

является следствием равенств (3.2.20), (9.6.2), (9.6.11). Следовательно, теорема верна для $n=k+1$. □

Теорема 9.6.4. *Решение системы дифференциальных уравнений* (9.6.2) *при начальном условии*

$$t=0 \qquad x = \begin{pmatrix} x^1 \\ ... \\ x^n \end{pmatrix} = c = \begin{pmatrix} c^1 \\ ... \\ c^n \end{pmatrix}$$

*имеет вид*

$$x = c^*{}_*e^{ta^*{}_*} \tag{9.6.12}$$

Доказательство. Согласно теореме 9.6.3 и определению ряда Тейлора в разделе [15]-4.2, ряд Тейлора для решения имеет вид

$$x = c^*{}_*\sum_{n=0}^{\infty}\frac{1}{n!}t^na^{n^*}{}_* \tag{9.6.13}$$

Равенство

$$x = c^*{}_*\sum_{n=0}^{\infty}\frac{1}{n!}(ta)^{n^*}{}_* \tag{9.6.14}$$

является следствием равенства (9.6.13). Равенство (9.6.12) является следствием равенств (8.5.1), (9.6.14). □

## 9.7. Дифференциальное уравнение $\frac{dx}{dt} = a^*{}_*x$

Пусть $A$ - банахова $D$-алгебра. Система дифференциальных уравнений

$$\begin{aligned} \frac{dx_1}{dt} &= a_1^1x_1 + ... + a_1^nx_n \\ &...... \\ \frac{dx_n}{dt} &= a_n^1x_1 + ... + a_n^nx_n \end{aligned} \tag{9.7.1}$$

где $a_j^i \in A$ и $x_i : R \to A$ - $A$-значная функция действительной переменной, называется однородной системой линейных дифференциальных уравнений.

Положим

$$x = \begin{pmatrix} x_1 & ... & x_n \end{pmatrix} \qquad \frac{dx}{dt} = \begin{pmatrix} \frac{dx_1}{dt} & ... & \frac{dx_n}{dt} \end{pmatrix}$$

$$a=\begin{pmatrix}a_1^1 & ... & a_n^1\\ ... & ... & ...\\ a_1^n & ... & a_n^n\end{pmatrix}$$

Тогда мы можем записать систему дифференциальных уравнений (9.7.1) в матричной форме

(9.7.2) $$\frac{dx}{dt}=a^*{}_*x$$

**9.7.1. Решение в виде экспоненты** $x=e^{bt}c$ **.** Мы будем искать решение системы дифференциальных уравнений (9.7.2) в виде экспоненты

(9.7.3) $$x=e^{bt}c=\begin{pmatrix}e^{bt}c_1 & ... & e^{bt}c_n\end{pmatrix}\qquad c=\begin{pmatrix}c_1 & ... & c_n\end{pmatrix}$$

Согласно теореме 9.1.6, равенство

(9.7.4) $$\frac{dx_i}{dt}=e^{bt}bc_i=be^{bt}c_i=bx_i$$

является следствием равенства (9.7.3). Равенство

(9.7.5) $$e^{bt}bc=a^*{}_*(e^{bt}c)$$

является следствием равенств (9.7.2), (9.7.4).

Теорема 9.7.1. *Пусть $b$ является $^*{}_*$-собственным значением матрицы $a$. Из условия*

(9.7.6) $$a_j^i\in A[b]\qquad i=1,...,n\quad j=1,...,n$$

*следует, что матрица отображений (9.7.3) является решением системы дифференциальных уравнений (9.7.2) для $^*{}_*$-собственной строки $c$.*

Доказательство. Согласно теореме 9.1.3, равенство

(9.7.7) $$e^{bt}bc=e^{bt}(a^*{}_*c)$$

является следствием равенства (9.7.5). Поскольку выражение $e^{bt}$, вообще говоря, отлично от 0, равенство

$$bc=a^*{}_*c$$

является следствием равенства (9.7.7). Согласно определению 5.6.2, $A$-число $b$ является $^*{}_*$-собственным значением матрицы $a$ и матрица $c$ является $^*{}_*$-собственной строкой матрицы $a$, соответствующим $^*{}_*$-собственному значению $b$. □

Теорема 9.7.2. *Пусть $b$ является $^*{}_*$-собственным значением матрицы $a$. Пусть элементы матрицы $a$ не удовлетворяют условию (9.7.6). Если элементы $^*{}_*$-собственной строки $c$ удовлетворяют условию*

(9.7.8) $$c^i\in A[b]\qquad i=1,...,n$$

*следует, что матрица отображений (9.7.3) является решением системы дифференциальных уравнений (9.7.2).*

Доказательство. Согласно теореме 9.1.3, равенство

$$bce^{bt} = (a^*{}_*c)e^{bt} \tag{9.7.9}$$

является следствием равенства (9.7.5). Поскольку выражение $e^{bt}$, вообще говоря, отлично от 0, равенство

$$bc = a^*{}_*c$$

является следствием равенства (9.7.9). Согласно определению 5.6.2, $A$-число $b$ является ${}^*{}_*$-собственным значением матрицы $a$ и матрица $c$ является ${}^*{}_*$-собственной строкой матрицы $a$, соответствующим ${}_*{}^*$-собственному значению $b$. □

Пусть элементы матрицы $a$ не удовлетворяют условию (9.3.6). Пусть элементы ${}_*{}^*$-собственной строки $c$ не удовлетворяют условию (9.7.8). Тогда матрица отображений (9.7.3) не является решением системы дифференциальных уравнений (9.7.2).

### 9.7.2. Метод последовательного дифференцирования.

Теорема 9.7.3. *Последовательно дифференцируя систему дифференциальных уравнений* (9.7.2) *мы получаем цепочку систем дифференциальных уравнений*

$$\frac{d^n x}{dt^n} = a^{n^*{}_*{}^*}{}_*x \tag{9.7.10}$$

Доказательство. Мы докажем теорему индукцией по $n$.

Согласно равенствам (3.2.21), (9.7.2), теорема верна для $n = 1$.

Допустим теорема верна для $n = k$

$$\frac{d^k x}{dt^k} = a^{k^*{}_*{}^*}{}_*x \tag{9.7.11}$$

Согласно определению [15]-4.1.4, равенство

$$\frac{d^{k+1}x}{dt^{k+1}} = \frac{d}{dt}\frac{d^k x}{dt^k} = \frac{d}{dt}(a^{k^*{}_*{}^*}{}_*x) = a^{k^*{}_*{}^*}{}_*\frac{dx}{dt} = a^{k^*{}_*{}^*}{}_*a^*{}_*x = a^{k+1^*{}_*{}^*}{}_*x$$

является следствием равенств (3.2.20), (9.7.2), (9.7.11). Следовательно, теорема верна для $n = k + 1$. □

Теорема 9.7.4. *Решение системы дифференциальных уравнений* (9.7.2) *при начальном условии*

$$t = 0 \qquad x = \begin{pmatrix} x^1 \\ \dots \\ x^n \end{pmatrix} = c = \begin{pmatrix} c^1 \\ \dots \\ c^n \end{pmatrix}$$

*имеет вид*

$$x = e^{ta^*{}_*{}^*}{}_*c \tag{9.7.12}$$

Доказательство. Согласно теореме 9.7.3 и определению ряда Тейлора в разделе [15]-4.2, ряд Тейлора для решения имеет вид

$$x = \sum_{n=0}^{\infty} \frac{1}{n!} t^n a^{n^*{}_*{}^*}{}_*c \tag{9.7.13}$$

Равенство

$$x = \sum_{n=0}^{\infty} \frac{1}{n!} (ta)^{n^*{}_*{}^*}{}_*c \tag{9.7.14}$$

является следствием равенства (9.7.13). Равенство (9.7.12) является следствием равенств (8.5.1), (9.7.14). □

## 9.8. Дифференциальное уравнение $\frac{dx}{dt} = x_*{}^*a$

Пусть $A$ - банахова $D$-алгебра. Система дифференциальных уравнений

$$\begin{aligned} \frac{dx_1}{dt} &= x_1 a_1^1 + ... + x_n a_1^n \\ &...... \\ \frac{dx_n}{dt} &= x_1 a_n^1 + ... + x_n a_n^n \end{aligned} \tag{9.8.1}$$

где $a_j^i \in A$ и $x^i : R \to A$ - $A$-значная функция действительной переменной, называется однородной системой линейных дифференциальных уравнений.

Положим

$$x = \begin{pmatrix} x_1 & ... & x_n \end{pmatrix} \quad \frac{dx}{dt} = \begin{pmatrix} \frac{dx_1}{dt} & ... & \frac{dx_n}{dt} \end{pmatrix}$$

$$a = \begin{pmatrix} a_1^1 & ... & a_n^1 \\ ... & ... & ... \\ a_1^n & ... & a_n^n \end{pmatrix}$$

Тогда мы можем записать систему дифференциальных уравнений (9.8.1) в матричной форме

$$\frac{dx}{dt} = x_*{}^*a \tag{9.8.2}$$

### 9.8.1. Решение в виде экспоненты $x = ce^{bt}$ .

Мы будем искать решение системы дифференциальных уравнений (9.8.2) в виде экспоненты

$$x = ce^{bt} = \begin{pmatrix} c_1 e^{bt} & ... & c_n e^{bt} \end{pmatrix} \quad c = \begin{pmatrix} c_1 & ... & c_n \end{pmatrix} \tag{9.8.3}$$

Согласно теореме 9.1.6, равенство

$$\frac{dx_i}{dt} = c_i e^{bt} b = c_i b e^{bt} = b x_i \tag{9.8.4}$$

является следствием равенства (9.8.3). Равенство

$$ce^{bt} b = cbe^{bt} = (ce^{bt})_*{}^*a \tag{9.8.5}$$

является следствием равенств (9.8.2), (9.8.4).

Теорема 9.8.1. *Пусть $b$ является ${}_*{}^*$-собственным значением матрицы $a$. Из условия*

$$a_j^i \in A[b] \quad i = 1, ..., n \quad j = 1, ..., n \tag{9.8.6}$$

*следует, что матрица отображений (9.8.3) является решением системы дифференциальных уравнений (9.8.2) для ${}_*{}^*$-собственной строки $c$.*

Доказательство. Согласно теореме 9.1.3, равенство

$$cbe^{bt} = (c_*{}^*a)e^{bt} \tag{9.8.7}$$

является следствием равенства (9.8.5). Поскольку выражение $e^{bt}$, вообще говоря, отлично от 0, равенство

$$cb = c_*{}^*a$$

является следствием равенства (9.8.7). Согласно определению 5.6.1, $A$-число $b$ является ${}_*{}^*$-собственным значением матрицы $a$ и матрица $c$ является ${}_*{}^*$-собственной строкой матрицы $a$, соответствующим ${}^*{}_*$-собственному значению $b$. □

Теорема 9.8.2. *Пусть $b$ является ${}_*{}^*$-собственным значением матрицы $a$. Пусть элементы матрицы $a$ не удовлетворяют условию (9.8.6). Если элементы ${}_*{}^*$-собственной строки $c$ удовлетворяют условию*

$$c^i \in A[b] \qquad i = 1, ..., n \tag{9.8.8}$$

*следует, что матрица отображений (9.8.3) является решением системы дифференциальных уравнений (9.8.2).*

Доказательство. Согласно теореме 9.1.3, равенство

$$e^{bt}cb = e^{bt}(c_*{}^*a) \tag{9.8.9}$$

является следствием равенства (9.8.5). Поскольку выражение $e^{bt}$, вообще говоря, отлично от 0, равенство

$$cb = c_*{}^*a$$

является следствием равенства (9.8.9). Согласно определению 5.6.1, $A$-число $b$ является ${}_*{}^*$-собственным значением матрицы $a$ и матрица $c$ является ${}_*{}^*$-собственной строкой матрицы $a$, соответствующим ${}^*{}_*$-собственному значению $b$. □

Пусть элементы матрицы $a$ не удовлетворяют условию (9.8.6). Пусть элементы ${}_*{}^*$-собственной строки $c$ не удовлетворяют условию (9.6.8). Тогда матрица отображений (9.6.3) не является решением системы дифференциальных уравнений (9.8.2).

### 9.8.2. Метод последовательного дифференцирования.

Теорема 9.8.3. *Последовательно дифференцируя систему дифференциальных уравнений (9.8.2) мы получаем цепочку систем дифференциальных уравнений*

$$\frac{d^n x}{dt^n} = x_*{}^*a^{n_*{}^*} \tag{9.8.10}$$

Доказательство. Мы докажем теорему индукцией по $n$.

Согласно равенствам (3.2.17), (9.8.2), теорема верна для $n = 1$.

Допустим теорема верна для $n = k$

$$\frac{d^k x}{dt^k} = x_*{}^*a^{k_*{}^*} \tag{9.8.11}$$

Согласно определению [15]-4.1.4, равенство

$$\frac{d^{k+1}x}{dt^{k+1}} = \frac{d}{dt}\frac{d^k x}{dt^k} = \frac{d}{dt}(x_*{}^*a^{k_*{}^*}) = \frac{dx}{dt}{}_*{}^*a^{k_*{}^*} = x_*{}^*a_*{}^*a^{k_*{}^*} = x_*{}^*a^{k+1_*{}^*}$$

является следствием равенств (3.2.16), (9.8.2), (9.8.11). Следовательно, теорема верна для $n = k + 1$. □

Теорема 9.8.4. *Решение системы дифференциальных уравнений* (9.8.2) *при начальном условии*

$$t = 0 \qquad x = \begin{pmatrix} x^1 \\ \dots \\ x^n \end{pmatrix} = c = \begin{pmatrix} c^1 \\ \dots \\ c^n \end{pmatrix}$$

*имеет вид*

$$x = c_*{}^* e^{ta_*{}^*} \tag{9.8.12}$$

Доказательство. Согласно теореме 9.8.3 и определению ряда Тейлора в разделе [15]-4.2, ряд Тейлора для решения имеет вид

$$x = c_*{}^* \sum_{n=0}^{\infty} \frac{1}{n!} t^n a^{n_*{}^*} \tag{9.8.13}$$

Равенство

$$x = c_*{}^* \sum_{n=0}^{\infty} \frac{1}{n!} (ta)^{n_*{}^*} \tag{9.8.14}$$

является следствием равенства (9.8.13). Равенство (9.8.12) является следствием равенств (8.5.1), (9.8.14). □

## 9.9. Краткий итог

Метод последовательного дифференцирования является важным методом решения дифференциальных уравнений, особенно в тех случаях, когда трудно найти другой метод решения. Однако ситуация меняется, когда мы рассматриваем дифференциальное уравнение над некоммутативной банаховой $D$-алгеброй.

Для дифференциального уравнения

$$\frac{dy}{dx} = 3x \otimes x \tag{9.9.1}$$

метод последовательного дифференцирования порождает отображение, которое не является решением дифференциального уравнения (9.9.1).

Согласно теореме 9.3.4, решение системы дифференциальных уравнений

$$\begin{aligned} \frac{dx^1}{dt} &= x^2 \\ \frac{dx^2}{dt} &= -x^1 \end{aligned} \tag{9.9.2}$$

которое удовлетворяет начальному условию

$$t = 0 \quad x^1 = 0 \quad x^2 = 1$$

имеет вид

$$\begin{aligned} x^1 &= \sin t \\ x^2 &= \cos t \end{aligned} \tag{9.9.3}$$

Согласно примеру 9.5.1, отображение

$$x^{1} = \frac{1}{2}(-i+j)(e^{it}-e^{jt}) \qquad (9.9.4)$$
$$x^{2} = \frac{1}{2}(-i+j)(ie^{it}-je^{jt})$$

также является решением системы дифференциальных уравнений (9.9.2) которое удовлетворяет начальному условию

$$t=0 \quad x^{1}=0 \quad x^{2}=1$$

Однако отображения (9.9.3), (9.9.4) имеют различные разложения в ряд Тейлора.

Из этих двух экстремально противоположных примеров следует, что мы должны быть осторожны с применением метода последовательного дифференцирования, чтобы решить дифференциальное уравнение над некоммутативной банаховой $D$-алгеброй.

Из рассуждения выше следует вопрос. Если мы имеем бесконечное семейство решений, как описать это семейство и как выбрать решение? Для ответа на этот вопрос хорошо бы иметь побольше примеров применения рассмотренной теории к решению различных задач. В частности, мне понравились примеры, описанные в разделе [17]-1.1.

Приложение А

# Сводка теорем

Пусть $D$ - полное коммутативное кольцо характеристики 0.

## А.1. Таблица производных

Теорема А.1.1. *Для любого* $b \in A$
$$\frac{db}{dx} = 0 \otimes 0$$

Доказательство. Теорема является следствием теоремы [15]-B.1.1. □

Теорема А.1.2.
$$\frac{dx}{dx} \circ dx = dx \qquad \frac{dx}{dx} = 1 \otimes 1 \tag{A.1.1}$$

Доказательство. Теорема является следствием теоремы [15]-B.1.2. □

Теорема А.1.3. *Для любых* $b,\ c \in A$
$$\left\{\begin{array}{c} \dfrac{dbf(x)c}{dx} = (b\otimes c)\circ\dfrac{df(x)}{dx} \\ \dfrac{dbf(x)c}{dx}\circ dx = b\left(\dfrac{df(x)}{dx}\circ dx\right)c \\ \\ \dfrac{d_{s\cdot 0}bf(x)c}{dx} = b\dfrac{d_{s\cdot 0}f(x)}{dx} \\ \\ \dfrac{d_{s\cdot 1}bf(x)c}{dx} = \dfrac{d_{s\cdot 1}f(x)}{dx}c \end{array}\right. \tag{A.1.2}$$

*Для любого* $F \in A \otimes A$,
$$\left\{\begin{array}{c} \dfrac{dF\circ f(x)}{dx} = F\circ\dfrac{df(x)}{dx} \\ \dfrac{dF\circ f(x)}{dx}\circ dx = F\circ\left(\dfrac{df(x)}{dx}\circ dx\right) \end{array}\right. \tag{A.1.3}$$

Доказательство. Теорема является следствием теоремы [15]-B.1.3. □

Теорема А.1.4. *Пусть*
$$f : A \to B$$
$$g : A \to B$$
*отображения банахова* $D$*-модуля* $A$ *в ассоциативную банахову* $D$*-алгебра* $B$. *Если производные* $\dfrac{df(x)}{dx}$, $\dfrac{dg(x)}{dx}$ *существуют, то существует производная*

$$\frac{d(f(x)+g(x))}{dx}$$

$$\frac{d(f(x)+g(x))}{dx}\circ dx=\frac{df(x)}{dx}\circ dx+\frac{dg(x)}{dx}\circ dx \tag{A.1.4}$$

$$\frac{d(f(x)+g(x))}{dx}=\frac{df(x)}{dx}+\frac{dg(x)}{dx} \tag{A.1.5}$$

*Если*

$$\frac{df(x)}{dx}=\frac{d_{s\cdot 0}f(x)}{dx}\otimes\frac{d_{s\cdot 1}f(x)}{dx} \tag{A.1.6}$$

$$\frac{dg(x)}{dx}=\frac{d_{t\cdot 0}g(x)}{dx}\otimes\frac{d_{t\cdot 1}g(x)}{dx} \tag{A.1.7}$$

*то*

$$\frac{d(f(x)+g(x))}{dx}=\frac{d_{s\cdot 0}f(x)}{dx}\otimes\frac{d_{s\cdot 1}f(x)}{dx}+\frac{d_{t\cdot 0}g(x)}{dx}\otimes\frac{d_{t\cdot 1}g(x)}{dx} \tag{A.1.8}$$

Доказательство. Теорема является следствием теоремы [15]-B.1.4. □

Теорема A.1.5. *Для любых $b$, $c\in A$*

$$\begin{cases}\frac{dbxc}{dx}=b\otimes c & \frac{dbxc}{dx}\circ dx=b\,dx\,c\\ \frac{d_{1\cdot 0}bxc}{dx}=b & \frac{d_{1\cdot 1}bxc}{dx}=c\end{cases} \tag{A.1.9}$$

Доказательство. Следствие теорем A.1.2, A.1.3, когда $f(x)=x$. □

Теорема A.1.6. *Пусть $f$ - линейное отображение*

$$f\circ x=(a_{s\cdot 0}\otimes a_{s\cdot 1})\circ x=a_{s\cdot 0}\ x\ a_{s\cdot 1}$$

*Тогда*

$$\frac{\partial f\circ x}{\partial x}=f$$

$$\frac{\partial f\circ x}{\partial x}\circ dx=f\circ dx$$

Доказательство. Следствие теорем A.1.4, A.1.5, [15]-3.3.13. □

Следствие A.1.7. *Для любого $b\in A$*

$$\begin{cases}\frac{d(xb-bx)}{dx}=1\otimes b-b\otimes 1\\ \frac{d(xb-bx)}{dx}\circ dx=dx\,b-b\,dx\\ \frac{d_{1\cdot 0}(xb-bx)}{dx}=1 & \frac{d_{1\cdot 1}(xb-bx)}{dx}=b\\ \frac{d_{2\cdot 0}(xb-bx)}{dx}=-b & \frac{d_{2\cdot 1}(xb-bx)}{dx}=1\end{cases}$$

□

Теорема A.1.8. *Пусть $D$ - полное коммутативное кольцо характеристики* 0. *Пусть $A$ - ассоциативная банахова $D$-алгебра. Тогда*[A.1]

$$\frac{dx^2}{dx}=x\otimes 1+1\otimes x \tag{A.1.10}$$

$$dx^2=x\,dx+dx\,x \tag{A.1.11}$$

$$\left\{\begin{array}{ll}\frac{d_{1\cdot 0}x^2}{dx}=x & \frac{d_{1\cdot 1}x^2}{dx}=1\\ \frac{d_{2\cdot 0}x^2}{dx}=1 & \frac{d_{2\cdot 1}x^2}{dx}=x\end{array}\right. \tag{A.1.12}$$

Доказательство. Теорема является следствием теоремы [15]-B.1.15. □

Теорема A.1.9. *Пусть $D$ - полное коммутативное кольцо характеристики* 0. *Пусть $A$ - ассоциативная банахова $D$-алгебра. Тогда*

$$\frac{dx^3}{dx}=x^2\otimes 1+x\otimes x+1\otimes x^2 \tag{A.1.13}$$

$$dx^3=x^2\,dx+x\,dx\,x+dx\,x^2 \tag{A.1.14}$$

Доказательство. Согласно теореме A.1.10,

$$\frac{dx^3}{dx}=\frac{dx^2x}{dx}=\frac{dx^2}{dx}x+x^2\frac{dx}{dx}=(x\otimes 1+1\otimes x)x+x^2(1\otimes 1) \tag{A.1.15}$$

Равенство (A.1.13) является следствием равенства (A.1.15). Равенство (A.1.14) является следствием равенства (A.1.13) и определения [15]-3.3.2. □

Теорема A.1.10. *Пусть $A$ - банаховый $D$-модуль. Пусть $B$ - банахова $D$-алгебра. Пусть $f$, $g$ - дифференцируемые отображения*

$$f:A\rightarrow B\quad g:A\rightarrow B$$

*Производная удовлетворяет соотношению*

$$\frac{df(x)g(x)}{dx}\circ dx=\left(\frac{df(x)}{dx}\circ dx\right)\ g(x)+f(x)\ \left(\frac{dg(x)}{dx}\circ dx\right)$$

$$\frac{df(x)g(x)}{dx}=\frac{df(x)}{dx}\ g(x)+f(x)\ \frac{dg(x)}{dx} \tag{A.1.16}$$

Доказательство. Теорема является следствием теоремы [15]-3.3.17. □

---

[A.1] Утверждение теоремы аналогично примеру VIII, [22], с. 451. Если произведение коммутативно, то равенство (A.1.10) принимает вид

$$dx^2\circ dx=2x\ dx$$

$$\frac{dx^2}{dx}=2x$$

Теорема A.1.11. *Пусть $A$ - банаховый $D$-модуль. Пусть $B$, $C$ - банаховые $D$-алгебры. Пусть $f$, $g$ - дифференцируемые отображения*

$$f : A \to B \quad g : A \to C$$

*Производная удовлетворяет соотношению*

$$\frac{df(x)\otimes g(x)}{dx}\circ a = \left(\frac{df(x)}{dx}\circ a\right)\otimes g(x) + f(x)\otimes\left(\frac{dg(x)}{dx}\circ a\right)$$

$$\frac{df(x)\otimes g(x)}{dx} = \frac{df(x)}{dx}\otimes g(x) + f(x)\otimes\frac{dg(x)}{dx}$$

Доказательство. Теорема является следствием теоремы [15]-3.3.20. □

Теорема A.1.12. *Пусть $D$-алгебра $A$ имеет единицу $e$. Пусть*

$$f : A \to A$$

*дифференцируемое отображение и*

$$P_0 : d \in D \to de_0 \in A$$

*вложение кольца $D$ в $D$-алгебру $A$. Тогда*

$$\frac{df(te)}{dt} = \frac{df(x)}{dx}\circ e \tag{A.1.17}$$

Доказательство. Согласно теореме [15]-3.3.23,

$$\frac{df(te)}{dt} = \frac{df(te)}{d\,te}\circ\frac{d\,te}{dt} \tag{A.1.18}$$

Равенство (A.1.17) является следствием равенства (A.1.18). □

Теорема A.1.13. *Пусть $D$ - полное коммутативное кольцо характеристики* 0. *Пусть $A$ - ассоциативная банаховая $D$-алгебра. Тогда*

$$\frac{d^n}{dx^n}\frac{d^m y}{dx^m} = \frac{d^{n+m}y}{dx^{n+m}} \tag{A.1.19}$$

Доказательство. Мы докажем теорему индукцией по $n$.

Теорема очевидна для $n = 0$. Для $n = 1$, теорема является следствием определения [15]-4.1.4.

Пусть терорема верна для $n = k$

$$\frac{d^k}{dx^k}\frac{d^m y}{dx^m} = \frac{d^{k+m}y}{dx^{k+m}} \tag{A.1.20}$$

Равенство

$$\frac{d^{k+1}}{dx^{k+1}}\frac{d^m y}{dx^m} = \frac{d}{dx}\left(\frac{d^k}{dx^k}\frac{d^m y}{dx^m}\right) = \frac{d}{dx}\frac{d^{k+m}y}{dx^{k+m}} = \frac{d^{k+1+m}y}{dx^{k+1+m}} \tag{A.1.21}$$

является следствием равенств [15]-(4.1.5), (A.1.20). Следовательно, терорема верна для $n = k + 1$. □

## А.2. Таблица интегралов

Теорема А.2.1. *Пусть отображение*

$$f : A \to B$$

*дифференцируемо. Тогда*

$$\int \frac{df(x)}{dx} \circ dx = f(x) + C \tag{А.2.1}$$

Доказательство. Теорема является следствием определения 6.1.1. □

Теорема А.2.2. *Пусть отображение*

$$f : A \to \mathcal{L}(D; A; B)$$

*интегрируемо. Тогда*

$$\frac{d}{dx} \int f(x) \circ dx = f(x) \tag{А.2.2}$$

Доказательство. Теорема является следствием определения 6.1.1. □

Теорема А.2.3.

$$\int (0 \otimes 0) \circ dx = C \tag{А.2.3}$$

Доказательство. Теорема является следствием теоремы А.1.1. □

Теорема А.2.4.

$$\int (f(x) + g(x)) \circ dx = \int f(x) \circ dx + \int g(x) \circ dx \tag{А.2.4}$$

Доказательство. Равенство (А.2.4) является следствием теорем А.1.4, А.2.1. □

Теорема А.2.5.

$$\int (f_{s\cdot 0} \otimes f_{s\cdot 1}) \circ dx = (f_{s\cdot 0} \otimes f_{s\cdot 1}) \circ x + C \tag{А.2.5}$$

$$\int f_{s\cdot 0}\, dx\, f_{s\cdot 1} = f_{s\cdot 0}\, x\, f_{s\cdot 1} + C \tag{А.2.6}$$

$$f_{s\cdot 0} \in A \qquad f_{s\cdot 1} \in A$$

Доказательство. Теорема является следствием теоремы [15]-B.2.2. □

Теорема А.2.6.

$$\int (1 \otimes x + x \otimes 1) \circ dx = x^2 + C \tag{А.2.7}$$

$$\int dx\, x + x\, dx = x^2 + C \tag{А.2.8}$$

Теорема А.2.7.

$$\int (1 \otimes x^2 + x \otimes x + x^2 \otimes 1) \circ dx = x^3 + C \tag{А.2.9}$$

$$\int dx\, x^2 + x\, dx\, x + x^2\, dx = x^3 + C \tag{А.2.10}$$

Доказательство. Теорема является следствием теоремы [15]-5.1.3. □

Теорема А.2.8.

$$\int (e^x \otimes 1 + 1 \otimes e^x) \circ dx = 2e^x + C \tag{A.2.11}$$

$$\int e^x\,dx + dx\,e^x = 2e^x + C \tag{A.2.12}$$

Доказательство. Теорема является следствием теоремы [15]-B.2.3. □

Теорема А.2.9.

$$\int (\sinh x \otimes 1 + 1 \otimes \sinh x) \circ dx = 2\cosh x + C \tag{A.2.13}$$

$$\int \sinh x\,dx + dx\,\sinh x = 2\cosh x + C \tag{A.2.14}$$

$$\int (\cosh x \otimes 1 + 1 \otimes \cosh x) \circ dx = 2\sinh x + C \tag{A.2.15}$$

$$\int \cosh x\,dx + dx\,\cosh x = 2\sinh x + C \tag{A.2.16}$$

Доказательство. Теорема является следствием теоремы [15]-B.2.4. □

Теорема А.2.10.

$$\int (\sin x \otimes 1 + 1 \otimes \sin x) \circ dx = -2\cos x + C \tag{A.2.17}$$

$$\int \sin x\,dx + dx\,\sin x = -2\cos x + C \tag{A.2.18}$$

$$\int (\cos x \otimes 1 + 1 \otimes \cos x) \circ dx = 2\sin x + C \tag{A.2.19}$$

$$\int \cos x\,dx + dx\,\cos x = 2\sin x + C \tag{A.2.20}$$

Доказательство. Теорема является следствием теоремы [15]-B.2.5. □

# Список литературы

# Предметный указатель

# Специальные символы и обозначения